\documentclass[11pt]{amsart}
\usepackage[margin=27mm]{geometry}
\usepackage{amsmath,amssymb,mathtools,bm,mathrsfs}
\usepackage{booktabs,array,tabularx,enumitem}
\usepackage[colorlinks=true,linkcolor=blue,citecolor=blue,urlcolor=blue,bookmarksnumbered=true]{hyperref}
\usepackage{microtype}
\usepackage{xcolor}
\newtheorem{theorem}{Theorem}[section]
\newtheorem{lemma}[theorem]{Lemma}
\newtheorem{proposition}[theorem]{Proposition}

\theoremstyle{definition}
\newtheorem{definition}[theorem]{Definition}

\newtheorem{remark}[theorem]{Remark}
\newtheorem*{unnumberedremark}{Remark}

\DeclareMathOperator{\supp}{supp}
\DeclareMathOperator{\diag}{diag}

\numberwithin{equation}{section}
\allowdisplaybreaks
\setlist{itemsep=2pt,topsep=4pt}
\title[Fixed-frequency and finite-time uniqueness in isotropic elasticity]{%
Fixed-Frequency and Finite-Time Uniqueness in Three-Dimensional Isotropic Elasticity}
\author[Huaian Diao]{Huaian Diao\textsuperscript{1}}
\author[Yang Dong]{Yang Dong\textsuperscript{1}}
\author[Umidjon Durdiev]{Umidjon Durdiev\textsuperscript{1,2,3}}
\author[Weisheng Zhou]{Weisheng Zhou\textsuperscript{1}}

\newcommand{\contactemails}{%
  \par\smallskip
  {\footnotesize \color{black}
  E-mail:\,
  \href{mailto:diao@jlu.edu.cn}{\textcolor{black}{diao@jlu.edu.cn}} (H. Diao);
  \href{mailto:1621284857@qq.com}{\textcolor{black}{1621284857@qq.com}} (Y. Dong);
  \href{mailto:umidjan93@list.ru}{\textcolor{black}{umidjan93@list.ru}} (U. Durdiev);
  \href{mailto:563323665zws@gmail.com}{\textcolor{black}{563323665zws@gmail.com}} (W. Zhou).
  }%
}

\makeatletter

\def\@setkeywords{%
  {\itshape Keywords:}\enspace
  \@keywords\@addpunct.
  \contactemails
}
\makeatother

\makeatletter
\@ifundefined{subjclassname@2020}{%
  \@namedef{subjclassname@2020}{%
    \textup{2020} Mathematics Subject Classification%
  }%
}{}
\makeatother

\makeatletter
\renewcommand{\@setauthors}{%
  \begin{center}
    \normalsize
    \def\and{\qquad}%
    \authors\par
    \medskip
    \small
    \textsuperscript{1}School of Mathematics, Jilin University,\\
    Changchun, P.\ R.\ of China\\[3pt]
    \textsuperscript{2}Bukhara State University, Uzbekistan\\[3pt]
    \textsuperscript{3}V.\ I.\ Romanovskiy Institute of Mathematics,\\
    Uzbekistan Academy of Science\\[5pt]
  
  \end{center}%
}
\makeatother

\hypersetup{pdftitle={Fixed-Frequency and Finite-Time Uniqueness in Three-Dimensional Isotropic Elasticity},
 pdfauthor={Huaian Diao, Yang Dong, Umidjon Durdiev, Weisheng Zhou}}
\date{}



\begin{document}
\begin{abstract}
We prove global uniqueness for two inverse boundary value problems
for three-dimensional isotropic elasticity on bounded domains with
$C^9$ boundary. The real Lam\'e parameters are of class $C^{16}$ and
satisfy strong convexity, and the positive mass density is of class
$C^{15}$. At one prescribed positive frequency where the Dirichlet
problem is uniquely solvable, the full displacement-to-traction map
determines the Lam\'e parameters and density. The proof recovers
boundary derivatives and uses three physical polarizations to obtain
matrix moment identities for the four-component reduction. These
identities give an exterior-normalized relative transport with
holomorphic dependence on the projective null direction. A finite-type
system then determines the Lam\'e parameters, and a limit of products
of elastic displacements determines the density by Fourier uniqueness.

For the time-dependent system with zero initial conditions and
time-independent coefficients, the full boundary response on a finite
interval determines the same coefficients when the observation time
exceeds twice the larger shear-wave filling time of the two candidate
media. We prove continuation of the boundary response using local
unique continuation and terminal states of the elastic wave equation.
A Laplace transform and analytic continuation transfer equality
of the causal responses to equality of the displacement-to-traction
maps at a common positive frequency. The time-harmonic theorem
then gives uniqueness.
\end{abstract}


\maketitle


\tableofcontents


\section{Introduction}\label{a28:setting}

\subsection{Elastic systems and inverse boundary problems}

We consider two inverse boundary problems for an isotropic elastic
medium. The time-harmonic problem uses the Dirichlet-to-Neumann map
at a prescribed positive frequency. The time-domain problem uses
the causal boundary response on a finite observation interval.
In each problem the unknowns are the Lam\'e parameters $\lambda,\mu$
and the mass density $\rho$.

Let $\Omega\subset\mathbb R^3$ be the bounded domain occupied by the
medium, with $C^9$ boundary.
We assume real Lam\'e parameters $\lambda,\mu\in C^{16}(\overline\Omega)$
and real mass density $\rho\in C^{15}(\overline\Omega)$.
Here $C^m(\overline\Omega)$ means that all Cartesian derivatives through
order $m$ extend continuously to $\overline\Omega$.
For a complex-valued displacement field
$u:\Omega\to\mathbb C^3$, the linearized strain tensor
and the Cauchy stress tensor are defined by
\begin{equation} \label{eq:stress}
 \varepsilon(u)
 =\tfrac12\bigl(\nabla u+(\nabla u)^T\bigr),\quad
 \sigma_{\lambda,\mu}(u) =2\mu\varepsilon(u)+\lambda(\operatorname{div}u)I_3.
\end{equation}
Here the superscript $T$ denotes matrix transposition, and $I_3$
is the $3\times3$ identity matrix. For vector fields $v$, we use the
gradient convention $(\nabla v)_{ij}=\partial_jv_i$.
The corresponding elasticity operator is given by
\begin{equation}\label{eq:elasticity-operator}
 L_{\lambda,\mu}u
 :=\operatorname{div}\sigma_{\lambda,\mu}(u).
\end{equation}
The coefficients are assumed to satisfy the strong convexity
and positivity conditions
\begin{equation}\label{a28:physical-assumptions}
 \mu>0,\qquad 3\lambda+2\mu>0,\qquad \rho>0
 \quad\text{on }\overline\Omega.
\end{equation}

For the time-harmonic problem, fix an angular frequency $\omega>0$.
The displacement field $u$ satisfies
\begin{equation}\label{en:elastic-system}
 L_{\lambda,\mu}u+\omega^2\rho u=0
 \quad\text{in }\Omega.
\end{equation}

Define the zeroth-order potential $W$ by
\begin{equation}\label{en:W-definition}
 W:=\omega^2\rho.
\end{equation}
The subsequent reduction employs $W$, whereas the inverse
boundary problem concerns the physical coefficient triple
$(\lambda,\mu,\rho)$. For the weak Dirichlet realization, write
\[
P_D^{\lambda,\mu,W}:
H_0^1(\Omega;\mathbb C^3)
\longrightarrow H^{-1}(\Omega;\mathbb C^3),
\qquad
P_D^{\lambda,\mu,W}v
:=(L_{\lambda,\mu}+WI_3)v.
\]
The nonresonance condition is the triviality of its kernel:
\begin{equation}\label{en:nonresonance}
\ker P_D^{\lambda,\mu,W}=\{0\}.
\end{equation}
Here the operator is understood in the weak sense, and
$H_0^1(\Omega;\mathbb C^3)$ incorporates the homogeneous
Dirichlet boundary condition.

Let
\[
\gamma_0:H^1(\Omega;\mathbb C^3)
\longrightarrow H^{1/2}(\partial\Omega;\mathbb C^3)
\]
be the Dirichlet trace operator. For each prescribed boundary
value $f\in H^{1/2}(\partial\Omega;\mathbb C^3)$, consider
\begin{equation}\label{en:dirichlet-problem}
\begin{cases}
(L_{\lambda,\mu}+WI_3)u=0
&\text{in }\Omega,\\
\gamma_0u=f
&\text{on }\partial\Omega.
\end{cases}
\end{equation}
Under the strong-convexity assumptions
\eqref{a28:physical-assumptions} and the nonresonance
condition~\eqref{en:nonresonance}, this problem is well posed:
for every $f\in H^{1/2}(\partial\Omega;\mathbb C^3)$, there
exists a unique weak solution $u_f\in H^1(\Omega;\mathbb C^3)$.
The well-posedness is established in
Subsection~\ref{an:tools}.
The outward unit normal to $\partial\Omega$ is denoted
by $\nu$. The Dirichlet-to-Neumann operator is defined by
\begin{equation}\label{en:dn-definition}
 \begin{aligned}
  \Lambda_{\lambda,\mu,\rho}^{\Omega}(\omega)
  &:H^{1/2}(\partial\Omega;\mathbb C^3)
    \longrightarrow
    H^{-1/2}(\partial\Omega;\mathbb C^3),\\
  \Lambda_{\lambda,\mu,\rho}^{\Omega}(\omega)f
  &=\left.\sigma_{\lambda,\mu}(u_f)\nu
    \right|_{\partial\Omega},
 \end{aligned}
\end{equation}
where the boundary traction is interpreted in the weak sense.

At the fixed frequency $\omega$, let $\mathcal A_\omega$
denote the class of all triples
\[
 (\lambda,\mu,\rho)\in
 C^{16}(\overline\Omega;\mathbb R)^2
 \times C^{15}(\overline\Omega;\mathbb R)
\]
satisfying the positivity conditions~\eqref{a28:physical-assumptions}
and the nonresonance condition \eqref{en:nonresonance}.  
The boundary map defines
\begin{equation}\label{en:inverse-uniqueness}
\begin{aligned}
\mathcal F_\omega:\mathcal A_\omega
&\longrightarrow
\mathcal L\bigl(H^{1/2}(\partial\Omega;\mathbb C^3),
H^{-1/2}(\partial\Omega;\mathbb C^3)\bigr),\\
(\lambda,\mu,\rho)
&\longmapsto \Lambda_{\lambda,\mu,\rho}^{\Omega}(\omega).
\end{aligned}
\end{equation}

The inverse boundary value problem is to determine
the coefficient triple $(\lambda,\mu,\rho)$ from
the fixed-frequency Dirichlet-to-Neumann map
$\Lambda_{\lambda,\mu,\rho}^{\Omega}(\omega)$.
With $\mathcal F_\omega$ defined in
\eqref{en:inverse-uniqueness}, the corresponding
uniqueness question is whether $\mathcal F_\omega$
is injective on $\mathcal A_\omega$. More precisely,
for any two admissible coefficient triples,
\[
\begin{aligned}
&\mathcal F_\omega(\lambda_1,\mu_1,\rho_1)
 =\mathcal F_\omega(\lambda_2,\mu_2,\rho_2) \qquad\text{only if}\qquad
(\lambda_1,\mu_1,\rho_1)
 =(\lambda_2,\mu_2,\rho_2)
 \quad\text{in }\overline\Omega.
\end{aligned}
\]
For $j=1,2$, we continue to write $W_j=\omega^2\rho_j$
for the zeroth-order potential associated with the fixed frequency.

We now consider the causal wave equation with time-independent
coefficients $\lambda,\mu,\rho$ satisfying~\eqref{a28:physical-assumptions}.
Let $T_{\rm obs}>0$ and take
$f\in C_c^\infty((0,T_{\rm obs});H^{1/2}(\partial\Omega;\mathbb C^3))$.
The causal displacement $u^f$ solves the initial-boundary value problem
\begin{equation}\label{td:wave}
 \begin{cases}
  \rho\,\partial_t^2u^f-L_{\lambda,\mu}u^f=0
    &\text{in }(0,T_{\rm obs})\times\Omega,\\
  \gamma_0u^f=f
    &\text{on }(0,T_{\rm obs})\times\partial\Omega,\\
  u^f|_{t=0}=\partial_tu^f|_{t=0}=0
    &\text{in }\Omega.
 \end{cases}
\end{equation}
Proposition~\ref{an:dynamic-wellposedness} constructs the unique
causal solution in the energy class and gives the response estimate
\eqref{td:response-bound}. With $\nu$ as above, define
\begin{equation}\label{td:response}
 \begin{aligned}
 R_{\lambda,\mu,\rho}^{\Omega,T_{\rm obs}}&:
 C_c^\infty\bigl((0,T_{\rm obs});H^{1/2}(\partial\Omega;\mathbb C^3)\bigr)\\
 &\hspace{4em}\longrightarrow
 \mathcal D'\bigl((0,T_{\rm obs});H^{-1/2}(\partial\Omega;\mathbb C^3)\bigr),\\
 R_{\lambda,\mu,\rho}^{\Omega,T_{\rm obs}}f
 &:=\left.\sigma_{\lambda,\mu}(u^f)\nu
              \right|_{(0,T_{\rm obs})\times\partial\Omega}.
 \end{aligned}
\end{equation}
The traction in~\eqref{td:response} is weak. For
$v\in H^1(\Omega;\mathbb C^3)$, its bilinear pairing is
\begin{equation}\label{td:weak-response}
 \bigl\langle (R_{\lambda,\mu,\rho}^{\Omega,T_{\rm obs}}f)(t),
              \gamma_0v\bigr\rangle
 =\int_\Omega
 \bigl(\sigma_{\lambda,\mu}(u^f(t)):\varepsilon(v)
       +\rho\,\partial_t^2u^f(t)\cdot v\bigr)\,dx
\end{equation}
in the sense of time distributions. Equation~\eqref{td:wave} makes
\eqref{td:weak-response} independent of the choice of $v$ with prescribed
trace. The input and measurement use every component of $\partial\Omega$.

Let $\mathcal A_{\rm dyn}$ be the class of triples in
$C^{16}(\overline\Omega;\mathbb R)^2\times C^{15}(\overline\Omega;\mathbb R)$
satisfying~\eqref{a28:physical-assumptions}. The time-domain inverse
problem is to determine $(\lambda,\mu,\rho)\in\mathcal A_{\rm dyn}$
from the operator
\begin{equation}\label{td:inverse-map}
 \mathcal F_{T_{\rm obs}}^{\rm dyn}(\lambda,\mu,\rho)
 :=R_{\lambda,\mu,\rho}^{\Omega,T_{\rm obs}}.
\end{equation}
The time sufficient for uniqueness is expressed in the shear metric.
For the $j$th candidate medium, set
\begin{equation}\label{td:filling-time}
 c_{S,j}:=\sqrt{\mu_j/\rho_j},\qquad
 g_{S,j}:=c_{S,j}^{-2}g_E,\qquad
 \tau_{S,j}:=\sup_{x\in\Omega} d_{g_{S,j}}(x,\partial\Omega),
\end{equation}
where $g_E$ is the Euclidean metric. The distance is the infimum of
$g_{S,j}$-lengths of piecewise smooth paths from $x$ to $\partial\Omega$
that remain in $\Omega$ before their boundary endpoint. Strong convexity
gives $\lambda_j+\mu_j>0$, so the pressure speed
$c_{P,j}:=\sqrt{(\lambda_j+2\mu_j)/\rho_j}$ satisfies $c_{P,j}>c_{S,j}$.
We call $\tau_{S,j}$ the shear-wave filling time and omit $j$ when
considering a single medium.

\subsection{Uniqueness theorems}
\begin{theorem}\label{a28:main}
Let $\Omega\subset\mathbb R^3$ be a bounded domain with
$C^9$ boundary, and fix $\omega>0$.
For $j=1,2$, let
\[
 \lambda_j,\mu_j\in C^{16}(\overline\Omega;\mathbb R),
 \qquad
 \rho_j\in C^{15}(\overline\Omega;\mathbb R),
\]
and suppose that the coefficients satisfy the strong convexity and positivity
conditions~\eqref{a28:physical-assumptions}.
For $j=1,2$, set
\[
P_{D,j}:=P_D^{\lambda_j,\mu_j,W_j}.
\]
Assume that both physical Dirichlet realizations are nonresonant:
\[
\ker P_{D,j}=\{0\},\qquad j=1,2.
\]
If the Dirichlet-to-Neumann maps satisfy
\begin{equation}\label{a28:DN}
 \Lambda_{\lambda_1,\mu_1,\rho_1}^{\Omega}(\omega)
 =
 \Lambda_{\lambda_2,\mu_2,\rho_2}^{\Omega}(\omega),
\end{equation}
then
\[
 \lambda_1=\lambda_2,\qquad
 \mu_1=\mu_2,\qquad
 \rho_1=\rho_2
 \quad\text{on }\overline\Omega.
\]
\end{theorem}

Thus $\mathcal F_\omega$ is injective on $\mathcal A_\omega$.
Proposition~\ref{bl:finite-jets} supplies the boundary derivatives used
in the proof of this theorem. Theorem~\ref{td:uniqueness} concerns
measurements in the time-domain model~\eqref{td:wave}.

\begin{theorem}[Finite-time full-boundary uniqueness]\label{td:uniqueness}
Let $\Omega\subset\mathbb R^3$ be a bounded domain with $C^9$ boundary.
For $j=1,2$, let $(\lambda_j,\mu_j,\rho_j)\in\mathcal A_{\rm dyn}$
and let $\tau_{S,j}$ be given by~\eqref{td:filling-time}.
Suppose the observation time and full boundary responses satisfy
\begin{equation}\label{td:data-equality}
 \begin{gathered}
 T_{\rm obs}>2\max\{\tau_{S,1},\tau_{S,2}\},\\
 R_{\lambda_1,\mu_1,\rho_1}^{\Omega,T_{\rm obs}}
 =
 R_{\lambda_2,\mu_2,\rho_2}^{\Omega,T_{\rm obs}}.
 \end{gathered}
\end{equation}
Then
\[
 \lambda_1=\lambda_2,\qquad \mu_1=\mu_2,\qquad
 \rho_1=\rho_2\quad\text{on }\overline\Omega.
\]
\end{theorem}

Theorem~\ref{td:uniqueness} gives injectivity of
$\mathcal F_{T_{\rm obs}}^{\rm dyn}$ on the class of triples with
$\tau_S<T_{\rm obs}/2$. The proof in Section~\ref{sec:dynamical}
selects a frequency $\omega_*$ with $\omega_*^2$ below the spectra
of both Dirichlet generators and applies Theorem~\ref{a28:main}
at that frequency.

\subsection{Related work}

Boundary determination for elastic systems has been studied
at both the level of the full boundary Taylor series and
through explicit reconstruction formulas. Nakamura and
Uhlmann~\cite{NU1995} determined the full boundary Taylor series
for isotropic media in dimensions $n\geq2$ and for generic
anisotropic tensors in dimension two. In the isotropic setting,
explicit reconstruction formulas for the Lam\'e
moduli and their higher-order normal derivatives were
obtained in~\cite{LN2017}. For the time-dependent isotropic
elastic wave equation, the boundary values and derivatives
of the Lam\'e parameters and density were recovered from
the dynamical Dirichlet-to-Neumann map in~\cite{dHNZ2017}.

For the static inverse boundary value problem in isotropic
elasticity, a global uniqueness result under a smallness
condition on the gradients of the shear moduli was established
in the corrected theorem of Nakamura and
Uhlmann~\cite[Theorem~1]{NU2003}.  Complex geometrical optics methods
and matrix-valued Cauchy--Riemann systems for isotropic
elasticity were developed by Eskin and Ralston
\cite{ER2002,ER2004}. Beretta, Francini, and Vessella~\cite{BFV2014}
proved uniqueness and Lipschitz stability for Lam\'e parameters
that are piecewise constant on a prescribed finite partition,
using local data on a flat boundary patch and chains of flat
interfaces. For smooth static isotropic elasticity,
Chen et al.~\cite{CJLT2026} obtained full-data uniqueness for
axisymmetric, separable, or one-direction quasianalytic shear
moduli $\mu$, under the corresponding domain hypotheses. The preprint~\cite{OpenAI2026}
establishes full-data uniqueness for the smooth static isotropic
problem in three dimensions. In the static setting, two of the
authors of the present paper established a
global uniqueness result for smooth Lam\'e parameters
using physical matrix moment identities and exact relative
transport~\cite{DD2026}. A subsequent work~\cite{DDStab2026}
established logarithmic interior stability on bounded, uniformly
strongly convex $C^{22}$ classes and stability estimates for
boundary derivatives. The boundary exponents are sharp on admissible
classes containing a constant pair satisfying the a priori
bounds strictly.

Density recovery for elastic systems has been studied
under various assumptions on the elastic coefficients
and the available boundary or scattering data.
For constant Lam\'e parameters, H\"ahner~\cite{Hahner2002}
proved uniqueness of an inhomogeneous scalar density from
far-field patterns for all incident plane waves at a single
frequency. For a homogeneous isotropic elastic background
and a symmetric matrix-valued density, Barcel\'o et al.~\cite{Barcelo2018}
proved uniqueness from boundary data and separately from
scattering data. In the time-harmonic setting,
simultaneous uniqueness and Lipschitz stability for
the Lam\'e parameters and density were established
in~\cite{BDFVZ2017} for coefficients that are
piecewise constant on a prescribed finite partition,
under suitable geometric and frequency conditions.
For the time-dependent elastic wave equation,
Zhai~\cite[Theorem~1]{Zhai2026} proved density uniqueness when
$c_P$ and $c_S$ are known, the shear metric satisfies the strictly
convex foliation condition, and the observation time exceeds
all shear geodesic lengths. Joint recovery of the Lam\'e parameters
and density follows under strictly convex foliation conditions
for both wave metrics and an observation time exceeding all geodesic
lengths in both metrics~\cite[Corollary~1]{Zhai2026}.

The dynamical argument in Section~\ref{sec:dynamical} follows the
boundary-control approach of Belishev and Lasiecka~\cite{BL2002}.
For full-boundary data, we prove the continuation step using
local unique continuation along shear-length paths and terminal
states that form a core for the Dirichlet generator. A full-time
Laplace transform and analytic continuation in the stationary
spectral parameter then give the displacement-to-traction map
at a common positive nonresonant frequency. The resulting
uniqueness statement uses observation times exceeding twice the
larger shear-wave filling time of the two media.

The present work establishes global uniqueness for the
finite-regularity isotropic triples $(\lambda,\mu,\rho)$ in
Theorem~\ref{a28:main} from the full displacement-to-traction map
at one fixed known positive frequency. The positive-frequency
map determines the mass density as well as the Lam\'e parameters. The static matrix-moment and relative-transport
framework of \cite{DD2026,DDStab2026} provides
a methodological starting point, while the fixed-frequency problem
requires a separate analysis. Equality of the stationary
Dirichlet-to-Neumann maps does not identify the corresponding
static maps, and the unknown potential $W=\omega^2\rho$ together
with its gradient $\nabla W$ enters the four-component reduction.
The reduced Green--Betti identity is therefore derived directly
from the physical fixed-frequency boundary data, with all
potential-dependent terms retained. Three physical polarizations
of the rank-three characteristic quotient yield the complete
$3\times3$ matrix moment identities. The resulting exact relative
transport, its determinant identity, and its holomorphic dependence
on the projective null direction lead to finite-type rigidity for
the Lam\'e parameters. A further family of physical elastic
solutions then recovers the density.

\subsection{Main analytical challenges}

The fixed-frequency inverse problem requires the Lam\'e parameters
to be determined while the mass density is still unknown. In the
time-harmonic system~\eqref{en:elastic-system}, the density appears
through the zeroth-order potential $W=\omega^2\rho$. Consequently,
the principal-symbol information alone does not separate all three
unknown coefficients. The proof addresses four related difficulties.

First, the four-component reduction~\eqref{bt:MW}, obtained from the
displacement representation~\eqref{bt:TS}, has principal part
$\Delta I_4$ but contains both $W$ and $\nabla W$ in its zeroth-order
coefficient matrix. The Green--Betti formula~\eqref{bt:defects} for
this reduced system produces additional boundary terms beyond the
pairing of the physical Dirichlet-to-Neumann maps. The matching
conditions~\eqref{bt:matching}, which follow from the boundary
determination result~\ref{bl:finite-jets}, eliminate these terms and yield
the exact identity~\eqref{bt:betti}.

Second, a complete $3\times3$ matrix of moment identities must be
realized by physical elastic solutions. The coefficient matrix in the leading four-component transport
equation~\eqref{pb27:two-sign-transport} has the common null
vector $(\theta,0)^T$. Factoring out its span leads to the
three-dimensional characteristic quotient described before
\eqref{i2:eq-frames}. The scaled lift~\eqref{i2:eq-scaled} provides
only two independent leading dual polarizations. A third polarization
requires the unscaled lift~\eqref{i2:eq-unscaled}, for which the
error in the unscaled physical lift involves both
$\|r_h\|_{H^1}$ and $h^{-1}\|r_h\|_{L^2}$. The same remainder is therefore constructed
so that
\[
\|r_h\|_{L^2(\Omega_0)}\le Ch^2,
\qquad
\|r_h\|_{H^1(\Omega_0)}\le Ch,
\]
as in~\eqref{pb27:two}. The transverse cutoff correction
\eqref{i2:eq-commutator} gives solutions with the prescribed leading
frame near the support of the transverse test. The quotient frames are
independent of the tests in the full matrix identities~\eqref{i2:eq-nine}.
The potential and potential-gradient differences enter only
through the zeroth-order matrix in~\eqref{bt:MW}. Under the
semiclassical scaling used to extract the matrix moments, their
contributions carry an additional factor $h$ and vanish in the
leading-order limit.

Third, the matrix moments must determine a relative transport
that is independent of local frame choices and normalized on a
common exterior region. The Cauchy transform formula
\eqref{if:intrinsic-formula} defines this matrix $H$. Its transport
equation and exterior normalization give the determinant identity
\[
\det H=\sqrt{\mu_2/\mu_1}>0,
\]
proved in~\eqref{if:detH}. Differentiating the transport equation at a fixed physical point $x$
gives holomorphic dependence on the projective null direction. The resulting
holomorphic block structure and transport equations yield the
finite-type system~\eqref{i4:closed-system}, which determines
$\mu$ and then $\lambda$ throughout the admissible coefficient
class.

Finally, after the Lam\'e parameters have been determined, the
physical Green identity~\eqref{pt:betti} isolates $W_1-W_2$.
A further family of exact elastic solutions is needed to obtain
the product limit~\eqref{pt:limit}. Using transverse tests to
remove the resulting weight, Fourier uniqueness then gives
$W_1=W_2$ and hence $\rho_1=\rho_2$.

For the time-domain problem, the finite response must first be related
to the full causal response. The observation time is controlled by the
shear distance~\eqref{td:filling-time}. In the duality argument, an
adjoint wave with zero Dirichlet trace and zero measured traction
extends by zero across $\partial\Omega$. For each interior point, the
strict time margin permits a family of timelike level surfaces along a
path of shear length less than the control time. Local unique continuation
for the shear and pressure parts yields density of the terminal
displacements. Recovering the generator from the Gram form of these
displacements requires density also in the graph norm of the generator.
A finite-interval Laplace transform retains the interior displacement
and velocity at the final time. The proof therefore first continues
the response and then transforms the full causal solution.

\subsection{The proof strategy}

Boundary determination gives derivatives of $\lambda$ and $\mu$
through order $16$, and of $W$ through order $14$, at $\partial\Omega$;
see~\ref{bl:finite-jets}. Proposition~\ref{ext:finite-common}
gives common $C^{16}$ Lam\'e extensions and separate $C^{15}$
extensions of $W_j$, with a common constant Lam\'e pair and
zero potentials sufficiently far away. The exact reduced
Green--Betti identity~\eqref{bt:betti} is applied on $\Omega$.
The boundary values of $\mu,W$ and the normal derivative of
$\mu$ cancel the additional boundary terms. After semiclassical
rescaling and passage to the limit on $\Omega$, the remaining
first-order difference depends only on the Lam\'e pair. Its
support allows the limiting integral to be taken over the
working domain.

On the extended domain, invertible transport frames are constructed
on prescribed bounded characteristic disks, with finite-order
differentiability in the local direction and transverse parameters.
The finite-regularity construction in Subsection~\ref{fr:interior-assembly}
gives a localized right inverse satisfying
\[
 h^{-1}\|\mathcal R^{(0)}_{h,\theta}F\|_{L^2(\Omega_0)}
 +\|\mathcal R^{(0)}_{h,\theta}F\|_{H^1(\Omega_0)}
 \le C\|F\|_{L^2(\Omega_0)}.
\]
After the first transport correction, the simultaneous remainder
bounds~\eqref{pb27:two} allow all three dual polarizations to be
realized by physical elastic solutions. The frames are selected
independently of the transverse and polynomial tests, and
Theorem~\ref{i2:thm-nine} gives the full matrix moment
identities~\eqref{i2:eq-nine}.

The moment identities yield a compactly supported Cauchy primitive
of the coefficient-difference matrix on each characteristic plane.
The formula~\eqref{if:intrinsic-formula} defines the relative
transport $H$, and the determinant identity~\eqref{if:detH}
proves its invertibility. Its normalization on the common exterior
region and holomorphic dependence on projective null directions
permit the block reduction leading to~\eqref{i4:closed-system}.
Uniqueness for the resulting homogeneous first-order system,
together with the exterior agreement, determines $\mu$ and then
$\lambda$.

With the Lam\'e parameters known, the physical integral identity
\eqref{pt:betti} and the product limit~\eqref{pt:limit} recover
$(W_1-W_2)/\mu$ through its Fourier transform. For the selected null phases, the cross terms involving the
longitudinal parts of the leading elastic fields tend to zero in
the product limit. A finite family of transverse test
functions removes the localization weight, yielding $W_1=W_2$
and hence $\rho_1=\rho_2$ because $\omega>0$ is known.

For the time-domain proof, fix a control time $\tau_{\rm c}$ with
$\max\{\tau_{S,1},\tau_{S,2}\}<\tau_{\rm c}<T_{\rm obs}/2$.
Lemma~\ref{td:dense-terminal} gives density of the terminal states
at $\tau_{\rm c}$ by local unique continuation. The response on
$(0,2\tau_{\rm c})$ determines the terminal Gram form
in~\eqref{td:terminal-gram}. In Proposition~\ref{td:response-continuation},
the terminal states form a core for each Dirichlet generator.
Equality of the Gram forms gives a unitary map between the state
spaces that sends corresponding terminal states to one another. The
core property and the generator action on terminal states identify
the full self-adjoint generators under this map as in
\eqref{td:generator-intertwining}. Their free evolutions then determine
the same tractions by the Green formula, which yields
\eqref{td:full-response}. The time bounds~\eqref{td:time-bound} justify
the full-time Laplace transform, which gives
\[
 \widetilde\Lambda_1^\Omega(-s^2)
   =\widetilde\Lambda_2^\Omega(-s^2),\qquad s>0,
\]
as in~\eqref{td:negative-parameter-equality}. Here
$\widetilde\Lambda_j^\Omega(z)$ is the stationary map at spectral
parameter $z$. Analytic continuation using the Dirichlet spectral
bounds~\eqref{td:spectral-gap} gives equality at
$z=\omega_*^2\in(0,\min\{\alpha_{D,1},\alpha_{D,2}\})$.
Theorem~\ref{a28:main} at $\omega_*$ gives the time-domain uniqueness.

\subsection{Conventions}\label{sec:notation}
All physical dot products, matrix transposes, and Green pairings are
complex bilinear. Hermitian products are used for Hilbert-space estimates.
The symbol $C$ in estimates denotes a constant
depending on fixed coefficients, domains, cutoffs, and chosen transport
frames; it may change from line to line. A frame denoted by $C$ is always
a matrix-valued function.

\subsection{Forward well-posedness and boundary maps}
\label{an:tools}

We establish well-posedness for both forward problems and give bounds
for their boundary maps. The elastic energy form provides a common
coercivity estimate for the time-harmonic Dirichlet problem and the
Dirichlet generator of the causal wave equation.

Let $V=H_0^1(\Omega;\mathbb C^3)$ and write $L=L_{\lambda,\mu}$.
With the convention of Subsection~\ref{sec:notation}, set
\[
 \mathfrak a(u,v)=\int_\Omega
 \left(2\mu\varepsilon(u):\overline{\varepsilon(v)}
 +\lambda\operatorname{div}u\,\overline{\operatorname{div}v}\right)dx,
 \qquad u,v\in H^1(\Omega;\mathbb C^3).
\]
Here $A:B=\sum_{i,j=1}^3A_{ij}B_{ij}$ is the bilinear matrix
contraction. For $v\in V$ the form represents $-L$:
$\mathfrak a(u,v)=\langle-Lu,v\rangle$, using the Hermitian duality.
Put $m_*=\min_{\overline\Omega}\min\{2\mu,3\lambda+2\mu\}>0$.
The deviatoric decomposition
$\operatorname{dev}A=A-\tfrac13(\operatorname{tr}A)I_3$ gives
\begin{equation}\label{an:elastic-coercivity}
\begin{aligned}
 \mathfrak a(v,v)
 &=\int_\Omega\left(2\mu|\operatorname{dev}\varepsilon(v)|^2
        +(\lambda+2\mu/3)|\operatorname{div}v|^2\right)dx\\
 &\ge m_*\|\varepsilon(v)\|_{L^2(\Omega)}^2
 \ge c_\Omega\|v\|_{H^1(\Omega)}^2,
 \qquad v\in V.
\end{aligned}
\end{equation}
The last inequality follows from Korn's and Poincar\'e's inequalities,
with $c_\Omega>0$ depending on $m_*$ and $\Omega$.
By Lax--Milgram, $-L:V\to H^{-1}(\Omega;\mathbb C^3)$ is a
bounded isomorphism.

We first consider the Dirichlet problem at the prescribed frequency
$\omega$.
\begin{proposition}\label{an:harmonic-wellposedness}
Assume \eqref{a28:physical-assumptions} and \eqref{en:nonresonance}
at the prescribed $\omega>0$. Then the physical weak Dirichlet
realization $P_D:V\to H^{-1}(\Omega;\mathbb C^3)$ is a bounded
isomorphism. For every $f\in H^{1/2}(\partial\Omega;\mathbb C^3)$,
problem~\eqref{en:dirichlet-problem} has a unique weak solution in
$H^1(\Omega;\mathbb C^3)$ and the map~\eqref{en:dn-definition} is
bounded from $H^{1/2}$ to $H^{-1/2}$.
\end{proposition}

\begin{proof}
For \eqref{en:dirichlet-problem}, write $P=L+WI_3$ and
$P_D=P_D^{\lambda,\mu,W}$ as in \eqref{en:nonresonance}.
Multiplication by $W\in L^\infty(\Omega)$ defines a compact operator from $H_0^1(\Omega;\mathbb C^3)$
to $H^{-1}(\Omega;\mathbb C^3)$, since it factors through the compact embedding into $L^2$. Thus,
\[
P_D=L(I+L^{-1}W)
\]
is Fredholm of index zero.
By the nonresonance condition~\eqref{en:nonresonance},
we have $\ker P_D=\{0\}$.
Consequently, $P_D$ is a bounded isomorphism. For boundary data $f\in H^{1/2}(\partial\Omega;\mathbb C^3)$,
choose a lifting $v\in H^1(\Omega;\mathbb C^3)$ satisfying
\[
\gamma_0 v=f,
\qquad
\|v\|_{H^1(\Omega)}
\le C\|f\|_{H^{1/2}(\partial\Omega)}.
\]
Set
\[
u=v-P_D^{-1}(Pv).
\]
Then $u$ is the unique weak solution of
\[
Pu=0\quad\text{in }\Omega,
\qquad
\gamma_0u=f\quad\text{on }\partial\Omega,
\]
and satisfies
\[
\|u\|_{H^1(\Omega)}
\le C\|f\|_{H^{1/2}(\partial\Omega)}.
\]
The independence of $u$ from the choice of lifting $v$
follows from uniqueness. For $g\in H^{1/2}(\partial\Omega;\mathbb C^3)$,
choose a lifting $v_g\in H^1(\Omega;\mathbb C^3)$ such that
\[
\gamma_0v_g=g,
\qquad
\|v_g\|_{H^1(\Omega)}
\le C\|g\|_{H^{1/2}(\partial\Omega)}.
\]
The map $\Lambda$ defined in~\eqref{en:dn-definition} admits the weak representation
\begin{equation}
\begin{aligned}
\langle\Lambda f,g\rangle
:=\int_\Omega\bigl[
&\lambda\operatorname{div}u\,
\operatorname{div}v_g\\
&+2\mu\varepsilon(u):\varepsilon(v_g)
-Wu\cdot v_g
\bigr]\,dx,
\end{aligned}
\label{en:dn-weak-form}
\end{equation}
where $u$ is the unique weak solution corresponding to
the boundary data $f$.

The weak formulation shows that this definition is
independent of the choice of the lifting $v_g$.
Moreover, the boundedness of the coefficients, together
with the stability estimate for $u$ and the lifting
estimate for $v_g$, yields
\[
|\langle\Lambda f,g\rangle|
\le C
\|f\|_{H^{1/2}(\partial\Omega)}
\|g\|_{H^{1/2}(\partial\Omega)}.
\]
Consequently,
\[
\Lambda:
H^{1/2}(\partial\Omega;\mathbb C^3)
\longrightarrow
H^{-1/2}(\partial\Omega;\mathbb C^3)
\]
is bounded.

\end{proof}

\begin{unnumberedremark}
The inverse $P_D^{-1}$ of Proposition~\ref{an:harmonic-wellposedness}
corrects approximate boundary layers to exact solutions with the same
trace. This gives the correction estimates~\eqref{bl:correct}
and~\eqref{bl:finite-exact-layer}, used in
Propositions~\ref{bl:layer} and~\ref{bl:finite-jets}.
The weak pairing~\eqref{en:dn-weak-form} gives the boundary integral
identity~\eqref{bl:alessandrini}.
\end{unnumberedremark}

For the time-domain problem, the density determines the Hilbert space
for the wave evolution. For the $j$th medium write
$L_j=L_{\lambda_j,\mu_j}$, set
$R_j^T=R_{\lambda_j,\mu_j,\rho_j}^{\Omega,T}$, and let
\[
 \mathcal H_j=L^2(\Omega;\rho_j\,dx;\mathbb C^3),\qquad
 (u,v)_{\mathcal H_j}=\int_\Omega\rho_j u\cdot\overline v\,dx .
\]
Write $\mathfrak a_j$ for the form $\mathfrak a$ with the $j$th
Lam\'e pair.

\begin{proposition}\label{an:dynamic-wellposedness}
Assume \eqref{a28:physical-assumptions} for the $j$th medium.
For $T>0$ and $f\in C_c^\infty((0,T);H^{1/2}(\partial\Omega;\mathbb C^3))$,
problem~\eqref{td:wave} on $(0,T)$ has a unique causal energy solution.
For any two such inputs $f,\widetilde f$, the weak response satisfies
\begin{equation}\label{td:response-bound}
 \|R_j^T(f-\widetilde f)\|_{L^\infty(0,T;H^{-1/2})}
 \le C_{T,j}\sum_{k=0}^3
       \|\partial_t^k(f-\widetilde f)\|_{L^\infty(0,T;H^{1/2})}.
\end{equation}
The boundary spaces in \eqref{td:response-bound} are those in
\eqref{td:response}.
\end{proposition}

\begin{proof}
Estimate~\eqref{an:elastic-coercivity} and the positive
upper and lower bounds on $\rho_j$ make $\mathfrak a_j$ a closed
positive form in $\mathcal H_j$ with domain $V$. Its self-adjoint
operator is
\begin{equation}\label{td:dirichlet-generator}
 \mathsf A_{D,j}=-\rho_j^{-1}L_j,\qquad
 D(\mathsf A_{D,j})=\{v\in V:\rho_j^{-1}L_jv\in\mathcal H_j\},
\end{equation}
where the domain condition is distributional. The identity
$(\mathsf A_{D,j}u,v)_{\mathcal H_j}=\mathfrak a_j(u,v)$ for
$u\in D(\mathsf A_{D,j})$ and $v\in V$, together with
\eqref{an:elastic-coercivity}, gives
\begin{equation}\label{td:elliptic-bound}
 \|v\|_{H^1(\Omega)}\le C_j\left(
 \|\mathsf A_{D,j}v\|_{\mathcal H_j}+\|v\|_{\mathcal H_j}\right),
 \qquad v\in D(\mathsf A_{D,j}).
\end{equation}
The compact embedding $V\hookrightarrow\mathcal H_j$ gives compact
resolvent. The same coercivity estimate gives
\begin{equation}\label{td:spectral-gap}
 \alpha_{D,j}=\inf_{0\ne v\in V}
 \frac{\mathfrak a_j(v,v)}{\|v\|_{\mathcal H_j}^2}>0,\qquad
 \operatorname{spec}(\mathsf A_{D,j})\subset[\alpha_{D,j},\infty).
\end{equation}

The trace theorem and \eqref{an:elastic-coercivity} give the bounded
static lift $G_j:H^{1/2}(\partial\Omega;\mathbb C^3)
\to H^1(\Omega;\mathbb C^3)$ characterized by
\[
 \gamma_0G_j g=g,\qquad \mathfrak a_j(G_j g,v)=0\quad(v\in V).
\]
For a smooth-in-time control $f$ vanishing near $t=0$, write
$u_j^f=G_j f+v_j$. Equation~\eqref{td:wave} reduces to
\begin{equation}\label{td:lifted-wave}
 \partial_t^2v_j+\mathsf A_{D,j}v_j=F_j,\qquad
 F_j=-G_j\partial_t^2f,\qquad v_j(0)=\partial_tv_j(0)=0.
\end{equation}
The spectral theorem gives the unique energy solution
\begin{equation}\label{td:causal-propagator}
 v_j(t)=\int_0^t
 \frac{\sin((t-s)\sqrt{\mathsf A_{D,j}})}{\sqrt{\mathsf A_{D,j}}}
 F_j(s)\,ds .
\end{equation}
In particular, setting
$E_j(v;t)=\tfrac12(\|\partial_tv(t)\|_{\mathcal H_j}^2+
\mathfrak a_j(v(t),v(t)))$, we have
\begin{equation}\label{td:energy-estimate}
 \frac{d}{dt}E_j(v_j;t)=\operatorname{Re}(F_j,\partial_tv_j)_{\mathcal H_j},
 \qquad
 E_j(v_j;t)^{1/2}\le\frac1{\sqrt2}\int_0^t\|F_j(s)\|_{\mathcal H_j}\,ds .
\end{equation}
The identity and bound follow first for finite spectral sums and
then by convergence in the energy norm. The formula
\eqref{td:causal-propagator} gives causality. Since $F_j$ vanishes
near $t=0$, differentiation gives the same formula for
$\partial_tv_j$ with forcing $\partial_tF_j$. Its energy estimate
controls $\partial_tv_j$ in $V$ and $\partial_t^2v_j$ in
$\mathcal H_j$. Thus $u_j^f$ solves \eqref{td:wave} with the
prescribed trace and initial data. For $f$ compactly supported in positive time,
$F_j,\partial_tF_j\in L^1((0,\infty);\mathcal H_j)$, so
\begin{equation}\label{td:time-bound}
 \sup_{t\ge0}
 \bigl(\|u_j^f(t)\|_{H^1(\Omega)}
       +\|\partial_tu_j^f(t)\|_{H^1(\Omega)}
       +\|\partial_t^2u_j^f(t)\|_{\mathcal H_j}\bigr)<\infty .
\end{equation}

On $(0,T)$ apply \eqref{td:energy-estimate} to the difference of the
zero-trace parts of the two solutions and to its time derivative.
Adding the static lift $G_j(f-\widetilde f)$ bounds the difference
of the physical solutions in $H^1(\Omega)$ and its second time
derivative in $\mathcal H_j$. The weak traction pairing
\eqref{td:weak-response} with a bounded trace lift then gives
\eqref{td:response-bound}.
\end{proof}

\begin{unnumberedremark}
The Dirichlet generator~\eqref{td:dirichlet-generator} and the
sine propagator in~\eqref{td:causal-propagator} enter the terminal
state argument of Lemma~\ref{td:dense-terminal}. The response
bound~\eqref{td:response-bound} extends the response equality from
separated controls to general controls in
Proposition~\ref{td:response-continuation}, and the uniform time
bound~\eqref{td:time-bound} justifies the full-time Laplace transform
in the proof of Theorem~\ref{td:uniqueness}.
\end{unnumberedremark}

\subsection{Organization of the paper}
Section~\ref{bl:section} determines the requisite finite boundary
derivatives and constructs exterior extensions.
Section~\ref{bt:section} derives the four-component reduced system
and its Green--Betti identity from the physical boundary map.
Sections~\ref{f27:section} and~\ref{r27:local-inverse} construct
the local transport frames and the right inverses used for the
complex geometrical optics solutions. Section~\ref{i2:sec}
constructs the first corrections and three physical polarizations,
fixes the quotient frames before the tests, and derives the full
matrix moment identities. Section~\ref{if:section} constructs the
normalized relative transport matrix, proves its holomorphic
dependence on direction, and determines the Lam\'e parameters by
a closed homogeneous system. Section~\ref{sec:potential} recovers
the density and completes the proof of Theorem~\ref{a28:main}.
Section~\ref{sec:dynamical} derives finite-time full-boundary uniqueness
for the dynamical system from the fixed-frequency theorem and the
response-continuation result.
Appendix~\ref{fr27:section} proves the free inverse and Fourier
operator estimates used in the solution construction.

\section{Boundary determination and a common exterior extension}
\label{bl:section}

This section first develops the boundary-layer formulas and
polarization matrices for smooth inputs. At each elastic order $j\ge2$,
the leading moment determines the $j$th normal derivatives of
$\lambda,\mu$ together with the $(j-2)$th normal derivative of $W$.
Proposition~\ref{bl:finite-jets} proves the finite-order determination
needed for Theorem~\ref{a28:main} on the physical domain with $C^9$
boundary. Proposition~\ref{ext:finite-common} then constructs the
exterior extensions for the interior argument.

In the smooth-input construction below, let
$\Omega\subset\mathbb R^3$ be a bounded domain with
$C^\infty$ boundary, possibly with several boundary components.
The coefficients $\lambda,\mu,W$ are real-valued and smooth
on $\overline\Omega$, and the Lam\'e parameters satisfy the
strong convexity conditions~\eqref{a28:physical-assumptions}.
For a fixed coefficient triple, write
\[
P:=L_{\lambda,\mu}+WI_3,
\qquad
P_D:=P_D^{\lambda,\mu,W}.
\]
When two coefficient triples are under consideration, we use
\[
P_j:=L_{\lambda_j,\mu_j}+W_jI_3,
\qquad
P_{D,j}:=P_D^{\lambda_j,\mu_j,W_j},
\qquad j=1,2.
\]
The nonresonance condition~\eqref{en:nonresonance} is imposed
on each physical coefficient triple. By the argument in
Subsection~\ref{an:tools}, the Dirichlet realization associated
with each such triple is a bounded isomorphism. In particular,
for a pair of triples, $P_{D,1}$ and $P_{D,2}$ are bounded
isomorphisms; for a construction involving one fixed triple,
$P_D$ denotes the corresponding realization.

Displacement components are expressed in a fixed Cartesian basis,
and bilinear pairings follow the convention of
Section~\ref{sec:notation}. The boundary-layer construction
and remainder estimates are formulated in the boundary normal
coordinates introduced below. Uniqueness for the frozen half-line
problem associated with~\eqref{bl:P0} follows from a Hermitian
energy identity. For each fixed finite expansion order, the
constants in the estimates for a single material may depend on
bounds for finitely many derivatives of the coefficients, the fixed geometry and cutoffs, and
$\|P_D^{-1}\|_{H^{-1}\to H_0^1}$. In estimates comparing two
materials, the dependence on the physical Dirichlet inverses is
through $\|P_{D,1}^{-1}\|_{H^{-1}\to H_0^1}$ and
$\|P_{D,2}^{-1}\|_{H^{-1}\to H_0^1}$. These constants are
independent of the scaling parameter $h$.

\subsection{Boundary coordinates and the scaled operator}

\begin{definition}
\label{bl:normal-coordinate-definition}

Let $\Phi_{\partial}:V\to\partial\Omega$ be a smooth boundary
parametrization, with $V\subset\mathbb R^2$ open and
$y=(y_1,y_2)$ the tangential coordinates.
Let $\nu$ denote the outward unit normal to $\partial\Omega$.
After restricting $V$, choose $r_0>0$ sufficiently small
so that the map
\begin{equation}\label{bl:coordinates}
 X(y,r)=\Phi_{\partial}(y)-r\nu(\Phi_{\partial}(y)),
 \qquad (y,r)\in V\times(-r_0,r_0),
\end{equation}
is a diffeomorphism onto an open neighborhood $U$ of
$\Phi_{\partial}(V)$. 

The inverse coordinates and the nearest boundary projection
associated with~\eqref{bl:coordinates} satisfy
\begin{equation}\label{bl:inverse-coordinates}
 \begin{aligned}
  X^{-1}(x)&=(y_1(x),y_2(x),r(x)),\\
  \pi_{\partial}(x)
  &=\Phi_{\partial}(y_1(x),y_2(x)),\\
  |x-\pi_{\partial}(x)|
  &=\operatorname{dist}(x,\partial\Omega)=|r(x)|,
  \qquad x\in U.
 \end{aligned}
\end{equation}
The normal coordinate is the signed distance
\begin{equation}\label{bl:signed-distance}
 r(x)=
 \begin{cases}
  \operatorname{dist}(x,\partial\Omega),
   & x\in U\cap\overline\Omega,\\
  -\operatorname{dist}(x,\partial\Omega),
   & x\in U\setminus\overline\Omega.
 \end{cases}
\end{equation}
The coordinates in~\eqref{bl:inverse-coordinates}, with the
sign convention~\eqref{bl:signed-distance}, are boundary
normal coordinates; see~\cite[Section~1]{PUZ2015}
and~\cite[Section~8]{DSFKSU2009}.
For the Euclidean metric, the normal geodesics are precisely
the straight normal segments in~\eqref{bl:coordinates}.

For any smooth physical coefficient $b=b(x)$, we use
the same symbol for its expression in the boundary
normal coordinates, defined by
\begin{equation}\label{bl:coordinate-fields}
\begin{aligned}
b(y,r)&:=b(X(y,r)),\\
n_{\mathrm{in}}(y,r)
&:=(\nabla_x r)(X(y,r)),\\
w^\alpha(y,r)
&:=(\nabla_x y_\alpha)(X(y,r)),
\qquad \alpha=1,2,\\
\mathcal X_i
&:=n_{\mathrm{in},i}(y,r)\partial_r
+\sum_{\alpha=1}^2
w_i^\alpha(y,r)\partial_{y_\alpha},
\qquad i=1,2,3.
\end{aligned}
\end{equation}
By the chain rule,
\begin{equation}
(\partial_{x_i}b)(X(y,r))
=\mathcal X_i b(y,r).
\end{equation}

The signed distance convention~\eqref{bl:signed-distance}
and the definitions in~\eqref{bl:coordinate-fields} imply
\begin{equation}\label{bl:normal-identities}
 \begin{aligned}
  n_{\mathrm{in}}(y,r)&=-\nu(\Phi_{\partial}(y)),\\
  |n_{\mathrm{in}}(y,r)|&=1,\qquad
  n_{\mathrm{in}}(y,r)\cdot w^\alpha(y,r)=0,
  \quad \alpha=1,2.
 \end{aligned}
\end{equation}
Thus, the vector $n_{\mathrm{in}}$ in~\eqref{bl:normal-identities}
is the inward unit normal to the parallel surfaces of constant
$r$, and the coordinate gradients $w^\alpha$ are tangent
to the parallel surfaces.

The Jacobian matrix of the coordinate map $X$ is
\[
D_{(y,r)}X(y,r)
=
\big(
\partial_{y_1}X(y,r),
\partial_{y_2}X(y,r),
\partial_r X(y,r)
\big).
\]

For a fixed compact patch $K\Subset V$ and $0<r_1<r_0$,
the volume density satisfies
\begin{equation}\label{bl:volume-jacobian}
 \begin{aligned}
  J(y,r)&:=|\det D_{(y,r)}X(y,r)|,
  \qquad dx=J(y,r)\,dy\,dr,\\
  0<c_J&\le J(y,r)\le C_J<\infty,
  \qquad (y,r)\in K\times[-r_1,r_1].
 \end{aligned}
\end{equation}
The constants $c_J$ and $C_J$ in~\eqref{bl:volume-jacobian}
depend on the fixed coordinate patch $K$ and the
normal coordinate interval $[-r_1,r_1]$.

Introduce the oscillation parameter, normal scale, stretched
normal coordinate, and phase sign by
\begin{equation}\label{bl:scaling}
 k\ge1,\qquad h=k^{-1},\qquad s_{\mathrm{bl}}=r/h,\qquad
 \mathfrak{s}\in\{+1,-1\}.
\end{equation}
The real parameter $k$ in~\eqref{bl:scaling} controls
tangential oscillation in the boundary data, $h$ is the scale
of normal localization, and $\mathfrak{s}$ determines the
sign of the tangential phase.
The measurement frequency $\omega$ remains fixed
as $k\to\infty$.
We refer to $(y,s_{\mathrm{bl}})=(y_1,y_2,r/h)$ as the
boundary-layer variables, where $y$ denotes the tangential
coordinates and $s_{\mathrm{bl}}$ is the stretched normal coordinate.

\end{definition}

\begin{lemma}
\label{bl:coordinate-calculus}
Let $U=U(y,s_{\mathrm{bl}})$ be a smooth vector-valued
function. For $X(y,hs_{\mathrm{bl}})$ in the coordinate neighborhood,
the conjugated Cartesian derivative satisfies
\begin{equation}\label{bl:D}
 \begin{aligned}
  &\left.
  e^{-i\mathfrak{s}y_1/h}\mathcal X_i
  \bigl(e^{i\mathfrak{s}y_1/h}U(y,r/h)\bigr)
  \right|_{r=hs_{\mathrm{bl}}}
 =h^{-1}(\mathcal D_{i,h}U)(y,s_{\mathrm{bl}}),\\
  &\mathcal D_{i,h}
  :=n_{\mathrm{in},i}(y,hs_{\mathrm{bl}})\partial_{s_{\mathrm{bl}}}
       +i\mathfrak{s}w_i^1(y,hs_{\mathrm{bl}})
       +h\sum_{\alpha=1}^2
          w_i^\alpha(y,hs_{\mathrm{bl}})\partial_{y_\alpha},
       \qquad i=1,2,3.
 \end{aligned}
\end{equation}

In the ordered product below, all coefficients and their derivatives
are evaluated at $(y,hs_{\mathrm{bl}})$. We have
\begin{align}
\mathcal D_{i,h}\mathcal D_{j,h}={}&n_{\mathrm{in},i}n_{\mathrm{in},j}\partial_{s_{\mathrm{bl}}}^2
 +h\sum_\alpha(n_{\mathrm{in},i}w_j^\alpha+w_i^\alpha n_{\mathrm{in},j})\partial_{s_{\mathrm{bl}}}\partial_{y_\alpha}
 +h^2\sum_{\alpha,\beta}w_i^\alpha w_j^\beta\partial_{y_\alpha}\partial_{y_\beta}
 \notag\\
 &+\big[i\mathfrak{s}(n_{\mathrm{in},i}w_j^1+w_i^1n_{\mathrm{in},j})+h\mathcal X_i(n_{\mathrm{in},j})\big]\partial_{s_{\mathrm{bl}}}
 \notag\\
 &+\sum_\beta\big[i\mathfrak{s} h(w_i^1w_j^\beta+w_i^\beta w_j^1)
                   +h^2\mathcal X_i(w_j^\beta)\big]\partial_{y_\beta}
 +\big[-w_i^1w_j^1+i\mathfrak{s} h\mathcal X_i(w_j^1)\big].
 \label{bl:ordered}
\end{align}

The derivatives of the coordinate fields are given by
\begin{equation}\label{bl:coordinate-derivatives}
\begin{aligned}
\mathcal X_i(n_{\mathrm{in},j})(y,r)
&=n_{\mathrm{in},i}(y,r)\partial_r n_{\mathrm{in},j}(y,r)
+\sum_{\beta=1}^2
w_i^\beta(y,r)\partial_{y_\beta}n_{\mathrm{in},j}(y,r)\\
&=(\partial_{x_i}\partial_{x_j}r)(X(y,r)),\\
\mathcal X_i(w_j^\alpha)(y,r)
&=n_{\mathrm{in},i}(y,r)\partial_r w_j^\alpha(y,r)
+\sum_{\beta=1}^2
w_i^\beta(y,r)\partial_{y_\beta}w_j^\alpha(y,r)\\
&=(\partial_{x_i}\partial_{x_j}y_\alpha)(X(y,r)).
\end{aligned}
\end{equation}

\end{lemma}

\begin{proof}
Apply $\mathcal X_i$ to the product
$e^{i\mathfrak{s}y_1/h}U(y,r/h)$. A derivative of the exponential
gives $i\mathfrak{s}h^{-1}w_i^1$, and differentiating $r/h$
gives $h^{-1}n_{\mathrm{in},i}$, which proves~\eqref{bl:D}.
In the composition, $\mathcal D_{i,h}$ acts on both the function
and the coefficients in $\mathcal D_{j,h}$. For a coefficient $b$
evaluated at $(y,hs_{\mathrm{bl}})$,
\[
 \left(n_{\mathrm{in},i}\partial_{s_{\mathrm{bl}}}
       +h\sum_\alpha w_i^\alpha\partial_{y_\alpha}\right)
 b(y,hs_{\mathrm{bl}})
 =h(\mathcal X_i b)(y,hs_{\mathrm{bl}}).
\]
The product rule and this identity give~\eqref{bl:ordered}.
The chain rule for the Cartesian coordinate gradients
$n_{\mathrm{in}}$ and $w^\alpha$ gives~\eqref{bl:coordinate-derivatives}.
\end{proof}

\begin{definition}
\label{bl:scaled-operator-definition}

To express the physical operator in the boundary-layer
variables, we define the scaled operator
$\mathcal P_{\mathfrak{s},h}$ by conjugating $P$ with
the factor $e^{i\mathfrak{s}y_1/h}$ and rescaling the
normal variable. For a smooth vector-valued function
$U=U(y,s_{\mathrm{bl}})$, define
\begin{equation*}
\begin{aligned}
(\mathcal P_{\mathfrak{s},h}U)(y,s_{\mathrm{bl}})
:=h^2e^{-i\mathfrak{s}y_1/h}
\left.
P\left(
e^{i\mathfrak{s}y_1(x)/h}
U\left(y(x),\frac{r(x)}{h}\right)
\right)
\right|_{x=X(y,hs_{\mathrm{bl}})}.
\end{aligned}
\end{equation*}
Here $y(x)$ and $r(x)$ are the boundary coordinate
functions defined in~\eqref{bl:inverse-coordinates}.
Thus $P$ acts on a function of the physical variable
$x$, whereas $\mathcal P_{\mathfrak{s},h}$ acts on the
vector-valued function $U=U(y,s_{\mathrm{bl}})$ in the boundary-layer
variables $(y,s_{\mathrm{bl}})$. By definition,
\begin{equation*}
\begin{aligned}
\left.
P\left(
e^{i\mathfrak{s}y_1(x)/h}
U\left(y(x),\frac{r(x)}{h}\right)
\right)\right|_{x=X(y,hs_{\mathrm{bl}})}
=h^{-2}e^{i\mathfrak{s}y_1/h}
(\mathcal P_{\mathfrak{s},h}U)(y,s_{\mathrm{bl}}).
\end{aligned}
\end{equation*}

\end{definition}

\begin{proposition}
\label{bl:scaled-expansion}
For $U$ as in the preceding definition, the scaled operator satisfies
\begin{align}
 (\mathcal P_{\mathfrak{s},h}U)_i={}&
 \mu\sum_j\mathcal D_{j,h}\mathcal D_{j,h}U_i
 +(\lambda+\mu)\sum_j
       \mathcal D_{i,h}\mathcal D_{j,h}U_j
 \notag\\
 &+h(\mathcal X_i(\lambda))
       \sum_j\mathcal D_{j,h}U_j
 +h\sum_j(\mathcal X_j(\mu))
       (\mathcal D_{j,h}U_i+\mathcal D_{i,h}U_j)
 +h^2WU_i.
 \label{bl:scaled}
\end{align}
All coefficients and their derivatives in
this formula are evaluated at $(y,hs_{\mathrm{bl}})$.

It has the representation
\begin{equation}\label{bl:A}
 \mathcal P_{\mathfrak{s},h}
 =\sum_{m=0}^2\sum_{|\beta|\le2}
 h^m A_{m,\beta}(y,hs_{\mathrm{bl}})
 \partial_{(y,s_{\mathrm{bl}})}^\beta.
\end{equation}
Here $A_{m,\beta}(y,hs_{\mathrm{bl}})$ are $3\times3$
matrix-valued coefficient functions, and
\[
\beta=(\beta_1,\beta_2,\beta_3)\in\mathbb N_0^3,
\qquad
|\beta|=\beta_1+\beta_2+\beta_3,
\]
where $\mathbb N_0=\{0,1,2,\ldots\}$.
We use the notation
\[
\partial_{(y,s_{\mathrm{bl}})}^\beta
=\partial_{y_1}^{\beta_1}
 \partial_{y_2}^{\beta_2}\partial_{s_{\mathrm{bl}}}^{\beta_3}.
\]
In~\eqref{bl:A}, each derivative acts only on $U$;
all terms arising when derivatives act on the
variable coefficients in~\eqref{bl:scaled} are
included in the coefficient matrices
$A_{m,\beta}$. Any coefficient that does not occur
in the expansion is understood to be zero.
In particular, whenever $A_{m,\beta}\ne0$, we have
\begin{equation}\label{bl:tangential-order}
\beta_1+\beta_2\le m.
\end{equation}

The matrices $A_{m,\beta}(y,r)$ are smooth
$3\times3$ matrix-valued functions determined by
$\lambda,\mu,W,n_{\mathrm{in}},w$, their first
Cartesian derivatives, and the fixed phase sign
$\mathfrak{s}$.

For each $q\ge0$, the coefficient of $h^q$ obtained by expanding
the coefficient functions at $r=0$ is
\begin{equation}\label{bl:Pq}
 \mathcal P_{\mathfrak{s}}^{(q)}
 =\sum_{m=0}^{\min(2,q)}\sum_{|\beta|\le2}
   \frac{s_{\mathrm{bl}}^{q-m}}{(q-m)!}
   \partial_r^{q-m}A_{m,\beta}(y,0)
   \partial_{(y,s_{\mathrm{bl}})}^\beta.
\end{equation}

In particular, $\mathcal P_{\mathfrak{s}}^{(0)}$ is the frozen
normal operator, and the explicit term involving $W$ in~\eqref{bl:scaled}
first contributes at order $h^2$.
\end{proposition}

\begin{proof}

In Cartesian components, the physical operator is
\[
 (Pu)_i=\mu\sum_j\partial_{x_j}^2u_i
 +(\lambda+\mu)\sum_j\partial_{x_i}\partial_{x_j}u_j
 +(\mathcal X_i(\lambda))\sum_j\partial_{x_j}u_j
 +\sum_j(\mathcal X_j(\mu))(\partial_{x_j}u_i+\partial_{x_i}u_j)+Wu_i.
\]

Using~\eqref{bl:D} in this Cartesian formula gives~\eqref{bl:scaled}.
Substituting~\eqref{bl:ordered} and collecting the terms of the same
power of $h$ and the same derivative of $U$ gives~\eqref{bl:A}.
Each tangential derivative of $U$ carries a factor $h$,
which gives~\eqref{bl:tangential-order}. To find the coefficient of
$h^q$, expand each $A_{m,\beta}(y,hs_{\mathrm{bl}})$ with
$m\le\min(2,q)$ to order $q-m$ in $r$. The term of Taylor order
$q-m$ gives~\eqref{bl:Pq}.
\end{proof}

The operators $\mathcal P_{\mathfrak{s}}^{(q)}$ with $q\ge1$
determine the successive profiles through~\eqref{bl:recursion}.

\begin{lemma}
\label{bl:taylor-remainder-lemma}

For a smooth scalar- or matrix-valued coefficient
$b(y,r)$ in the boundary normal coordinates and an
integer $N\ge0$, let
\begin{equation}\label{bl:taylor}
 \begin{aligned}
 b(y,hs_{\mathrm{bl}})&=\sum_{q=0}^N\frac{h^q s_{\mathrm{bl}}^q}{q!}\partial_r^qb(y,0)
                   +h^{N+1}R_{N,h}[b](y,s_{\mathrm{bl}}),\\
 R_{N,h}[b](y,s_{\mathrm{bl}})&=
 \frac{s_{\mathrm{bl}}^{N+1}}{N!}\int_0^1(1-s)^N
                          \partial_r^{N+1}b(y,shs_{\mathrm{bl}})\,ds.
 \end{aligned}
\end{equation}

For a fixed finite derivative order $S$ and a fixed compact
coordinate patch, if the required $y$- and $r$-derivatives of $b$
are bounded on the collar, then for $|\gamma|\le S$,
$0<h\le1$, and $0\le hs_{\mathrm{bl}}\le r_0$,
\begin{equation}\label{bl:taylor-bound}
\left|\partial_{(y,s_{\mathrm{bl}})}^\gamma R_{N,h}[b](y,s_{\mathrm{bl}})\right|
\le C_{N,S}(1+s_{\mathrm{bl}})^{N+1}.
\end{equation}

\end{lemma}

\begin{proof}
Formula~\eqref{bl:taylor} is Taylor's formula with integral
remainder in the normal coordinate. Here
$\gamma=(\gamma_1,\gamma_2,\gamma_3)$ and
$\alpha=(\gamma_1,\gamma_2)$.

Differentiation
with respect to $y$ gives
\[
\partial_y^\alpha R_{N,h}[b]
=R_{N,h}[\partial_y^\alpha b].
\]
For $0\le\gamma_3\le N$, differentiation with
respect to $s_{\mathrm{bl}}$ gives
\[
\partial_y^\alpha\partial_{s_{\mathrm{bl}}}^{\gamma_3}R_{N,h}[b]
=R_{N-\gamma_3,h}
 [\partial_y^\alpha\partial_r^{\gamma_3}b].
\]
Using~\eqref{bl:taylor}, we obtain
\[
\begin{aligned}
&\left|\partial_y^\alpha\partial_{s_{\mathrm{bl}}}^{\gamma_3}
R_{N,h}[b](y,s_{\mathrm{bl}})\right|\le C_{N,S}s_{\mathrm{bl}}^{N-\gamma_3+1}
\le C_{N,S}(1+s_{\mathrm{bl}})^{N+1}.
\end{aligned}
\]
For $\gamma_3\ge N+1$, we have
\[
\partial_y^\alpha\partial_{s_{\mathrm{bl}}}^{\gamma_3}R_{N,h}[b]
=h^{\gamma_3-N-1}
\partial_y^\alpha\partial_r^{\gamma_3}b(y,hs_{\mathrm{bl}}),
\]
and hence
\[
\left|\partial_y^\alpha\partial_{s_{\mathrm{bl}}}^{\gamma_3}
R_{N,h}[b](y,s_{\mathrm{bl}})\right|
\le C_{N,S}h^{\gamma_3-N-1}
\le C_{N,S}(1+s_{\mathrm{bl}})^{N+1}.
\]

Combining the two cases proves~\eqref{bl:taylor-bound}.
\end{proof}

\begin{proposition}
\label{bl:finite-expansion}
Fix a truncation order $N_{\mathrm{bl}}\ge0$ and define the
coefficients and the remainder operator by
\begin{equation}\label{bl:remainder-coefficients}
\begin{aligned}
 B^{(N_{\mathrm{bl}})}_{m,\beta,h}(y,s_{\mathrm{bl}})
 &=
 \begin{cases}
 R_{N_{\mathrm{bl}}-m,h}[A_{m,\beta}](y,s_{\mathrm{bl}}),
 &0\le m\le\min(2,N_{\mathrm{bl}}),\\
 h^{m-N_{\mathrm{bl}}-1}A_{m,\beta}(y,hs_{\mathrm{bl}}),
 &N_{\mathrm{bl}}<m\le2,
 \end{cases}\\
 \mathfrak R^{\mathrm{bl}}_{\mathfrak{s},N_{\mathrm{bl}},h}
 &=
 \sum_{m=0}^2\sum_{|\beta|\le2}
 B^{(N_{\mathrm{bl}})}_{m,\beta,h}(y,s_{\mathrm{bl}})
 \partial_{(y,s_{\mathrm{bl}})}^\beta.
\end{aligned}
\end{equation}

The operator $\mathfrak R^{\mathrm{bl}}_{\mathfrak{s},N_{\mathrm{bl}},h}$
has differential order at most two and satisfies the exact identity
\begin{equation}\label{bl:remainder}
 \mathcal P_{\mathfrak{s},h}
 =\sum_{q=0}^{N_{\mathrm{bl}}}
 h^q\mathcal P_{\mathfrak{s}}^{(q)}
 +h^{N_{\mathrm{bl}}+1}
 \mathfrak R^{\mathrm{bl}}_{\mathfrak{s},N_{\mathrm{bl}},h}.
\end{equation}

On each fixed compact coordinate patch, the remainder
coefficients and their derivatives of any fixed finite order have at
most polynomial growth in $s_{\mathrm{bl}}$, uniformly for
$0<h\le1$ and $0\le hs_{\mathrm{bl}}\le r_0$.
\end{proposition}

\begin{proof}

For a truncation order $N_{\mathrm{bl}}\ge0$, we
substitute~\eqref{bl:taylor} into~\eqref{bl:A}.
For $0\le m\le\min(2,N_{\mathrm{bl}})$, this gives
\[
\begin{aligned}
h^m A_{m,\beta}(y,hs_{\mathrm{bl}})
={}&\sum_{p=0}^{N_{\mathrm{bl}}-m}
\frac{h^{m+p}s_{\mathrm{bl}}^p}{p!}
\partial_r^p A_{m,\beta}(y,0)+h^{N_{\mathrm{bl}}+1}
R_{N_{\mathrm{bl}}-m,h}[A_{m,\beta}](y,s_{\mathrm{bl}}).
\end{aligned}
\]
For $N_{\mathrm{bl}}<m\le2$, the term already
contains a power of $h$ greater than
$N_{\mathrm{bl}}$, and we simply write
\[
h^m A_{m,\beta}(y,hs_{\mathrm{bl}})
=h^{N_{\mathrm{bl}}+1}
h^{m-N_{\mathrm{bl}}-1}A_{m,\beta}(y,hs_{\mathrm{bl}}).
\]
Thus, after collecting the terms of powers
$h^0,\ldots,h^{N_{\mathrm{bl}}}$, we factor out
$h^{N_{\mathrm{bl}}+1}$ from all remaining terms.

The remaining coefficients are precisely those in
~\eqref{bl:remainder-coefficients}. Since $|\beta|\le2$, the
remainder has differential order at most two. This proves~\eqref{bl:remainder}.

For $m\le N_{\mathrm{bl}}$, estimate~\eqref{bl:taylor-bound}
gives polynomial bounds in $s_{\mathrm{bl}}$ for
$B^{(N_{\mathrm{bl}})}_{m,\beta,h}$ and all its
derivatives of any fixed finite order. For
$m>N_{\mathrm{bl}}$, the same bounds follow from
the smoothness of $A_{m,\beta}$ and $0<h\le1$.
Consequently, the coefficients of
$\mathfrak R^{\mathrm{bl}}_{\mathfrak{s},N_{\mathrm{bl}},h}$
and all their derivatives of any fixed finite order
have at most polynomial growth in $s_{\mathrm{bl}}$, uniformly
for $0<h\le1$ and $0\le hs_{\mathrm{bl}}\le r_0$.

\end{proof}

Since the profiles constructed below are polynomials
in $s_{\mathrm{bl}}$ multiplied by $e^{-c_{\mathrm{bl}}(y)s_{\mathrm{bl}}}$,
with $c_{\mathrm{bl}}(y)\ge c_0>0$, for any
$0<c<c_0$ and any fixed integer $M\ge0$,
\[
(1+s_{\mathrm{bl}})^M e^{-c_{\mathrm{bl}}(y)s_{\mathrm{bl}}}
\le C_{M,c}e^{-cs_{\mathrm{bl}}}.
\]
The same type of bound holds for all derivatives
of these profiles. This exponential decay, together
with the polynomial bounds for the remainder
coefficients, is used to establish the residual
bound~\eqref{bl:residual}.

\subsection{The frozen half-line operator}

In this subsection, $n_{\mathrm{in}}$ is evaluated
at $r=0$. Set
\[
\tau_{\mathrm{tan}}=w^1(y,0),\qquad
c_{\mathrm{bl}}=|\tau_{\mathrm{tan}}|,\qquad
e_{\mathrm{bl}}=\frac{\tau_{\mathrm{tan}}}{c_{\mathrm{bl}}}.
\]
Thus, $e_{\mathrm{bl}}$ is the unit tangent vector
in the direction of $\tau_{\mathrm{tan}}$.
By~\eqref{bl:normal-identities},
\[
n_{\mathrm{in}}\cdot\tau_{\mathrm{tan}}=0,
\qquad |n_{\mathrm{in}}|=1,
\]
and $c_{\mathrm{bl}}\ge c_0>0$ on the compact patch.
Writing $\lambda_b=\lambda(y,0)$ and
$\mu_b=\mu(y,0)$, the leading operator is
\begin{equation}\label{bl:P0}
 \mathcal P_{\mathfrak{s}}^{(0)}
 =\mu_b(\partial_{s_{\mathrm{bl}}}^2-c_{\mathrm{bl}}^2)I_3
 +(\lambda_b+\mu_b)(i\mathfrak{s}\tau_{\mathrm{tan}}+n_{\mathrm{in}}\partial_{s_{\mathrm{bl}}})
                       (i\mathfrak{s}\tau_{\mathrm{tan}}+n_{\mathrm{in}}\partial_{s_{\mathrm{bl}}})^T.
\end{equation}

The matrix $F_{\mathfrak{s}}$ is given by
\begin{equation}\label{bl:frozen-profile}
 \begin{aligned}
 a_{\mathrm{bl}}&=\frac{\lambda_b+\mu_b}{\lambda_b+3\mu_b}\in(1/7,1),\\
 N_{\mathfrak{s}}&=a_{\mathrm{bl}}\big[e_{\mathrm{bl}}\otimes e_{\mathrm{bl}}-n_{\mathrm{in}}\otimes n_{\mathrm{in}}
             +i\mathfrak{s}(e_{\mathrm{bl}}\otimes n_{\mathrm{in}}+n_{\mathrm{in}}\otimes e_{\mathrm{bl}})\big],\\
 F_{\mathfrak{s}}(y,s_{\mathrm{bl}})&=e^{-c_{\mathrm{bl}}s_{\mathrm{bl}}}(I_3-c_{\mathrm{bl}}s_{\mathrm{bl}} N_{\mathfrak{s}}).
 \end{aligned}
\end{equation}
Here $(a\otimes b)_{ij}=a_ib_j$ denotes the tensor product of vectors.

\begin{lemma}
\label{bl:frozen}
Let $d_p\in\mathbb N_0$ and
\[
 p(y,s_{\mathrm{bl}})=\sum_{j=0}^{d_p}p_j(y)s_{\mathrm{bl}}^j,
 \qquad p_j\in C^\infty(K;\mathbb C^3),
\]
where $K\Subset V$ is the fixed compact parameter patch.
For each $y\in K$, the problem
\[
 \mathcal P_{\mathfrak{s}}^{(0)}V(y,s_{\mathrm{bl}})
 =e^{-c_{\mathrm{bl}}(y)s_{\mathrm{bl}}}p(y,s_{\mathrm{bl}}),\qquad s_{\mathrm{bl}}>0,
 \qquad V(y,0)=0,
\]
has a unique solution decaying as $s_{\mathrm{bl}}\to\infty$.
Moreover,
\[
 V(y,s_{\mathrm{bl}})=e^{-c_{\mathrm{bl}}(y)s_{\mathrm{bl}}}q(y,s_{\mathrm{bl}}),
 \qquad \deg_{s_{\mathrm{bl}}} q(y,\cdot)\le d_p+2,
\]
where $q$ is a vector polynomial in $s_{\mathrm{bl}}$ with smooth
coefficients in $y$. For every multi-index
$\alpha\in\mathbb N_0^2$ and integer $j\ge0$, there
exist constants $C_{\alpha,j}>0$ and
$M_{\alpha,j}\in\mathbb N_0$ such that
\[
 |\partial_y^\alpha\partial_{s_{\mathrm{bl}}}^j V(y,s_{\mathrm{bl}})|
 \le C_{\alpha,j}(1+s_{\mathrm{bl}})^{M_{\alpha,j}}e^{-c_0s_{\mathrm{bl}}},
 \qquad y\in K,\quad s_{\mathrm{bl}}\ge0.
\]
The constants depend on finitely many seminorms of
the data and the frozen coefficients on $K$.
\end{lemma}
\begin{proof}
Regard $\partial_{s_{\mathrm{bl}}}$ as a commuting indeterminate $D$, and put
$d=D^2-c_{\mathrm{bl}}^2$ and $z=i\mathfrak{s}\tau_{\mathrm{tan}}+n_{\mathrm{in}}D$. The identity $z^Tz=d$ gives
\[
 \det(\mu_bdI_3+(\lambda_b+\mu_b)zz^T)
 =\mu_b^2(\lambda_b+2\mu_b)d^3.
\]
This polynomial identity holds also at $d=0$. To check the homogeneous
profile, set $v_{\mathfrak{s}}=i\mathfrak{s} e_{\mathrm{bl}}-n_{\mathrm{in}}$. Then
$v_{\mathfrak{s}}^Tv_{\mathfrak{s}}=0$, $N_{\mathfrak{s}}=-a_{\mathrm{bl}}v_{\mathfrak{s}} v_{\mathfrak{s}}^T$ and
$n_{\mathrm{in}}^TN_{\mathfrak{s}}=a_{\mathrm{bl}}v_{\mathfrak{s}}^T$, so $N_{\mathfrak{s}}^2=0$. Direct differentiation gives
\[
 \begin{aligned}
 (D^2-c_{\mathrm{bl}}^2)F_{\mathfrak{s}}&=2c_{\mathrm{bl}}^2e^{-c_{\mathrm{bl}}s_{\mathrm{bl}}}N_{\mathfrak{s}},\\
 z(D)z(D)^TF_{\mathfrak{s}}&=-c_{\mathrm{bl}}^2(1-a_{\mathrm{bl}})a_{\mathrm{bl}}^{-1}e^{-c_{\mathrm{bl}}s_{\mathrm{bl}}}N_{\mathfrak{s}}.
 \end{aligned}
\]
The identity
$2\mu_b a_{\mathrm{bl}}
=(\lambda_b+\mu_b)(1-a_{\mathrm{bl}})$
gives $\mathcal P_{\mathfrak{s}}^{(0)}F_{\mathfrak{s}}=0$.
Moreover, $F_{\mathfrak{s}}(y,0)=I_3$.

The polynomial matrix
\begin{equation}\label{bl:quadratic-inverse}
 \begin{aligned}
 \mathcal B_{\mathfrak{s}}(D)&=(\lambda_b+2\mu_b)(D^2-c_{\mathrm{bl}}^2)I_3
                  -(\lambda_b+\mu_b)zz^T,\\
 \mathcal P_{\mathfrak{s}}^{(0)}\mathcal B_{\mathfrak{s}}(D)
     &=\mu_b(\lambda_b+2\mu_b)(D^2-c_{\mathrm{bl}}^2)^2I_3
 \end{aligned}
\end{equation}
reduces the construction of a right inverse to two scalar integrations. Indeed,
$(zz^T)^2=dzz^T$, so the coefficient of $dzz^T$ in the product is
\[
 (\lambda_b+\mu_b)
 [-\mu_b+(\lambda_b+2\mu_b)-(\lambda_b+\mu_b)]=0.
\]
Thus \eqref{bl:quadratic-inverse} holds in the ring of polynomial
matrices, including at characteristic roots.

Set $C_b=\mu_b(\lambda_b+2\mu_b)>0$, and take two componentwise
antiderivatives with zero integration constants to obtain $Q_*$ with
$\partial_{s_{\mathrm{bl}}}^2Q_*=p/C_b$. On polynomials of degree at most $d_p+2$ set
\begin{equation}\label{bl:polynomial-inverse}
 J_{c_{\mathrm{bl}}}=-\frac1{2c_{\mathrm{bl}}}\sum_{q=0}^{d_p+2}
                 \left(\frac{\partial_{s_{\mathrm{bl}}}}{2c_{\mathrm{bl}}}\right)^q,
 \qquad Q=J_{c_{\mathrm{bl}}}^2Q_*,\qquad
 V_{\rm part}=\mathcal B_{\mathfrak{s}}(\partial_{s_{\mathrm{bl}}})(e^{-c_{\mathrm{bl}}s_{\mathrm{bl}}}Q).
\end{equation}
Differentiation is nilpotent on this space, and therefore
$(\partial_{s_{\mathrm{bl}}}-2c_{\mathrm{bl}})J_{c_{\mathrm{bl}}}=I$. Its commutation with $\partial_{s_{\mathrm{bl}}}$ gives
$C_b\partial_{s_{\mathrm{bl}}}^2(\partial_{s_{\mathrm{bl}}}-2c_{\mathrm{bl}})^2Q=p$. Conjugating
$(\partial_{s_{\mathrm{bl}}}^2-c_{\mathrm{bl}}^2)^2$ by $e^{-c_{\mathrm{bl}}s_{\mathrm{bl}}}$ in
\eqref{bl:quadratic-inverse} now yields
$\mathcal P_{\mathfrak{s}}^{(0)}V_{\rm part}=e^{-c_{\mathrm{bl}}s_{\mathrm{bl}}}p$. Hence
\[
 V=V_{\rm part}-F_{\mathfrak{s}} V_{\rm part}(0)
\]
is the required zero-trace solution. The polynomial degree of
$e^{c_{\mathrm{bl}}s_{\mathrm{bl}}}V_{\rm part}$ is at most $d_p+2$, and that of the subtracted
term is at most one. All operations in~\eqref{bl:polynomial-inverse} are finite,
with denominators involving only $c_{\mathrm{bl}}$ and
$C_b=\mu_b(\lambda_b+2\mu_b)$, which are bounded away
from zero on the compact parameter patch.
The coefficients of the resulting polynomial are smooth
functions of $y$. Differentiation in $y$ preserves this
structure, except that derivatives falling on
$e^{-c_{\mathrm{bl}}(y)s_{\mathrm{bl}}}$ introduce additional polynomial
factors in $s_{\mathrm{bl}}$. Consequently, for every multi-index
$\alpha\in\mathbb N_0^2$ and integer $j\ge0$, there exist
constants $C_{\alpha,j}>0$ and $M_{\alpha,j}\in\mathbb N_0$
such that
\[
 |\partial_y^\alpha\partial_{s_{\mathrm{bl}}}^j V(y,s_{\mathrm{bl}})|
 \le C_{\alpha,j}(1+s_{\mathrm{bl}})^{M_{\alpha,j}}e^{-c_0s_{\mathrm{bl}}},
 \qquad y\in K,\quad s_{\mathrm{bl}}\ge0.
\]
The estimates remain uniform when the degree of the
forcing polynomial drops.

For uniqueness, the adjugate identity applied to the homogeneous
system shows that every component satisfies $(D^2-c_{\mathrm{bl}}^2)^3V_i=0$.
A solution tending to zero at infinity is consequently a polynomial
times $e^{-c_{\mathrm{bl}}s_{\mathrm{bl}}}$ and has square-integrable derivatives. If it has zero
trace, integration by parts is justified and gives
\[
 -\int_0^\infty(\mathcal P_{\mathfrak{s}}^{(0)}V)\cdot\overline V\,d s_{\mathrm{bl}}
 =\mu_b\int_0^\infty(|V'|^2+c_{\mathrm{bl}}^2|V|^2)\,d s_{\mathrm{bl}}
 +(\lambda_b+\mu_b)\int_0^\infty
             |i\mathfrak{s}\tau_{\mathrm{tan}}\cdot V+n_{\mathrm{in}}\cdot V'|^2\,d s_{\mathrm{bl}}.
\]
Strong convexity implies $\lambda_b+\mu_b>0$, so $V=0$.
The same argument applies to homogeneous $H^1$ solutions with zero trace.
\end{proof}

\subsection{Finite boundary layers and exact solutions}

\begin{proposition}\label{bl:layer}
Fix an integer $N_{\mathrm{bl}}\ge0$, an amplitude
$A\in C_c^\infty(V;\mathbb C^3)$, and a phase sign
$\mathfrak{s}\in\{+1,-1\}$. Let $\eta(r)$ be a fixed smooth
normal cutoff, compactly supported in the coordinate collar
and equal to one near $r=0$. Then, for all sufficiently small
$h>0$, there exist profiles
$U_0,\ldots,U_{N_{\mathrm{bl}}}$, independent of $h$, such
that each $U_\ell(y,s_{\mathrm{bl}})$ is a polynomial in $s_{\mathrm{bl}}$ multiplied by
$e^{-c_{\mathrm{bl}}(y)s_{\mathrm{bl}}}$ and is supported in a fixed compact
set $K\Subset V$ containing $\operatorname{supp}A$. The field
\[
 u_{h,N_{\mathrm{bl}}}^{\mathrm{app}}(X(y,r))
 =e^{i\mathfrak{s}y_1/h}\eta(r)
   \sum_{\ell=0}^{N_{\mathrm{bl}}}h^\ell
   U_\ell(y,r/h)
\]
extends by zero to a smooth field on $\Omega$ and has boundary
trace $e^{i\mathfrak{s}y_1/h}A(y)$.

There exists an exact weak solution
$u_h\in H^1(\Omega;\mathbb C^3)$ of $Pu_h=0$ with the same
boundary trace, and
\begin{align}
 \|Pu_{h,N_{\mathrm{bl}}}^{\mathrm{app}}\|_{L^2(\Omega)}
 &\le C_{N_{\mathrm{bl}}}
       h^{N_{\mathrm{bl}}-1/2},
 \label{bl:residual}\\
 \|u_{h,N_{\mathrm{bl}}}^{\mathrm{app}}\|_{L^2(\Omega)}
 &\le C_{N_{\mathrm{bl}}}h^{1/2},\qquad
 \|u_{h,N_{\mathrm{bl}}}^{\mathrm{app}}\|_{H^1(\Omega)}
 \le C_{N_{\mathrm{bl}}}h^{-1/2},
 \label{bl:appnorm}\\
 \|u_h-u_{h,N_{\mathrm{bl}}}^{\mathrm{app}}\|_{H^1(\Omega)}
 &\le C_{N_{\mathrm{bl}}}
       h^{N_{\mathrm{bl}}-1/2}.
 \label{bl:correct}
\end{align}
The constants are independent of $h$.

More precisely, for each fixed finite integer $S\ge0$, the
profiles and all derivatives required in the construction
are controlled by the $C^{2N_{\mathrm{bl}}+S+4}$ norm of
$A$, the coefficient bounds
\[
 \max_{\substack{|\alpha|\le2N_{\mathrm{bl}}+S+4\\
                  0\le b\le N_{\mathrm{bl}}+S+2}}
 \left\|
 \partial_y^\alpha\partial_r^b
 (\lambda,\mu,W,n_{\mathrm{in}},w)
 \right\|_{L^\infty}
\]
on the fixed compact collar, finitely many cutoff norms,
positive lower bounds for
$c_{\mathrm{bl}},\mu,\lambda+2\mu,\lambda+3\mu$,
and upper and positive lower bounds for the coordinate
Jacobian. The bound for the difference
$u_h-u_{h,N_{\mathrm{bl}}}^{\mathrm{app}}$ also depends on
$\|P_D^{-1}\|_{H^{-1}\to H^1_0}$.
\end{proposition}

\begin{proof}
Using the profile in~\eqref{bl:frozen-profile}, set
$U_0=F_{\mathfrak{s}}A$. By the defining properties of the
frozen profile,
\[
 \mathcal P_{\mathfrak{s}}^{(0)}U_0=0,
 \qquad U_0(y,0)=A(y).
\]
For the finite sum
$\sum_{\ell=0}^{N_{\mathrm{bl}}}h^\ell U_\ell$, the expansion
\eqref{bl:remainder} shows that the coefficient of $h^\ell$
is
$\sum_{q=0}^{\ell}\mathcal P_{\mathfrak{s}}^{(q)}U_{\ell-q}$
for $\ell\le N_{\mathrm{bl}}$. We cancel these coefficients
successively by imposing
\begin{equation}\label{bl:recursion}
 \mathcal P_{\mathfrak{s}}^{(0)}U_\ell
 =-\sum_{q=1}^{\ell}
       \mathcal P_{\mathfrak{s}}^{(q)}U_{\ell-q},
 \qquad U_\ell(y,0)=0,\quad
 1\le\ell\le N_{\mathrm{bl}}.
\end{equation}
Every right-hand side is a polynomial in $s_{\mathrm{bl}}$ multiplied by
$e^{-c_{\mathrm{bl}}s_{\mathrm{bl}}}$. Hence
Lemma~\ref{bl:frozen} applies at each step and produces the
required decaying profile while preserving the prescribed
tangential support.

The bound~\eqref{bl:tangential-order} gives the sharper estimate
\begin{equation}\label{bl:profile-degree}
 \deg(e^{c_{\mathrm{bl}}s_{\mathrm{bl}}}U_\ell)\le1+3\ell.
\end{equation}
Indeed, a term of $\mathcal P_{\mathfrak{s}}^{(q)}$ contains
$s_{\mathrm{bl}}^{q-m}$ and at most $m$ derivatives with respect to
the tangential variables $y$.
Its polynomial degree therefore increases by at most $q$,
including derivatives falling on
$e^{-c_{\mathrm{bl}}(y)s_{\mathrm{bl}}}$. The inverse in
Lemma~\ref{bl:frozen} adds at most two to the polynomial
degree. Thus, assuming the estimate for lower indices, the
$q$th source in the recursion produces degree at most
\[
 1+3(\ell-q)+q+2\le1+3\ell,
 \qquad q\ge1,
\]
which proves \eqref{bl:profile-degree} by induction.

The coefficients in \eqref{bl:Pq} involve at most $q+1$
normal derivatives and one tangential derivative of the
original coefficients, while each recursion applies at most
two further tangential derivatives. The derivative bounds displayed above therefore control all profiles with
$\ell\le N_{\mathrm{bl}}$ and the Taylor remainders. The
additional $S+4$ derivatives cover the final operator
application and the prescribed differentiations. The frozen inverse uses finite algebraic operations with positive
denominators, so the required normal derivatives remain within the
ranges of the derivative bounds stated above.

On the region where $\eta=1$, the recursion
\eqref{bl:recursion} cancels all powers
$h^0,\ldots,h^{N_{\mathrm{bl}}}$. Let
$p_{N_{\mathrm{bl}}}\in\mathbb N_0$ be an integer large enough to bound the
polynomial factors arising from the Taylor remainders
and the finitely many profiles at this truncation order.
Taylor's formula
\eqref{bl:taylor} gives
\[
 \left|
 \mathcal P_{\mathfrak{s},h}
 \sum_{\ell=0}^{N_{\mathrm{bl}}}h^\ell U_\ell
 \right|
 \le C_{N_{\mathrm{bl}}}
 h^{N_{\mathrm{bl}}+1}
 (1+s_{\mathrm{bl}})^{p_{N_{\mathrm{bl}}}}e^{-c_0s_{\mathrm{bl}}}.
\]
Multiplication by $h^{-2}$ shows that the physical residual
is bounded by
\[
 C_{N_{\mathrm{bl}}}h^{N_{\mathrm{bl}}-1}
 (1+s_{\mathrm{bl}})^{p_{N_{\mathrm{bl}}}}e^{-c_0s_{\mathrm{bl}}}.
\]
Since the volume form is $J(y,r)\,dy\,dr$, with $dr=h\,d s_{\mathrm{bl}}$
and $J$ bounded above and below, integration in the stretched
normal variable gives \eqref{bl:residual}.

On $\operatorname{supp}\eta'$, one has $r\ge r_*>0$.
Consequently, every cutoff commutator is bounded by
$C h^{-M}e^{-c_0r_*/h}$ for some fixed $M\ge0$, and hence
decays faster than any prescribed power of $h$. The same
integration, with a factor at most $h^{-1}$ for each physical
derivative, proves \eqref{bl:appnorm}.

By Section~\ref{an:tools}, the physical Dirichlet realization
$P_D$ has a bounded inverse
\[
 P_D^{-1}:H^{-1}(\Omega;\mathbb C^3)
 \longrightarrow H^1_0(\Omega;\mathbb C^3).
\]
Define
\[
 u_h=u_{h,N_{\mathrm{bl}}}^{\mathrm{app}}
      -P_D^{-1}
       \bigl(Pu_{h,N_{\mathrm{bl}}}^{\mathrm{app}}\bigr).
\]
The correction belongs to $H^1_0(\Omega;\mathbb C^3)$, so
it preserves the boundary trace. Moreover,
$Pu_h=0$ in the weak sense. Finally,
\[
\begin{aligned}
 \|u_h-u_{h,N_{\mathrm{bl}}}^{\mathrm{app}}\|_{H^1(\Omega)}
 &\le
 \|P_D^{-1}\|_{H^{-1}\to H^1_0}
 \|Pu_{h,N_{\mathrm{bl}}}^{\mathrm{app}}\|_{H^{-1}(\Omega)}
 \\
 &\le C
 \|Pu_{h,N_{\mathrm{bl}}}^{\mathrm{app}}\|_{L^2(\Omega)}
 \le C_{N_{\mathrm{bl}}}
 h^{N_{\mathrm{bl}}-1/2}.
\end{aligned}
\]
This proves \eqref{bl:correct} and completes the proof.
\end{proof}

\subsection{Recovery at the first two elastic orders}

Consider two materials satisfying the preceding assumptions
and having equal full boundary maps~\eqref{en:dn-definition}.
For either material, the weak representation
\eqref{en:dn-weak-form} expresses the boundary pairing in
terms of the corresponding interior solutions. We first use
the leading boundary response to recover the boundary values
of the Lam\'e parameters, and then use the bilinear identity
between the two materials to recover their first normal
derivatives.

At the boundary point under consideration, choose coordinates
so that
\[
 e_{\mathrm{bl}}=e_1,\qquad
 n_{\mathrm{in}}=e_3,\qquad
 \nu=-e_3,\qquad c_{\mathrm{bl}}=1.
\]
In these frozen coordinates, write
$N_\pm=N_{\mathfrak{s}}|_{\mathfrak{s}=\pm1}$,
where $N_{\mathfrak{s}}$ is defined in~\eqref{bl:frozen-profile}.
Thus
\[
 N_\pm=a_{\mathrm{bl}}
 \begin{pmatrix}
 1&0&\pm i\\
 0&0&0\\
 \pm i&0&-1
 \end{pmatrix}.
\]

For each material separately, write $\Lambda$ for its physical
Dirichlet-to-Neumann map and let
$F_\pm(y,s_{\mathrm{bl}})=F_{\mathfrak{s}}(y,s_{\mathrm{bl}})|_{\mathfrak{s}=\pm1}$
be the frozen profiles in~\eqref{bl:frozen-profile}.
With $\lambda_b=\lambda(y,0)$, $\mu_b=\mu(y,0)$ and all
coordinate vectors evaluated at $(y,0)$, define
\[
\begin{aligned}
Z_+(y):={}&-\big[\mu_b I_3+(\lambda_b+\mu_b)
 n_{\mathrm{in}}\otimes n_{\mathrm{in}}\big]
 \partial_{s_{\mathrm{bl}}}F_+(y,0)\\
&-i\lambda_b n_{\mathrm{in}}\otimes\tau_{\mathrm{tan}}
 -i\mu_b\tau_{\mathrm{tan}}\otimes n_{\mathrm{in}}.
\end{aligned}
\]
This matrix gives the frozen outward traction at $s_{\mathrm{bl}}=0$
when tangential differentiation is $i\tau_{\mathrm{tan}}$
and normal differentiation is $n_{\mathrm{in}}\partial_{s_{\mathrm{bl}}}$;
the outward normal is $-n_{\mathrm{in}}$.

Fix $N_{\mathrm{bl}}\ge1$ and
$A_+,A_-\in C_c^\infty(V;\mathbb C^3)$ independently of $h$.
Let
\[
 f_{\pm,h}(X(y,0))=e^{\pm iy_1/h}A_\pm(y),
\]
extended by zero outside the boundary patch, and let
$u_{\pm,h}$ and $v_{\pm,h}$ be the exact and approximate
solutions supplied by Proposition~\ref{bl:layer} for these
boundary data. That proposition gives
\[
 \|v_{\pm,h}\|_{H^1(\Omega)}\le Ch^{-1/2},\qquad
 \|u_{\pm,h}-v_{\pm,h}\|_{H^1(\Omega)}
 \le Ch^{N_{\mathrm{bl}}-1/2}.
\]
Use $u_{-,h}$ as the lifting in~\eqref{en:dn-weak-form}.
Replacing both exact solutions by $v_{\pm,h}$ in this
volume pairing changes the normalized boundary pairing
$h\langle\Lambda f_{+,h},f_{-,h}\rangle$ by at most
\[
 Ch\big(h^{N_{\mathrm{bl}}-1/2}h^{-1/2}
       +h^{2N_{\mathrm{bl}}-1}\big)
 =C\big(h^{N_{\mathrm{bl}}}+h^{2N_{\mathrm{bl}}}\big).
\]
This estimate uses the global $H^1$ bounds. The normalized
potential integral for the approximate layers is $O(h^2)$
by~\eqref{bl:appnorm}. Formula~\eqref{bl:D}, the finite
polynomial--exponential profiles, and $r=hs_{\mathrm{bl}}$ show that the
remaining normalized approximate pairing differs from its
frozen leading elastic integral by $O(h)$.

For each fixed $y$, integrate that frozen elastic integral
by parts in $s_{\mathrm{bl}}$, pairing $F_+(y,s_{\mathrm{bl}})A_+(y)$ with
$F_-(y,s_{\mathrm{bl}})A_-(y)$. The opposite tangential signs give the
bilinear pairing for the frozen equation, whose interior
term vanishes by~\eqref{bl:P0} and
$\mathcal P_{+1}^{(0)}F_+=0$.
Exponential decay removes the endpoint at infinity, while
$F_\pm(y,0)=I_3$ makes the endpoint at $s_{\mathrm{bl}}=0$ equal to
$A_-(y)^TZ_+(y)A_+(y)$. Since $J(y,0)\,dy$ is the
boundary surface measure, it follows that
\[
 \lim_{h\to0}h\langle\Lambda f_{+,h},f_{-,h}\rangle
 =\int_V J(y,0)A_-(y)^TZ_+(y)A_+(y)\,dy.
\]
The limit is taken for each fixed pair of amplitudes before
pointwise localization. At the chosen boundary point $y_0$,
use the normalized frame above and write $Z_+=Z_+(y_0)$.
Then $\partial_{s_{\mathrm{bl}}}F_+(y_0,0)=-(I_3+N_+)$, and direct evaluation
gives

\begin{equation}\label{bl:Z}
 Z_+=
 \begin{pmatrix}
 d_{\mathrm{bd}}&0&-ib_*\\
 0&\mu&0\\
 ib_*&0&d_{\mathrm{bd}}
 \end{pmatrix},
 \qquad
 d_{\mathrm{bd}}=
 \frac{2\mu(\lambda+2\mu)}{\lambda+3\mu},
 \qquad
 b_*=\frac{2\mu^2}{\lambda+3\mu}.
\end{equation}
Equivalently,
\[
 Z_+=A_0(I_3+N_+)-iR,
 \qquad
 A_0=\operatorname{diag}(\mu,\mu,\lambda+2\mu),
 \qquad
 R=
 \begin{pmatrix}
 0&0&\mu\\
 0&0&0\\
 \lambda&0&0
 \end{pmatrix}.
\]
Reversing the phase gives $Z_-=Z_+^T$. Since the two
boundary maps agree, their leading traction matrices agree.
The entries of $Z_+$ therefore determine the boundary values:
\[
 \mu=(Z_+)_{22},\qquad
 \lambda=
 \frac{3d_{\mathrm{bd}}\mu-4\mu^2}
 {2\mu-d_{\mathrm{bd}}},
 \qquad
 2\mu-d_{\mathrm{bd}}=b_*>0.
\]
Thus the boundary values of $\lambda$ and $\mu$ agree for
the two materials.

Define the boundary Taylor coefficients by
\begin{equation}\label{bl:boundary-jets}
\begin{aligned}
L_j^{\mathrm{bd}}(y)
&:=\partial_r^j(\lambda_1-\lambda_2)(y,0),
&&j\ge1,\\
M_j^{\mathrm{bd}}(y)
&:=\partial_r^j(\mu_1-\mu_2)(y,0),
&&j\ge1,\\
V_j^{\mathrm{bd}}(y)
&:=\partial_r^{j-2}(W_1-W_2)(y,0),
&&j\ge2.
\end{aligned}
\end{equation}
In~\eqref{bl:boundary-jets}, $j$ is the elastic Taylor order.
At $j=1$ the boundary Taylor coefficients are $L_1^{\mathrm{bd}}$
and $M_1^{\mathrm{bd}}$. For the two media write
$\delta_{12}b:=b_1-b_2$.

Define the bilinear form
\[
\mathfrak B(v_1,v_2)
:=\int_\Omega\left[
 \delta_{12}\lambda\operatorname{div}v_1\operatorname{div}v_2
 +2\delta_{12}\mu\varepsilon(v_1):\varepsilon(v_2)
 -\delta_{12}W\,v_1\cdot v_2
 \right]\,dx.
\]
Applying~\eqref{en:dn-weak-form} to solutions $u_1,u_2$
for the two materials and using equality of the boundary
maps gives the bilinear identity
\begin{equation}\label{bl:alessandrini}
\mathfrak B(u_1,u_2)=0.
\end{equation}

Since the boundary values of the Lam\'e parameters agree,
Taylor's formula~\eqref{bl:taylor} gives
\[
 \delta_{12}\lambda
 =rL_1^{\mathrm{bd}}+O(r^2),
 \qquad
 \delta_{12}\mu
 =rM_1^{\mathrm{bd}}+O(r^2).
\]
Fix the boundary point $y_0$ and suppress the dependence
on this point in the following formulas. Take the exact solutions
$u_{1,h},u_{2,h}$ from Proposition~\ref{bl:layer} with phases $+,-$
and $N_{\mathrm{bl}}\ge3$. First take the high-frequency limit
for a fixed tangential cutoff, and then localize at $y_0$.
By \eqref{bl:appnorm} and \eqref{bl:correct},
\[
\left|\int_\Omega\delta_{12}W\,u_{1,h}\cdot u_{2,h}\,dx\right|
 \le \|\delta_{12}W\|_{L^\infty(\Omega)}
      \|u_{1,h}\|_{L^2(\Omega)}\|u_{2,h}\|_{L^2(\Omega)}
 =O(h).
\]
The potential first contributes to the leading boundary moment at
elastic order $j=2$, through $V_2^{\mathrm{bd}}$ in
\eqref{bl:leading-moment}.
Rows of the moment matrix correspond to the negative-phase
testing polarization, and columns to the positive-phase
input polarization.

For $p=1,2,3$, define the polynomial gradient matrices
\begin{equation}\label{bl:H1-gradient}
 G_\pm^p=
 \big[\,\pm i\widetilde F_\pm^p,\ 0,\
                   ((s_{\mathrm{bl}}-1)N_\pm-I_3)^p\,\big],
 \qquad
 \widetilde F_\pm=I_3-s_{\mathrm{bl}} N_\pm.
\end{equation}
Here the superscript $p$ denotes the $p$th column. We write
$\operatorname{sym}A=\tfrac12(A+A^T)$ for the symmetric part of a matrix.
For $p,q\in\{1,2,3\}$, introduce the bilinear products
\begin{equation}\label{bl:polarization-products}
\begin{aligned}
D_{pq}(s_{\mathrm{bl}})
&:=\operatorname{tr}G_-^p(s_{\mathrm{bl}})\,
  \operatorname{tr}G_+^q(s_{\mathrm{bl}}),\\
S_{pq}(s_{\mathrm{bl}})
&:=2\operatorname{sym}G_-^p(s_{\mathrm{bl}}):
       \operatorname{sym}G_+^q(s_{\mathrm{bl}}),\\
Q_{pq}(s_{\mathrm{bl}})
&:=\widetilde F_-^p(s_{\mathrm{bl}})\cdot
       \widetilde F_+^q(s_{\mathrm{bl}}).
\end{aligned}
\end{equation}
For the diagonal moments write
$D_p:=D_{pp}$, $S_p:=S_{pp}$ and $Q_p:=Q_{pp}$.

To identify the first-order boundary moment, set $r=hs_{\mathrm{bl}}$.
At the chosen boundary point $y_0$, suppress the dependence
on $y_0$ in the following formulas. Taylor's formula gives
\[
\begin{aligned}
\delta_{12}\lambda(y_0,hs_{\mathrm{bl}})
 &=hs_{\mathrm{bl}} L_1^{\mathrm{bd}}+O(h^2s_{\mathrm{bl}}^2),\\
\delta_{12}\mu(y_0,hs_{\mathrm{bl}})
 &=hs_{\mathrm{bl}} M_1^{\mathrm{bd}}+O(h^2s_{\mathrm{bl}}^2).
\end{aligned}
\]
Substituting these expansions and the leading gradient
matrices into the elastic terms of~\eqref{bl:alessandrini},
the principal integrand, before multiplication by the
volume Jacobian and the common localization factors, is
\[
h^{-1}s_{\mathrm{bl}} e^{-2s_{\mathrm{bl}}}
\left[
L_1^{\mathrm{bd}}D_{pq}(s_{\mathrm{bl}})
+M_1^{\mathrm{bd}}S_{pq}(s_{\mathrm{bl}})
\right].
\]
The opposite phases cancel. Under the change of variables
$r=hs_{\mathrm{bl}}$, the volume element is
\[
dx=J(y,hs_{\mathrm{bl}})\,dy\,h\,d s_{\mathrm{bl}}.
\]
Thus the factor $h^{-1}$ is canceled by the factor $h$
in the normal integration measure. The error estimates in
Proposition~\ref{bl:layer} justify the high-frequency limit.
After localization at $y_0$ and division by the positive
Jacobian and localization factors, the remaining normal
integral is the first-order moment.

Accordingly, define the first-order boundary moment matrix
$H_1$ by
\begin{equation}\label{bl:H1-moment}
 (H_1)_{pq}
 =\int_0^\infty s_{\mathrm{bl}} e^{-2s_{\mathrm{bl}}}
 \left[
 L_1^{\mathrm{bd}}D_{pq}(s_{\mathrm{bl}})
 +M_1^{\mathrm{bd}}S_{pq}(s_{\mathrm{bl}})
 \right]d s_{\mathrm{bl}},
 \qquad p,q=1,2,3.
\end{equation}
Since~\eqref{bl:alessandrini} vanishes for every admissible
choice of boundary data, its high-frequency limit and
localization give
\[
(H_1(y_0))_{pq}=0,\qquad p,q=1,2,3.
\]
Hence $H_1(y_0)=0$. We now compute its entries explicitly.

Using~\eqref{bl:H1-gradient} and
\[
 \int_0^\infty s_{\mathrm{bl}}^se^{-2s_{\mathrm{bl}}}\,d s_{\mathrm{bl}}
 =\frac{s!}{2^{s+1}},
 \qquad s\in\mathbb N_0,
\]
gives
\[
 H_{1,22}=\frac12M_1^{\mathrm{bd}},
 \qquad H_{1,31}=-H_{1,13},
\]
and
\begin{align*}
 -iH_{1,13}
 &=\tfrac14\big[
 (1-a_{\mathrm{bl}})^2L_1^{\mathrm{bd}}
 +(7a_{\mathrm{bl}}^2-2a_{\mathrm{bl}}-1)
 M_1^{\mathrm{bd}}\big],\\
 H_{1,11}
 &=\tfrac14\big[
 (1-a_{\mathrm{bl}})^2L_1^{\mathrm{bd}}
 +(7a_{\mathrm{bl}}^2-6a_{\mathrm{bl}}+3)
 M_1^{\mathrm{bd}}\big],\\
 H_{1,33}
 &=\tfrac14\big[
 (1-a_{\mathrm{bl}})^2L_1^{\mathrm{bd}}
 +(7a_{\mathrm{bl}}^2+2a_{\mathrm{bl}}+3)
 M_1^{\mathrm{bd}}\big].
\end{align*}
The other entries are zero. Using
\[
 \lambda=\mu\frac{3a_{\mathrm{bl}}-1}{1-a_{\mathrm{bl}}},
\]
the expression for $-iH_{1,13}$ reduces to
\[
 b=
 \frac{\mu^2L_1^{\mathrm{bd}}
       +(\lambda^2-2\mu^2)M_1^{\mathrm{bd}}}
      {(\lambda+3\mu)^2}.
\]
The diagonal entries then give the explicit matrix
\begin{equation}\label{bl:H1}
 H_1=
 \begin{pmatrix}
 M_1^{\mathrm{bd}}(1-a_{\mathrm{bl}})+b&0&ib\\
 0&M_1^{\mathrm{bd}}/2&0\\
 -ib&0&M_1^{\mathrm{bd}}(1+a_{\mathrm{bl}})+b
 \end{pmatrix}.
\end{equation}
Since $H_1(y_0)=0$, its $(2,2)$ entry gives
\[
 \frac12M_1^{\mathrm{bd}}=0,
 \qquad M_1^{\mathrm{bd}}=0.
\]
Then $H_{1,13}=0$ and the formula for $-iH_{1,13}$ imply
\[
 (1-a_{\mathrm{bl}})^2L_1^{\mathrm{bd}}=0.
\]
Since $a_{\mathrm{bl}}<1$, we conclude that
$L_1^{\mathrm{bd}}=0$.

\begin{lemma}
\label{bl:frozen-moments}
Fix $a_{\mathrm{bl}}\in(0,1)$ and an integer $j\ge2$.
Use the frozen matrices $N_\pm,\widetilde F_\pm,G_\pm^p$
defined in the preceding subsection, with
$c_{\mathrm{bl}}=1$ and
$e_{\mathrm{bl}}=e_1$, $n_{\mathrm{in}}=e_3$.
For the diagonal products, write
\[
 D_p:=D_{pp},\qquad
 S_p:=S_{pp},\qquad
 Q_p:=Q_{pp},\qquad p=1,2,3,
\]
where $D_{pq},S_{pq},Q_{pq}$ are defined in
\eqref{bl:polarization-products}.

At a fixed boundary point $y$, regard
$L_j^{\mathrm{bd}}(y)$, $M_j^{\mathrm{bd}}(y)$, and
$V_j^{\mathrm{bd}}(y)$, defined in~\eqref{bl:boundary-jets},
as scalars, and define
\begin{equation}\label{bl:diagonal-moments}
 E_{pp}=2^{j+3}\int_0^\infty e^{-2s_{\mathrm{bl}}}
 \left[\frac{s_{\mathrm{bl}}^j}{j!}(L_j^{\mathrm{bd}}D_p+M_j^{\mathrm{bd}}S_p)
       -\frac{s_{\mathrm{bl}}^{j-2}}{(j-2)!}V_j^{\mathrm{bd}}Q_p\right]d s_{\mathrm{bl}},
 \qquad p=1,2,3.
\end{equation}
Then
\begin{align}
 E_{22}&=8M_j^{\mathrm{bd}}-16V_j^{\mathrm{bd}},\notag\\
 E_{33}-E_{11}
 &=32a_{\mathrm{bl}}[jM_j^{\mathrm{bd}}-(j-1)V_j^{\mathrm{bd}}],
 \notag\\
 E_{11}
 &=4(1-a_{\mathrm{bl}})^2L_j^{\mathrm{bd}}
   +B_jM_j^{\mathrm{bd}}-C_jV_j^{\mathrm{bd}},
 \label{bl:moments}\\
 B_j&=8j(j+1)a_{\mathrm{bl}}^2-16j a_{\mathrm{bl}}
       +12a_{\mathrm{bl}}^2-8a_{\mathrm{bl}}+12,\notag\\
 C_j&=8j(j-1)a_{\mathrm{bl}}^2-16(j-1)a_{\mathrm{bl}}+16.
 \notag
\end{align}
In particular, $E_{11}=E_{22}=E_{33}=0$ implies
$L_j^{\mathrm{bd}}=M_j^{\mathrm{bd}}=V_j^{\mathrm{bd}}=0$.
\end{lemma}

\begin{proof}
The bilinear products of the matrices in
\eqref{bl:H1-gradient} yield
\begin{align*}
 (D_1,D_2,D_3)
 &=((1-a_{\mathrm{bl}})^2,\,0,\,(1-a_{\mathrm{bl}})^2),\\
 (S_1,S_2,S_3)
 &=\big(
 8a_{\mathrm{bl}}^2s_{\mathrm{bl}}^2-8a_{\mathrm{bl}}(a_{\mathrm{bl}}+1)s_{\mathrm{bl}}
       +3a_{\mathrm{bl}}^2+2a_{\mathrm{bl}}+3,\ 2,\\
 &\hspace{7em}
 8a_{\mathrm{bl}}^2s_{\mathrm{bl}}^2-8a_{\mathrm{bl}}(a_{\mathrm{bl}}-1)s_{\mathrm{bl}}
       +3a_{\mathrm{bl}}^2-6a_{\mathrm{bl}}+3
 \big),\\
 (Q_1,Q_2,Q_3)
 &=(1-2a_{\mathrm{bl}}s_{\mathrm{bl}}+2a_{\mathrm{bl}}^2s_{\mathrm{bl}}^2,\ 1,\
     1+2a_{\mathrm{bl}}s_{\mathrm{bl}}+2a_{\mathrm{bl}}^2s_{\mathrm{bl}}^2).
\end{align*}
For example,
\[
\operatorname{tr}G_+^1=i(1-a_{\mathrm{bl}}),
\qquad
\operatorname{tr}G_-^1=-i(1-a_{\mathrm{bl}}),
\]
and both traces vanish for $p=2$. For $p=3$ both traces
equal $a_{\mathrm{bl}}-1$. In particular,
\[
 S_3-S_1=16a_{\mathrm{bl}}s_{\mathrm{bl}}-8a_{\mathrm{bl}},
 \qquad
 Q_3-Q_1=4a_{\mathrm{bl}}s_{\mathrm{bl}}.
\]

The factorial identity
\begin{equation}\label{bl:factorial-moment}
 \frac1{n!}\int_0^\infty s_{\mathrm{bl}}^{n+q}e^{-2s_{\mathrm{bl}}}\,d s_{\mathrm{bl}}
 =\frac{(n+q)!}{n!\,2^{n+q+1}},
 \qquad n,q\ge0,
\end{equation}
now gives
\[
 E_{22}
 =2^{j+3}\left(
 \frac{2M_j^{\mathrm{bd}}}{2^{j+1}}
 -\frac{V_j^{\mathrm{bd}}}{2^{j-1}}
 \right)
 =8M_j^{\mathrm{bd}}-16V_j^{\mathrm{bd}},
\]
and, by the two displayed differences,
\[
\begin{aligned}
 E_{33}-E_{11}
 &=2^{j+3}\bigg[
 \left(
 \frac{16a_{\mathrm{bl}}(j+1)}{2^{j+2}}
 -\frac{8a_{\mathrm{bl}}}{2^{j+1}}
 \right)M_j^{\mathrm{bd}}\\
 &\hspace{5em}
 -\frac{4a_{\mathrm{bl}}(j-1)}{2^j}V_j^{\mathrm{bd}}
 \bigg]\\
 &=32a_{\mathrm{bl}}
 \big[jM_j^{\mathrm{bd}}-(j-1)V_j^{\mathrm{bd}}\big].
\end{aligned}
\]

The first diagonal polarization gives the separate
coefficients. Here $[X]E_{11}$ denotes the coefficient
of the formal variable $X$ in $E_{11}$, regarded as a
linear expression in
$L_j^{\mathrm{bd}},M_j^{\mathrm{bd}},V_j^{\mathrm{bd}}$:
\begin{align*}
 [L_j^{\mathrm{bd}}]E_{11}
 &=4(1-a_{\mathrm{bl}})^2,\\
 [M_j^{\mathrm{bd}}]E_{11}
 &=8a_{\mathrm{bl}}^2(j+1)(j+2)
   -16a_{\mathrm{bl}}(a_{\mathrm{bl}}+1)(j+1)\\
 &\qquad
   +4(3a_{\mathrm{bl}}^2+2a_{\mathrm{bl}}+3)
   =B_j,\\
 -[V_j^{\mathrm{bd}}]E_{11}
 &=8a_{\mathrm{bl}}^2j(j-1)
   -16a_{\mathrm{bl}}(j-1)+16
   =C_j.
\end{align*}
This proves~\eqref{bl:moments} by finite algebra for every
integer $j\ge2$, including $j=2$.

If all three moments vanish, the first identity yields
\[
 M_j^{\mathrm{bd}}=2V_j^{\mathrm{bd}}.
\]
The second then gives
\[
 32a_{\mathrm{bl}}(j+1)V_j^{\mathrm{bd}}=0,
\]
so $V_j^{\mathrm{bd}}=M_j^{\mathrm{bd}}=0$. The third
identity now gives
\[
 4(1-a_{\mathrm{bl}})^2L_j^{\mathrm{bd}}=0.
\]
Since $a_{\mathrm{bl}}\in(0,1)$, all three unknowns vanish.
\end{proof}

\subsection{Determination of Boundary Derivatives of Arbitrary Order}

Recall from~\eqref{bl:frozen-profile} that the frozen
decaying profile is
\[
F_{\mathfrak{s}}(y,s_{\mathrm{bl}})
=e^{-c_{\mathrm{bl}}(y)s_{\mathrm{bl}}}
 \bigl(I_3-c_{\mathrm{bl}}(y)s_{\mathrm{bl}} N_{\mathfrak{s}}(y)\bigr),
\qquad \mathfrak{s}\in\{+,-\}.
\]
For fixed boundary amplitudes
$A_\pm\in C_c^\infty(V;\mathbb C^3)$, define the frozen
displacement profiles and their leading gradient matrices by
\begin{equation}\label{bl:leading-fields}
 \mathcal U_{\mathfrak{s}}=F_{\mathfrak{s}}A_{\mathfrak{s}},
 \qquad
 \mathcal G_{\mathfrak{s}}
 =i\mathfrak{s}\mathcal U_{\mathfrak{s}}\otimes\tau_{\mathrm{tan}}
  +(\partial_{s_{\mathrm{bl}}}\mathcal U_{\mathfrak{s}})\otimes n_{\mathrm{in}},
 \qquad \mathfrak{s}\in\{+,-\}.
\end{equation}
The signs $+$ and $-$ are identified with $+1$ and $-1$,
respectively. The vectors $\tau_{\mathrm{tan}}$ and
$n_{\mathrm{in}}$ are evaluated at $(y,0)$ as in~\eqref{bl:P0}.

At a fixed boundary point $y_0$, use the normalized frame
\[
c_{\mathrm{bl}}(y_0)=1,\qquad
\tau_{\mathrm{tan}}(y_0)=e_1,\qquad
n_{\mathrm{in}}(y_0)=e_3.
\]
Then
\[
F_\pm(y_0,s_{\mathrm{bl}})=e^{-s_{\mathrm{bl}}}\widetilde F_\pm(s_{\mathrm{bl}}),
\qquad
\widetilde F_\pm(s_{\mathrm{bl}})=I_3-s_{\mathrm{bl}} N_\pm.
\]
Writing $A_\pm(y_0)=\sum_{p=1}^3A_{\pm,p}(y_0)e_p$,
the general profiles satisfy
\[
\begin{aligned}
\mathcal U_\pm(y_0,s_{\mathrm{bl}})
 &=e^{-s_{\mathrm{bl}}}\sum_{p=1}^3
 A_{\pm,p}(y_0)\widetilde F_\pm^p(s_{\mathrm{bl}}),\\
\mathcal G_\pm(y_0,s_{\mathrm{bl}})
 &=e^{-s_{\mathrm{bl}}}\sum_{p=1}^3
 A_{\pm,p}(y_0)G_\pm^p(s_{\mathrm{bl}}).
\end{aligned}
\]
Here $\widetilde F_\pm^p$ and $G_\pm^p$ are the columns
defined in~\eqref{bl:H1-gradient}. These identities hold
at the chosen point in the normalized frozen frame; the
general profiles in~\eqref{bl:leading-fields} are used
on the full boundary patch.

For each integer $j\ge2$, define the leading boundary
moment functional associated with $A_+$ and $A_-$ by
\begin{equation}\label{bl:leading-moment}
\begin{split}
\mathfrak M_j
=\int_V\int_0^\infty J(y,0)\bigg[
&\frac{s_{\mathrm{bl}}^j}{j!}\big\{
 L_j^{\mathrm{bd}}(y)\operatorname{tr}\mathcal G_+
                   \operatorname{tr}\mathcal G_-
 +2M_j^{\mathrm{bd}}(y)\operatorname{sym}\mathcal G_+:
                   \operatorname{sym}\mathcal G_-\big\}\\
&-\frac{s_{\mathrm{bl}}^{j-2}}{(j-2)!}V_j^{\mathrm{bd}}(y)
  \mathcal U_+\cdot\mathcal U_-\bigg] \,d s_{\mathrm{bl}}\,dy.
\end{split}
\end{equation}
Here $\mathcal U_\pm$ and $\mathcal G_\pm$ are the
profiles defined in~\eqref{bl:leading-fields} and
are evaluated at $(y,s_{\mathrm{bl}})$.

\begin{lemma}
\label{bl:order-separation}
Fix an integer $j\ge2$, a compact boundary patch, and
amplitudes $A_+,A_-\in C_c^\infty(V;\mathbb C^3)$.
Suppose the two materials have common boundary values of
$\lambda,\mu$ and, uniformly on the collar,
\[
\begin{aligned}
 \delta_{12}\lambda
 &=\frac{r^j}{j!}L_j^{\mathrm{bd}}+O(r^{j+1}),\\
 \delta_{12}\mu
 &=\frac{r^j}{j!}M_j^{\mathrm{bd}}+O(r^{j+1}),\\
 \delta_{12}W
 &=\frac{r^{j-2}}{(j-2)!}V_j^{\mathrm{bd}}+O(r^{j-1}),
\end{aligned}
\]
where the boundary Taylor coefficients are defined in
\eqref{bl:boundary-jets}. Let $N_{\mathrm{bl}}\ge j+1$,
and let $u_{1,h},u_{2,h}$ be the exact solutions of
Proposition~\ref{bl:layer} with phases $+,-$,
amplitudes $A_+,A_-$, and a common normal cutoff.
Then, with $\mathfrak B$ and $\mathfrak M_j$ defined above,
\begin{equation}\label{bl:order-separation-bound}
 \left|
 h^{1-j}\mathfrak B(u_{1,h},u_{2,h})-\mathfrak M_j
 \right|
 \le C\left(
 h+h^{N_{\mathrm{bl}}-j}
  +h^{2N_{\mathrm{bl}}-j}
 \right),\qquad 0<h\le h_0\le1.
\end{equation}
Here $h_0$ is chosen so that both layers are defined.
For each fixed $j$, $N_{\mathrm{bl}}$, and pair of amplitudes,
the constant depends on the corresponding bounds for finitely many coefficient
derivatives, geometry, cutoffs, and the two physical Dirichlet
inverse norms
\[
\|P_{D,1}^{-1}\|_{H^{-1}\to H_0^1},
\qquad
\|P_{D,2}^{-1}\|_{H^{-1}\to H_0^1},
\]
and is independent of $h$.
\end{lemma}

\begin{proof}
Write $v_{1,h},v_{2,h}$ for the corresponding approximate
layers. For $\mathfrak{s}=+$ use $v_{1,h}$ and for
$\mathfrak{s}=-$ use $v_{2,h}$, and set in the collar
\[
\begin{aligned}
 \widetilde v_{\mathfrak{s},h}(y,s_{\mathrm{bl}})
 &=e^{-i\mathfrak{s}y_1/h}
   v_{\mathfrak{s},h}(X(y,hs_{\mathrm{bl}})),\\
 \widetilde G_{\mathfrak{s},h}(y,s_{\mathrm{bl}})
 &=he^{-i\mathfrak{s}y_1/h}
   \nabla_xv_{\mathfrak{s},h}(X(y,hs_{\mathrm{bl}})).
\end{aligned}
\]
The notation $v_{+,h}=v_{1,h}$ and $v_{-,h}=v_{2,h}$
is used in these two definitions only. Formula~\eqref{bl:D},
the finite collection of profiles, and Taylor's formula for
$n_{\mathrm{in}},w$ give, on $0\le hs_{\mathrm{bl}}\le r_0$,
\begin{equation}\label{bl:profile-error}
 \begin{aligned}
 |\widetilde v_{\mathfrak{s},h}-\mathcal U_{\mathfrak{s}}|
  +|\widetilde G_{\mathfrak{s},h}-\mathcal G_{\mathfrak{s}}|
       &\le Ch(1+s_{\mathrm{bl}})^Ne^{-c_0s_{\mathrm{bl}}},\\
 |\widetilde v_{\mathfrak{s},h}|+|\widetilde G_{\mathfrak{s},h}|
  +|\mathcal U_{\mathfrak{s}}|+|\mathcal G_{\mathfrak{s}}|
       &\le C(1+s_{\mathrm{bl}})^Ne^{-c_0s_{\mathrm{bl}}}.
 \end{aligned}
\end{equation}
Here $N$ is finite and may depend on $N_{\mathrm{bl}}$.
To verify the first line, all profiles of index at least
one have a factor $h$; the tangential derivatives of the
profiles in the scaled gradient appear with a factor $h$;
and
\[
\begin{aligned}
n_{\mathrm{in}}(y,hs_{\mathrm{bl}})-n_{\mathrm{in}}(y,0)&=O(hs_{\mathrm{bl}}),\\
w^1(y,hs_{\mathrm{bl}})-\tau_{\mathrm{tan}}(y)&=O(hs_{\mathrm{bl}}).
\end{aligned}
\]
The normal cutoff satisfies $|\eta(hs_{\mathrm{bl}})-1|\le Chs_{\mathrm{bl}}$, while
its derivative in the scaled gradient appears as
$h\eta'(hs_{\mathrm{bl}})$. These bounds prove the first line also in
the cutoff transition region. The second line follows from
the same finite profile bounds.

After the substitution $r=hs_{\mathrm{bl}}$, the scaled elastic
coefficient factors are
\[
\begin{aligned}
 h^{-j}\delta_{12}\lambda(y,hs_{\mathrm{bl}})
 &=\frac{s_{\mathrm{bl}}^j}{j!}L_j^{\mathrm{bd}}(y)+O(hs_{\mathrm{bl}}^{j+1}),\\
 h^{-j}\delta_{12}\mu(y,hs_{\mathrm{bl}})
 &=\frac{s_{\mathrm{bl}}^j}{j!}M_j^{\mathrm{bd}}(y)+O(hs_{\mathrm{bl}}^{j+1}),
\end{aligned}
\]
and the potential factor is
\[
 h^{2-j}\delta_{12}W(y,hs_{\mathrm{bl}})
 =\frac{s_{\mathrm{bl}}^{j-2}}{(j-2)!}V_j^{\mathrm{bd}}(y)
   +O(hs_{\mathrm{bl}}^{j-1}).
\]
Moreover, $J(y,hs_{\mathrm{bl}})-J(y,0)=O(hs_{\mathrm{bl}})$. The signs of the two
phases are opposite, so their product is exactly one.
Using~\eqref{bl:profile-error} in each bilinear product
bounds the difference between the scaled approximate-layer
integrand and the integrand in~\eqref{bl:leading-moment}
by $Ch(1+s_{\mathrm{bl}})^{N'}e^{-2c_0s_{\mathrm{bl}}}$, for a fixed finite $N'$.
The leading integral over $s_{\mathrm{bl}}>r_0/h$ is smaller than every
prescribed power of $h$, since its integrand has
polynomial-exponential decay. Integration yields
\[
 |\,
 h^{1-j}\mathfrak B(v_{1,h},v_{2,h})-\mathfrak M_j
 \,|\le Ch.
\]

Finally, boundedness of the coefficient differences gives
\[
 |\mathfrak B(v_1,v_2)|
 \le C\|v_1\|_{H^1}\|v_2\|_{H^1}.
\]
Expanding the difference of the two bilinear forms and
using~\eqref{bl:appnorm}--\eqref{bl:correct} gives
\[
\begin{aligned}
 h^{1-j}|\mathfrak B(u_{1,h},u_{2,h})
                     -\mathfrak B(v_{1,h},v_{2,h})|
 &\le Ch^{1-j}
 \big(h^{N_{\mathrm{bl}}-1/2}h^{-1/2}
       +h^{2N_{\mathrm{bl}}-1}\big)\\
 &=C\big(h^{N_{\mathrm{bl}}-j}
          +h^{2N_{\mathrm{bl}}-j}\big).
\end{aligned}
\]
This estimate holds on all of $\Omega$, including the
region outside the coordinate collar. Combining the last
two bounds proves~\eqref{bl:order-separation-bound}.
\end{proof}

\begin{theorem}
\label{bl:jets}
Let the two coefficient triples satisfy the standing
smoothness and strong-convexity assumptions, and suppose
that the corresponding Dirichlet problems satisfy the
nonresonance condition~\eqref{en:nonresonance}.
Assume that their full displacement-to-outward-traction
maps defined in~\eqref{en:dn-definition} agree on
$\partial\Omega$. Then, in every inward normal coordinate
chart, for every $y$ on the boundary patch and every
$q\in\mathbb N_0$,
\[
\partial_r^q(\lambda_1-\lambda_2)(y,0)
=
\partial_r^q(\mu_1-\mu_2)(y,0)
=
\partial_r^q(W_1-W_2)(y,0)
=0.
\]
Consequently, for every tangential multi-index
$\alpha\in\mathbb N_0^2$,
\[
\begin{aligned}
\partial_y^\alpha\partial_r^q(\lambda_1-\lambda_2)(y,0)
&=0,\\
\partial_y^\alpha\partial_r^q(\mu_1-\mu_2)(y,0)
&=0,\\
\partial_y^\alpha\partial_r^q(W_1-W_2)(y,0)
&=0.
\end{aligned}
\]
If $W_j=\omega^2\rho_j$ for a fixed known
$\omega\ne0$, then the same boundary-derivatives equality holds
for $\rho_1$ and $\rho_2$.
\end{theorem}

\begin{proof}
Equations~\eqref{bl:Z} and~\eqref{bl:H1} recover the zeroth and
first normal derivatives of the Lam\'e parameters. Let $j\ge2$ be the
elastic Taylor order. Assume inductively that the Lam\'e derivatives
agree for $0\le q<j$ and the potential derivatives agree for
$0\le q<j-2$. Use the notation
$L_j^{\mathrm{bd}},M_j^{\mathrm{bd}},V_j^{\mathrm{bd}}$
from~\eqref{bl:boundary-jets} for the current induction step.
Taylor's formula~\eqref{bl:taylor} then gives, uniformly on the fixed
compact patch,
\begin{align*}
 \delta_{12}\lambda(y,r)&=\frac{r^j}{j!}L_j^{\mathrm{bd}}(y)+O(r^{j+1}),\quad
 \delta_{12}\mu(y,r)=\frac{r^j}{j!}M_j^{\mathrm{bd}}(y)+O(r^{j+1}),\\
 \delta_{12}W(y,r)&=\frac{r^{j-2}}{(j-2)!}V_j^{\mathrm{bd}}(y)+O(r^{j-1}).
\end{align*}
For $j=2$, the assumption on lower-order potential derivatives is empty.
Fix $N_{\mathrm{bl}}\ge j+1$. Lemma~\ref{bl:order-separation} applies to every pair of
fixed tangential amplitudes. In particular, its global correction bound
can still be written with $k=h^{-1}$ as
\begin{equation}\label{bl:error}
 Ck^{j-1}\big(k^{1/2-N_{\mathrm{bl}}}k^{1/2}+k^{1-2N_{\mathrm{bl}}}\big)
 =C(k^{j-N_{\mathrm{bl}}}+k^{j-2N_{\mathrm{bl}}})=o(1).
\end{equation}
The additional approximate-layer error is $O(h)$ by
\eqref{bl:order-separation-bound}. Thus the identity
\eqref{bl:alessandrini} gives $\mathfrak M_j=0$ for every fixed pair
of amplitudes.

Choose $A_+=A_-=\chi e_p$ with a real $\chi\in C_c^\infty(V)$.
The leading functional is the integral of $\chi^2J(y,0)$ times a
continuous pointwise diagonal moment. If that moment were nonzero at
some point, multiplication by a fixed complex number would make its real
part strictly positive on a smaller neighborhood. A nonzero real cutoff
supported there would contradict $\mathfrak M_j=0$. Hence each pointwise
moment is zero. For each fixed cutoff, the limit $h\to0$ is taken first;
pointwise localization is then performed on the resulting moment identity.

At a fixed boundary point $y_0$, apply the preceding
localization in a Cartesian basis with
$e_{\mathrm{bl}}=e_1$ and $n_{\mathrm{in}}=e_3$.
In each localized normal moment set
$\widetilde s_{\mathrm{bl}}=c_{\mathrm{bl}}(y_0)s_{\mathrm{bl}}$.
The rescaling multiplies each leading gradient by
$c_{\mathrm{bl}}(y_0)$ and leaves each leading displacement unchanged. Thus both elastic
terms acquire the factor
$c_{\mathrm{bl}}(y_0)^{2-j-1}=c_{\mathrm{bl}}(y_0)^{1-j}$,
and the potential term acquires the same factor
$c_{\mathrm{bl}}(y_0)^{-(j-2)-1}=c_{\mathrm{bl}}(y_0)^{1-j}$.
Divide by this common positive factor and write $s_{\mathrm{bl}}$ for the normalized integration variable $\widetilde s_{\mathrm{bl}}$.
We may therefore use $c_{\mathrm{bl}}=1$ in the following
frozen moment calculation without changing the boundary
Taylor coefficients.

The limiting gradient and displacement factors are precisely
$G_\pm^p$ and $\widetilde F_\pm^p$ in
Lemma~\ref{bl:frozen-moments}: the rescaled normal derivative is
$(\partial_{s_{\mathrm{bl}}}-1)\widetilde F_\pm=((s_{\mathrm{bl}}-1)N_\pm-I_3)$, and
the rescaled tangential derivative is $\pm i\widetilde F_\pm$.
The localization argument therefore gives
$E_{11}=E_{22}=E_{33}=0$, where $E_{pp}$ is defined in
\eqref{bl:diagonal-moments} using the boundary Taylor
coefficients from~\eqref{bl:boundary-jets}.
Lemma~\ref{bl:frozen-moments} proves that this finite linear system
has only the zero solution. This determines the normal derivatives of the Lam\'e parameters through order $j$ and those of $W$ through order $j-2$, completing the induction.
In particular, the $q$th normal derivative of $W$ is recovered at
elastic order $j=q+2$. For each fixed order $j$, the induction uses a finite truncation
$N_{\mathrm{bl}}\ge j+1$, and the limit $h\to0$ is taken at that order.

Tangential differentiation of the recovered boundary identities gives every mixed derivative and hence the full boundary derivatives, because in each boundary chart every Cartesian differential operator $\partial_x^\beta$ is a finite sum of operators $\partial_y^\alpha\partial_r^q$ with smooth coefficients. Applying the construction on every compact patch and every boundary component proves the assertion. The density follows by division by $\omega^2$.
\end{proof}

\subsection{Boundary determination with finite regularity}

\begin{proposition}\label{bl:finite-jets}
Let $\Omega$ be bounded with $C^9$ boundary. Let two material triples
satisfy the strong convexity, positivity, and nonresonance hypotheses
of Theorem~\ref{a28:main}, with
$\lambda_j,\mu_j\in C^{16}(\overline\Omega)$ and
$W_j\in C^{15}(\overline\Omega)$, $j=1,2$.
If their physical Dirichlet-to-Neumann maps agree, then
\[
 \left.\partial_x^\alpha\delta_{12}\lambda\right|_{\partial\Omega}
 =\left.\partial_x^\alpha\delta_{12}\mu\right|_{\partial\Omega}=0
 \quad (|\alpha|\le16),\qquad
 \left.\partial_x^\alpha\delta_{12}W\right|_{\partial\Omega}=0
 \quad (|\alpha|\le14).
\]
\end{proposition}

\begin{proof}
Work in the coordinates~\eqref{bl:coordinates} on a relatively compact
boundary patch. The amplitudes and the normal cutoff $\eta$ are fixed.
Write $\Phi_{\partial}\in C^9$ and $n_{\mathrm{in}}\in C^8$. The identity
\[
 D_{(y,r)}X=
 \bigl(\partial_{y_1}\Phi_{\partial}+r\partial_{y_1}n_{\mathrm{in}},       \partial_{y_2}\Phi_{\partial}+r\partial_{y_2}n_{\mathrm{in}},       n_{\mathrm{in}}\bigr)
\]
gives the chart derivative counts for the first two rows
$(w^\alpha)^T$ of $(D_{(y,r)}X)^{-1}$. With $s$ tangential derivatives,
$w^\alpha(y,0)$ uses at most $s+1$ derivatives of $\Phi_{\partial}$.
A positive $r$-derivative uses $s+2$, since
$\partial_rD_{(y,r)}X$ contains $\partial_y n_{\mathrm{in}}$.
The terms $\mathcal X_i(n_{\mathrm{in},j})$ and
$\mathcal X_i(w_j^\alpha)$ use one additional derivative and carry a
factor of $h$ in~\eqref{bl:ordered}--\eqref{bl:scaled}. For the
coefficients in~\eqref{bl:Pq}, consider the chart factors of
$\partial_y^\alpha\partial_r^{q-m}A_{m,\beta}(y,0)$ with
$|\alpha|=s$. When $q=m$, these factors use at most $s+1$
derivatives of $\Phi_{\partial}$ for $m=0$ and $s+2$ for $m\ge1$.
When $q>m$, the corresponding bounds are $s+2$ and $s+3$.
In all cases the chart factors use at most $q+s+1$ derivatives
of $\Phi_{\partial}$, and the material factors use at most $q+s+1$
Cartesian derivatives. Thus
\begin{equation}\label{bl:finite-chart-count}
 \partial_r^{q-m}A_{m,\beta}(\cdot,0)\in C^s(V)
 \qquad\text{if }q+s+1\le9.
\end{equation}
For the Taylor expansion~\eqref{bl:remainder} with
$N_{\mathrm{bl}}=8$, the remainder uses
$\partial_r^{9-m}A_{m,\beta}$ on the collar. These derivatives exist
because $(D_{(y,r)}X)^{-1}$ is smooth in $r$ at fixed $y$ and the
material coefficients are differentiated along straight normal segments.

Construct $U_0=F_{\mathfrak{s}}A$ and apply~\eqref{bl:recursion} through
$\ell=8$. Consider $s$ tangential derivatives of a term of
$\mathcal P_{\mathfrak{s}}^{(q)}U_{\ell-q}$, with $\sigma$ derivatives
falling on its coefficient. By~\eqref{bl:tangential-order}, the
derivative orders satisfy
\begin{equation}\label{bl:finite-profile-count}
 q+\sigma+1\le \ell+s+1\le9,\qquad
 (\ell-q)+(s-\sigma+\beta_1+\beta_2)
 \le\ell+s\le8 .
\end{equation}
The first inequality gives the coefficient regularity
by~\eqref{bl:finite-chart-count}, and the second bounds the tangential
derivative order on the preceding profile. The frozen
inverse of Lemma~\ref{bl:frozen} acts in $s_{\mathrm{bl}}$ and depends
smoothly on the frozen parameters. By induction, $U_\ell$ has
continuous tangential derivatives through order $8-\ell$, with the
polynomial--exponential bounds and tangential support of the
construction in Proposition~\ref{bl:layer}.

Extend the profiles by zero in $y$ from a compact patch with a fixed
margin in $V$. For $\ell=7,8$, let $U_{\ell,h}$ be the convolution
in $y$ with an even smooth compactly supported unit-mass kernel at
scale $h$. Put $U_{\ell,h}=U_\ell$ for $0\le\ell\le6$. These
profiles have trace zero at $s_{\mathrm{bl}}=0$ for $\ell\ge1$.
Uniformly in a supremum norm with a fixed exponential weight
in $s_{\mathrm{bl}}$ smaller than the profile decay rate,
\begin{equation}\label{bl:finite-mollification}
 \begin{aligned}
  \partial_y^\alpha(U_{\ell,h}-U_\ell)
     &=o(h^{8-\ell-|\alpha|}) && (|\alpha|\le8-\ell),\\
  \partial_y^\alpha U_{\ell,h}
     &=o(h^{8-\ell-|\alpha|}) && (8-\ell<|\alpha|\le2),
 \end{aligned}
 \qquad \ell=7,8 .
\end{equation}
The same estimates hold after up to two $s_{\mathrm{bl}}$ derivatives.
The estimate $U_{7,h}-U_7=o(h)$ uses the vanishing first moment of
the mollifier and uniform continuity of $\partial_yU_7$. For the
derivative estimates on continuous functions, differentiate the kernel
and subtract the value at the center, using its zero integral. Apply
this argument
to $\partial_y U_7$ for its second-derivative bound. Uniform continuity
in the weighted norm follows from the polynomial--exponential form
on the compact patch.

Define
\[
 v_h(X(y,r))=e^{i\mathfrak{s}y_1/h}\eta(r)
               \sum_{\ell=0}^8 h^\ell U_{\ell,h}(y,s_{\mathrm{bl}}).
\]
In the expansions~\eqref{bl:Pq} and~\eqref{bl:remainder}, the terms
of total degree $\ell+q\le8$ for the unmollified profiles exist by
$\beta_1+\beta_2\le m\le q$. Their coefficients vanish by
\eqref{bl:recursion} and $\mathcal P_{\mathfrak{s}}^{(0)}U_0=0$.
For $\ell=7,8$ and $|\beta|\le2$, set $s=\beta_1+\beta_2$.
By~\eqref{bl:tangential-order}, the power $h^{\ell+q}$ of
the term has $q\ge s$. Applying~\eqref{bl:finite-mollification}
in the same weighted norm gives
\[
 \begin{aligned}
  h^{\ell+q}\partial_{(y,s_{\mathrm{bl}})}^\beta
      (U_{\ell,h}-U_\ell)
     &=o(h^8) && (\ell+s\le8),\\
  h^{\ell+q}\partial_{(y,s_{\mathrm{bl}})}^\beta U_{\ell,h}
     &=o(h^8) && (\ell+s>8).
 \end{aligned}
\]
The higher terms using existing derivatives are $O(h^9)$ with a
polynomial--exponential envelope.
The Taylor remainder in~\eqref{bl:remainder} is bounded by $O(h^9)$
with the same type of envelope. In fact, for $|\beta|\le2$,
\[
 h^\ell\partial_{(y,s_{\mathrm{bl}})}^\beta U_{\ell,h}=O(1)
\]
in the weighted supremum norm, by~\eqref{bl:finite-profile-count}
and~\eqref{bl:finite-mollification}, while the remainder coefficients
grow at most polynomially in $s_{\mathrm{bl}}$.
Cutoff commutators are supported where $r$ is bounded away from zero
and have exponential decay in $h$. After decreasing the decay
exponent $c>0$ to absorb the polynomial factors and undoing the
$h^2$ scaling in~\eqref{bl:scaled}, we obtain
\begin{equation}\label{bl:finite-layer-residual}
 |Pv_h(X(y,r))|\le o(h^6)e^{-cr/h},\qquad
 \|rPv_h\|_{L^2}^2\le o(h^{12})\int_0^{r_0}r^2e^{-2cr/h}\,dr
 =o(h^{15}).
\end{equation}
The factor $J(y,r)$ and the area of the compact support patch are
included in the second bound. The prescribed trace follows from the
profile traces. Formula~\eqref{bl:D} and~\eqref{bl:finite-mollification}
give $O(h)$ for all scaled displacement and gradient contributions
after the leading profile, with an exponential envelope. This
proves~\eqref{bl:profile-error} and
\[
 \|v_h\|_{L^2(\Omega)}=O(h^{1/2}),\qquad
 \|v_h\|_{H^1(\Omega)}=O(h^{-1/2})
\]
as in~\eqref{bl:appnorm}.

For $\phi\in H_0^1(\Omega;\mathbb C^3)$, the normal Hardy inequality
and the positive upper and lower Jacobian bounds give
$\|\phi/r\|_{L^2}\le C\|\phi\|_{H^1(\Omega)}$ on the collar.
Using~\eqref{bl:finite-layer-residual} in the dual pairing,
\[
 |\langle Pv_h,\phi\rangle|
 \le \|rPv_h\|_{L^2}\|\phi/r\|_{L^2},
 \qquad
 \|Pv_h\|_{H^{-1}(\Omega)}=o(h^{15/2}).
\]
The physical Dirichlet inverse now gives
\begin{equation}\label{bl:finite-exact-layer}
 u_h=v_h-P_D^{-1}(Pv_h),\qquad
 Pu_h=0,\qquad u_h-v_h\in H_0^1(\Omega;\mathbb C^3),\qquad
 \|u_h-v_h\|_{H^1(\Omega)}=o(h^{15/2}).
\end{equation}

For either bounded weak pairing used at elastic orders zero and one,
replacing the exact fields by the approximate layers changes the
unscaled pairing by at most
$C(o(h^{15/2})O(h^{-1/2})+o(h^{15}))=o(h^7)$.
Thus the leading limits used in~\eqref{bl:Z} and~\eqref{bl:H1} apply
with these layers. The frozen boundary matrices are continuous on
the patch. They give
\[
 \delta_{12}\lambda(y,0)=\delta_{12}\mu(y,0)=0,\qquad
 L_1^{\mathrm{bd}}(y)=M_1^{\mathrm{bd}}(y)=0.
\]

Let $2\le j\le16$. Assume that the normal traces of the Lam\'e
differences vanish below order $j$ and those of $\delta_{12}W$ below
order $j-2$. Taylor expansion on the normal segments gives, with
remainders uniform on the compact patch,
\begin{equation}\label{bl:finite-jet-expansion}
\begin{aligned}
 \delta_{12}\lambda(y,r)&=\frac{r^j}{j!}L_j^{\mathrm{bd}}(y)+o(r^j),\\
 \delta_{12}\mu(y,r)&=\frac{r^j}{j!}M_j^{\mathrm{bd}}(y)+o(r^j),\\
 \delta_{12}W(y,r)&=\frac{r^{j-2}}{(j-2)!}V_j^{\mathrm{bd}}(y)+o(r^{j-2}).
\end{aligned}
\end{equation}
Use the layers $v_{1,h},v_{2,h}$ for the two materials with opposite
phase signs. By~\eqref{bl:finite-jet-expansion},
\eqref{bl:profile-error} and $r=hs_{\mathrm{bl}}$, dominated
convergence gives
\begin{equation}\label{bl:finite-moment-limit}
 \lim_{h\to0}h^{1-j}\mathfrak B(v_{1,h},v_{2,h})=\mathfrak M_j .
\end{equation}
The moment on the right is~\eqref{bl:leading-moment}. The same bounds
give, for $i=1,2$,
\[
 \|\delta_{12}\lambda\,\nabla v_{i,h}\|_{L^2}
 +\|\delta_{12}\mu\,\nabla v_{i,h}\|_{L^2}=O(h^{j-1/2}),\qquad
 \|\delta_{12}W\,v_{i,h}\|_{L^2}=O(h^{j-3/2}).
\]
With~\eqref{bl:finite-exact-layer}, the terms with one correction in
the scaled elastic and potential pairings are bounded respectively by
$h^{1-j}O(h^{j-1/2})o(h^{15/2})=o(h^8)$ and
$h^{1-j}O(h^{j-3/2})o(h^{15/2})=o(h^7)$. Boundedness of
$\mathfrak B$ on $H^1\times H^1$ gives the two-correction bound
$h^{1-j}o(h^{15})=o(h^{16-j})$. Consequently
\begin{equation}\label{bl:finite-pairing-error}
 h^{1-j}\bigl[
 \mathfrak B(u_{1,h},u_{2,h})-\mathfrak B(v_{1,h},v_{2,h})
 \bigr]=o(h^8)+o(h^7)+o(h^{16-j}) .
\end{equation}
This tends to zero also at $j=16$. By
\eqref{bl:alessandrini}, \eqref{bl:finite-moment-limit}
and~\eqref{bl:finite-pairing-error}, $\mathfrak M_j=0$ for every
fixed pair of amplitudes.

Fix a boundary point $y_0$ and choose a constant orthonormal frame
with $e_1=e_{\mathrm{bl}}(y_0)$ and $e_3=n_{\mathrm{in}}(y_0)$.
For $A_+=A_-=\chi e_p$, with $\chi\in C_c^\infty(V;\mathbb R)$,
the identity $\mathfrak M_j=0$ tests the continuous diagonal moment
against $\chi^2J(y,0)$. Since $J(y,0)>0$, localization as in the proof
of Theorem~\ref{bl:jets} gives vanishing of the moment at $y_0$. Set
$\widetilde s_{\mathrm{bl}}=c_{\mathrm{bl}}(y_0)s_{\mathrm{bl}}$.
Then the elastic weights, including the two leading gradients and
the measure, scale by $c_{\mathrm{bl}}(y_0)^{2-j-1}$; the potential
weight and measure scale by $c_{\mathrm{bl}}(y_0)^{-(j-2)-1}$.
These are the same positive factor. After division, the moments
are those in~\eqref{bl:diagonal-moments}, so
\[
 E_{11}=E_{22}=E_{33}=0,\qquad
 L_j^{\mathrm{bd}}=M_j^{\mathrm{bd}}=V_j^{\mathrm{bd}}=0
\]
by Lemma~\ref{bl:frozen-moments}. This proves the vanishing normal
traces through orders $16,16,14$ for $\lambda,\mu,W$, respectively.

Finally, write $f$ for any one of the three coefficient differences.
Suppose that the Cartesian derivatives through order $q-1$ have
zero trace. Differentiating the order-$q-1$ tensor along a boundary
curve $\gamma$ with $\gamma(0)=x$ and $\gamma'(0)=\tau$ gives, for
constant vectors $\xi_1,\ldots,\xi_{q-1}$,
\[
 D_x^q f(x)[\tau,\xi_1,\ldots,\xi_{q-1}]
 =\left.\frac{d}{da}
 D_x^{q-1}f(\gamma(a))[\xi_1,\ldots,\xi_{q-1}]\right|_{a=0}=0 .
\]
This tangential chain rule follows by applying the fundamental
theorem of calculus along parallel interior curves in a graph chart
and taking their limit, using continuity of the derivative traces.
By symmetry only the pure normal component remains. The affine
dependence on $r$ in~\eqref{bl:coordinates} gives
\[
 D_x^q f(x)[n_{\mathrm{in}},\ldots,n_{\mathrm{in}}]
   =\partial_r^q f(y,0)=0,
 \qquad x=X(y,0).
\]
Induction on $q$ and a finite cover of $\partial\Omega$ complete
the proof.
\end{proof}

\subsection{Common exterior extension for smooth inputs}\label{sec:extension}

By Theorem~\ref{bl:jets}, the two coefficient triples
have identical boundary Taylor data. In every inward normal
coordinate chart,
\[
\begin{aligned}
\partial_y^\alpha\partial_r^q(\lambda_1-\lambda_2)(y,0)
&=0,\\
\partial_y^\alpha\partial_r^q(\mu_1-\mu_2)(y,0)
&=0,\\
\partial_y^\alpha\partial_r^q(W_1-W_2)(y,0)
&=0,
\qquad
\alpha\in\mathbb N_0^2,\quad q\in\mathbb N_0.
\end{aligned}
\]
We use these identities to construct smooth extensions
that agree outside $\overline\Omega$. Common exterior
extensions are also used in inverse boundary problems;
see~\cite[Section~6]{DSFKSU2009}. We first establish a
zero-extension result for functions whose derivatives of all orders vanish on the boundary.

\begin{lemma}
\label{ext:zero-extension}
Let $f\in C^\infty(\overline\Omega)$ satisfy
\[
\left.\partial^\alpha f\right|_{\partial\Omega}=0
\qquad\text{for every }\alpha\in\mathbb N_0^3.
\]
Define
\[
(\mathcal E_0f)(x)=
\begin{cases}
f(x),&x\in\Omega,\\
0,&x\notin\Omega.
\end{cases}
\]
Then
\[
\mathcal E_0f\in C_c^\infty(\mathbb R^3),
\qquad
\operatorname{supp}(\mathcal E_0f)
\subset\overline\Omega.
\]
\end{lemma}

\begin{proof}
In each chart~\eqref{bl:coordinates}, Taylor's formula
in the normal variable gives, for every multi-index
$\alpha$ and integer $N\ge1$,
\[
(\partial^\alpha f)(X(y,r))
=\frac{r^N}{(N-1)!}
\int_0^1(1-s)^{N-1}
\left.\partial_\rho^N
\big[(\partial^\alpha f)(X(y,\rho))\big]\right|_{\rho=sr}\,ds.
\]
Consequently, on every compact subpatch,
\[
|(\partial^\alpha f)(X(y,r))|
\le C_{\alpha,N}r^N.
\]
Thus every derivative has zero trace on the boundary.
Integration by parts and induction on $|\alpha|$ yield
\[
\partial^\alpha(\mathcal E_0f)
=\mathcal E_0(\partial^\alpha f)
\quad\text{in }\mathcal D'(\mathbb R^3).
\]
The right-hand side is continuous across the boundary
for every $\alpha$, so $\mathcal E_0f$ is smooth.
A finite chart cover handles all boundary components,
and boundedness of $\Omega$ gives compact support.
\end{proof}

\begin{proposition}
\label{ext:common}
Let $(\lambda_j,\mu_j,W_j)$, $j=1,2$, be the two smooth
coefficient triples satisfying the standing strong-convexity
assumptions. Suppose that they have identical boundary
Taylor data on every boundary component. Then they admit
smooth extensions to $\mathbb R^3$ that agree on
$\mathbb R^3\setminus\overline\Omega$, preserve strong
convexity of the Lam\'e pairs, and equal the same constant
strongly convex Lam\'e pair with zero potentials outside
a compact set.
\end{proposition}

\begin{proof}
We first construct a smooth extension operator using the
classical reflection construction underlying Seeley's
theorem~\cite{Seeley1964}. For extension operators with
stronger mapping properties, see~\cite{LuYao2024}.

Put $b_j=2^j$ for $j\ge0$. For $0\le j\le N$, define
\begin{equation}\label{ext:reflection-coefficients}
a_{j,N}
=\prod_{\substack{0\le k\le N\\k\ne j}}
\frac{1+b_k}{b_k-b_j},
\qquad a_{j,N}=0\quad(j>N).
\end{equation}
Lagrange interpolation at $-b_0,\ldots,-b_N$,
evaluated at $1$, gives
\[
\sum_{j=0}^N a_{j,N}(-b_j)^p=1,
\qquad 0\le p\le N.
\]

For each fixed $j\ge0$, define
\[
a_j:=\lim_{N\to\infty}a_{j,N}.
\]
The limit exists because the factors with $k>j$ are
\[
\frac{1+2^{-k}}{1-2^{j-k}}.
\]
Splitting the product at $j$ gives
\[
|a_{j,N}|\le C2^{-j(j+1)/2},
\qquad
C=
\prod_{k\ge0}(1+2^{-k})
\prod_{\ell\ge1}(1-2^{-\ell})^{-2}
<\infty.
\]
For every fixed integer $p\ge0$, the bound is summable
after multiplication by $b_j^p$. Dominated convergence
therefore gives
\begin{equation}\label{ext:reflection-moments}
\sum_{j\ge0}a_j(-b_j)^p=1,
\qquad
\sum_{j\ge0}|a_j|b_j^p<\infty
\quad(p\ge0).
\end{equation}

Use the inward coordinates~\eqref{bl:coordinates}.
Choose a finite smooth partition of unity near
$\overline\Omega$, with its interior term compactly
supported in $\Omega$ and each boundary term $\eta$
supported strictly inside a smaller coordinate box.
For a boundary term, write
\[
F(y,r)=(\eta f)(X(y,r)),\qquad r\ge0.
\]
Choose a fixed smooth normal cutoff $\chi_{\mathrm n}$,
equal to one near $0$ and supported in $[0,\varepsilon)$
on the nonnegative half-axis, where $\varepsilon$ is
smaller than the collar width. Set
\[
g(y,s)=\chi_{\mathrm n}(s)F(y,s),\qquad s\ge0,
\]
and extend $g$ by zero for $s\ge\varepsilon$. Define
\begin{equation}\label{ext:reflection-extension}
(E_{\rm loc}F)(y,r)=
\begin{cases}
F(y,r),&r\ge0,\\[2mm]
\displaystyle\sum_{j\ge0}a_jg(y,-b_jr),&r<0.
\end{cases}
\end{equation}

For each fixed $r<0$, the sum is finite. By
\eqref{ext:reflection-moments}, every differentiated series
converges uniformly up to $r=0$ on compact tangential
patches, and
\[
\begin{aligned}
\lim_{r\uparrow0}
\partial_y^\alpha\partial_r^p(E_{\rm loc}F)(y,r)
&=\sum_{j\ge0}a_j(-b_j)^p
\partial_y^\alpha\partial_r^pF(y,0)\\
&=\partial_y^\alpha\partial_r^pF(y,0).
\end{aligned}
\]
Thus the two pieces join smoothly. Since $b_j\ge1$,
the exterior piece vanishes for $-r\ge\varepsilon$.
The fixed tangential supports and a larger coordinate box
permit smooth extension by zero at all chart edges.

All partitions and cutoffs are fixed independently of
$f$ and of the derivative order. Adding the finitely many
chartwise extensions and the zero extension of the interior
term gives a real linear operator
\begin{equation}\label{ext:smooth-operator}
E:C^\infty(\overline\Omega;\mathbb R)
\longrightarrow C^\infty(\mathbb R^3;\mathbb R),
\qquad
(Ef)|_{\overline\Omega}=f.
\end{equation}

Put
\[
\lambda_1^e=E\lambda_1,\qquad
\mu_1^e=E\mu_1,\qquad
W_1^e=EW_1.
\]
Strong convexity and compactness of $\overline\Omega$
give a relatively compact open neighborhood
$U\supset\overline\Omega$ and a constant $c_*>0$ such that
\[
\mu_1^e\ge c_*,
\qquad
3\lambda_1^e+2\mu_1^e\ge c_*
\quad\text{on }U.
\]
Choose $\theta\in C_c^\infty(U)$ with
$0\le\theta\le1$ and $\theta=1$ near
$\overline\Omega$. Choose constants $\lambda_*,\mu_*$
such that
\[
\mu_*>0,\qquad 3\lambda_*+2\mu_*>0.
\]
Define
\begin{equation}\label{ext:first-extension}
\begin{aligned}
\widetilde\lambda_1
 &=\theta\lambda_1^e+(1-\theta)\lambda_*,\\
\widetilde\mu_1
 &=\theta\mu_1^e+(1-\theta)\mu_*,\\
\widetilde W_1&=\theta W_1^e.
\end{aligned}
\end{equation}
These functions agree with the original coefficients on
$\overline\Omega$. Moreover,
\[
\begin{aligned}
\widetilde\mu_1
&\ge\min(c_*,\mu_*)>0,\\
3\widetilde\lambda_1+2\widetilde\mu_1
&\ge\min(c_*,3\lambda_*+2\mu_*)>0.
\end{aligned}
\]
Thus the cutoff interpolation preserves strong convexity,
and
\[
(\widetilde\lambda_1,\widetilde\mu_1,\widetilde W_1)
=(\lambda_*,\mu_*,0)
\quad\text{outside }\operatorname{supp}\theta.
\]

The matching boundary Taylor data and the smooth collar coordinates give
\[
 \left.\partial^\alpha(f_2-f_1)\right|_{\partial\Omega}=0,
 \qquad f\in\{\lambda,\mu,W\},\quad \alpha\in\mathbb N_0^3.
\]
Lemma~\ref{ext:zero-extension} therefore permits the definition
\begin{equation}\label{ext:second-extension}
\widetilde f_2
=\widetilde f_1+\mathcal E_0(f_2-f_1).
\end{equation}
Consequently,
\[
\widetilde f_j|_{\overline\Omega}=f_j,\qquad j=1,2,
\]
and
\[
\widetilde f_2=\widetilde f_1
\quad\text{on }\mathbb R^3\setminus\Omega,
\qquad
\operatorname{supp}(\widetilde f_2-\widetilde f_1)
\subset\overline\Omega.
\]
The second Lam\'e pair is strongly convex in $\Omega$
by hypothesis, and outside $\Omega$ it agrees with the
strongly convex first extension. Smoothness across the
boundary follows from~\eqref{ext:second-extension}.
Hence both extended triples have the same exterior values
and equal $(\lambda_*,\mu_*,0)$ outside a compact set.
Their differences, together with all derivatives of those
differences, are supported in $\overline\Omega$. Thus
integrals involving these differences may be evaluated on
any larger construction domain.
\end{proof}

In the subsequent auxiliary constructions, we use
$\lambda_j,\mu_j,W_j$ also for the smooth extensions
$\widetilde\lambda_j,\widetilde\mu_j,\widetilde W_j$ from
Proposition~\ref{ext:common}; their restrictions to
$\Omega$ are the original physical coefficients.
All boundary maps and Dirichlet problems continue to
refer to the original domain $\Omega$.

\subsection{Finite-order exterior extensions}
\label{sec:finite-extension}

\begin{proposition}\label{ext:finite-common}
Let $\partial\Omega$ be of class $C^9$, let
$\lambda_j,\mu_j\in C^{16}(\overline\Omega;\mathbb R)$ satisfy
$\mu_j>0$ and $3\lambda_j+2\mu_j>0$ on $\overline\Omega$, and let
$W_j\in C^{15}(\overline\Omega;\mathbb R)$ for $j=1,2$.
Suppose that the Cartesian derivative traces of $\lambda_1-\lambda_2$
and $\mu_1-\mu_2$ vanish through order $16$. Then the Lam\'e pairs
have $C^{16}(\mathbb R^3)$ extensions that agree on
$\mathbb R^3\setminus\overline\Omega$, are strongly convex
everywhere, and equal the same strongly convex constant pair
$(\lambda_*,\mu_*)$ outside a compact set. Each $W_j$ has a
$C_c^{15}(\mathbb R^3)$ extension.
\end{proposition}

\begin{proof}
Fix $f\in C^m(\overline\Omega;\mathbb R)$, where $m=15$ or $16$.
Write $P_x^mf$ for its Cartesian Taylor polynomial at $x$:
\[
 P_x^mf(z)=\sum_{|\beta|\le m}
             \frac{\partial^\beta f(x)}{\beta!}(z-x)^\beta .
\]
We verify the compatibility of the corresponding Whitney data. In a
Lipschitz graph chart $\Omega=\{x_3>\varphi(y)\}$, two sufficiently
close interior points $x=(y,x_3)$ and $x'=(y',x'_3)$ can be joined
by vertical segments to the height
$\max\{x_3,x'_3\}+(1+\operatorname{Lip}\varphi)|x-x'|$ and
a horizontal segment at that height. For points in a smaller chart
and sufficiently close to each other the segments remain inside
the chart. A finite cover of the boundary by such charts, and
straight segments away from the boundary, give a polygonal path
$\gamma$ in $\Omega$ with
\begin{equation}\label{ext:finite-short-path}
 \gamma(0)=x,\qquad \gamma(l)=x',\qquad
 l\le C|x-x'|,\qquad
 |\gamma(\tau)-x|\le C|x-x'|\quad(0\le\tau\le l),
\end{equation}
where $\tau$ is Euclidean arclength and $C$ is uniform for close pairs.

Set $g_x=f-P_x^mf$, and let $\eta_m$ be the common modulus of
continuity on $\overline\Omega$ of the order-$m$ derivatives of $f$.
For $|\alpha|\le m$ the derivatives of $g_x$ vanish at $x$, and at
order $m$ we have
\[
 \partial^\alpha g_x(\gamma(\tau))
   =\partial^\alpha f(\gamma(\tau))-\partial^\alpha f(x),
   \qquad |\alpha|=m.
\]
For $|\alpha|<m$, integration on the segments of $\gamma$ gives
\[
 \partial^\alpha g_x(\gamma(\tau))
   =\int_0^\tau
      \nabla\partial^\alpha g_x(\gamma(v))\cdot\dot\gamma(v)\,dv .
\]
Applying this identity successively from order $m-1$ downwards and
using \eqref{ext:finite-short-path} gives the uniform remainder bounds
\begin{equation}\label{ext:finite-whitney-jets}
 \bigl|\partial^\alpha f(x')-\partial^\alpha P_x^mf(x')\bigr|
 \le C|x-x'|^{m-|\alpha|}\eta_m(C|x-x'|)
 =o\bigl(|x-x'|^{m-|\alpha|}\bigr),
 \qquad |\alpha|\le m.
\end{equation}
The same bounds hold for endpoints in $\overline\Omega$ by continuity,
enlarging $C$ if necessary. By \eqref{ext:finite-whitney-jets},
Whitney's extension theorem~\cite[Theorem~I]{Whitney1934} gives
a real extension $f^e\in C^m(\mathbb R^3)$ with the prescribed
Cartesian derivatives through order $m$ on $\overline\Omega$.

For $f=\lambda_2-\lambda_1$ or $f=\mu_2-\mu_1$, use the zero
extension $\mathcal E_0f$ defined in Lemma~\ref{ext:zero-extension}.
If $|\alpha|<16$ and $\phi\in C_c^\infty(\mathbb R^3)$, the trace
condition on $\partial^\alpha f$ gives
\[
 \bigl\langle
   \partial_i\mathcal E_0(\partial^\alpha f),\phi\bigr\rangle
 =-\int_\Omega\partial^\alpha f\,\partial_i\phi\,dx
 =\int_\Omega\partial^{\alpha+e_i} f\,\phi\,dx .
\]
Induction on $|\alpha|$ yields
\begin{equation}\label{ext:finite-zero-derivatives}
 \partial^\alpha(\mathcal E_0f)
      =\mathcal E_0(\partial^\alpha f)
      \quad\text{in }\mathcal D'(\mathbb R^3),\qquad |\alpha|\le16 .
\end{equation}
Each right-hand side in \eqref{ext:finite-zero-derivatives} is
continuous across $\partial\Omega$ by the corresponding vanishing
trace. Thus $\mathcal E_0f\in C^{16}(\mathbb R^3)$ and
$\operatorname{supp}\mathcal E_0f\subset\overline\Omega$.

Apply the extension construction to $\lambda_1,\mu_1$. By compactness
and continuity, there are an open neighborhood $U$ of
$\overline\Omega$ and $c_*>0$ such that
\[
 \mu_1^e\ge c_*,\qquad 3\lambda_1^e+2\mu_1^e\ge c_*
 \quad\text{on }U.
\]
Take $\theta\in C_c^\infty(U)$ with $0\le\theta\le1$ and
$\theta=1$ near $\overline\Omega$, and choose a strongly convex
constant pair $(\lambda_*,\mu_*)$. Define
\begin{equation}\label{ext:finite-lame-pairs}
 \begin{aligned}
 \widetilde\lambda_1&=\theta\lambda_1^e+(1-\theta)\lambda_*,
 &\widetilde\mu_1&=\theta\mu_1^e+(1-\theta)\mu_*,\\
 \widetilde\lambda_2&=\widetilde\lambda_1+\mathcal E_0(\lambda_2-\lambda_1),
 &\widetilde\mu_2&=\widetilde\mu_1+\mathcal E_0(\mu_2-\mu_1).
 \end{aligned}
\end{equation}
The first pair in \eqref{ext:finite-lame-pairs} satisfies
\[
 \widetilde\mu_1\ge\min\{c_*,\mu_*\}>0,\qquad
 3\widetilde\lambda_1+2\widetilde\mu_1
       \ge\min\{c_*,3\lambda_*+2\mu_*\}>0 .
\]
By \eqref{ext:finite-zero-derivatives}, the second pair is $C^{16}$.
It equals $(\lambda_2,\mu_2)$ on $\overline\Omega$ and the first
pair on $\mathbb R^3\setminus\overline\Omega$, so it is strongly
convex everywhere. Both pairs equal $(\lambda_*,\mu_*)$ outside
$\operatorname{supp}\theta$.

Finally extend each $W_j$ individually to $W_j^e$ by
\eqref{ext:finite-whitney-jets} and Whitney's theorem. Set
$\widetilde W_j=\theta W_j^e$. Then
$\widetilde W_j|_{\overline\Omega}=W_j$ and
$\operatorname{supp}\widetilde W_j\subset\operatorname{supp}\theta$,
so each $\widetilde W_j$ belongs to $C_c^{15}(\mathbb R^3)$.
\end{proof}

For these extensions the coefficients of $M_j^{W_j}-\Delta I_4$
in~\eqref{bt:MW} are of class $C^{14}$. The first-order
coefficients depend only on the Lam\'e pair and thus have the
same exterior values for the two media. The multiplication
matrices can differ in the exterior of $\overline\Omega$.
The affine background conjugation preceding Proposition~\ref{pb27:rightinverse}
uses the common constant pair and the vanishing of both
potentials sufficiently far away.

\begin{remark}\label{ext:finite-betti}
The reduced identity~\eqref{bt:betti} is applied on $\Omega$
itself. The boundary matching needed for it is
$\mu_1=\mu_2$, $\partial_\nu\mu_1=\partial_\nu\mu_2$, and
$W_1=W_2$ on $\partial\Omega$, as in~\eqref{bt:matching}.
Its derivation remains valid for the finite-order coefficients
just described. Indeed, the factorization~\eqref{bt:factorization}
is a differential identity for coefficients with the indicated
finite derivatives. The additional boundary cancellation uses
the tangential derivative of the trace and the normal-normal
Hessian form of Lemma~\ref{bt:hessian}. This only requires a
$C^2$ boundary. For application to the constructed solutions,
fix the CGO parameter and consider an exact reduced solution
in $H^1$ on an open neighborhood of $\overline\Omega$.
Its equation expresses its Laplacian in terms of first and
zeroth derivatives multiplied by $C^{14}$ coefficients.
Interior elliptic regularity and iteration give $H^8_{\rm loc}$
there, supplying the classical derivatives required by the
physical lifts and the boundary integration by parts. Each
lift is therefore a weak physical solution on $\Omega$, and
its traction agrees with the weak Dirichlet-to-Neumann traction.
The bootstrap takes place in the open neighborhood of
$\overline\Omega$.

The compact support of the first-order difference allows the
leading moment integral to be extended after its limit is
taken on $\Omega$. The exact integral involving $D^W$ is kept
on $\Omega$ while taking that limit.
\end{remark}

\section{The reduced Green--Betti identity}
\label{bt:section}
In this section, we derive a Green--Betti identity for the four-component
system~\eqref{bt:MW} using only the Dirichlet-to-Neumann
maps. We first compute the boundary terms introduced by the
reduction, including the contribution of the auxiliary scalar
component, and then prove the cancellation of the corresponding
boundary terms under the matching conditions~\eqref{bt:matching},
which follow from the boundary determination in
Section~\ref{bl:section}. The resulting identity~\eqref{bt:betti}
expresses an interior pairing involving the difference of the
reduced operators in terms of the difference of the
boundary maps. The identity is applied on the original domain
$\Omega$. The matrix-moment argument uses the exterior coefficient
identity~\eqref{bt:status}.

For the reduced operators $M_j^{W_j}$ defined
in~\eqref{bt:MW}, set
\begin{equation}\label{bt:operator-difference}
 D^W:=M_1^{W_1}-M_2^{W_2}.
\end{equation}
The order of subtraction in~\eqref{bt:operator-difference}
fixes the signs in the reduced Green--Betti
identity~\eqref{bt:defects}. Throughout the section, $\nu$
denotes the outward unit normal to $\partial\Omega$, and
we use the bilinear pairing and transpose convention
fixed in Section~\ref{sec:notation}.

\subsection{The four-component reduction}
\label{bt:reduction}

Let $\Omega\subset\mathbb R^3$ be a bounded domain with smooth
boundary. For $j=1,2$, consider smooth real coefficient triples
$(\lambda_j,\mu_j,W_j)$ on $\overline\Omega$ satisfying
$\mu_j>0$ and $3\lambda_j+2\mu_j>0$. Set
\begin{equation}
 \Sigma_j=\lambda_j+2\mu_j,\quad
 \kappa_j=\frac{\lambda_j+\mu_j}{\Sigma_j},\quad
 a_j=\mu_j^{-1},\quad c_j=\mu_j^{-1/2},\quad
 \ell_j=\log\mu_j,\quad \mathbf m_j=\nabla\ell_j.
 \label{bt:notation}
\end{equation}
The strong convexity conditions imply
$\Sigma_j>0$, $1/4<\kappa_j<1$, and
$1-\kappa_j=\mu_j/\Sigma_j$. Equivalently,
$\lambda_j=\mu_j(2\kappa_j-1)/(1-\kappa_j)$.
Let $\nu$ denote the outward unit normal to $\partial\Omega$.
With the convention $(\nabla u)_{ik}=\partial_k u_i$,
define the stress, elasticity operator, and traction by
\begin{equation}
 \begin{aligned}
  \sigma_j(u)
  &=2\mu_j\varepsilon(u)
    +\lambda_j\operatorname{div}u\,I_3,\\
  L_ju
  &:=L_{\lambda_j,\mu_j}u
   =\operatorname{div}\sigma_j(u),\qquad
  t_j(u)=\sigma_j(u)\nu.
 \end{aligned}
 \label{bt:stress}
\end{equation}
The full elastic operator associated with the coefficient
triple $(\lambda_j,\mu_j,W_j)$ is $L_j+W_jI_3$.

For three-component vector fields $v,p$ and scalar fields $f,s$,
write $z=(v,f)^T$ and introduce the displacement map $T_j$
and the source map $S_j$:
\begin{equation}
 T_jz=c_jv+a_j\nabla f-f\nabla a_j,\qquad
 S_j(p,s)=c_j^{-1}p+
          \nabla\!\left(\frac{s}{1-\kappa_j}\right).
 \label{bt:TS}
\end{equation}
For $j=1,2$, define the coefficients
and the four-component reduced operator by
\begin{align}
 F_{\mu_j}
 &=2c_j\bigl\{
     \nabla^2\ell_j-\mathbf m_j\otimes\mathbf m_j
     +(|\mathbf m_j|^2-\Delta\ell_j)I_3
   \bigr\},
 \notag\\
 V_{\mu_j}
 &=-\nabla^2\ell_j
   -\tfrac12\mathbf m_j\otimes\mathbf m_j
   -\bigl(
      \tfrac12\Delta\ell_j
      +\tfrac14|\mathbf m_j|^2
    \bigr)I_3,
 \notag\\
 q_{\mu_j}
 &=\Delta\ell_j-|\mathbf m_j|^2,\qquad
 d_j=\kappa_j\sqrt{\mu_j},\qquad
 e_j^{\mathrm{red}}=\tfrac12(3\kappa_j-2)\sqrt{\mu_j},
 \notag\\
 M_j^0(v,f)
 &=\begin{pmatrix}
 \Delta v-\mathbf m_j\operatorname{div}v
       +F_{\mu_j}\nabla f+V_{\mu_j}v
       +F_{\mu_j}\mathbf m_j f\\
 \Delta f+d_j\operatorname{div}v
       -e_j^{\mathrm{red}}\,\mathbf m_j\cdot v+q_{\mu_j}f
 \end{pmatrix},
 \notag\\
 \mathcal W_j
 &=\begin{pmatrix}
 (W_j/\mu_j)I_3
 &\mu_j^{-3/2}(2W_j\mathbf m_j-\nabla W_j)\\
 0&(1-\kappa_j)W_j/\mu_j
 \end{pmatrix},
 \qquad
 M_j^{W_j}=M_j^0+\mathcal W_j.
 \label{bt:MW}
\end{align}
The operator $M_j^{W_j}$ has principal part $\Delta I_4$.
The potential $W_j$ and the gradient $\nabla W_j$ enter only
through the $4\times4$ zeroth-order multiplication matrix
$\mathcal W_j$.
For every integer $k\ge0$, smoothness of the coefficients
gives the bounded mappings
\[
 \begin{aligned}
 T_j,S_j&:H^{k+1}(\Omega;\mathbb C^4)
             \longrightarrow H^k(\Omega;\mathbb C^3),\\
 M_j^{W_j}&:H^{k+2}(\Omega;\mathbb C^4)
             \longrightarrow H^k(\Omega;\mathbb C^4).
 \end{aligned}
\]

Lemma~\ref{bt:scalarlemma} establishes the
factorization~\eqref{bt:factorization}, allowing the CGO
solutions constructed for $M^W$ to generate  elastic
displacements through $T$. The scalar-row
identity~\eqref{bt:scalar} identifies the fourth-component
boundary contribution in~\eqref{bt:defects}, including the
density-dependent term. Under the matching
conditions~\eqref{bt:matching}, the fourth-component contribution
cancels the traction discrepancy, yielding the reduced
Green--Betti identity~\eqref{bt:betti} in terms of the 
Dirichlet-to-Neumann maps.

\begin{lemma}
\label{bt:scalarlemma}
Suppose $\Omega\subset\mathbb R^3$ is a bounded domain with
$C^\infty$ boundary and, for $j=1,2$,
$\lambda_j,\mu_j,W_j\in C^\infty(\overline\Omega;\mathbb R)$
satisfy $\mu_j>0$ and $3\lambda_j+2\mu_j>0$ on $\overline\Omega$.
The notation in~\eqref{bt:notation} is used for each material.
The elasticity operator $L_j=L_{\lambda_j,\mu_j}$, the maps $T_j,S_j$,
and the reduced operator $M_j^{W_j}$ are defined
in~\eqref{bt:stress}, \eqref{bt:TS}, and~\eqref{bt:MW},
respectively. The formal transpose $S_j^T$ is defined using
the complex-bilinear $L^2$ pairing on compactly supported
smooth fields.

For every $z=(v,f)^T\in C^\infty(\overline\Omega;\mathbb C^4)$
and $u\in C^\infty(\overline\Omega;\mathbb C^3)$, the
factorization and transpose identities hold pointwise
in $\Omega$ for $j=1,2$:
\begin{equation}
 \begin{aligned}
  (L_j+W_jI_3)(T_jz)&=S_j(M_j^{W_j}z),\\
  S_j^Tu&=
  \begin{pmatrix}
   \sqrt{\mu_j}\,u\\
   -(1-\kappa_j)^{-1}\operatorname{div}u
  \end{pmatrix}.
 \end{aligned}
 \label{bt:factorization}
\end{equation}
The fourth component of $M_j^{W_j}z$ satisfies
\begin{equation}
 \begin{aligned}
  \frac{(M_j^{W_j}z)_4}{1-\kappa_j}
  &=\Sigma_j\operatorname{div}(T_jz)
    -\sqrt{\mu_j}\,\operatorname{div}v\\
  &\quad+\tfrac32\sqrt{\mu_j}\,\mathbf m_j\cdot v
    +\frac{W_j}{\mu_j}f.
 \end{aligned}
 \label{bt:scalar}
\end{equation}
Consequently, the displacement map $T_j$ sends every smooth
solution $z$ of $M_j^{W_j}z=0$ in $\Omega$ to an elastic
displacement $u=T_jz$ satisfying $(L_j+W_jI_3)u=0$ in $\Omega$.
\end{lemma}

\begin{proof}
Fix $j\in\{1,2\}$. The Cartesian indices $i,k$ range from $1$ to $3$,
and repeated Cartesian indices are summed.
We first verify the factorization~\eqref{bt:factorization} by treating
the scalar and vector contributions to $T_jz$ in~\eqref{bt:TS}
separately. For the scalar contribution, define
$b_j=a_j\nabla f-f\nabla a_j$ and
$p_{f,j}=\Delta f-(\Delta a_j/a_j)f=\Delta f+q_{\mu_j}f$,
where the last equality follows from~\eqref{bt:notation}
and~\eqref{bt:MW}.
The product rule gives
$\varepsilon(b_j)=a_j\nabla^2f-f\nabla^2a_j$ and
$\operatorname{div}b_j=a_jp_{f,j}$.
With $\partial_{ik}=\partial_i\partial_k$,
expansion of $\partial_k(\sigma_j(b_j))_{ik}$
according to~\eqref{bt:stress} yields
\begin{equation}
 \begin{aligned}
  \bigl(L_jb_j-\nabla(a_j\Sigma_jp_{f,j})\bigr)_i
  &=2\mu_j(\Delta a_j\,\delta_{ik}-\partial_{ik}a_j)
       \partial_k f\\
  &\quad
    +2\bigl((\partial_i\mu_j)\Delta a_j
            -(\partial_k\mu_j)\partial_{ik}a_j\bigr)f.
 \end{aligned}
 \label{bt:scalar-index}
\end{equation}
The coefficient relations
$\nabla\mu_j=-\mu_j^2\nabla a_j$,
$\nabla a_j=-a_j\mathbf m_j$,
$\Delta a_jI_3-\nabla^2a_j=(c_j/2)F_{\mu_j}$, and
$a_j\Sigma_j=(1-\kappa_j)^{-1}$,
obtained from~\eqref{bt:notation} and~\eqref{bt:MW},
reduce~\eqref{bt:scalar-index} to
\begin{equation}
 L_jb_j
 =\nabla\!\left(
     \frac{\Delta f+q_{\mu_j}f}{1-\kappa_j}
   \right)
  +c_j^{-1}F_{\mu_j}(\nabla f+\mathbf m_jf).
 \label{bt:scalarfactor}
\end{equation}

For the vector contribution, define
$d_v=\operatorname{div}v$ and
$\Psi_{v,j}=c_j(\lambda_j+\mu_j)d_v
 +(\lambda_j-\mu_j)\nabla c_j\cdot v$.
Subtracting $\nabla\Psi_{v,j}$ from $L_j(c_jv)$ cancels the term
$c_j(\lambda_j+\mu_j)\nabla d_v$ and all terms involving
$\nabla\lambda_j$.
The relation $2\mu_j\nabla c_j+c_j\nabla\mu_j=0$,
following from~\eqref{bt:notation}, reduces the
remaining terms to
\begin{equation}
 \begin{aligned}
  \bigl(L_j(c_jv)-\nabla\Psi_{v,j}\bigr)_i
  &=c_j^{-1}\bigl(\Delta v_i-(\mathbf m_j)_id_v\bigr)
    +(Z_j)_{ik}v_k,\\
  (Z_j)_{ik}
  &=2\mu_j\partial_{ik}c_j
    +(\partial_i\mu_j)(\partial_kc_j)
    +(\partial_k\mu_j)(\partial_ic_j)\\
  &\quad
    +(\mu_j\Delta c_j+\nabla\mu_j\cdot\nabla c_j)\delta_{ik}.
 \end{aligned}
 \label{bt:vector-index}
\end{equation}
The differentiation formulas
$\partial_ic_j=-c_j(\mathbf m_j)_i/2$ and
$\partial_{ik}c_j
 =c_j((\mathbf m_j)_i(\mathbf m_j)_k/4-\partial_{ik}\ell_j/2)$,
obtained from~\eqref{bt:notation}, give
$c_jZ_j=V_{\mu_j}$ by~\eqref{bt:MW}.
The definitions of $d_j$ and $e_j^{\mathrm{red}}$ in~\eqref{bt:MW} also give
$\Psi_{v,j}=(d_jd_v-e_j^{\mathrm{red}}\mathbf m_j\cdot v)/(1-\kappa_j)$.
Consequently, equation~\eqref{bt:vector-index} becomes
\begin{equation}
 \begin{aligned}
  L_j(c_jv)
  &=c_j^{-1}(\Delta v-\mathbf m_jd_v+V_{\mu_j}v)\\
  &\quad
    +\nabla\!\left(
       \frac{d_jd_v-e_j^{\mathrm{red}}\mathbf m_j\cdot v}{1-\kappa_j}
     \right).
 \end{aligned}
 \label{bt:vectorfactor}
\end{equation}
Adding~\eqref{bt:scalarfactor} and~\eqref{bt:vectorfactor},
with $T_j,S_j$ defined in~\eqref{bt:TS}
and $M_j^0$ defined in~\eqref{bt:MW}, yields
$L_jT_jz=S_jM_j^0z$.

For the potential contribution, the relation
$\nabla a_j=-a_j\mathbf m_j$ and the definitions
in~\eqref{bt:TS} and~\eqref{bt:MW} give
\begin{equation}
 \begin{aligned}
  S_j(\mathcal W_jz)
  &=W_jc_jv+a_j(2W_j\mathbf m_j-\nabla W_j)f
    +\nabla(W_ja_jf)\\
  &=W_j(c_jv+a_j\nabla f-f\nabla a_j)
   =W_jT_jz.
 \end{aligned}
 \label{bt:potentialfactor}
\end{equation}
Hence $(L_j+W_jI_3)T_jz=S_j(M_j^{W_j}z)$.
Integration by parts against compactly supported test fields
in the source map~\eqref{bt:TS}, with the complex-bilinear
$L^2$ pairing, gives
$S_j^Tu=(c_j^{-1}u,-(1-\kappa_j)^{-1}\operatorname{div}u)^T$,
completing the verification of~\eqref{bt:factorization}.

To establish the scalar identity~\eqref{bt:scalar}, differentiation
of the displacement map~\eqref{bt:TS} gives
\begin{equation}
 \operatorname{div}(T_jz)
 =c_jd_v-\tfrac12c_j\mathbf m_j\cdot v
  +a_j(\Delta f+q_{\mu_j}f).
 \label{bt:divT}
\end{equation}
Substituting~\eqref{bt:divT} into the right-hand side
of~\eqref{bt:scalar} and using
$\Sigma_ja_j=(1-\kappa_j)^{-1}$ from~\eqref{bt:notation},
together with the definitions of $d_j$ and $e_j^{\mathrm{red}}$
in~\eqref{bt:MW}, we obtain
\begin{equation}
 \frac{\Delta f+d_jd_v-e_j^{\mathrm{red}}\mathbf m_j\cdot v+q_{\mu_j}f}
      {1-\kappa_j}
 +\frac{W_j}{\mu_j}f
 =\frac{(M_j^{W_j}z)_4}{1-\kappa_j},
 \label{bt:scalarverification}
\end{equation}
which proves the scalar identity~\eqref{bt:scalar}.
\end{proof}

The maps $T_j$ and $S_j$ in~\eqref{bt:TS} give rise to the
displacement discrepancy, traction discrepancy, and auxiliary
scalar term in the reduced Green--Betti identity~\eqref{bt:defects}.
Subsections~\ref{bt:total} and~\ref{bt:weak} derive these terms
and show how the matching conditions~\eqref{bt:matching} cancel them.

\subsection{Boundary terms in the reduced Green--Betti identity}
\label{bt:total}
To derive the reduced Green--Betti identity~\eqref{bt:defects},
we first introduce a reduced solution
$z=(v,f)^T\in C^\infty(\overline\Omega;\mathbb C^4)$
for the first material and an elastic displacement
$u_2\in C^\infty(\overline\Omega;\mathbb C^3)$
for the second material, satisfying
\[
 M_1^{W_1}z=0,\qquad
 (L_2+W_2I_3)u_2=0
 \quad\text{in }\Omega.
\]
With $\gamma_0$ denoting the Dirichlet trace on $\partial\Omega$,
the associated fields and boundary traces are defined by
\begin{equation}
 \begin{gathered}
  u_1=T_1z,\qquad V=T_2z,\qquad w=M_2^{W_2}z,\\
  y_2=S_2^Tu_2,\qquad
  F_i=\gamma_0u_i\quad(i=1,2).
 \end{gathered}
 \label{bt:fields}
\end{equation}

The following proposition relates the reduced volume pairing
$\int_\Omega y_2^TD^Wz\,dx$ to the physical boundary pairing
associated with $\Lambda_2^\Omega-\Lambda_1^\Omega$.
The identity~\eqref{bt:defects} identifies the additional
contributions from the displacement trace discrepancy,
the traction discrepancy, and the scalar component $w_4$.
The formula is derived directly from the physical boundary maps and records all discrepancies before the matching conditions~\eqref{bt:matching} are imposed.
\begin{proposition}
\label{bt:totalprop}
Assume that the coefficient triples
$(\lambda_j,\mu_j,W_j)\in C^\infty(\overline\Omega;\mathbb R^3)$,
$j=1,2$, satisfy $\mu_j>0$ and $3\lambda_j+2\mu_j>0$
on $\overline\Omega$, and that the corresponding 
dirichlet problems satisfy the nonresonance
condition~\eqref{en:nonresonance}. For $z=(v,f)^T\in C^\infty(\overline\Omega;\mathbb C^4)$
and $u_2\in C^\infty(\overline\Omega;\mathbb C^3)$ satisfying
\begin{equation}
 M_1^{W_1}z=0,\qquad
 (L_2+W_2I_3)u_2=0
 \quad\text{in }\Omega,
 \label{bt:solution-equations}
\end{equation}
the fields and boundary traces defined in~\eqref{bt:fields}
satisfy the Green--Betti identity
\begin{align}
 \int_\Omega y_2^TD^Wz\,dx
 &=\langle(\Lambda_2^\Omega-\Lambda_1^\Omega)F_1,F_2\rangle
   +\langle\Lambda_2^\Omega F_2,\gamma_0V-F_1\rangle
   \notag\\
 &\quad-\langle t_2(V)-t_1(u_1),F_2\rangle
   +\int_{\partial\Omega}
      \frac{w_4}{1-\kappa_2}(u_2\cdot\nu)\,dS.
 \label{bt:defects}
\end{align}
 In the above formulation, the scalar $w_4=(M_2^{W_2}z)_4$ represents the fourth component of $w=M_2^{W_2}z$. The traction operators $t_j$ and the Dirichlet-to-Neumann maps $\Lambda_j^\Omega$ for $j=1,2$ are defined in~\eqref{bt:stress} and~\eqref{en:dn-definition}, respectively. The operator difference $D^W$ is introduced in~\eqref{bt:operator-difference}. All boundary pairings are understood as complex bilinear forms.
Identity~\eqref{bt:defects} holds under the stated coefficient assumptions,
prior to imposing the boundary matching conditions~\eqref{bt:matching}.
\end{proposition}
\begin{proof}
Integration by parts in the gradient term of the source map
$S_2$ defined in~\eqref{bt:TS} gives
\begin{equation}
 \int_\Omega u_2^TS_2w\,dx
 -\int_\Omega(S_2^Tu_2)^Tw\,dx
 =\int_{\partial\Omega}
   \frac{w_4}{1-\kappa_2}(u_2\cdot\nu)\,dS.
 \label{bt:Sgreen}
\end{equation}

The constitutive law~\eqref{bt:stress} implies
\begin{equation}
 \begin{aligned}
  \nabla u_2:\sigma_2(V)
  &=2\mu_2\varepsilon(u_2):\varepsilon(V)
    +\lambda_2(\operatorname{div}u_2)(\operatorname{div}V)\\
  &=\nabla V:\sigma_2(u_2).
 \end{aligned}
 \label{bt:stress-symmetry}
\end{equation}
Integration by parts, using~\eqref{bt:stress-symmetry}
and the cancellation of the scalar potential terms,
therefore yields the  Green identity
\begin{equation}
 \begin{aligned}
  &\int_\Omega
   \bigl[
    u_2^T(L_2+W_2I_3)V
    -V^T(L_2+W_2I_3)u_2
   \bigr]\,dx\\
  &\qquad
   =\langle t_2(V),F_2\rangle
    -\langle t_2(u_2),\gamma_0V\rangle.
 \end{aligned}
 \label{bt:physicalgreen}
\end{equation}
For each material, application of the  Green identity
to two homogeneous solutions with smooth boundary data
$\varphi,\psi$ gives
\begin{equation}
 \langle\Lambda_j^\Omega\varphi,\psi\rangle
 =\langle\Lambda_j^\Omega\psi,\varphi\rangle,
 \qquad j=1,2.
 \label{bt:bilinear-reciprocity}
\end{equation}
The complex-bilinear reciprocity relation extends to all
$\varphi,\psi\in H^{1/2}(\partial\Omega;\mathbb C^3)$
by density of smooth boundary data and boundedness of
$\Lambda_j^\Omega$.

The factorization identity~\eqref{bt:factorization}
and the field definitions~\eqref{bt:fields} give
$(L_2+W_2I_3)V=S_2w$ and $y_2=S_2^Tu_2$,
while $M_1^{W_1}z=0$ gives $D^Wz=-w$.
Substitution of the field relations into~\eqref{bt:Sgreen}
and~\eqref{bt:physicalgreen}, together with
$(L_2+W_2I_3)u_2=0$ and
$t_2(u_2)=\Lambda_2^\Omega F_2$, yields
\begin{equation}
 \begin{aligned}
  \int_\Omega y_2^TD^Wz\,dx
  &=-\langle t_2(V),F_2\rangle
    +\langle\Lambda_2^\Omega F_2,\gamma_0V\rangle\\
  &\quad
    +\int_{\partial\Omega}
      \frac{w_4}{1-\kappa_2}(u_2\cdot\nu)\,dS.
 \end{aligned}
 \label{bt:green-intermediate}
\end{equation}
To isolate the physical Dirichlet-to-Neumann difference,
we write
$\gamma_0V=F_1+(\gamma_0V-F_1)$ and
$t_2(V)=t_1(u_1)+(t_2(V)-t_1(u_1))$.
Substituting the decompositions into~\eqref{bt:green-intermediate}
and using $t_1(u_1)=\Lambda_1^\Omega F_1$
together with the reciprocity relation~\eqref{bt:bilinear-reciprocity},
we obtain~\eqref{bt:defects}.

The discrepancy from the physical Dirichlet-to-Neumann
pairing therefore admits the boundary representation
\begin{equation}
 \begin{aligned}
  &\int_\Omega y_2^TD^Wz\,dx
   -\langle(\Lambda_2^\Omega-\Lambda_1^\Omega)F_1,F_2\rangle\\
  &\quad
   =\int_{\partial\Omega}\bigg[
      t_2(u_2)\cdot(V-u_1)
      -\bigl(t_2(V)-t_1(u_1)\bigr)\cdot u_2\\
  &\hspace{9em}
      +\frac{w_4}{1-\kappa_2}(u_2\cdot\nu)
     \bigg]\,dS.
 \end{aligned}
 \label{bt:exact-boundary-defect}
\end{equation}
\end{proof}
\subsection{Elimination of the boundary terms}
\label{bt:weak}

Under the matching conditions~\eqref{bt:matching}, the boundary
terms in~\eqref{bt:defects} cancel, so equality of the
Dirichlet-to-Neumann maps implies that the reduced volume pairing
vanishes. To compute the traction discrepancy~\eqref{bt:tractiondiff},
we apply the boundary Hessian formula~\eqref{bt:normalhessian}
of Lemma~\ref{bt:hessian} to $a_2-a_1$, with $a_j$ defined
in~\eqref{bt:notation}.
\begin{lemma}
\label{bt:hessian}
Let $\Omega\subset\mathbb R^3$ be a domain with $C^2$ boundary,
and let $\nu$ denote the outward unit normal to $\partial\Omega$.
Suppose that $b\in C^2(\overline\Omega;\mathbb R)$ satisfies
\begin{equation}
 b=0,\qquad \nabla b=0
 \quad\text{on }\partial\Omega.
 \label{bt:vanishing-first-jet}
\end{equation}
Then the boundary Hessian of $b$ has only a normal--normal
component:
\begin{equation}
 \nabla^2b(x)=\eta(x)\,\nu(x)\otimes\nu(x),
 \qquad x\in\partial\Omega,
 \label{bt:normalhessian}
\end{equation}
where $\eta(x):=\nu(x)^T\nabla^2b(x)\nu(x)$.
\end{lemma}

\begin{proof}
Fix $x\in\partial\Omega$ and $\tau\in T_x(\partial\Omega)$.
Tangential differentiation of the gradient condition
in~\eqref{bt:vanishing-first-jet} gives
\[
 \nabla^2b(x)\tau=0.
\]
The symmetry of $\nabla^2b(x)$ therefore implies
\eqref{bt:normalhessian}.
\end{proof}

Proposition~\ref{bt:weakprop} identifies the cancellation of the boundary terms in~\eqref{bt:defects} under~\eqref{bt:matching} and computes the residual boundary contribution when the matching of $W_j$ is relaxed.
\begin{proposition}
\label{bt:weakprop}
Under the coefficient assumptions of
Proposition~\ref{bt:totalprop}, let $z=(v,f)^T$ and $u_2$
be smooth solutions of~\eqref{bt:solution-equations},
with the associated fields and boundary traces defined
by~\eqref{bt:fields}. Suppose that the boundary matching conditions
\begin{equation}
 \mu_1=\mu_2,\qquad
 \partial_\nu\mu_1=\partial_\nu\mu_2,\qquad
 W_1=W_2
 \quad\text{on }\partial\Omega
 \label{bt:matching}
\end{equation}
hold. Then the reduced Green--Betti identity
\begin{equation}
 \int_\Omega y_2^TD^Wz\,dx
 =\langle
   (\Lambda_2^\Omega-\Lambda_1^\Omega)F_1,F_2
  \rangle
 \label{bt:betti}
\end{equation}
holds.

If the boundary condition $W_1=W_2$ on $\partial\Omega$ is dropped while
the two matching conditions for the shear modulus in~\eqref{bt:matching}
are retained, the right-hand side of~\eqref{bt:betti} acquires the boundary defect term
\begin{equation}
 \int_{\partial\Omega}
 \frac{W_2-W_1}{\mu_1}\,f(u_2\cdot\nu)\,dS.
 \label{bt:Wdefect}
\end{equation}
The matching conditions for this cancellation concern the traces of $\mu_j$,
$\partial_\nu\mu_j$, and $W_j$, as specified in~\eqref{bt:matching}.
\end{proposition}

\begin{proof}
Under the shear-modulus matching conditions in~\eqref{bt:matching},
differentiating $\mu_1-\mu_2=0$ along tangent directions and using
the equality $\partial_\nu\mu_1=\partial_\nu\mu_2$, we obtain
$\nabla\mu_1=\nabla\mu_2$ on $\partial\Omega$.
For a coefficient $b$ carrying a material index, use
$\delta_{21}b:=b_2-b_1$, the opposite subtraction order to
$\delta_{12}b$.
Consequently, the coefficient differences
$\delta_{21}a:=a_2-a_1$ and
$\delta_{21}c:=c_2-c_1$, with $a_j,c_j$ defined
in~\eqref{bt:notation}, satisfy
\[
 \delta_{21}a=\delta_{21}c=0,\qquad
 \nabla\delta_{21}a=\nabla\delta_{21}c=0
 \quad\text{on }\partial\Omega.
\]
Applying Lemma~\ref{bt:hessian}, we find that the Hessian of
$\delta_{21}a$ assumes the normal form~\eqref{bt:normalhessian}:
\[
 \nabla^2\delta_{21}a=\eta_\mu\nu\otimes\nu,
 \qquad
 \eta_\mu:=\nu^T(\nabla^2\delta_{21}a)\nu
 \quad\text{on }\partial\Omega.
\]

From the displacement representation~\eqref{bt:TS} and the
field definitions~\eqref{bt:fields}, the difference of
displacements satisfies
\begin{equation}
 V-u_1
 =\delta_{21}c\,v+\delta_{21}a\nabla f
  -f\nabla\delta_{21}a.
 \label{bt:dispdiff}
\end{equation}
Since $\delta_{21}a$, $\delta_{21}c$, and the corresponding
first-order derivatives vanish on $\partial\Omega$,
formula~\eqref{bt:dispdiff} implies
$\gamma_0V=\gamma_0u_1=F_1$.
Differentiating~\eqref{bt:dispdiff} and evaluating the result
on $\partial\Omega$, we obtain
\[
 \nabla(V-u_1)=-f\eta_\mu\nu\otimes\nu,
 \qquad
 \operatorname{div}(V-u_1)=-f\eta_\mu
 \quad\text{on }\partial\Omega.
\]
By~\eqref{bt:matching}, the equality $\mu_1=\mu_2$
holds on $\partial\Omega$.
Evaluating the traction difference on $\partial\Omega$
via~\eqref{bt:stress}, we deduce
\begin{align}
 t_2(V)-t_1(u_1)
 &=\bigl\{
     2\mu_1\varepsilon(V-u_1)
     +(\lambda_2\operatorname{div}V
       -\lambda_1\operatorname{div}u_1)I_3
   \bigr\}\nu
   \notag\\
 &=\bigl\{
     (\lambda_2-\lambda_1)\operatorname{div}u_1
     -(\lambda_2+2\mu_1)f\eta_\mu
   \bigr\}\nu
   \quad\text{on }\partial\Omega.
 \label{bt:tractiondiff}
\end{align}

The matching conditions for the first-order derivatives
in~\eqref{bt:matching} also imply that
$\mathbf m_1=\mathbf m_2$ on $\partial\Omega$.
Subtracting the scalar identities~\eqref{bt:scalar}
for $j=1,2$ and using $\mu_1=\mu_2$ on $\partial\Omega$
eliminates the terms involving $\operatorname{div}v$
and $\mathbf m_j\cdot v$.
Since $(M_1^{W_1}z)_4=0$ by~\eqref{bt:solution-equations},
evaluating the resulting difference on $\partial\Omega$ yields
\begin{align}
 \frac{w_4}{1-\kappa_2}
 &=\Sigma_2\operatorname{div}V
   -\Sigma_1\operatorname{div}u_1
   +\frac{W_2-W_1}{\mu_1}f
   \notag\\
 &=(\lambda_2-\lambda_1)\operatorname{div}u_1
   -(\lambda_2+2\mu_1)f\eta_\mu
   +\frac{W_2-W_1}{\mu_1}f
   \quad\text{on }\partial\Omega.
 \label{bt:fourthdiff}
\end{align}
Comparing~\eqref{bt:tractiondiff} and~\eqref{bt:fourthdiff},
the contributions involving
$(\lambda_2-\lambda_1)\operatorname{div}u_1$
and $\eta_\mu$ cancel identically, yielding
\begin{equation}
 \begin{aligned}
  &\frac{w_4}{1-\kappa_2}(u_2\cdot\nu)
   -\bigl(t_2(V)-t_1(u_1)\bigr)\cdot u_2\\
  &\qquad
   =\frac{W_2-W_1}{\mu_1}f(u_2\cdot\nu)
   \quad\text{on }\partial\Omega.
 \end{aligned}
 \label{bt:explicit-cancellation}
\end{equation}
Because $\gamma_0V=F_1$, the boundary defect
in~\eqref{bt:exact-boundary-defect} reduces, by
\eqref{bt:explicit-cancellation}, to the boundary
integral~\eqref{bt:Wdefect}.
When $W_1=W_2$ on $\partial\Omega$, the boundary defect
term vanishes and the reduced Green--Betti
identity~\eqref{bt:betti} follows.
\end{proof}

For the extensions in Propositions~\ref{ext:common}
and~\ref{ext:finite-common}, \eqref{bt:MW} gives
\begin{equation}\label{bt:status}
 \begin{aligned}
 D^W&=(M_1^0-M_2^0)+(\mathcal W_1-\mathcal W_2),\\
 M_1^0-M_2^0&=0
       &&\text{on }\mathbb R^3\setminus\overline\Omega .
 \end{aligned}
\end{equation}
For the smooth common extensions, the multiplication matrices
$\mathcal W_j$ also agree outside $\overline\Omega$. With the
separate finite-order extensions of $W_j$, $\mathcal W_1-\mathcal W_2$ may have
exterior support. The matrix-moment argument in
Theorem~\ref{i2:thm-nine} and Subsection~\ref{fr:interior-assembly}
takes the scaled limit of \eqref{bt:betti} on $\Omega$. The limiting
integrand contains the first-order part of $M_1^0-M_2^0$.
By \eqref{bt:status}, its coefficients are supported in
$\overline\Omega$, and the integral extends to the working domain.

\section{Transport frames on the plane}\label{f27:section}
The construction of complex geometrical optics solutions
requires invertible transport frames with smooth dependence
on the direction and transverse variables. The relevant
transport equations arise from the conjugated matrix operator
\begin{equation}\label{a27:generic}
 \begin{aligned}
  P_\delta
  &=(D+\delta)^2
    +2\sum_{k=1}^3 A_k(x)(D_k+\delta_k)+q(x),\\
  A_k,q
  &\in C_c^\infty(\mathbb R^3;M_m(\mathbb C)),
  \qquad k=1,2,3.
 \end{aligned}
\end{equation}
In~\eqref{a27:generic}, $m$ denotes the matrix size,
$D=-i\nabla$ is the Fourier differential operator, with the pairing
convention of Section~\ref{sec:notation}, and the complex phase vector is
$\delta=\tau(\mathbf v_0+i\mathbf v_1)$ for real
orthonormal vectors $\mathbf v_0,\mathbf v_1$ and
$\tau\ge1$. This phase vector is distinct from the ordered coefficient
differences $\delta_{12}b$ and $\delta_{21}b$. The construction applies to arbitrary $m$,
with $m=2,3,4$ required in the elasticity applications.
The relation between $P_\delta$ and the reduced elastic
operator is established in Section~\ref{pb27:section}.

In planar coordinates, the leading transport equations
take the form~\eqref{f27:frameeq}, and a smooth solution
matrix $C$ with values in $GL(m,\mathbb C)$ is called
an invertible transport frame.
Theorem~\ref{f27:main} establishes global frames for
compactly supported smooth planar coefficients and
jointly smooth local families on prescribed disks.
The global construction combines weighted Cauchy estimates
with Fredholm theory. Polynomial corrections enforce
compatibility with the finite-dimensional cokernel,
producing a global solution matrix. Local holomorphic
factorizations then permit the successive removal of
determinant zeros, yielding an invertible frame for every
compactly supported smooth planar coefficient.

The global frame supplies the reference solution for
the local parameter construction. The right-inverse
identity~\eqref{f27:rightinverse} yields the
Neumann-series representation~\eqref{f27:param-neumann},
and interior regularity gives joint smoothness in the
planar variables and real parameters. The mixed
seminorms in~\eqref{f27:mixed} measure the finite-order
derivatives of the frame and inverse needed in the
subsequent operator estimates.
For each fixed smooth coefficient family, these bounds are
uniform on compact parameter subcharts.

The planar construction supplies the two-, three-, and
four-component transport frames in
Proposition~\ref{tb27:F-supply}.
For the subsequent local right inverse,
Proposition~\ref{tb27:budget} specifies the derivative
orders required when the frames are expressed in moving
plane coordinates. The geometry, coefficients, cutoffs,
local parameter charts, and frames are fixed before
the limit $\tau\to\infty$ is taken. The estimate constants
depend on the fixed data and the selected finite frame
seminorms, and are independent of $\tau$.

The planar notation used in the construction is
\begin{equation}\label{f27:planar-notation}
 \begin{gathered}
  z=x_1+ix_2,\qquad
  D_R=\{z\in\mathbb C:|z|<R\},\qquad
  \langle z\rangle=(1+|z|^2)^{1/2},\\
  \partial_{\bar z}
  =\tfrac12(\partial_{x_1}+i\partial_{x_2}),
  \qquad
  \partial_z
  =\tfrac12(\partial_{x_1}-i\partial_{x_2}).
 \end{gathered}
\end{equation}
The auxiliary Hilbert inner products are linear in the
first variable. With the planar notation
in~\eqref{f27:planar-notation}, the Cauchy operators on
compactly supported smooth scalar, vector, or matrix
fields are defined componentwise by
\begin{equation}\label{f27:cauchy-operators}
 \mathscr C f(z)
 =\frac1\pi\int_{\mathbb C}\frac{f(w)}{z-w}\,dA(w),
 \qquad
 \mathscr C^* f(z)
 =-\frac1\pi\int_{\mathbb C}
       \frac{f(w)}{\bar z-\bar w}\,dA(w).
\end{equation}
In~\eqref{f27:cauchy-operators}, $dA$ denotes Euclidean
area measure, and the superscript $*$ denotes the
unweighted formal Hermitian adjoint. In the weighted
cokernel system~\eqref{f27:cokernel-eq}, $B^*$ denotes
the Hermitian transpose of the matrix coefficient $B$.
The weighted Hilbert adjoint is defined separately
in Lemma~\ref{f27:cokernel}.

\begin{theorem}
\label{f27:main}
Let $m\ge1$ be an integer.
\begin{enumerate}
\item
For every $B\in C_c^\infty(\mathbb C;M_m(\mathbb C))$,
there exists a frame
$C\in C^\infty(\mathbb C;GL(m,\mathbb C))$
satisfying the transport equation
\begin{equation}
 \partial_{\bar z}C=BC.
 \label{f27:frameeq}
\end{equation}
Moreover, the frame $C$ can be chosen to satisfy
\[
 \det C
 =c\,\exp\!\bigl(\mathscr C(\operatorname{tr}B)\bigr)
\]
for some constant $c\in\mathbb C\setminus\{0\}$, where
$\mathscr C$ is the Cauchy operator defined
in~\eqref{f27:cauchy-operators}.

\item
Let $0<L<R_{\mathrm{fr}}$, let $Y$ be a finite-dimensional
smooth real manifold, and let
\[
 B\in C^\infty
 \bigl(D_{R_{\mathrm{fr}}}\times Y;M_m(\mathbb C)\bigr).
\]
Here $R_{\mathrm{fr}}$ is the outer radius of the planar disk on which
the coefficient family is prescribed.
For every $y_0\in Y$, there exist a radius
$L'\in(L,R_{\mathrm{fr}})$, an open coordinate
neighborhood $U\subset Y$ of $y_0$, and a frame
\[
 C\in C^\infty
 \bigl(D_{L'}\times U;GL(m,\mathbb C)\bigr)
\]
such that $C(\cdot,y)$ satisfies~\eqref{f27:frameeq}
in $D_{L'}$ for every $y\in U$.

For every compact set $K\Subset U$ and every
integer $J\ge0$, the frame $C$ and the inverse
matrix $C^{-1}$ satisfy
\begin{equation}
 \max_{|\alpha|+|\beta|\le J}
 \sup_{\overline D_L\times K}
 \left(
  \left|\partial_x^\alpha\partial_y^\beta C\right|
  +\left|\partial_x^\alpha\partial_y^\beta C^{-1}\right|
 \right)<\infty,
 \label{f27:mixed}
\end{equation}
where derivatives are taken with respect to
$x=(\operatorname{Re}z,\operatorname{Im}z)$ and the
real coordinates on $U$.

For every prescribed function
\[
 g\in C^\infty
 \bigl(D_{L'}\times U;\mathbb C\setminus\{0\}\bigr)
\]
satisfying $\partial_{\bar z}g=(\operatorname{tr}B)g$,
the local frame $C$ can additionally be chosen to satisfy
$\det C=g$ on $D_{L'}\times U$.
Under the determinant normalization $\det C=g$, the mixed
derivative bound~\eqref{f27:mixed} continues to hold.

In particular, every compact subset of $Y$ admits a finite cover
by coordinate neighborhoods carrying local frame families.
\end{enumerate}
\end{theorem}

\subsection{The Cauchy operator and local uniqueness}

For the operators in \eqref{f27:cauchy-operators}, integration by parts outside a disk of radius $\varepsilon$, followed by $\varepsilon\downarrow0$, gives the distributional identity
$\partial_{\bar z}(1/(\pi z))=\delta_0$. Consequently,
\begin{equation}\label{f27:cauchy-identities}
 \partial_{\bar z}\mathscr{C} f=f,
 \qquad -\partial_z\mathscr{C}^* f=f.
\end{equation}
The Fourier multiplier of $\partial_z\mathscr{C}$ is
\[
 \frac{i\xi_1+\xi_2}{i\xi_1-\xi_2}\quad(\xi\ne0),
\]
which has modulus one. Plancherel's theorem therefore gives
$\|\partial_z\mathscr{C} f\|_2=\|f\|_2$.
On fixed bounded disks, local integrability of $|z-w|^{-1}$ and Young's inequality control
the local $L^2$ norm of $\mathscr{C}f$. Density extends these identities to compactly supported $L^2$ functions.
In particular, $\mathscr{C}:L^2(D_R)\to H^1(D_R)$ is bounded when the source is extended by zero.

We also record a derivative estimate for compact sources on a fixed
disk. Fix $S>0$ and $\chi\in C_c^\infty(D_S)$. For
$0\le k\le6$ and $v\in C^k(\overline D_S)$, the product $\chi v$
has a $C^k$ zero extension to $\mathbb C$, and
\begin{equation}\label{f27:uniform-cauchy-bound}
 \partial_z^\alpha\mathscr C(\chi v)
   =\mathscr C\bigl(\partial_z^\alpha(\chi v)\bigr),
 \qquad
 \|\mathscr C(\chi v)\|_{C^k(\overline D_S)}
   \le C_{S,k,\chi}\|v\|_{C^k(\overline D_S)},
 \quad |\alpha|\le k.
\end{equation}
Here $\partial_z^\alpha$ denotes derivatives in the two real leaf
coordinates. Distributional integration by parts transfers each
derivative to the compact source. The absolute integral of
$|z-w|^{-1}$ over $D_S$ is uniformly bounded for $z\in\overline D_S$,
and translation continuity of this locally integrable kernel makes
the resulting convolutions continuous. Thus the derivative identities
hold classically and give the displayed norm bound. The same argument
applies to mixed leaf and real parameter derivatives on a fixed compact
parameter set: parameter differentiation acts on the source, and the
fixed compact leaf support justifies differentiation under the integral.
In particular, the Cauchy operator acts separately at each transverse
parameter and preserves compact support in that parameter.

\begin{lemma}\label{f27:regularity}
If $u\in L^2_{\mathrm{loc}}(V)$ satisfies
$\partial_{\bar z}u=Bu+f$ with smooth $B,f$, then $u$ is smooth.
The same holds for $-\partial_z u=B^*u+f$.
On each relatively compact subdisk of $V$, every finite Sobolev norm is controlled by
$\|u\|_2$ on a larger subdisk and finitely many seminorms of $B$ and $f$.
\end{lemma}
\begin{proof}
For compactly supported $v$, the symbol of $\partial_{\bar z}$ has modulus $|\xi|/2$, so
\[
 \|v\|_{H^{k+1}}\le C_k
       (\|\partial_{\bar z}v\|_{H^k}+\|v\|_2).
\]
Apply this estimate to nested cutoffs of $u$. The equation first gives $H^1_{\rm loc}$ regularity;
the product rule then gives higher Sobolev regularity by induction. Sobolev embedding proves the assertion and the estimates.
Complex conjugation gives the result for the second operator.
\end{proof}

\begin{lemma}\label{f27:local-gauge}
Every smooth $B$ admits a smooth invertible transport frame near each point.
On a connected planar open set, a solution $u\in L^2_{\rm loc}$ of
$\partial_{\bar z}u=Bu$ that vanishes on a nonempty open subset vanishes identically.
The same uniqueness statement holds for $-\partial_z u=B^*u$.
\end{lemma}
\begin{proof}
Translate the point to the origin and choose a smooth cutoff supported in $D_{2r}$.
Let $B_{\rm cut}$ be the resulting coefficient, equal to $B$ on $D_r$.
Since
\[
 \sup_z\frac1\pi\int_{D_{2r}}|z-w|^{-1}\,dA(w)\le Cr,
\]
the equation $d=I+\mathscr{C}(B_{\rm cut}d)$ is a contraction on bounded continuous matrices for small $r$,
and $\|d-I\|_\infty<1/2$ can be imposed. The Cauchy transform of a bounded compactly supported function is continuous,
as follows from translation continuity of the locally integrable kernel.
Thus $d$ is invertible everywhere. Lemma~\ref{f27:regularity} makes it smooth, and on $D_r$ it satisfies
$\partial_{\bar z}d=Bd$.
For a local solution $u$, regularity and the equation imply
$\partial_{\bar z}(d^{-1}u)=0$. Hence $d^{-1}u$ is holomorphic.
The identity theorem propagates vanishing from one small disk to the next throughout the connected open set.
Complex conjugation reduces the last equation to the same argument with coefficient $-B^T$.
\end{proof}

\subsection{Weighted estimates and compactness}

Fix $1/2<s<1$ and set
\[
 L^2_t=\{f:\langle z\rangle^t f\in L^2\},\qquad
 \mathcal H=L^2_{-s}(\mathbb C;\mathbb C^m),\qquad
 \varpi=\langle z\rangle^{-2s}.
\]
We use $\|f\|_{L^2_t}=\|\langle z\rangle^t f\|_{L^2(\mathbb C)}$ and
$(u,v)_{\mathcal H}=\int_{\mathbb C}u\cdot\overline v\,\varpi\,dA$.
Multiplication by $B$ followed by the Cauchy operator defines the compact operator used in the Fredholm argument.
\begin{lemma}\label{f27:weighted}
The operators $\mathscr{C},\mathscr{C}^*:L^2_s\to L^2_{-s}$ are bounded.
For every compactly supported smooth matrix $B$,
$\mathcal K=\mathscr{C} B:\mathcal H\to\mathcal H$ is compact.
\end{lemma}
\begin{proof}
After conjugation by the weights, the absolute kernel is
$\pi^{-1}\langle z\rangle^{-s}|z-w|^{-1}\langle w\rangle^{-s}$.
Use $a(z)=\langle z\rangle^{-1}$ in Schur's test.
For $R=|z|\ge2$, split
$I(z)=\int|z-w|^{-1}\langle w\rangle^{-1-s}\,dA(w)$
into the regions $|w|\le R/2$, $|z-w|\le R/2$, their complement with $|w|\le2R$,
and $|w|\ge2R$. The respective bounds are
\[
 CR^{-1}\int_0^{R/2}r\langle r\rangle^{-1-s}\,dr,\quad
 CR^{-1-s}\int_0^{R/2}dr,\quad
 CR^{-1}R^{-1-s}R^2,\quad
 C\int_{2R}^\infty r^{-1-s}\,dr.
\]
Each is at most $CR^{-s}$. The integral is uniformly bounded for bounded $z$.
Both Schur quotients are thus bounded by $C\langle z\rangle^{1-2s}$, which is bounded because $s>1/2$.
The absolute kernel of $\mathscr{C}^*$ is identical.

Suppose $\operatorname{supp}B\subset D_{R_0}$ with $R_0\ge1$.
A bounded set in $\mathcal H$ is mapped by multiplication by $B$ to a bounded set in $L^2$ on a fixed compact set.
The Fourier identity and the weighted estimate bound $\mathscr{C}Bu$ in
$H^1$ on every fixed disk. Rellich's theorem gives local $L^2$ compactness. For the tail,
when $|z|>2R_0$ we have
\[
 |\mathscr{C}Bu(z)|\le C|z|^{-1}\|u\|_{\mathcal H},\qquad
 \|\mathscr{C}Bu\|_{\mathcal H(|z|>R)}^2
       \le CR^{-2s}\|u\|_{\mathcal H}^2.
\]
The second bound follows by integrating $r^{-1-2s}$ from $R$ to infinity.
Uniform smallness of the tail and compactness on a sufficiently large fixed disk prove compactness in $\mathcal H$.
\end{proof}

\subsection{The weighted cokernel and polynomial corrections}

\begin{lemma}\label{f27:cokernel}
Let $(I-\mathcal K)^*_{\mathcal H}$ denote the Hilbert adjoint in the weighted inner product just defined. If $h\in\ker(I-\mathcal K)^*_{\mathcal H}$, then
$\phi=\mathscr{C}^*(\varpi h)$ is well defined, belongs to $\mathcal H$, and satisfies
\begin{equation}\label{f27:cokernel-eq}
 \varpi h=B^*\phi,\qquad
 \phi=\mathscr{C}^*(B^*\phi),\qquad
 -\partial_z\phi=B^*\phi.
\end{equation}
In particular, $B^*\phi$ is compactly supported and $\phi$ is smooth. For every compactly supported $f\in L^2$,
\begin{equation}\label{f27:pairing}
 (\mathscr{C}f,h)_{\mathcal H}=(f,\phi)_{L^2}.
\end{equation}
\end{lemma}
\begin{proof}
The identity $\|\varpi h\|_{L^2_s}=\|h\|_{\mathcal H}$ permits the weighted estimate.
Interchanging the kernels gives~\eqref{f27:pairing} first for smooth compactly supported data;
weighted boundedness and density give the general case. The cokernel condition yields, for every $u\in\mathcal H$,
\[
 (u,h)_{\mathcal H}=(\mathscr{C}Bu,h)_{\mathcal H}
       =(Bu,\phi)_{L^2}=(u,B^*\phi)_{L^2}.
\]
It follows that $\varpi h=B^*\phi$.
Substituting this into the definition of $\phi$ and using~\eqref{f27:cauchy-identities} gives the remaining identities.
Smoothness follows from Lemma~\ref{f27:regularity}.
\end{proof}

\begin{lemma}\label{f27:polynomials}
Let $h_1,\ldots,h_r$ be a basis of $\ker(I-\mathcal K)^*_{\mathcal H}$, the finite-dimensional cokernel supplied by Lemma~\ref{f27:weighted}. For each $j$, set $\phi_j=\mathscr C^*(\varpi h_j)$ as in Lemma~\ref{f27:cokernel}, and define
\[
 \ell_j(P)=(BP,\phi_j)_{L^2},\qquad P\in\mathbb C[z]^m.
\]
Then $P\mapsto(\ell_1(P),\ldots,\ell_r(P))$ maps onto $\mathbb C^r$.
Its restriction to polynomials of degree at most $N$ is onto for some finite $N\ge0$.
\end{lemma}
\begin{proof}
Otherwise, a nonzero cokernel element $h$ would have a corresponding $\phi$ such that
$(BP,\phi)=0$ for every holomorphic vector polynomial $P$.
The coefficients are conjugated according to the convention that the first variable is linear. Put $f=B^*\phi$.
For the standard coordinate vectors $e_a\in\mathbb C^m$, choose $P=z^ke_a$ and take complex conjugates to obtain
\[
 \int_{\mathbb C}\bar w^k f(w)\,dA(w)=0\quad(k=0,1,\ldots).
\]
The function $f$ is compactly supported and belongs to $L^1$. Outside a disk containing its support, the geometric series converges absolutely and uniformly, so
\[
 \phi(z)=-\frac1\pi\sum_{k=0}^\infty\bar z^{-k-1}
                       \int\bar w^k f(w)\,dA(w)=0.
\]
Uniqueness in Lemma~\ref{f27:local-gauge}, applied to $-\partial_z\phi=B^*\phi$,
gives $\phi=0$, and~\eqref{f27:cokernel-eq} then gives $h=0$, a contradiction.
Since the target is finite dimensional, finitely many polynomial images already span it; their maximum degree gives $N$.
If $r=0$, take $N=0$.
\end{proof}

\begin{proof}[Proof of Theorem~\ref{f27:main}, part~1]
By Lemma~\ref{f27:weighted}, $I-\mathcal K$ is Fredholm of index zero. Choose the polynomial degree $N$ supplied by Lemma~\ref{f27:polynomials} and put $d=N+1$. That lemma gives, for each standard vector $e_i\in\mathbb C^m$,
a polynomial $Q_i$ of degree at most $N$ such that
\[
 P_i=z^de_i+Q_i,\qquad \ell_j(P_i)=0\quad(1\le j\le r).
\]
By~\eqref{f27:pairing}, $\mathscr{C}(BP_i)$ is orthogonal to the weighted cokernel.
The Fredholm alternative gives $v_i\in\mathcal H$ satisfying
\[
 (I-\mathcal K)v_i=\mathscr{C}(BP_i).
\]
Thus $C_i=P_i+v_i$ solves $\partial_{\bar z}C_i=BC_i$, with
$v_i=\mathscr{C}(BC_i)$. Local regularity makes $C_i$ smooth, and expansion outside the source support gives
$v_i=O(1/z)$. Hence $C=(C_1,\ldots,C_m)$ satisfies
\begin{equation}\label{f27:leading}
 C=z^dI+O(|z|^{d-1}),\qquad
 \det C=z^{md}(1+O(|z|^{-1}))\quad(|z|\to\infty).
\end{equation}
In particular, its determinant is not identically zero.

Multilinearity of the determinant, including at singular matrices, gives
\[
 \partial_{\bar z}\det C=(\operatorname{tr}B)\det C.
\]
Set $q=\mathscr{C}(\operatorname{tr}B)$. Then
$F=e^{-q}\det C$ is entire. Since $q=O(1/z)$, formula~\eqref{f27:leading} shows that
$F$ has polynomial growth with leading term $z^{md}$. Cauchy's estimates make it
a nonzero polynomial of degree $md$, whose zeros have finite total multiplicity.

At a zero $z_0$, choose a local frame $d_0$ from Lemma~\ref{f27:local-gauge} and write
$C=d_0H$ with $H$ holomorphic. Extend a nonzero $v\in\ker H(z_0)$ to the first column
of a constant $T\in GL(m,\mathbb C)$. The first column of $HT$ vanishes at $z_0$.
Divide the first column of the global matrix $CT$ by $z-z_0$ and retain the other columns.
The local expression $d_0HT$ gives a smooth removable extension at $z_0$,
because every component of its holomorphic first column is divisible by $z-z_0$.
The extended matrix solves the same equation, while its modified determinant $e^{-q}\det C$ is the previous one multiplied by a nonzero constant and divided by $z-z_0$.
The total zero multiplicity decreases by one. Repeating finitely many times yields $\det C=ce^q$ with $c\ne0$,
and hence a smooth invertible frame on the plane.
\end{proof}

\subsection{Frames on a fixed disk and smooth real parameters}

\begin{proof}[Proof of Theorem~\ref{f27:main}, part~2]
Choose $L<L'<R<R_{\mathrm{fr}}$, with finitely many intermediate radii as needed.
Let
$\chi\in C_c^\infty(D_{R_{\mathrm{fr}}})$ equal one near $\overline D_R$, and set
$\widetilde B_0=\chi B(\cdot,y_0)$, extended by zero outside $D_{R_{\mathrm{fr}}}$.

Part~1 gives a global frame $C_0$ for $\widetilde B_0$.
On $D_R$, define
\[
 S_0f=C_0\mathscr{C}(\mathbf1_{D_R}C_0^{-1}f).
\]
The local estimate following \eqref{f27:cauchy-identities} gives a bounded map $S_0:L^2(D_R)\to H^1(D_R)$, applied columnwise to matrix fields. Since $\partial_{\bar z}C_0=B(\cdot,y_0)C_0$ on $D_R$, the product rule gives
\begin{equation}\label{f27:rightinverse}
 (\partial_{\bar z}-B(\cdot,y_0))S_0f=f\quad\hbox{in }D_R.
\end{equation}
Here the source is extended by zero as an $L^2$ function.

Write $B_y=B(\cdot,y)$ and $b_y=B_y-B_{y_0}$, where $y$ ranges over a local real coordinate chart near $y_0$.
Smoothness on the product of $\overline D_R$ and a smaller parameter chart gives
\[
 A_y=S_0 b_y:H^1(D_R)\to H^1(D_R),\qquad
 \|A_y\|\le C\|b_y\|_\infty
\]
with norm less than $1/2$ near $y_0$. Therefore
\begin{equation}\label{f27:param-neumann}
 C_y=C_0+S_0(b_yC_y),\qquad C_y=(I-A_y)^{-1}C_0
\end{equation}
defines an $H^1$ family. By~\eqref{f27:rightinverse}, it satisfies
$\partial_{\bar z}C_y=B_yC_y$ on all of $D_R$.
The map $y\mapsto A_y$ is smooth in operator norm; differentiation of the bounded inverse,
or the Neumann series, gives all real parameter derivatives of $C_y$ in $H^1(D_R)$, continuously in $y$.
The Neumann bound is controlled by $\|B_y-B_{y_0}\|_{L^\infty}$;
the reference coefficient $B(\cdot,y_0)$ is fixed within the given
smooth coefficient class.
More quantitatively, $\|A_y\|\leq1/2$ and
$C_y-C_0=(I-A_y)^{-1}S_0(b_yC_0)$ give
\begin{equation}\label{f27:param-stability}
 \|C_y-C_0\|_{H^1(D_R)}
 \leq 2\|S_0\|_{L^2(D_R)\to H^1(D_R)}
          \|B_y-B_{y_0}\|_{L^\infty(D_R)}\|C_0\|_{L^2(D_R)}.
\end{equation}
For each fixed reference coefficient $B(\cdot,y_0)$,
\eqref{f27:param-stability} controls the dependence of $C_y$
on $B(\cdot,y)$ in a neighborhood of $y_0$.

To obtain uniform invertibility, put $W_y=C_y-C_0$.
Estimate~\eqref{f27:param-stability} gives $W_y\to0$ in $H^1(D_R)$ and
\[
 \partial_{\bar z}W_y=B_yW_y+b_yC_0.
\]
On an intermediate disk $D_{R_1}\Subset D_R$, localized Fourier estimates and the product rule give
\[
 \|W_y\|_{H^2(D_{R_1})}
 \le C\big(\|W_y\|_{H^1(D_R)}
             +\|b_y\|_{C^1(\overline D_R)}\big)\longrightarrow0,
\]
where the $C^1$ norms of $B_y$ are uniformly bounded after restricting the chart.
For $L'<R_1$, the two-dimensional Sobolev embedding gives $C_y\to C_0$ in
$C^0(\overline D_{L'})$.
Let $\sigma_0=\min_{\overline D_{L'}}\sigma_{\min}(C_0)>0$.
Restrict the parameter chart so that $\sup_{\overline D_{L'}}|C_y-C_0|<\sigma_0/2$.
All $C_y$ are then invertible on this fixed disk. In the Euclidean
operator norm the singular-value inequality
$\sigma_{\min}(C_y)\geq\sigma_{\min}(C_0)-\|C_y-C_0\|$
therefore yields the quantitative uniform bound
\begin{equation}\label{f27:param-inverse-bound}
 \sup_{y\in U}\sup_{z\in\overline D_{L'}}\|C_y(z)^{-1}\|
 \leq 2\sigma_0^{-1},
\end{equation}
after shrinking $U$ if needed. This is the bound on the inverse frame
used when compact parameter subcharts are selected later.

For every parameter multi-index $\beta$, the $H^1$ derivatives already constructed satisfy
\[
 \partial_{\bar z}\partial_y^\beta C_y
 =\sum_{\gamma\le\beta}\binom{\beta}{\gamma}
       (\partial_y^\gamma B_y)(\partial_y^{\beta-\gamma}C_y).
\]
Fix this parameter neighborhood after the invertibility step.
Writing $u_\beta=\partial_y^\beta C_y$, an induction on the Sobolev order in the leaf variable $z$ on nested subdisks bounds all $|\beta|\le q$, for arbitrary finite $q,k$, in $H^k$ by bounds for finitely many derivatives of the coefficients and the previously obtained norms on larger disks.
Applying the same estimates to $u_\beta(y)-u_\beta(y_*)$ gives continuity in the interior $H^k$ spaces. The $H^1$ identity
\[
 u_\beta(y+te_j)-u_\beta(y)
   =\int_0^t u_{\beta+e_j}(y+se_j)\,ds
\]
therefore holds in $H^k$ as well, by continuous $H^k$-valued integration and injectivity of the embedding. Apply the preceding estimates through order $q+1$ to cover the derivative in the integrand.
Thus the initially chosen family $C_y$ is smooth with values in every interior $H^k$ space.
Sobolev embedding for sufficiently large $k$ gives joint continuity of all mixed derivatives.
Hence $C$ is jointly smooth. The smooth map $A\mapsto A^{-1}$ on $GL(m,\mathbb C)$
gives joint smoothness of $C^{-1}$; equivalently, differentiate $\partial C^{-1}=-C^{-1}(\partial C)C^{-1}$.
Boundedness of continuous derivatives on $\overline D_L\times K$ proves~\eqref{f27:mixed}.
A finite subcover treats compact parameter sets, pairing each frame with its own inverse.

If $g$ is prescribed, $a=g/\det C_y$ is jointly smooth, nonvanishing, and holomorphic in $z$.
Multiplying the first column by $a$ gives determinant $g$.
The mixed seminorm bounds for the new frame and its inverse follow from smoothness and compactness.
\end{proof}

\begin{remark}\label{f27:Fexit}
Choose smooth leaf coordinates in a local real chart of directions. The two-, three-, and four-component systems below then have the form above.
For a transport field $\partial_{x_1}+i\partial_{x_2}=2\partial_{\bar z}$,
divide its coefficient by two. Direction and transverse coordinates together form the real parameter $y$.
For a prescribed disk and fixed material, choose the charts and
frames independently of the transverse test $b$, the leaf polynomial
$g$, and the CGO parameter in Section~\ref{i2:sec}. The mixed
seminorm bounds for the frame and its inverse are uniform on compact
subcharts. A finite collection of such subcharts covers any compact
parameter set.
\end{remark}

\subsection{Frames for the transport systems}\label{tb27:section}
We apply Theorem~\ref{f27:main} to the reduced system, its dual and quotient
systems, and the moving directions in the local right inverse. Use the
first-order representation~\eqref{pb27:MW} and set
$A_\theta^{\rm phys}=\sum_{k=1}^3\theta_kA_k^{\rm phys}$ and
$\mathsf T_\theta=2D_\theta+A_\theta^{\rm phys}$, where
$D_\theta=\theta\cdot\nabla$.
Thus $\mathsf T_\theta$ is the leading four-component transport
operator.

\begin{proposition}\label{tb27:F-supply}
Fix a sufficiently large leaf disk $\overline D_L$ and a center for the direction and transverse parameters. The two-, three-, and four-component systems used in Sections~\ref{i2:sec} and~\ref{if:section} admit smooth invertible frames on this disk. Each frame and its inverse are jointly smooth in a local real parameter chart and have all prescribed finite mixed seminorms. They can be fixed before all tests $b,g$ and all small parameters.
\end{proposition}
\begin{proof}
Choose a smooth local orthonormal frame $\theta=e_1+ie_2$, $\nu_\theta=e_1\times e_2$, and write
$x=z_1e_1+z_2e_2+t\nu_\theta$, $z=z_1+iz_2$.
Here $\nu_\theta$ is the real unit normal to the oriented
$\theta$-plane, distinct from the outward boundary normal $\nu$.
Then $D_\theta=2\partial_{\bar z}$.
The four-component leading transport equation $\mathsf T_\theta C=0$ becomes
\[
 \partial_{\bar z}C=-\tfrac14 A_\theta^{\rm phys}C.
\]
The dual equation $(2D_\theta-(A_\theta^{\rm phys})^T)Y=0$ has planar coefficient
$+(A_\theta^{\rm phys})^T/4$.
The three-component quotient coefficient $\mathsf P$ is the matrix induced by $A_\theta^{\rm phys}/2$ on $\mathbb C^4/\mathbb C(\theta,0)^T$, in local quotient coordinates independent of $x$; its explicit formula is given in Lemma~\ref{i2:lem-fixed}. Thus $D_\theta X+\mathsf{P}X=0$ has planar coefficient $-\mathsf{P}/2$. The two-component system in Lemma~\ref{i2:lem-f} is treated in the same way. All these coefficients are jointly smooth on any fixed larger disk and local real parameter chart.

Choose $L<L_1<R_{\mathrm{fr}}$ and a fixed cutoff in the leaf coordinates, equal to one near $D_{L_1}$ and compactly supported in $D_{R_{\mathrm{fr}}}$. Multiplication by this cutoff and zero extension give smooth parameter families with common compact support, agreeing with the original systems on the required disk.
Theorem~\ref{f27:main} gives local frames and their inverses for these families.
The frames are independent of $b,g,h$ and apply to every fixed test whose transverse support lies in the chosen chart. A finite subcover treats compact parameter sets. Proposition~\ref{if:binding} subsequently glues the relative fields by uniqueness on their common open shell.
\end{proof}

\begin{proposition}\label{tb27:budget}
For the near-characteristic symbols of the local right inverse, mixed derivatives of each frame and its inverse through order ten in the leaf and transverse variables and order two in the direction variables suffice for eight $x$ derivatives and two $\xi$ derivatives of the symbols. Theorem~\ref{f27:main} with total order $J=12$ supplies these bounds.
\end{proposition}
\begin{proof}
On the fixed compact support in $x$, write the moving leaf coordinates as
$Z(x,\eta)=(z_1,z_2,t)$, with local real direction coordinates $\eta$. Frequency enters through
$\eta=\eta(\xi,\tau)$ and the cutoff functions, and in the near region
$\partial_\xi^k\eta=O(\tau^{-k})$, $k=1,2$.
A frequency derivative acts by the chain rule through
\[
 (\partial_Z C)\partial_\xi Z
       +(\partial_\eta C)\partial_\xi\eta
\]
and adds at most one Cartesian derivative or one derivative in the direction parameters. After eight Cartesian derivatives and two $\xi$ derivatives, ten Cartesian derivatives and two direction derivatives therefore suffice. Derivatives of the frequency cutoffs are estimated in the local right inverse construction.
On the fixed compact subcharts and disk, $J=12$ in Theorem~\ref{f27:main} bounds all the required terms. Compactness also controls the fixed derivatives of the transverse weights and coordinate changes.
\end{proof}

\section{A local right inverse for variable coefficients}\label{r27:local-inverse}

This section constructs a source-independent linear right inverse for
the variable-coefficient conjugated operator $P_\delta$
in~\eqref{r27:generic-P}. For each sufficiently large $\tau$ and
each direction in the fixed compact parameter family, we select a
linear operator $R_\delta$ independently of the source. It satisfies
\[
 P_\delta R_\delta f=f\quad\text{in }\Omega_1,
 \qquad
 \tau\|R_\delta f\|_{L^2(\Omega_0)}
 +\|R_\delta f\|_{H^1(\Omega_0)}
 \le C\|f\|_{L^2(\Omega_0)}
\]
for every $f\in L^2(\Omega_0;\mathbb C^m)$ supported in
$K_{\rm src}$.
The construction combines the weighted free inverse $E_\delta$
with local planar transport frames $C_a$. A finite partition of the
direction parameters and a near--far decomposition in frequency
space produce a parametrix: in the near-frequency region the
transport equations cancel the leading first-order contribution,
whereas in the far-frequency region the free symbol is separated
from its characteristic set. The resulting localized error satisfies
\[
\|\mathcal T_{\tau,\gamma}\|_{L^2\to L^2}
\le C(\tau^{-\gamma}+\tau^{\gamma-1}),
\qquad 0<\gamma<1,
\]
as quantified in~\eqref{r27:error-budget}. A Neumann-series
correction gives the exact right inverse. The resulting inverse
yields exact CGO solutions for the reduced elastic system in
Subsection~\ref{pb27:section}, with simultaneous $L^2$ and $H^1$
bounds for a single remainder.

\subsection{Geometry and coefficient bounds}\label{r27:inverse-inputs}
Use $D=-i\nabla$, with the pairing convention of Section~\ref{sec:notation};
let $m$ be the matrix size, and write
\begin{equation}\label{r27:generic-P}
 \begin{gathered}
 P_\delta=(D+\delta)^2+2\sum_{j=1}^3A_j(x)(D_j+\delta_j)+q(x),\\
 \delta=\tau(\mathbf v_0+i{\mathbf v_1}),\quad \tau\geq1,\quad
 |\mathbf v_0|=|{\mathbf v_1}|=1,\quad\mathbf v_0\cdot{\mathbf v_1}=0.
 \end{gathered}
\end{equation}
The matrices $A_j,q$ are fixed smooth $m\times m$ coefficients. For $w\in\mathbb C^3$, write $A_w=\sum_jw_jA_j$; in expressions such as $w\cdot(A b-i\nabla_xb)$, this means $\sum_jw_j(A_jb-i\partial_jb)$. Extend and cut the coefficients off in $\Omega_0$, preserving them near the cutoffs in $x$ used below. Their supports may extend beyond the source ball.

First fix nested domains, a source ball, and a source cutoff:
\begin{equation}\label{r27:geometry}
\begin{aligned}
 &\overline\Omega\Subset\Omega_{\rm CGO}
   \Subset\operatorname{int}K_{\rm src},\qquad
   K_{\rm src}\Subset\Omega_1\Subset\Omega_0,\\
 &\chi_0\in C_c^\infty(\operatorname{int}K_{\rm src}),\qquad
   \chi_0=1\quad\text{near }\overline\Omega_{\rm CGO}.
\end{aligned}
\end{equation}
The domains may be taken as concentric balls, with $K_{\rm src}$
a closed ball. The compact inclusion in~\eqref{r27:geometry} gives
the cutoff $\chi_0$. Fix $\psi_2,\psi_3,\chi_{\rm sp}\in
C_c^\infty(\Omega_0)$ equal to one near $\overline\Omega_1$,
$\operatorname{supp}\psi_2$, and $\operatorname{supp}\psi_3$,
respectively. These domains and cutoffs are independent of $\tau$ and
the source.

All weighted spaces in this section are on $\mathbb R^3$ with values in
$\mathbb C^m$; pairings follow Section~\ref{sec:notation}. Put
$\langle x\rangle=(1+|x|^2)^{1/2}$ and, for $s\in\mathbb R$ and
$j=0,1$, define
\[
 \|u\|_{H^j_s}^2=\sum_{|\alpha|\le j}
 \|\langle x\rangle^s\partial^\alpha u\|_{L^2}^2,
 \qquad L^2_s=H^0_s.
\]
The weight exponent $s$ is distinct from the boundary-layer phase sign
$\mathfrak{s}$.
The free inverse is supplied by Theorem~\ref{fr27:free}.
\begin{description}
\item[The free inverse $\mathsf E$.]
Let $E_\delta$ be the selected free inverse from Theorem~\ref{fr27:free}.
Its defining identity and weighted estimate are~\eqref{fr27:right-identity}
and~\eqref{fr27:free-estimate}; both hold uniformly on the compact direction
families under consideration.
\item[Transport frames $\mathsf F$.]
Choose $R_* >\sup_{x\in\Omega_0}|x|$ and $R_*<L<M$. Near each parameter center for ordered orthonormal pairs $(\widehat{\mathbf v},{\mathbf v_1})$ and $|t|\leq R_*$, there is a smooth invertible frame $C_a$, defined near the fixed disk $D_L$, such that
\begin{equation}\label{r27:frame-input}
 i(\widehat{\mathbf v}+i{\mathbf v_1})\cdot\nabla_x C_a=A_{\widehat{\mathbf v}+i{\mathbf v_1}}C_a,
 \quad x=z_1\widehat{\mathbf v}+z_2{\mathbf v_1}+t(\widehat{\mathbf v}\times{\mathbf v_1}),\quad z=z_1+iz_2.
\end{equation}
In these coordinates the frame equation becomes
\[
 \partial_{\bar z}C_a=-\tfrac i2 A_{\widehat{\mathbf v}+i\mathbf v_1}C_a,
 \qquad 2\partial_{\bar z}
   =(\widehat{\mathbf v}+i\mathbf v_1)\cdot\nabla_x.
\]
Part~2 of Theorem~\ref{f27:main} applies with $R_{\mathrm{fr}}=M$
and $B=-\tfrac i2 A_{\widehat{\mathbf v}+i\mathbf v_1}$.
It gives the joint smoothness and the mixed bounds below for
each frame and its inverse.
\end{description}

Cover the compact parameter set by finitely many charts, with weights supported strictly inside them. Let $\eta$ be local real coordinates for the ordered orthonormal pair. The following bounds suffice:
\begin{equation}\label{r27:frame-budget}
 \max_a\max_{|\alpha|\leq10,\ |\beta|\leq2}
 \sup\left(\|\partial_{(z_1,z_2,t)}^\alpha\partial_\eta^\beta C_a\|
 +\|\partial_{(z_1,z_2,t)}^\alpha\partial_\eta^\beta C_a^{-1}\|\right)<\infty
\end{equation}
on the chosen closed subcharts and $|z|\leq L$, together with the corresponding $\partial_t^{\leq10}\partial_\eta^{\leq2}$ bounds on the weights, six Cartesian derivatives of the coefficients, eight Cartesian derivatives of the cutoffs in $x$, and two derivatives of the frequency cutoffs. Compactness controls coordinate changes. The two extra Cartesian derivatives accommodate frequency differentiation of the moving leaf coordinates.

\begin{proposition}\label{r27:local-right-inverse}
Fix the preceding data and $0<\gamma<1$. For each sufficiently large $\tau$ and each chosen direction, there is a linear operator $R_\delta$, selected independently of the source, such that every $f\in L^2(\Omega_0;\mathbb C^m)$ with $\operatorname{supp}f\subset K_{\rm src}$ satisfies
\begin{equation}\label{r27:generic-exit}
 P_\delta R_\delta f=f\quad\text{in }\Omega_1,\qquad
 \tau\|R_\delta f\|_{L^2(\Omega_0)}+\|R_\delta f\|_{H^1(\Omega_0)}
 \leq C\|f\|_2.
\end{equation}
More precisely, the threshold $\tau_0$ and constant $C$ may be chosen
from $\gamma$, the matrix size, the fixed nested domains and cutoffs
in $x$ in~\eqref{r27:geometry}, the free-inverse constant in
\eqref{fr27:free-estimate}, the coefficient and cutoff seminorms specified
above, the bounds for the frames and their inverses in
\eqref{r27:frame-budget}, and the finitely many chart weights. Once these
data are fixed, $\tau_0$ and $C$ are independent of the source $f$,
$\tau\ge\tau_0$, and the direction in the selected compact family.
\end{proposition}

\subsection{Fourier kernel estimates and weighted localization}\label{r27:kernel-estimates}
For Fourier analysis we use
$\widehat u(\xi)=\int_{\mathbb R^3}e^{-ix\cdot\xi}u(x)\,dx$,
with inverse-transform factor $(2\pi)^{-3}$.
For a symbol, $\widehat a(v,\xi)$ denotes Fourier transformation in $x$ alone. Left quantization and the finite symbol seminorms are
\[
 \operatorname{Op}(a)u(x)=(2\pi)^{-3}\int e^{ix\cdot\xi}a(x,\xi)\widehat u(\xi)\,d\xi,
 \qquad
 s_{N,r}(a)=\max_{|\alpha|\leq N,|\beta|\leq r}
 \|\partial_x^\alpha\partial_\xi^\beta a\|_\infty.
\]
Lemma~\ref{fr27:schur} gives the unweighted bounds needed here. For symbols with fixed compact support in $x$, its proof by integration by parts with $(1-\Delta_x)^3$ gives
$\|\widehat a(\eta-\xi,\xi)\|\leq C s_{6,0}(a)\langle\eta-\xi\rangle^{-6}$.
Schur's test and $\langle\eta\rangle/\langle\xi\rangle\leq C\langle\eta-\xi\rangle$ imply $\operatorname{Op}(a):H^j\to H^j$, $j=0,1$, with norm controlled by $s_{6,0}(a)$.

For symbols $a,c$ with compact support in $x$, Proposition~\ref{fr27:composition}
gives the composition-defect estimate~\eqref{fr27:composition-bound}; the
exact kernel and its matrix order are recorded in~\eqref{fr27:exact-defect-kernel}.
Thus if $s_{6,0}(\partial_{\xi_k}a)=O(\epsilon)$ for $k=1,2,3$ and
$s_{6,0}(c)=O(1)$, the defect is $O(\epsilon)$ on $H^j$, $j=0,1$.
The pure Fourier-multiplier case is included in the same proposition.

The weighted localization estimate follows directly from Proposition~\ref{fr27:weight-exchange}:
the exact identity is~\eqref{fr27:exact-weight}, and the resulting bound is
~\eqref{fr27:weighted-local}. The proposition's extension applies to the
noncompact outputs of $E_\delta$ as well.

\subsection{Frequency cutoffs and paired symbols}\label{r27:near-far-symbols}
Fix $\gamma\in(0,1)$. Choose real smooth weights $w_a(\widehat{\mathbf v},{\mathbf v_1},t)$ supported inside the frame charts, with $\sum_aw_a^2=1$ on the compact parameter set. Since their dependence on $x$ is only through $t$,
$(\widehat{\mathbf v}+i{\mathbf v_1})\cdot\nabla_xw_a=0$. Each frame is paired with its own inverse.

Write
\[
 \xi_{\mathbf v_1}=\xi\cdot{\mathbf v_1},\quad \xi'=\xi-\xi_{\mathbf v_1}{\mathbf v_1},\quad
 r=|\xi'+\tau\mathbf v_0|,\quad \widehat{\mathbf v}=(\xi'+\tau\mathbf v_0)/r.
\]
The real unit vector $\widehat{\mathbf v}$ is the moving direction
used to select the near-frequency transport frame.
Choose a fixed $\varrho\in C_c^\infty(\mathbb R)$ with $0\leq\varrho\leq1$, equal to one on $[-1,1]$ and zero outside $[-2,2]$, and define
\begin{equation}\label{r27:window}
 t_{\tau,\gamma}=\varrho(\xi_{\mathbf v_1}/\tau^\gamma)
 \varrho((r^2-\tau^2)/\tau^{1+\gamma}),\qquad
 \chi_N=\sin(\pi t_{\tau,\gamma}/2),\quad
 \chi_F=\cos(\pi t_{\tau,\gamma}/2).
\end{equation}
For sufficiently large $\tau$, on $\operatorname{supp}\chi_N$,
\begin{equation}\label{r27:near-geometry}
 r\asymp\tau,\quad |r-\tau|+|\xi_{\mathbf v_1}|\leq C\tau^\gamma,\quad |\xi|\leq4\tau.
\end{equation}
The radial bound follows by division by $r+\tau\asymp\tau$. The near cutoff vanishes near $r=0$, so the moving direction is smooth on its support. Cutoff derivatives of order $k=1,2$ are $O(\tau^{-k\gamma})$.
To verify the angular derivatives entering the frame symbols, set
$\Pi=I_3-{\mathbf v_1}{\mathbf v_1}^T$ and
$v(\xi)=\Pi\xi+\tau\mathbf v_0=r\widehat{\mathbf v}$.
For $r>0$ direct differentiation gives
\begin{equation}\label{r27:moving-direction-derivatives}
 \partial_{\xi_j}r=\widehat{\mathbf v}\cdot\Pi e_j,\qquad
 \partial_{\xi_j}\widehat{\mathbf v}
 =r^{-1}(I_3-\widehat{\mathbf v}\widehat{\mathbf v}^T)\Pi e_j.
\end{equation}
Since $r\asymp\tau$ on the near support, repeated differentiation
of~\eqref{r27:moving-direction-derivatives} gives
$\partial_\xi^k\widehat{\mathbf v}=O(\tau^{-k})$ for $k=1,2$.
The corresponding changes in the leaf coordinates are also
$O(\tau^{-k})$ on the fixed compact set in $x$.
For the cutoff calculation, with $s=(r^2-\tau^2)/\tau^{1+\gamma}$,
the full second radial derivative is
\[
 \partial_r^2\varrho(s)=4r^2\tau^{-2-2\gamma}\varrho''(s)
      +2\tau^{-1-\gamma}\varrho'(s)=O(\tau^{-2\gamma})
\]
on the transition region $r\asymp\tau$. The same bound holds for the sine and cosine cutoff profiles by the chain rule. Direction and moving-coordinate derivatives of order $k$ are $O(\tau^{-k})$.

Let the positive indices $a$ run over the finitely many chosen frame charts, and reserve $a=0$ for the far-frequency term. For $a\geq1$, set
\begin{equation}\label{r27:paired-symbols}
 b_a=\chi_{\rm sp}\chi_Nw_aC_a,\qquad
 d_a=\chi_{\rm sp}\chi_Nw_aC_a^{-1},\qquad
 b_0=d_0=\chi_F I_m.
\end{equation}
Extend the near symbols by zero outside their parameter charts; the strictly interior supports of the weights ensure smoothness. The chain rule and~\eqref{r27:frame-budget} give
\begin{equation}\label{r27:symbol-budget}
 \|\partial_x^\alpha\partial_\xi^\beta b_a\|_\infty+
 \|\partial_x^\alpha\partial_\xi^\beta d_a\|_\infty
 \leq C\tau^{-\gamma|\beta|},\qquad |\alpha|\leq8,\quad |\beta|\leq2.
\end{equation}
The index-zero symbols are pure multipliers. On $\operatorname{supp}\psi_3$ we have the exact identities
\begin{equation}\label{r27:paired-identity}
 \sum_{a\geq0}b_ad_a=I_m,
 \qquad (\widehat{\mathbf v}+i{\mathbf v_1})\cdot(A b_a-i\nabla_xb_a)=0\quad(a\geq1).
\end{equation}
The second identity uses $\chi_{\rm sp}=1$ on a neighborhood of this support.

For any $g\in L^2(\mathbb R^3;\mathbb C^m)$, define
\begin{equation}\label{r27:parametrix}
 H_ag=\psi_3\operatorname{Op}(d_a)g,
 \qquad Q_{\tau,\gamma}g=\sum_{a\geq0}\operatorname{Op}(b_a)E_\delta H_ag.
\end{equation}
The operators $H_a:L^2\to L^2_1$ are uniformly bounded and have output supported in $\operatorname{supp}\psi_3$. Indeed, the ordinary Fourier--Schur estimate and
$\sup_{\supp\psi_3}\langle x\rangle<\infty$ give, for every near index,
\begin{equation}\label{r27:source-localization}
 \|H_ag\|_{L^2_1}
 \leq\sup_{\supp\psi_3}\langle x\rangle
       \|\psi_3\|_\infty\, C s_{6,0}(d_a)\|g\|_2
 \leq C_{\rm src}\|g\|_2.
\end{equation}
For $a=0$, the identical conclusion follows from Plancherel for
the bounded multiplier $d_0=\chi_FI_m$, followed by multiplication
by $\psi_3$. These bounds hold for arbitrary $L^2$ sources $g$; the cutoff in $x$ is applied to the \emph{output} of $\operatorname{Op}(d_a)$.
Proposition~\ref{fr27:composition} and its multiplier case give
\begin{equation}\label{r27:S-small}
 S=\sum_a\operatorname{Op}(b_a)\psi_3\operatorname{Op}(d_a)-\psi_3I_m,
 \qquad \|S\|_{L^2\to L^2}\leq C\tau^{-\gamma}.
\end{equation}
Here the right symbol is the compactly supported $\psi_3d_a$, and $\sum_a b_a\psi_3d_a=\psi_3I_m$.

The differential composition in the parametrix follows from
the exact left-quantization identities
\begin{align}
 (D_j+\delta_j)\operatorname{Op}(b)
    &=\operatorname{Op}\bigl((\xi_j+\delta_j)b-i\partial_j b\bigr),\notag\\
 (D+\delta)^2\operatorname{Op}(b)
    &=\operatorname{Op}\bigl(p_\delta b
           -2i(\xi+\delta)\cdot\nabla_x b-\Delta_x b\bigr).
 \label{r27:exact-symbol-product}
\end{align}
Here the differential operators act on the output variable $x$;
the symbols $b$ are matrix-valued and are never commuted across the
matrices $A_j$.  In particular, multiplying the first identity by
$2A_j(x)$ on the left and then adding $q(x)\operatorname{Op}(b)$
gives, preserving matrix order,
\begin{align}
 P_\delta\operatorname{Op}(b_a)E_\delta
 &=\operatorname{Op}(b_a)+\operatorname{Op}(r_{0,a})E_\delta
                       +\operatorname{Op}(r_{1,a})E_\delta,\label{r27:differential-identity}\\
 r_{0,a}&=-\Delta_xb_a-2i\sum_jA_j\partial_jb_a+qb_a,\notag\\
 r_{1,a}&=2\sum_j(\xi_j+\delta_j)(A_jb_a-i\partial_jb_a).\notag
\end{align}
Define the differential residual operator by
\[
 E_{\mathrm{err}}=\sum_{a\ge0}\bigl[\operatorname{Op}(r_{0,a})E_\delta H_a+
                         \operatorname{Op}(r_{1,a})E_\delta H_a\bigr].
\]
Thus $S$ in \eqref{r27:S-small} is the composition defect and $E_{\mathrm{err}}$ is the residual from differentiating the paired symbols. Substituting \eqref{r27:differential-identity} into \eqref{r27:parametrix} gives
\begin{equation}\label{r27:PQ}
 P_\delta Q_{\tau,\gamma}=\psi_3I_m+S+E_{\mathrm{err}}.
\end{equation}
For fixed $\tau$, the identity first holds on smooth sources and extends by the local estimates below and continuity in distributions.

\subsection{Near and far residual estimates}\label{r27:residual-norm}
The seminorms $s_{6,2}(\psi_3r_{0,a})$ are uniformly bounded; six derivatives of $-\Delta b_a$ require eight Cartesian derivatives of $b_a$. Applying~\eqref{fr27:free-estimate} followed by~\eqref{fr27:weighted-local} gives
\begin{equation}\label{r27:r0-bound}
 \|\operatorname{Op}(\psi_3r_{0,a})E_\delta H_a\|_{L^2\to L^2}
 \leq C\tau^{-1}.
\end{equation}
In the near region, the second identity in~\eqref{r27:paired-identity} cancels $\tau(\widehat{\mathbf v}+i{\mathbf v_1})$, leaving exactly
\begin{equation}\label{r27:remaining-vector}
 \xi+\delta-\tau(\widehat{\mathbf v}+i{\mathbf v_1})=(r-\tau)\widehat{\mathbf v}+\xi_{\mathbf v_1}{\mathbf v_1}.
\end{equation}
On the support of the near symbols, put
$R(\xi)=(r-\tau)\widehat{\mathbf v}+\xi_{\mathbf v_1}\mathbf v_1$.
The estimates~\eqref{r27:near-geometry}, $r\asymp\tau$, and
$\partial_{\xi}\widehat{\mathbf v}=O(\tau^{-1})$ give
\[
 |R|\le C\tau^\gamma,\qquad
 |\partial_{\xi}R|\le C,\qquad
 |\partial_{\xi}^2R|\le C\tau^{-1}.
\]
Here $|\partial_{\xi}^qR|$ denotes the maximum of the
component derivatives of exact order $q$. On $\operatorname{supp}\psi_3$, where
the transport cancellation is exact, the near residual is
\[
 r_{1,a}=2R(\xi)\cdot(A b_a-i\nabla_xb_a),\qquad a\ge1.
\]
Every frequency derivative of order at most two falls on $R$ or on
a near symbol. Applying~\eqref{r27:symbol-budget} and the product
rule, for $|\alpha|\le6$, yields
\[
 \|\partial_x^\alpha\partial_\xi^\beta(\psi_3r_{1,a})\|_\infty
 \leq C\tau^{\gamma-\gamma|\beta|},\qquad |\alpha|\leq6,\quad |\beta|\leq2,
 \quad a\geq1.
\]
Hence $s_{6,2}(\tau^{-\gamma}\psi_3r_{1,a})=O(1)$. The free estimate followed by weighted localization yields
\begin{equation}\label{r27:near-residual}
 \|\operatorname{Op}(\psi_3r_{1,a})E_\delta H_a\|_{L^2\to L^2}
 \leq C\tau^{\gamma-1},\qquad a\geq1.
\end{equation}
This estimate differentiates only the smooth symbols.

For the far region, the identity
\[
 p_\delta=\xi_{\mathbf v_1}^2+r^2-\tau^2+2i\tau\xi_{\mathbf v_1}
\]
gives, on $\operatorname{supp}\chi_F$,
\begin{equation}\label{r27:far-gap}
 |p_\delta(\xi)|\geq
 \begin{cases}
 \tfrac12\tau^{1+\gamma},& |\xi|\leq4\tau,\\
 (1-1/\sqrt2)|\xi|^2,& |\xi|\geq4\tau.
 \end{cases}
\end{equation}
On the far support, either $|\xi_{\mathbf v_1}|\geq\tau^\gamma$ or $|r^2-\tau^2|\geq\tau^{1+\gamma}$. The first case follows from the imaginary part. Otherwise $\xi_{\mathbf v_1}^2\leq\tau^{2\gamma}$, so the real part has modulus at least $\tau^{1+\gamma}-\tau^{2\gamma}\geq\tfrac12\tau^{1+\gamma}$. For high frequencies use $|2\delta\cdot\xi|\leq2\sqrt2\tau|\xi|$.

Set
\begin{equation}\label{r27:far-symbol}
 a_F(x,\xi)=\frac{2\sum_jA_j(x)(\xi_j+\delta_j)\chi_F(\xi)}{p_\delta(\xi)}.
\end{equation}
The far cutoff vanishes near the characteristic set, so the symbol extends smoothly by zero there. For $|\xi|\leq4\tau$, direct differentiation yields
\[
 \partial_{\xi_j}p_\delta=2(\xi_j+\delta_j),\qquad
 \partial_{\xi_j}\partial_{\xi_k}p_\delta=2\delta_{jk}.
\]
On the far support in this bounded-frequency region,
\eqref{r27:far-gap} and $|\partial_\xi p_\delta|\le C\tau$
give, for $|\beta|\le2$,
\[
 |\partial_\xi^\beta(p_\delta^{-1})|
 \le C\tau^{-1-\gamma-\gamma|\beta|}.
\]
Indeed the first derivative is
$-(\partial_\xi p_\delta)p_\delta^{-2}$, while the second
is a sum of $(\partial_\xi p_\delta)^2p_\delta^{-3}$
and $(\partial_\xi^2p_\delta)p_\delta^{-2}$; the latter is
bounded by the displayed power because $0<\gamma<1$.
The numerator and its first two frequency derivatives, including the
cutoff, have bounds $C\tau^{1-\gamma|\beta|}$.
Consequently, the product rule gives
\begin{equation}\label{r27:far-seminorm}
 \|\partial_x^\alpha\partial_\xi^\beta(\psi_3a_F)\|_\infty
 \leq C\tau^{-\gamma(1+|\beta|)},\qquad |\alpha|\leq6,\quad|\beta|\leq2.
\end{equation}
For $|\xi|\ge4\tau$, one has $\chi_F=1$ and
$|p_\delta|\ge(1-1/\sqrt2)|\xi|^2$; direct differentiation
gives $|\partial_\xi^\beta a_F|\le C|\xi|^{-1-|\beta|}$,
which is bounded by the right side of~\eqref{r27:far-seminorm}
for $0<\gamma<1$. The bounds hold in every transition region, and
all frequency derivatives in the far region fall on the smooth symbol
factors, while the singular multiplier is used at zeroth order on the
characteristic set. Since the free inverse is a multiplier,
$\operatorname{Op}(r_{1,0})E_\delta=\operatorname{Op}(a_F)$ exactly. The ordinary kernel estimate bounds its localization composed with $H_0$ by $O(\tau^{-\gamma})$.

Combining the estimates gives the localized error on $L^2(\mathbb R^3)$:
\begin{equation}\label{r27:error-budget}
 \mathcal T_{\tau,\gamma}=\psi_3(S+E_{\mathrm{err}}),\qquad
 \|\mathcal T_{\tau,\gamma}\|_{L^2\to L^2}
 \leq C\bigl(\tau^{-\gamma}+\tau^{\gamma-1}\bigr).
\end{equation}
The derivatives in the transition region of $\chi_{\rm sp}$ disappear after localization because $\chi_{\rm sp}=1$ near $\operatorname{supp}\psi_3$. Both terms tend to zero for $0<\gamma<1$; the choice $\gamma=1/2$ optimizes their common decay rate.

For each index put $a_{2,a}=\psi_2b_a$. This symbol has compact support in $x$, even for the far index $a=0$, and
$\sup_\tau s_{6,2}(a_{2,a})<\infty$ by
\eqref{r27:symbol-budget} and the cutoff estimates. Applying
\eqref{fr27:weighted-local} to $a_{2,a}$, the \emph{same} free field
$E_\delta H_ag$ satisfies
\[
 \tau\|\operatorname{Op}(a_{2,a})E_\delta H_ag\|_2
 +\|\operatorname{Op}(a_{2,a})E_\delta H_ag\|_{H^1}
 \leq C\|H_ag\|_{L^2_1}
 \leq C\|g\|_2.
\]
The last step uses~\eqref{r27:source-localization} (or Plancherel
for $a=0$). As the parameter cover has a fixed finite number of
indices and $\psi_2\operatorname{Op}(b_a)=\operatorname{Op}(\psi_2b_a)$,
summing proves
\begin{equation}\label{r27:Q-two-norms}
 \tau\|\psi_2 Q_{\tau,\gamma}g\|_2+
 \|\psi_2 Q_{\tau,\gamma}g\|_{H^1}\leq C\|g\|_2.
\end{equation}
In particular, the $H^1$ bound uses the weighted derivative estimate
for the selected free inverse and is controlled by $\|g\|_{L^2}$ for
arbitrary $L^2$ sources.

\subsection{The Neumann inverse and support properties}\label{r27:selected-inverse}
Let $C_{\rm err}$ be a constant for which the right side of
\eqref{r27:error-budget} is at most
$C_{\rm err}(\tau^{-\gamma}+\tau^{\gamma-1})$ on the selected
compact direction family. Choose $\tau_0$ so that
\begin{equation}\label{r27:neumann-threshold}
 C_{\rm err}(\tau_0^{-\gamma}+\tau_0^{\gamma-1})\le\tfrac12.
\end{equation}
Both exponents are negative, so the same inequality holds for every
$\tau\ge\tau_0$. Define on $L^2(\mathbb R^3)$
\begin{equation}\label{r27:selected-R}
 \mathcal G_{\tau,\gamma}=(I+\mathcal T_{\tau,\gamma})^{-1}
 =\sum_{k=0}^\infty(-\mathcal T_{\tau,\gamma})^k,
 \qquad R_\delta=\psi_2Q_{\tau,\gamma}\mathcal G_{\tau,\gamma}.
\end{equation}
Here $\mathcal G=\mathcal G_{\tau,\gamma}$ and
$\mathcal T=\mathcal T_{\tau,\gamma}$. Estimate
\eqref{r27:neumann-threshold} makes the series converge in operator
norm, with $\|\mathcal G\|_{2\to2}\le2$. For a source $f$,
extended by zero, set $g=\mathcal Gf$. Then
$\|g\|_2\le2\|f\|_2$. Because $\mathcal T$ has output supported
in $\operatorname{supp}\psi_3$ and
$\operatorname{supp}f\subset K_{\rm src}\subset\{\psi_3=1\}$,
the identity $g=f-\mathcal Tg$ implies
\begin{equation}\label{r27:supports}
 \begin{split}
 &\operatorname{supp}f\subset K_{\rm src},\qquad
 \operatorname{supp}g\subset\operatorname{supp}\psi_3,\qquad
 \operatorname{supp}H_ag\subset\operatorname{supp}\psi_3,\\
 &E_\delta H_ag\text{ may have a noncompact tail},\qquad
 \operatorname{supp}R_\delta f\subset\operatorname{supp}\psi_2\Subset\Omega_0.
 \end{split}
\end{equation}
The construction localizes the final output while retaining the free tail in the preceding estimates.

Since $\psi_2=\psi_3=1$ in a neighborhood of $\Omega_1$,
formula~\eqref{r27:PQ} and $g+\mathcal Tg=f$ give
\[
 P_\delta R_\delta f
 =g+(S+E_{\mathrm{err}})g=g+\mathcal Tg=f\qquad\text{in }\Omega_1.
\]
With $C_{\rm para}$ denoting the constant in
\eqref{r27:Q-two-norms}, the field $R_\delta f=\psi_2Q_{\tau,\gamma}g$
satisfies
\[
 \tau\|R_\delta f\|_2+\|R_\delta f\|_{H^1}
 \le C_{\rm para}\|g\|_2\le2C_{\rm para}\|f\|_2.
\]
The displayed equation holds on all of $\Omega_1$, whereas the norms
are taken on $\Omega_0$. The free outputs $E_\delta H_ag$ retain their possible noncompact
tails, while the final cutoff $\psi_2$ localizes the right-inverse
output inside $\Omega_0$. Since the frames, compact subcharts, and cutoffs were fixed
before $f$ and $\tau$ were chosen, these constants depend only on the
finite data stated in Proposition~\ref{r27:local-right-inverse}.

\subsection{Complex geometrical optics solutions for the reduced elastic system}\label{pb27:section}
Set
\[
 c=\mu^{-1/2},\quad \kappa=\frac{\lambda+\mu}{\lambda+2\mu},
 \quad \mathbf{m}=\nabla\log\mu,\quad d=\kappa\sqrt\mu,
\quad e^{\mathrm{red}}=\tfrac12(3\kappa-2)\sqrt\mu.
\]
The reduced system is
\begin{equation}\label{pb27:MW}
 M^W=\Delta I_4+\sum_{k=1}^3A_k^{\rm phys}\partial_k+V^W,
 \qquad
 A_k^{\rm phys}=\begin{pmatrix}-\mathbf{m} e_k^T&F_\mu e_k\\d e_k^T&0\end{pmatrix},
\end{equation}
Here $e_1,e_2,e_3$ are the standard coordinate vectors in $\mathbb R^3$. The matrices $F_\mu,V_\mu$ and the scalars $e^{\mathrm{red}},q_\mu$ are defined in \eqref{bt:MW}. Collecting its zeroth-order terms gives
\[
 V^W=\begin{pmatrix}
 V_\mu+(W/\mu)I_3&F_\mu\mathbf m+\mu^{-3/2}(2W\mathbf m-\nabla W)\\
 -e^{\mathrm{red}}\mathbf m^T&q_\mu+(1-\kappa)W/\mu
 \end{pmatrix}.
\]
Thus $V^W$ is a $4\times4$ multiplication matrix, including the
potential and its gradient. The extensions in
Propositions~\ref{ext:common} and~\ref{ext:finite-common} supply
the constant background Lam\'e pair and $W=0$ sufficiently far away.

\begin{lemma}\label{pb27:background}
Let $\lambda_*,\mu_*$ be the constant exterior Lam\'e coefficients
from Proposition~\ref{ext:common} or~\ref{ext:finite-common}, and put $\kappa_*=(\lambda_*+\mu_*)/(\lambda_*+2\mu_*)$ and $d_*=\kappa_*\sqrt{\mu_*}$. With $e_4=(0,0,0,1)^T\in\mathbb C^4$, define
\[
 A_{k,*}=d_*e_4(e_k^T,0),\quad
 N(x)=\sum_kx_kA_{k,*},\quad
 \mathcal B_{\rm bg}(x)=I_4-\tfrac12N(x).
\]
Then $N^2=0$, $\det\mathcal B_{\rm bg}=1$, $\mathcal B_{\rm bg}^{-1}=I_4+N/2$, and
\[
 M_*:=\Delta I_4+\sum_kA_{k,*}\partial_k,
 \qquad \mathcal B_{\rm bg}^{-1}M_*\mathcal B_{\rm bg}=\Delta I_4.
\]
If $A_k^{\rm phys}-A_{k,*}$ and $V^W$ in~\eqref{pb27:MW} are smooth and compactly supported, then all lower-order coefficients of
\begin{equation}\label{pb27:compact}
 \widetilde M^W=\mathcal B_{\rm bg}^{-1}M^W\mathcal B_{\rm bg}
 =\Delta I_4+\sum_k\widetilde A_k\partial_k+\widetilde V^W
\end{equation}
are smooth and compactly supported. On each fixed working ball, $\mathcal B_{\rm bg}^{\pm1}$ and all their finite-order derivatives are bounded independently of $h$.
\end{lemma}
\begin{proof}
Since $A_{j,*}A_{k,*}=0$ for all $j,k$, the assertions about $N$ and the inverse follow. The product rule gives
\[
 \widetilde A_k=\mathcal B_{\rm bg}^{-1}(2\partial_k\mathcal B_{\rm bg}+A_k^{\rm phys}\mathcal B_{\rm bg}),
 \quad
 \widetilde V^W=\mathcal B_{\rm bg}^{-1}
 (\Delta\mathcal B_{\rm bg}+\sum_k A_k^{\rm phys}\partial_k\mathcal B_{\rm bg}+V^W\mathcal B_{\rm bg}).
\]
In the background region, $\partial_k\mathcal B_{\rm bg}=-A_{k,*}/2$ and $\Delta\mathcal B_{\rm bg}=0$, and both parentheses vanish. Thus the conjugated coefficients are supported in the original nonconstant region. The remaining bounds follow from the affine formulas for $\mathcal B_{\rm bg}^{\pm1}$.
\end{proof}

\begin{proposition}\label{pb27:rightinverse}
For the operator in~\eqref{a27:generic}, let $R_\delta$ be the source-independent
linear right inverse selected in Proposition~\ref{r27:local-right-inverse};
for sources supported in $K_{\rm src}$ it satisfies~\eqref{r27:generic-exit}.
For $\theta=a+ib$ with real orthonormal $a,b$, make the substitution
\begin{equation}\label{pb27:signs}
 A_k=\frac{\widetilde A_k}{2i},\quad q=-\widetilde V^W,
 \quad\delta=\frac{i\theta}{h}=\frac{-b+ia}{h},\quad\tau=h^{-1}.
\end{equation}
Then
\[
 \mathcal R_{h,\theta}^{(0)}F
 :=-\mathcal B_{\rm bg} R_\delta(\mathcal B_{\rm bg}^{-1}F)
\]
defines a linear map from $L^2(\Omega_0;\mathbb C^4)$ sources supported in $K_{\rm src}$ to $H^1(\Omega_0;\mathbb C^4)$, for $0<h\le\tau_0^{-1}$. For every such $F$,
\begin{equation}\label{pb27:physicalR}
 (e^{\theta\cdot x/h}M^We^{-\theta\cdot x/h})\mathcal R_{h,\theta}^{(0)}F=F
 \quad\text{in }\Omega_1,
 \qquad h^{-1}\|\mathcal R_{h,\theta}^{(0)}F\|_{L^2(\Omega_0)}
 +\|\mathcal R_{h,\theta}^{(0)}F\|_{H^1(\Omega_0)}\le C'\|F\|_{L^2(\Omega_0)}.
\end{equation}
\end{proposition}
\begin{proof}
Since $\theta\cdot\theta=0$, the convention $D=-i\nabla$ and the
bilinear pairing of Section~\ref{sec:notation} give
\[
 P_\delta=-e^{\theta\cdot x/h}\widetilde M^We^{-\theta\cdot x/h}.
\]
The scalar exponential commutes with $\mathcal B_{\rm bg}$. Multiplication by $\mathcal B_{\rm bg}^{-1}$ preserves source support, and multiplication by $\mathcal B_{\rm bg}^{\pm1}$ is bounded on the fixed $L^2$ and $H^1$ domains independently of $h$. Both estimates therefore hold for the same output. Replacing $\theta$ by $-\theta$ gives the opposite phase, with uniform constants on compact direction subcharts.
\end{proof}

For the moving directions in the local right inverse, write $w=\widehat{\mathbf v}+i{\mathbf v_1}$ and let $z_g$ be its leaf coordinate. For any complex direction $w$, use the contractions $A_w=\sum_kw_kA_k$ and $\widetilde A_w=\sum_kw_k\widetilde A_k$. Here $A_k$ is from \eqref{a27:generic} and $\widetilde A_k$ from \eqref{pb27:compact}. The frame equation is
$i w\cdot\nabla C=A_w C$. In the leaf coordinate $z_g$, the
identity $w\cdot\nabla=2\partial_{\bar z_g}$ gives
\begin{equation}\label{pb27:moving-frame-equation}
 \partial_{\bar z_g}C=-\tfrac i2 A_w C .
\end{equation}
Substituting~\eqref{pb27:signs} into \eqref{pb27:moving-frame-equation}
gives the coefficient
\[
 B_g=-\tfrac i2 A_w=-\tfrac14\widetilde A_w.
\]
At a fixed phase $\theta=a+ib$, with real orthonormal $a,b$, let $z_{\rm phys}=a\cdot x+ib\cdot x$. The central moving direction is $w_0=i\theta=-b+ia$.
Thus $z_g=-b\cdot x+ia\cdot x=i z_{\rm phys}$, so
$\partial_{\bar z_g}=i\partial_{\bar z_{\rm phys}}$, while
$B_g=-i\widetilde A_\theta/4$.
Dividing the equation by $i$ gives, in the phase coordinates,
$\partial_{\bar z_{\rm phys}}C=-\widetilde A_\theta C/4$.
Thus the moving-direction frame is a leading transport frame for the conjugated operator $\widetilde M^W$. Multiplication by $\mathcal B_{\rm bg}$ gives a leading transport frame for $M^W$, since
\[
 (2D_\theta+A_\theta^{\rm phys})(\mathcal B_{\rm bg} C)
 =\mathcal B_{\rm bg}(2D_\theta+\widetilde A_\theta)C=0.
\]

\begin{remark}\label{pb27:paths}
If a null path $\theta(h)$ has real and imaginary parts of common length $a(h)>0$, bounded above and below, apply the result with $\widehat\theta(h)=\theta(h)/a(h)$ and $\widehat h=h/a(h)$. The estimates retain their stated orders in $h$. The paths $\vartheta_j(h)$ in Section~\ref{sec:potential} are normalized and lie in fixed compact subcharts near their limiting directions; the uniform constants above apply to them.
\end{remark}

\begin{proposition}\label{pb27:remainder}
Fix $\chi_0\in C_c^\infty(\operatorname{int}K_{\rm src})$ equal
to one near $\overline\Omega_{\rm CGO}$.  Choose a fixed open set $U$
such that
$\operatorname{supp}\chi_0\cup\overline\Omega_{\rm CGO}
\subset U\Subset\Omega_1$.
Let $\mathfrak{s}\in\{+1,-1\}$ and let $\theta$ be a normalized null
direction in one of the fixed compact direction families.  Suppose
smooth amplitudes near $U$ satisfy
\begin{equation}\label{pb27:two-sign-transport}
 \mathsf T_\theta a_0=0,\qquad
 \mathsf T_\theta a_1=\mathfrak{s} M^Wa_0,
 \qquad \mathsf T_\theta=2\theta\cdot\nabla+A_\theta^{\rm phys}.
\end{equation}
Assume $\|a_1\|_{H^2(U)}$ is uniformly bounded along the selected
parameter family.  Put
\begin{equation}\label{pb27:two-sign-correction}
 F=\chi_0M^Wa_1,\qquad
 r_{h,\mathfrak{s}}=-\mathcal R_{h,\mathfrak{s}\theta}^{(0)}(hF),\qquad
 Z_{h,\mathfrak{s}}=e^{-\mathfrak{s}\theta\cdot x/h}
                 (a_0+ha_1+r_{h,\mathfrak{s}}).
\end{equation}
Then $M^WZ_{h,\mathfrak{s}}=0$ in $\Omega_{\rm CGO}$, and this single
remainder satisfies
\begin{equation}\label{pb27:two}
 \|r_{h,\mathfrak{s}}\|_{L^2(\Omega_0)}\le Ch^2,
 \qquad
 \|r_{h,\mathfrak{s}}\|_{H^1(\Omega_0)}\le Ch.
\end{equation}
The constants are uniform for either sign on the fixed compact
parameter family.  The case $\mathfrak{s}=+1$ is the negative physical phase
$e^{-\theta\cdot x/h}$; the case $\mathfrak{s}=-1$ is the positive phase
$e^{\theta\cdot x/h}$.
\end{proposition}
\begin{proof}
For $\mathfrak{s}=\pm1$, the null identity $\theta\cdot\theta=0$
and the product rule give the exact conjugation formula
\begin{equation}\label{pb27:two-sign-conjugation}
 e^{\mathfrak{s}\theta\cdot x/h}M^We^{-\mathfrak{s}\theta\cdot x/h}
     =M^W-\frac{\mathfrak{s}}{h}\mathsf T_\theta.
\end{equation}
The transport conditions~\eqref{pb27:two-sign-transport} imply
\begin{equation}\label{pb27:two-sign-residual}
 \left(M^W-\frac{\mathfrak{s}}{h}\mathsf T_\theta\right)
                       (a_0+ha_1)
 =M^Wa_0+hM^Wa_1-\mathfrak{s}\mathsf T_\theta a_1
 =hM^Wa_1.
\end{equation}
The compact support of $\chi_0$ in $U$ permits extension of $F$
by zero.  Since $M^W:H^2(U)\to L^2(U)$ is bounded,
\[
 \|F\|_2\le C\|a_1\|_{H^2(U)}.
\]
Proposition~\ref{pb27:rightinverse}, applied with the normalized
direction $\mathfrak{s}\theta$, shows that
\[
 \left(M^W-\frac{\mathfrak{s}}{h}\mathsf T_\theta\right)
                     r_{h,\mathfrak{s}}=-h\chi_0M^Wa_1
                     \quad\text{in }\Omega_1.
\]
Adding this identity to \eqref{pb27:two-sign-residual} gives
\[
 e^{\mathfrak{s}\theta\cdot x/h}M^WZ_{h,\mathfrak{s}}
 =h(1-\chi_0)M^Wa_1\quad\text{in }U.
\]
Since $\chi_0=1$ near $\overline\Omega_{\rm CGO}$, this yields
$M^WZ_{h,\mathfrak{s}}=0$ in $\Omega_{\rm CGO}$. The same
right-inverse output in
\eqref{pb27:two-sign-correction} obeys, by
\eqref{pb27:physicalR},
\[
 h^{-1}\|r_{h,\mathfrak{s}}\|_{L^2(\Omega_0)}
 +\|r_{h,\mathfrak{s}}\|_{H^1(\Omega_0)}
 \le Ch\|F\|_2,
\]
which implies both estimates in~\eqref{pb27:two}.
The norm of the inverse and the bound for $F$ are uniform over the
chosen compact families, including their negatives.  For each fixed
$h$, local elliptic regularity makes $Z_{h,\mathfrak{s}}$ smooth in the
target whenever the coefficients and amplitudes are smooth.
\end{proof}

For the fixed smooth material, the simultaneous estimates
\eqref{pb27:two} give
\[
 C\big(h\|a_1\|_{H^1(U)}
      +\|r_h\|_{H^1(\Omega_0)}
      +h^{-1}\|r_h\|_{L^2(\Omega_0)}\big)=O(h).
\]
Lemma~\ref{i2:lem-shell} applies this bound to the first-order
representation of the unscaled third column. The estimate holds
throughout the target, including all transverse cutoff transitions,
because both norms refer to the same remainder and the lift identity
is applied to the complete exact solution.

\section{The full matrix moment identities}
\label{i2:sec}

This section derives the complete $3\times3$ matrix of polynomial
moment identities from the reduced Green--Betti identity
\eqref{i2:eq-betti}. The primal and dual quotient frames are realized
by physical elastic solutions with compatible high-frequency
normalizations. We first construct the transport corrections
(Subsection~\ref{i2:preliminaries}) and derive the asymptotics of the
canonical lift (Subsection~\ref{i2:canonical-lift}).
Subsection~\ref{i2:polarizations} constructs two leading columns
using scaled lifts and a third using an unscaled lift with
$J_\theta A_0^S=0$. The determinant computation in
Lemma~\ref{i2:lem-fixed} proves invertibility independently of the
scalar entries of the third correction. The quotient frames are
fixed before the polynomial and transverse tests are chosen;
Subsection~\ref{i2:moment-extraction} then extracts the leafwise
matrix identities of Theorem~\ref{i2:thm-nine}.

The first two physical testing columns are normalized by $-h$,
whereas the third is unscaled. Both limits therefore rely on the
same remainder satisfying
\[
\|r_h\|_{L^2}=O(h^2),
\qquad
\|r_h\|_{H^1}=O(h),
\]
as stated in~\eqref{pb27:two}. All pairings and transposes below
are complex-bilinear, with the convention of
Section~\ref{sec:notation}.

\subsection{Preliminaries and first-correction estimates}
\label{i2:preliminaries}

\begin{definition}\label{i2:def-input}

The two materials are described by the coefficient triples
\[
(\lambda_j,\mu_j,W_j),\qquad j=1,2,
\]
on the common physical domain $\Omega$.
For formulas involving one fixed material, we suppress the
index $j$ and write
\[
\lambda,\mu,W,\kappa,a,c,\ell,\mathbf m,d,e^{\mathrm{red}},
q_\mu,F_\mu,V_\mu,M^W,\mathcal W
\]
for the corresponding coefficients and derived quantities,
with the meanings fixed in~\eqref{en:W-definition},
\eqref{bt:notation}, and~\eqref{bt:MW}.
Write $M=M^0$ for the zero-potential part of the
reduced operator.

Use the leaf coordinates $z,x_\perp$, the directional derivative
$D_\theta$, and the leading transport operator $\mathsf T_\theta$
fixed in Section~\ref{tb27:section}. The coefficient contraction
$A_\theta^{\rm phys}$ is determined by~\eqref{pb27:MW}.
The leafwise right inverse $\mathcal C_\theta$ of $D_\theta$
is defined by the Cauchy formula
\[
 (\mathcal C_\theta H)(z,t)=\frac1{2\pi}\int_{\mathbb C}
                         \frac{H(w,t)}{z-w}\,dA(w),\qquad t=x_\perp.
\]
Here $H$ has compact support in the leaf variable and $dA$ is planar area
measure. In applications a leaf cutoff is included in $H$. Thus
$D_\theta\mathcal C_\theta H=H$ on the working disk, and the operator
leaves the transverse parameter $t$ fixed.
\end{definition}

\begin{proposition}\label{a28:moment-inputs}
Suppose the boundary and materials are smooth, and that both
physical Dirichlet realizations are nonresonant and the physical
Dirichlet-to-Neumann maps agree. The following constructions and
identities are available. The finite-order counterparts used for
Theorem~\ref{a28:main} are verified in
Subsection~\ref{fr:interior-assembly} on the same fixed working sets.
\begin{enumerate}
\item Theorem~\ref{f27:main} and Proposition~\ref{tb27:F-supply} supply smooth
invertible transport frames on the prescribed leaf disk, with all required
finite derivatives of each frame and its inverse bounded on compact parameter
subcharts.
\item Fix the nested domains and source set in~\eqref{r27:geometry}. For each material, use the selected right inverse
$\mathcal R_{h,\theta}^{(0)}$ of Proposition~\ref{pb27:rightinverse};
its equation and simultaneous estimate are given in~\eqref{pb27:physicalR}
for sources supported in $K_{\rm src}$. The phase $e^{-\theta\cdot x/h}$
uses this operator, and the opposite phase uses the operator selected at
$-\theta$, as in Proposition~\ref{pb27:remainder}. The exact solutions
are used on $\Omega_{\rm CGO}$, with
$\overline\Omega\Subset\Omega_{\rm CGO}\Subset\Omega_1$.
\item Theorem~\ref{bl:jets} and Proposition~\ref{ext:common} give smooth common
extensions of the two coefficient sets, with differences supported in
$\overline\Omega$. Their Dirichlet-to-Neumann maps on the original domain
agree. For every exact reduced solution $M_1^{W_1}z_1=0$ and elastic solution
$(L_2+W_2I)u_2=0$, the field $y_2=S_2^Tu_2$ satisfies
\begin{equation}\label{i2:eq-betti}
 \int_\Omega y_2^T(M_1^{W_1}-M_2^{W_2})z_1\,dx=0.
\end{equation}
\end{enumerate}
The source cutoff, target domain, and transport region are fixed nested
neighborhoods. Constants may depend on each fixed test amplitude and its
polynomial degree.

For the later use of Proposition~\ref{if:actualH}, fix its complete
working leaf disk $D_L$, with $L>R$, and the bounded transverse
parameter ranges in advance. In choosing the concentric auxiliary
domains in~\eqref{r27:geometry}, take $U$ large enough that
\[
\{x:z(x)\in\overline D_L,\ x_\perp\in\overline I\}
\Subset U
\]
for every selected compact direction chart and local transverse
interval $I$ within these prescribed transverse ranges. All these
cylinders lie in a fixed ball, and the
coefficients have already been extended to $\mathbb R^3$, so this
choice is possible before constructing the right inverses and frames.
Choose the leaf cutoffs to equal one near $\pi_\theta(\overline U)$
and the finite transverse partitions to sum to one near
$\pi_\perp(\overline U)$, as in Lemma~\ref{a31:first-correction}.
The uncut frame equations constructed in Lemma~\ref{i2:lem-fixed}
then hold on the whole of $D_L\times I$, before the test-dependent
transverse localization is performed.

\end{proposition}
\begin{proof}
Theorem~\ref{bl:jets} and Proposition~\ref{ext:common} give the common smooth
extension and compact support of all coefficient differences. The boundary
value of $W$ is recovered at $j=2$ in Theorem~\ref{bl:jets}. Together with the
value and first normal derivative of $\mu$, this gives the matching conditions
of Proposition~\ref{bt:weakprop}, and hence~\eqref{i2:eq-betti} from the
original boundary data.

Theorem~\ref{f27:main} and Proposition~\ref{tb27:F-supply} supply the
required frames and their inverses, jointly smoothly in local real parameters.
Use the bounds for the leaf primitive $\mathcal C_\theta$ of
Definition~\ref{i2:def-input} in Lemma~\ref{a31:first-correction}.
The disk, chart, frames, and primitive are independent of $b,g,h$.
For each $b$ supported in the chart and each $g$, the transverse
cutoffs, amplitudes, and corrections may depend on $b,g$ and
are fixed as $h\to0$.

Use the selected right inverse of Proposition~\ref{pb27:rightinverse}.
Choose the first correction on the common leaf disk by
Lemma~\ref{a31:first-correction}.
The test amplitudes have compact transverse support inside a local
frame chart. The amplitudes and their corrections are smoothly extended
by zero in the transverse variable, so both transport equations hold
on the whole fixed set $U$. For a leading amplitude already defined on
all source cylinders, the finite-cover construction in that lemma
gives the same conclusion.
Proposition~\ref{pb27:remainder}, applied
to $h\chi_0M^WA_1$, gives a single remainder with $L^2=O(h^2)$ and
$H^1=O(h)$. Since $\chi_0=1$ near all of $\Omega_{\rm CGO}$ and
$\overline\Omega\Subset\Omega_{\rm CGO}$, elliptic regularity makes the
resulting exact field smooth near both sides of $\partial\Omega$. Its
restriction to $\Omega$ therefore satisfies the Betti identity. Use
$-\theta$ for the positive phase. All constants may depend on the fixed test.
\end{proof}

By Proposition~\ref{pb27:remainder}, an exact CGO solution
of the form~\eqref{pb27:two-sign-correction} is obtained
once the transport conditions~\eqref{pb27:two-sign-transport}
are satisfied. The next lemma uses the
leafwise Cauchy right inverse to construct the first correction $a_1$
for a prescribed leading amplitude $a_0$ and gives the $H^2(U)$ bound
required in Proposition~\ref{pb27:remainder}.

\begin{lemma}
\label{a31:first-correction}
Let $U$ be the fixed open set in
Proposition~\ref{pb27:remainder}, and let
$\pi_\theta(x)=z$ and $\pi_\perp(x)=x_\perp$ be the leaf and
transverse coordinates from Definition~\ref{i2:def-input}.
Choose a relatively compact transverse interval $I_\perp$ in a
local parameter chart and radii $L_{\rm src}<L_{\rm fr}$ such that
\[
\pi_\theta(\overline U)\Subset D_{L_{\rm src}}
\Subset D_{L_{\rm fr}}.
\]
Define the source cylinder
\begin{equation}\label{a31:source-cylinder}
\mathcal Q_{\rm src}
:=\{x\in\mathbb R^3:
 z(x)\in D_{L_{\rm src}},\ x_\perp\in I_\perp\}.
\end{equation}
Let $a_0$ be smooth on a neighborhood of
$\overline{\mathcal Q}_{\rm src}$ and satisfy
$\mathsf T_\theta a_0=0$ near
$\overline U\cap\overline{\mathcal Q}_{\rm src}$.
Choose the invertible four-component transport frame $C_\theta$
supplied by Theorem~\ref{f27:main} and
Proposition~\ref{tb27:F-supply}, smooth on a neighborhood of
$\overline D_{L_{\rm fr}}\times\overline I_\perp$, with
\[
\mathsf T_\theta C_\theta=0,\qquad \det C_\theta\ne0.
\]
Let $\mathcal C_\theta$ be the leafwise Cauchy right inverse from
Definition~\ref{i2:def-input}, and choose
\[
\chi_{\rm leaf}\in C_c^\infty(D_{L_{\rm src}}),\qquad
\chi_{\rm leaf}=1\quad\text{near }\pi_\theta(\overline U),
\qquad \chi_\perp\in C_c^\infty(I_\perp).
\]
Set
\[
\widetilde a_0(x)=\chi_\perp(x_\perp)a_0(x),
\]
with smooth transverse zero extension outside $I_\perp$. For
$\mathfrak s\in\{+1,-1\}$, with the phase convention of
Proposition~\ref{pb27:remainder}, define on the local cylinder
\begin{equation}\label{a31:first-correction-formula}
 a_1=\frac{\mathfrak s}{2}
 C_\theta\mathcal C_\theta
 \left(\chi_{\rm leaf}C_\theta^{-1}
 M^W\widetilde a_0\right),
\end{equation}
and extend $a_1$ by zero in the transverse variable. Then
\[
\mathsf T_\theta\widetilde a_0=0,\qquad
\mathsf T_\theta a_1=\mathfrak s M^W\widetilde a_0
\quad\text{on all of }U.
\]
Thus $(\widetilde a_0,a_1)$ satisfies
\eqref{pb27:two-sign-transport} on $U$.
If $a_0$ already has compact transverse support in
$J\Subset I_\perp$ and $\chi_\perp=1$ near $J$, then
$\widetilde a_0=a_0$, understood with its smooth zero extension,
and the same conclusion holds for $(a_0,a_1)$.

For each fixed smooth amplitude family, $a_1$ has an $H^2(U)$ bound
uniform on a compact parameter subchart with fixed positive support
margins, uniformly in $h$ for fixed coefficients, amplitudes,
frames, disks, and cutoffs.

If a leading amplitude $a_0$ is already smooth on neighborhoods of
all source cylinders in a finite transverse cover of
$\pi_\perp(\overline U)$ and satisfies $\mathsf T_\theta a_0=0$
near $\overline U$, the local construction can be summed over a
transverse partition of unity. It then gives a smooth first correction
for this uncut $a_0$ on all of $U$, with the same type of $H^2(U)$
bound on a sufficiently small compact direction subchart.
\end{lemma}

\begin{proof}
Since $D_\theta\chi_\perp=0$, multiplication by $\chi_\perp$
preserves the first transport equation wherever $a_0$ is defined.
The product rule and $\mathsf T_\theta C_\theta=0$ give
\[
\mathsf T_\theta(C_\theta v)=2C_\theta D_\theta v.
\]
Using $D_\theta\mathcal C_\theta H=H$ from
Definition~\ref{i2:def-input}, we obtain on
$U\cap\{x_\perp\in I_\perp\}$
\begin{align*}
\mathsf T_\theta a_1
&=\mathfrak s C_\theta D_\theta\mathcal C_\theta
 \left(\chi_{\rm leaf}C_\theta^{-1}M^W\widetilde a_0\right)\\
&=\mathfrak s\chi_{\rm leaf}M^W\widetilde a_0
 =\mathfrak s M^W\widetilde a_0.
\end{align*}
The last equality uses the leaf cutoff, which equals one on the
projection of $U$. The source $M^W(\chi_\perp a_0)$ retains the
terms in which derivatives fall on $\chi_\perp$.

The source $M^W\widetilde a_0$ is supported transversely in
$\operatorname{supp}\chi_\perp$. Since $\mathcal C_\theta$ acts
separately at each value of $x_\perp$, we have
\[
\operatorname{supp}_{x_\perp}(M^W\widetilde a_0),\quad
\operatorname{supp}_{x_\perp}a_1
\ \subseteq\ \operatorname{supp}\chi_\perp\Subset I_\perp.
\]
Both amplitudes and their derivatives vanish in a neighborhood of
the ends of $I_\perp$. Their zero extensions are smooth, and the
transport equations hold also on the rest of $U$, where both sides
vanish locally. The derivative-transfer estimate
\eqref{f27:uniform-cauchy-bound}, the fixed support margins, and the
smooth bounds for $C_\theta$ and $C_\theta^{-1}$ give the asserted
$H^2(U)$ bound. Parameter derivatives obey the same estimates.

For the final assertion, choose finitely many local transverse
intervals $I_{\perp,\alpha}$ and cutoffs
$\chi_{\perp,\alpha}\in C_c^\infty(I_{\perp,\alpha})$ such that
\[
\sum_\alpha\chi_{\perp,\alpha}=1
\quad\text{near }\pi_\perp(\overline U).
\]
Apply the local construction with the corresponding frame
$C_{\theta,\alpha}$ to each $\chi_{\perp,\alpha}a_0$, and denote
the zero-extended correction by $a_{1,\alpha}$. Their finite sum
$a_1=\sum_\alpha a_{1,\alpha}$ satisfies
\[
\mathsf T_\theta a_1
=\mathfrak s\sum_\alpha M^W(\chi_{\perp,\alpha}a_0)
=\mathfrak s M^W\!\left(\sum_\alpha\chi_{\perp,\alpha}a_0\right)
=\mathfrak s M^Wa_0
\quad\text{on }U.
\]
The derivatives of the transverse cutoffs cancel in this sum.
Compactness gives the finite cover. After shrinking the direction
neighborhood, the same cover, cutoffs, and leaf disks retain positive
margins around all relevant projections. Smooth local frame families
and a finite sum then give joint smoothness and the stated uniform
bound. The cover, cutoffs, disks, and frames are independent of
subsequent test functions.
\end{proof}

\subsection{Canonical lifts and high-frequency asymptotics}
\label{i2:canonical-lift}

The Green--Betti identity~\eqref{i2:eq-betti} involves the
dual testing field $y_2=S_2^Tu_2$, where $u_2$ is an elastic
solution. To construct this field from an exact CGO solution
of the reduced system, let $z_h^-$ be the solution constructed
in~\eqref{pb27:two-sign-correction}, and set
\[
u_h:=Tz_h^-,
\qquad
y_h:=S^Tu_h=S^TTz_h^-.
\]
Thus $u_h$ is the associated elastic displacement and $y_h$
is the corresponding dual field. This motivates the definition
\[
K_{\mathrm{lift}}:=S^TT,
\qquad y_h=K_{\mathrm{lift}}z_h^-.
\]
The operator $K_{\mathrm{lift}}$ contains second derivatives,
whereas the remainder estimates in~\eqref{pb27:two} are stated in
$L^2$ and $H^1$. The following lemma rewrites it as
$K_{\mathrm{lift}}=\widehat K_{\mathrm{lift}}^W
-E_{\mathrm{lift}}M^W$. On exact solutions of the reduced system,
the term involving $M^W$ vanishes, and the first-order operator
$\widehat K_{\mathrm{lift}}^W$ yields the required high-frequency
asymptotics of $y_h$.

\begin{lemma}
\label{i2:canonical-lift-identities}
With $K_{\mathrm{lift}}$ defined above, let
$e_4=(0,0,0,1)^T$. For $z=(v,f)^T$, define
\begin{align}
K_{\mathrm{lift}}(v,f)
&=
\begin{pmatrix}
v+c\nabla f+cf\mathbf m\\[1mm]
-\dfrac{
c\operatorname{div}v-\frac12c\mathbf m\cdot v
+a(\Delta f+q_\mu f)}
{1-\kappa}
\end{pmatrix},
\label{i2:eq-k}\\
\widehat K_{\mathrm{lift}}^W(v,f)
&=
\begin{pmatrix}
v+c\nabla f+cf\mathbf m\\[1mm]
-c\operatorname{div}v+\frac32c\mathbf m\cdot v
+W\mu^{-2}f
\end{pmatrix},
\qquad
E_{\mathrm{lift}}
=\frac{e_4e_4^T}{\mu(1-\kappa)},
\nonumber\\
K_{\mathrm{lift}}
&=\widehat K_{\mathrm{lift}}^W
-E_{\mathrm{lift}}M^W,
\label{i2:eq-khat}\\
N_\theta(v,f)
&=(-cf\theta,c\theta\cdot v)^T,
\qquad
J_\theta=-N_\theta-E_{\mathrm{lift}}\mathsf T_\theta,
\label{i2:eq-j}\\
J_\theta(v,f)
&=
\left(
cf\theta,
-\frac{c\theta\cdot v+2aD_\theta f}{1-\kappa}
\right)^T.
\nonumber
\end{align}
For every exact solution $z$ satisfying
$M^Wz=0$ in $\Omega_{\mathrm{CGO}}$,
\begin{equation}
K_{\mathrm{lift}}z
=\widehat K_{\mathrm{lift}}^Wz
\quad\text{in }\Omega_{\mathrm{CGO}}.
\label{i2:eq-k-exact}
\end{equation}
Moreover, the conjugated lift satisfies
\begin{equation}
e^{\theta\cdot x/h}K_{\mathrm{lift}}
e^{-\theta\cdot x/h}
=K_{\mathrm{lift}}-h^{-1}J_\theta.
\label{i2:eq-lift-conjugation}
\end{equation}
For amplitudes $a_0,a_1$ satisfying
\eqref{pb27:two-sign-transport} with $\mathfrak s=+1$,
\begin{equation}
-N_\theta a_0=J_\theta a_0,
\qquad
\widehat K_{\mathrm{lift}}^W a_0+N_\theta a_1
=K_{\mathrm{lift}}a_0-J_\theta a_1.
\label{i2:eq-amp-shell}
\end{equation}
\end{lemma}

\begin{proof}
The reduction maps $T,S$ and the formal bilinear transpose
$S^T$ are given in~\eqref{bt:TS} and
Lemma~\ref{bt:scalarlemma}. Differentiating $T(v,f)$ gives
\[
\operatorname{div}T(v,f)
=c\operatorname{div}v
-\frac12c\mathbf m\cdot v
+a(\Delta f+q_\mu f).
\]
Applying $S^T$ yields~\eqref{i2:eq-k}.

To prove~\eqref{i2:eq-khat}, compare the fourth components
of the operators. The coefficient of $\operatorname{div}v$
in $\widehat K_{\mathrm{lift}}^W-E_{\mathrm{lift}}M^W$ is
\[
-c-\frac{\kappa c}{1-\kappa}
=-\frac{c}{1-\kappa},
\]
and that of $\mathbf m\cdot v$ is
\[
\frac32c+\frac{(3\kappa-2)c}{2(1-\kappa)}
=\frac{c}{2(1-\kappa)}.
\]
The terms involving $W\mu^{-2}f$ cancel. The remaining
coefficients agree with~\eqref{i2:eq-k}, proving the operator
identity~\eqref{i2:eq-khat}. Equation~\eqref{i2:eq-k-exact}
then follows from $M^Wz=0$.

Direct differentiation of the first-order operator gives
\begin{equation}
e^{\theta\cdot x/h}\widehat K_{\mathrm{lift}}^W
e^{-\theta\cdot x/h}
=\widehat K_{\mathrm{lift}}^W+h^{-1}N_\theta.
\label{i2:eq-khat-conjugation}
\end{equation}
The conjugation identity for $M^W$ is
\eqref{pb27:two-sign-conjugation} with $\mathfrak s=+1$.
Consequently, using~\eqref{i2:eq-khat},
\begin{align*}
e^{\theta\cdot x/h}K_{\mathrm{lift}}
e^{-\theta\cdot x/h}
&=\widehat K_{\mathrm{lift}}^W+h^{-1}N_\theta
-E_{\mathrm{lift}}
\bigl(M^W-h^{-1}\mathsf T_\theta\bigr)\\
&=K_{\mathrm{lift}}
+h^{-1}\bigl(N_\theta+E_{\mathrm{lift}}\mathsf T_\theta\bigr)\\
&=K_{\mathrm{lift}}-h^{-1}J_\theta,
\end{align*}
which proves~\eqref{i2:eq-lift-conjugation}.

Finally, the transport conditions
\eqref{pb27:two-sign-transport} with $\mathfrak s=+1$ give
\[
\mathsf T_\theta a_0=0,
\qquad
\mathsf T_\theta a_1=M^Wa_0.
\]
By the definition of $J_\theta$,
\[
J_\theta a_0
=-N_\theta a_0-E_{\mathrm{lift}}\mathsf T_\theta a_0
=-N_\theta a_0.
\]
Similarly,~\eqref{i2:eq-khat} gives
\begin{align*}
K_{\mathrm{lift}}a_0-J_\theta a_1
&=\widehat K_{\mathrm{lift}}^W a_0
-E_{\mathrm{lift}}M^Wa_0
+N_\theta a_1
+E_{\mathrm{lift}}\mathsf T_\theta a_1\\
&=\widehat K_{\mathrm{lift}}^W a_0+N_\theta a_1.
\end{align*}
This proves~\eqref{i2:eq-amp-shell}.
\end{proof}

The exact-solution identity~\eqref{i2:eq-k-exact}, the
conjugation formula~\eqref{i2:eq-khat-conjugation}, and
the amplitude identities~\eqref{i2:eq-amp-shell} from
Lemma~\ref{i2:canonical-lift-identities}, together with
the remainder estimates~\eqref{pb27:two}, yield the
high-frequency asymptotics of
$y_h=K_{\mathrm{lift}}z_h^-$.

\begin{lemma}
\label{i2:lem-shell}
Let $a_0,a_1$ satisfy
\eqref{pb27:two-sign-transport} with $\mathfrak s=+1$,
and let $z_h^-$ be the exact CGO solution furnished by
Proposition~\ref{pb27:remainder}, with remainder $r_h$.
Then the scaled lift satisfies
\begin{equation}
e^{\theta\cdot x/h}
\bigl(-hK_{\mathrm{lift}}z_h^-\bigr)
=J_\theta a_0
+O_{L^2(\Omega_{\mathrm{CGO}})}(h).
\label{i2:eq-scaled}
\end{equation}
If $J_\theta a_0=0$, then the unscaled lift satisfies
\begin{equation}
e^{\theta\cdot x/h}K_{\mathrm{lift}}z_h^-
=K_{\mathrm{lift}}a_0-J_\theta a_1
+O_{L^2(\Omega_{\mathrm{CGO}})}(h).
\label{i2:eq-unscaled}
\end{equation}
Both estimates hold for the same remainder $r_h$.
\end{lemma}

\begin{proof}
By~\eqref{i2:eq-k-exact} and the first-order conjugation
identity~\eqref{i2:eq-khat-conjugation},
\[
\begin{aligned}
e^{\theta\cdot x/h}K_{\mathrm{lift}}z_h^-
={}&h^{-1}N_\theta a_0
+\widehat K_{\mathrm{lift}}^W a_0
+N_\theta a_1\\
&+h\widehat K_{\mathrm{lift}}^W a_1
+\widehat K_{\mathrm{lift}}^W r_h
+h^{-1}N_\theta r_h.
\end{aligned}
\]
The amplitude identities~\eqref{i2:eq-amp-shell} give
\[
-N_\theta a_0=J_\theta a_0,
\qquad
\widehat K_{\mathrm{lift}}^W a_0+N_\theta a_1
=K_{\mathrm{lift}}a_0-J_\theta a_1.
\]
Consequently,
\[
\begin{aligned}
e^{\theta\cdot x/h}K_{\mathrm{lift}}z_h^-
={}&-h^{-1}J_\theta a_0
+K_{\mathrm{lift}}a_0-J_\theta a_1\\
&+h\widehat K_{\mathrm{lift}}^W a_1
+\widehat K_{\mathrm{lift}}^W r_h
+h^{-1}N_\theta r_h.
\end{aligned}
\]

The remainder estimates~\eqref{pb27:two} and the fact that
$\widehat K_{\mathrm{lift}}^W$ is first order imply
\[
\begin{aligned}
\|h\widehat K_{\mathrm{lift}}^W a_1\|_{L^2}
&=O(h),\\
\|\widehat K_{\mathrm{lift}}^W r_h\|_{L^2}
&\le C\|r_h\|_{H^1}=O(h),\\
\|h^{-1}N_\theta r_h\|_{L^2}
&\le Ch^{-1}\|r_h\|_{L^2}=O(h).
\end{aligned}
\]
Thus
\[
e^{\theta\cdot x/h}K_{\mathrm{lift}}z_h^-
=-h^{-1}J_\theta a_0
+K_{\mathrm{lift}}a_0-J_\theta a_1
+O_{L^2(\Omega_{\mathrm{CGO}})}(h).
\]
Multiplication by $-h$ gives~\eqref{i2:eq-scaled}. If
$J_\theta a_0=0$, the leading $h^{-1}$ term vanishes and the
expansion gives~\eqref{i2:eq-unscaled}. Both estimates follow from
this expansion of the lift of $z_h^-$ with remainder $r_h$.
\end{proof}

\subsection{Construction of the three polarization columns}
\label{i2:polarizations}

We construct the first two leading polarization columns from the
scaled physical lifts and then construct a third column from an
unscaled lift. The latter is chosen so that its leading amplitude
lies in the kernel of $J_\theta$, allowing the second asymptotic
formula in~\eqref{i2:eq-unscaled} to be used.

\begin{lemma}\label{i2:lem-f}
Fix the leaf disks as in Lemma~\ref{a31:first-correction}.
At each direction and transverse-parameter centre, choose a common
local real parameter chart on which the two- and four-component
frames of part~2 of Theorem~\ref{f27:main} and
Proposition~\ref{tb27:F-supply} are available on a neighborhood
of $\overline D_{L_{\rm fr}}$. Fix $I_\perp$ relatively compact
in its transverse-parameter interval. These choices precede all
tests. The uncut frames and amplitudes below are defined on this
chart; the four-component transport and lift identities are
asserted on $U\cap\{x_\perp\in I_\perp\}$.
Consider the two-component system
\begin{equation}
D_\theta f=-\frac12\mu\kappa R,
\qquad
D_\theta R=fD_\theta^2a,
\qquad a=\mu^{-1}.
\label{i2:eq-fr-system}
\end{equation}

On this fixed common local parameter chart, apply part~2 of
Theorem~\ref{f27:main} and Proposition~\ref{tb27:F-supply}
on the prescribed larger leaf disk. They supply a smooth
invertible matrix on a neighborhood of
$\overline D_{L_{\rm fr}}\times\overline I_\perp$
\[
\Phi=
\begin{pmatrix}
f_1&f_2\\
R_1&R_2
\end{pmatrix}
\]
satisfying
\[
D_\theta\Phi=
\begin{pmatrix}
0&-\frac12\mu\kappa\\
D_\theta^2a&0
\end{pmatrix}\Phi,
\qquad
\det\Phi=1.
\]
The normalization $\det\Phi=1$ is obtained in the proof.

For $j=1,2$, define
\begin{equation}
Y_j^F=
\begin{pmatrix}
cf_j\theta\\
-R_j
\end{pmatrix}.
\label{i2:eq-yf}
\end{equation}
Then
\[
(2D_\theta-(A_\theta^{\rm phys})^T)Y_j^F=0.
\]

Moreover, there exist smooth local amplitudes $A_{0,j}^F$ satisfying
\[
\mathsf T_\theta A_{0,j}^F=0,\qquad
J_\theta A_{0,j}^F=Y_j^F
\quad\text{on }U\cap\{x_\perp\in I_\perp\}.
\]
For every $\chi_\perp\in C_c^\infty(I_\perp)$, there are smooth
first corrections $A_{1,j}^F$, depending on this cutoff, such that
\[
\mathsf T_\theta A_{1,j}^F=M^W(\chi_\perp A_{0,j}^F)
\quad\text{on all of }U,
\]
where the localized amplitudes and corrections are extended by zero
in the transverse variable. Let $z_{h,j}^-$ be the exact CGO solution
furnished by Proposition~\ref{pb27:remainder}, with
$\mathfrak s=+1$, $a_0=\chi_\perp A_{0,j}^F$, and
$a_1=A_{1,j}^F$. Then
\begin{equation}
e^{\theta\cdot x/h}
\bigl(-hK_{\mathrm{lift}}z_{h,j}^-\bigr)
=\chi_\perp Y_j^F+O_{L^2(\Omega_{\mathrm{CGO}})}(h),
\qquad j=1,2.
\label{i2:eq-yf-realization}
\end{equation}
In particular, the leading column is $Y_j^F$ wherever
$\chi_\perp=1$.

\end{lemma}

\begin{proof}
The coefficient matrix in~\eqref{i2:eq-fr-system} has trace
zero. Hence the frame equation gives
\[
D_\theta\det\Phi
=
\operatorname{tr}
\begin{pmatrix}
0&-\frac12\mu\kappa\\
D_\theta^2a&0
\end{pmatrix}
\det\Phi=0.
\]
The chosen frame is invertible on the same neighborhood of
$\overline D_{L_{\rm fr}}\times\overline I_\perp$.
Its determinant is nonzero and holomorphic in the leaf variable there.
Multiplying its first column by $(\det\Phi)^{-1}$ on this whole
neighborhood preserves the frame equation, remains jointly smooth
in the local real parameters, and normalizes $\det\Phi=1$ there.

For the normalized null direction, put
\[
a_\theta=\frac{\bar\theta}{2},
\qquad \theta\cdot a_\theta=1.
\]
For $j=1,2$, define
\[
H_j
=
2D_\theta(c^{-1}R_j a_\theta)
-\mathbf m c^{-1}R_j
+F_\mu\theta f_j.
\]

Using the coefficient identity
\begin{equation}
D_\theta^2a=-\frac c2\theta^TF_\mu\theta,
\label{i2:eq-d2a}
\end{equation}
together with $D_\theta c=-\frac c2\theta\cdot\mathbf m$
and~\eqref{i2:eq-fr-system}, we obtain
\begin{align*}
\theta\cdot H_j
&=2D_\theta(c^{-1}R_j)
-c^{-1}R_j\,\theta\cdot\mathbf m
+f_j\theta^TF_\mu\theta\\
&=2c^{-1}D_\theta R_j
+f_j\theta^TF_\mu\theta\\
&=2c^{-1}f_jD_\theta^2a
+f_j\theta^TF_\mu\theta\\
&=-f_j\theta^TF_\mu\theta
+f_j\theta^TF_\mu\theta\\
&=0.
\end{align*}

Using the leafwise Cauchy right inverse and cutoff fixed in
Lemma~\ref{a31:first-correction}, set
\[
w_j=-\frac12\mathcal C_\theta
\bigl(\chi_{\rm leaf}H_j\bigr),
\qquad
A_{0,j}^F=
\begin{pmatrix}
c^{-1}R_j a_\theta+w_j\\
f_j
\end{pmatrix}.
\]

The preceding calculation of $\theta\cdot H_j$ holds throughout
$D_{L_{\rm fr}}\times I_\perp$, where the chosen columns of $\Phi$
satisfy~\eqref{i2:eq-fr-system}.
Since $\operatorname{supp}\chi_{\rm leaf}\Subset D_{L_{\rm src}}
\Subset D_{L_{\rm fr}}$ and $\theta$ is independent of \(x\),
\[
\theta\cdot w_j
=-\frac12\mathcal C_\theta
 \bigl(\chi_{\rm leaf}\,\theta\cdot H_j\bigr)=0.
\]

Moreover,
\[
2D_\theta w_j=-H_j
\]
on the working disk, where $\chi_{\rm leaf}=1$.
Consequently,
\[
\theta\cdot(c^{-1}R_j a_\theta+w_j)=c^{-1}R_j.
\]
The vector and scalar components of
$\mathsf T_\theta A_{0,j}^F$ therefore reduce to
\[
\begin{aligned}
&H_j+2D_\theta w_j=0,\\
&2D_\theta f_j+d c^{-1}R_j
=-\mu\kappa R_j+\kappa\mu R_j=0.
\end{aligned}
\]
Thus $\mathsf T_\theta A_{0,j}^F=0$.

The explicit formula for $J_\theta$ in~\eqref{i2:eq-j}
and the first equation of~\eqref{i2:eq-fr-system} give
\[
\begin{aligned}
J_\theta A_{0,j}^F
&=
\begin{pmatrix}
cf_j\theta\\[1mm]
-\dfrac{c\theta\cdot(c^{-1}R_j a_\theta+w_j)
+2aD_\theta f_j}{1-\kappa}
\end{pmatrix}\\
&=
\begin{pmatrix}
cf_j\theta\\[1mm]
-\dfrac{R_j-\kappa R_j}{1-\kappa}
\end{pmatrix}
=Y_j^F.
\end{aligned}
\]

It remains to verify the dual transport equation. For
$Y=(p,s)^T$, the equation
$(2D_\theta-(A_\theta^{\rm phys})^T)Y=0$ is equivalent to
\begin{equation}
\begin{aligned}
2D_\theta p+\theta(\mathbf m\cdot p)-d\theta s&=0,\\
2D_\theta s-\theta^TF_\mu p&=0.
\end{aligned}
\label{i2:eq-dualrows}
\end{equation}
Substituting $p=cf_j\theta$ and $s=-R_j$ reduces the two
rows to
\[
\theta(2cD_\theta f_j+dR_j)=0,
\qquad
-2(D_\theta R_j-f_jD_\theta^2a)=0,
\]
respectively. Both vanish by~\eqref{i2:eq-fr-system}.
Hence $(2D_\theta-(A_\theta^{\rm phys})^T)Y_j^F=0$.

Finally, apply the local construction of
Lemma~\ref{a31:first-correction} to $a_0=A_{0,j}^F$, the
chosen $\chi_\perp$, and $\mathfrak s=+1$. Denote the resulting
zero-extended first correction by $A_{1,j}^F$. Then
\[
\mathsf T_\theta A_{1,j}^F=M^W(\chi_\perp A_{0,j}^F)
\quad\text{on all of }U,
\qquad A_{1,j}^F\in H^2(U).
\]
Proposition~\ref{pb27:remainder}, with
$a_0=\chi_\perp A_{0,j}^F$ and $a_1=A_{1,j}^F$, gives the
exact CGO solution $z_{h,j}^-$. Since $D_\theta\chi_\perp=0$
and the differential part of $J_\theta$ contains only $D_\theta$,
\[
J_\theta(\chi_\perp A_{0,j}^F)=\chi_\perp Y_j^F.
\]
The scaled lift formula~\eqref{i2:eq-scaled} now proves
\eqref{i2:eq-yf-realization} throughout the target domain.

\end{proof}

The first two columns $Y_1^F,Y_2^F$ are constructed in
Lemma~\ref{i2:lem-f}. We now construct a third column
$Y^S$ by choosing a leading amplitude $A_0^S$ independent of $x$
and satisfying
\[
\mathsf T_\theta A_0^S=0,
\qquad J_\theta A_0^S=0.
\]
The second identity allows us to use the unscaled lift
formula~\eqref{i2:eq-unscaled}. We will also show that the
resulting leading column $Y^S$ is independent of the
potential $W$.

\begin{proposition}
\label{i2:prop-s}
Represent the fixed null direction $\theta$ by a nonzero
spinor $(s_0,s_1)\in\mathbb C^2$:
\[
\theta=
\begin{pmatrix}
s_0^2-s_1^2\\
i(s_0^2+s_1^2)\\
2s_0s_1
\end{pmatrix},
\qquad
Q=\frac12
\begin{pmatrix}
-s_1&-s_0\\
is_1&-is_0\\
s_0&-s_1
\end{pmatrix},
\qquad
n_{\mathrm{sp}}=
\begin{pmatrix}-s_1\\s_0\end{pmatrix}.
\]
Then
\[
Qn_{\mathrm{sp}}=-\frac{\theta}{2},
\qquad
\operatorname{ran}Q=\theta^\perp.
\]
On the two spinor charts choose
\[
 \xi_q^{(0)}=(s_0^{-1},0)^T\quad(s_0\ne0),\qquad
 \xi_q^{(1)}=(0,s_1^{-1})^T\quad(s_1\ne0).
\]
On their overlap,
\[
 \xi_q^{(0)}-\xi_q^{(1)}
   =-\frac{n_{\mathrm{sp}}}{s_0s_1},
 \qquad
 Q\xi_q^{(0)}-Q\xi_q^{(1)}=\frac{\theta}{2s_0s_1}.
\]
Fix either chart and write $\xi_q$ for the corresponding vector
and $q_\theta:=Q\xi_q$. Thus
\[
\theta\cdot q_\theta=0,
\qquad
\det(n_{\mathrm{sp}},\xi_q)=-1.
\]
The vector $q_\theta$ is independent of $x$. Define
\[
A_0^S:=\begin{pmatrix}q_\theta\\0\end{pmatrix}.
\]
Then
\[
\mathsf T_\theta A_0^S=0,
\qquad J_\theta A_0^S=0.
\]

Choose a smooth first correction
\[
A_{1,0}^S=\begin{pmatrix}v_S\\f_S\end{pmatrix}
\]
satisfying
\[
\mathsf T_\theta A_{1,0}^S=MA_0^S,
\]
where $M$ is the zero-potential operator defined by
$M=M^0$ in~\eqref{bt:MW}. Set
\[
R_S:=c\theta\cdot v_S.
\]
For the full operator $M^W$, there exists a first correction
$A_{1,W}^S$ satisfying
\[
\mathsf T_\theta A_{1,W}^S=M^WA_0^S.
\]
The third column is
\begin{equation}
Y^S=
\begin{pmatrix}
q_\theta-cf_S\theta\\[1mm]
R_S+\dfrac32c\,\mathbf m\cdot q_\theta
\end{pmatrix},
\qquad
(2D_\theta-(A_\theta^{\rm phys})^T)Y^S=0.
\label{i2:eq-ys}
\end{equation}
Moreover, for the exact CGO solution $z_h^-$ constructed
from $A_0^S,A_{1,W}^S$ by
Proposition~\ref{pb27:remainder},
\[
e^{\theta\cdot x/h}K_{\mathrm{lift}}z_h^-
=Y^S+O_{L^2(\Omega_{\mathrm{CGO}})}(h).
\]
The leading column $Y^S$ is independent of $W$.
\end{proposition}

\begin{proof}

The amplitude $A_0^S$ is independent of \(x\) and is defined on every
source cylinder. Apply the finite-cover construction of
Lemma~\ref{a31:first-correction} with $a_0=A_0^S$, $W=0$, and
$\mathfrak s=+1$. It supplies a smooth correction $A_{1,0}^S$
on the whole working set, with
\[
\mathsf T_\theta A_{1,0}^S
=\sum_\alpha M(\chi_{\perp,\alpha}A_0^S)
=M A_0^S
\quad\text{on }U.
\]
The local four-component frames, transverse partition, and resulting
static correction are fixed before any later transverse test is
chosen. The construction gives the $H^2(U)$ bound and smooth local
parameter dependence required by Proposition~\ref{pb27:remainder}.

Writing $A_{1,0}^S=(v_S,f_S)^T$ and
$R_S=c\theta\cdot v_S$, the scalar row and the contraction
of the vector row with $\theta$ give
\begin{align}
D_\theta f_S
&=-\frac12\mu\kappa R_S
-\frac12e^{\mathrm{red}}\,\mathbf m\cdot q_\theta,
\label{i2:eq-sfirst}\\
D_\theta R_S
&=f_SD_\theta^2a
+\frac c2\theta^TV_\mu q_\theta.
\label{i2:eq-ssecond}
\end{align}
In the second identity, differentiating $c$ in
$D_\theta(c\theta\cdot v_S)$ cancels the term involving
$\mathbf m$ from the contracted vector equation.

We next account for the potential. By~\eqref{bt:MW},
\[
M^W=M+\mathcal W,
\qquad
\mathcal W=
\begin{pmatrix}
(W/\mu)I_3&
\mu^{-3/2}(2W\mathbf m-\nabla W)\\
0&(1-\kappa)W/\mu
\end{pmatrix}.
\]
Since $A_0^S=(q_\theta,0)^T$, it follows that
\[
\mathcal WA_0^S
=\begin{pmatrix}(W/\mu)q_\theta\\0\end{pmatrix}.
\]
On each parallel leaf meeting $U$, take $\chi_{\rm leaf}=1$ on a
disk containing $\pi_\theta(\overline U)$ as in
Proposition~\ref{a28:moment-inputs}. Use the same cutoff on those
leaves and define
\[
\beta=\frac12\mathcal C_\theta
\left(\chi_{\rm leaf}\frac W\mu\right),
\qquad
A_{1,W}^S
=A_{1,0}^S+\begin{pmatrix}\beta q_\theta\\0\end{pmatrix}.
\]
On the working disk,
\[
D_\theta\beta=\frac{W}{2\mu}.
\]
Since $D_\theta q_\theta=0$ and
$\theta\cdot q_\theta=0$, the definition
$\mathsf T_\theta=2D_\theta+A_\theta^{\rm phys}$ gives
\[
\begin{aligned}
\mathsf T_\theta
\begin{pmatrix}\beta q_\theta\\0\end{pmatrix}
&=\begin{pmatrix}2D_\theta\beta\,q_\theta\\0\end{pmatrix}\\
&=\begin{pmatrix}(W/\mu)q_\theta\\0\end{pmatrix}
=\mathcal WA_0^S.
\end{aligned}
\]
Consequently,
\[
\mathsf T_\theta A_{1,W}^S
=MA_0^S+\mathcal WA_0^S
=M^WA_0^S.
\]
Furthermore, the explicit formula for $J_\theta$
in~\eqref{i2:eq-j} and $\theta\cdot q_\theta=0$ imply
\[
J_\theta
\begin{pmatrix}\beta q_\theta\\0\end{pmatrix}=0,
\qquad
J_\theta A_{1,W}^S=J_\theta A_{1,0}^S.
\]

We now compute the unscaled lift. Formula~\eqref{i2:eq-k}
gives
\[
K_{\mathrm{lift}}A_0^S
=\begin{pmatrix}
q_\theta\\[1mm]
\dfrac{c\,\mathbf m\cdot q_\theta}{2(1-\kappa)}
\end{pmatrix}.
\]
By~\eqref{i2:eq-j},
\[
J_\theta A_{1,0}^S
=\begin{pmatrix}
cf_S\theta\\[1mm]
-\dfrac{R_S+2aD_\theta f_S}{1-\kappa}
\end{pmatrix}.
\]
Thus
\[
\begin{aligned}
K_{\mathrm{lift}}A_0^S-J_\theta A_{1,W}^S
&=K_{\mathrm{lift}}A_0^S-J_\theta A_{1,0}^S\\
&=\begin{pmatrix}
q_\theta-cf_S\theta\\[1mm]
\dfrac{\frac12c\,\mathbf m\cdot q_\theta
+R_S+2aD_\theta f_S}{1-\kappa}
\end{pmatrix}.
\end{aligned}
\]
Substituting~\eqref{i2:eq-sfirst} and the coefficient
identity $ae^{\mathrm{red}}=(3\kappa-2)c/2$ from~\eqref{bt:notation},
the last component simplifies to
\[
\frac{\frac12c\,\mathbf m\cdot q_\theta
+R_S+2aD_\theta f_S}{1-\kappa}
=R_S+\frac32c\,\mathbf m\cdot q_\theta.
\]
Hence
\[
K_{\mathrm{lift}}A_0^S-J_\theta A_{1,W}^S=Y^S.
\]
Since $J_\theta A_0^S=0$, the unscaled lift formula
\eqref{i2:eq-unscaled} yields
\[
e^{\theta\cdot x/h}K_{\mathrm{lift}}z_h^-
=Y^S+O_{L^2(\Omega_{\mathrm{CGO}})}(h).
\]
The expression for $Y^S$ depends on $q_\theta,f_S,R_S$
and the Lam\'e coefficients, so its leading term is independent
of $W$.

It remains to verify the equation in~\eqref{i2:eq-ys}.
Using the two rows of the dual equation
\eqref{i2:eq-dualrows}, substitute
\[
p=q_\theta-cf_S\theta,
\qquad
s=R_S+\frac32c\,\mathbf m\cdot q_\theta.
\]
Because $D_\theta q_\theta=0$ and
$\theta\cdot q_\theta=0$, the vector row becomes
\[
\theta\left[
-2cD_\theta f_S
+\left(1-\frac32\kappa\right)\mathbf m\cdot q_\theta
-dR_S
\right].
\]
This vanishes by~\eqref{i2:eq-sfirst} and the coefficient
identities $ce^{\mathrm{red}}=(3\kappa-2)/2$ and
$d=\kappa\sqrt\mu$ with $c=\mu^{-1/2}$.

The scalar row becomes
\[
2D_\theta R_S
+3D_\theta(c\,\mathbf m\cdot q_\theta)
-\theta^TF_\mu q_\theta
+cf_S\theta^TF_\mu\theta.
\]
Using~\eqref{i2:eq-ssecond} and~\eqref{i2:eq-d2a},
the terms containing $f_S$ cancel, leaving
\[
\theta^T
\left[cV_\mu+3(\nabla(c\mathbf m))^T-F_\mu\right]
q_\theta.
\]
With $\mathbf m=-2\nabla c/c$, direct expansion gives
\begin{align*}
F_\mu
&=-4\nabla^2c
-4c^{-1}\nabla c\otimes\nabla c
+4(\Delta c+c^{-1}|\nabla c|^2)I_3,\\
cV_\mu
&=2\nabla^2c
-4c^{-1}\nabla c\otimes\nabla c
+(\Delta c-2c^{-1}|\nabla c|^2)I_3,\\
3(\nabla(c\mathbf m))^T
&=-6\nabla^2c.
\end{align*}
Consequently,
\begin{equation}
F_\mu-cV_\mu-3(\nabla(c\mathbf m))^T
=(3\Delta c+6c^{-1}|\nabla c|^2)I_3.
\label{i2:eq-tensor-id}
\end{equation}
The scalar row is therefore a scalar multiple of
$\theta\cdot q_\theta$, which is zero. Both rows of the
dual equation hold, completing the proof.
\end{proof}

\subsection{\texorpdfstring{Extraction of the matrix moments}{Extraction of the matrix moments}}
\label{i2:moment-extraction}

The three leading columns are now assembled into fixed primal and
dual quotient frames. We next show that these frames can be chosen
before the polynomial and transverse tests, so that the same matrix
field occurs in every moment identity.

\begin{lemma}\label{i2:lem-fixed}
At each direction and transverse-parameter centre, choose
the common local chart of Lemma~\ref{i2:lem-f}, intersected with
local charts for the primal quotient and four-component frames
on the prescribed larger leaf disk $D_{L_{\rm fr}}$.
Fix an open interval $I$ whose closure is contained in $I_\perp$
and in the transverse-parameter range of this common chart.
The common chart and its interval $I$ are independent of the tests.
For material $j$, let $A_{j,\theta}^{\rm phys}$ denote the matrix $A_\theta^{\rm phys}$
evaluated at its coefficients. On this fixed chart, choose invertible $3\times3$ quotient frames $X_1,Y_2$ satisfying
\begin{equation}\label{i2:eq-frames}
 D_{\theta}X_1+\mathsf{P}_1X_1=0,\qquad D_{\theta}Y_2-\mathsf{P}_2^TY_2=0,\qquad\det Y_2=2c_2\ne0.
\end{equation}
Here $U_\theta=\operatorname{diag}(Q,1)$ is a $4\times3$ matrix. The map $U_\theta^T:\mathbb C^4\to\mathbb C^3$ has kernel $\mathbb C(\theta,0)^T$ and gives coordinates on the characteristic quotient. Choose the spatially constant right inverse
$V_\theta=\overline U_\theta(U_\theta^T\overline U_\theta)^{-1}$,
which is smooth in the local direction parameters because $U_\theta$
has full column rank, and set
$\mathsf{P}_j=U_\theta^TA_{j,\theta}^{\rm phys}V_\theta/2$.
Since $A_{j,\theta}^{\rm phys}(\theta,0)^T=0$, this $3\times3$ quotient coefficient is independent of the choice of $V_\theta$. The frame $X_1$ is primal and $Y_2$ is dual; the transpose in the second equation of~\eqref{i2:eq-frames} expresses their bilinear duality.
For every subsequently chosen $b\in C_c^\infty(I)$, all three columns can be realized by elastic solutions with
$O_{L^2}(h)$ error throughout the target, and with the same leading matrix $U_\theta Y_2$ near $\operatorname{supp}b$.
\end{lemma}
\begin{proof}

On the fixed common chart, apply part~2 of
Theorem~\ref{f27:main} and Proposition~\ref{tb27:F-supply}
on the prescribed larger leaf disk. Choose $X_1$ on a neighborhood of
$\overline D_{L_{\rm fr}}\times\overline I$, satisfying its quotient
transport equation there. Put $Z_1=V_\theta X_1+e_\theta\alpha$,
where $e_\theta=(\theta,0)^T$. Throughout this larger source cylinder,
\[
U_\theta^T\mathsf T_{1,\theta}(V_\theta X_1)
=2(D_\theta X_1+\mathsf P_1X_1)=0.
\]
Thus $\mathsf T_{1,\theta}(V_\theta X_1)=e_\theta f$ there,
with $f$ a smooth row of length three on the whole Cauchy source.
Choose $\alpha=-\frac12\mathcal C_\theta(\chi_{\rm leaf}f)$
componentwise. Since $A_{1,\theta}^{\rm phys}e_\theta=0$,
this gives $\mathsf T_{1,\theta}Z_1=0$ and $U_\theta^TZ_1=X_1$
on $U\cap\{x_\perp\in I\}$, where $\chi_{\rm leaf}=1$.

Lemma~\ref{i2:lem-f} and Proposition~\ref{i2:prop-s} give three dual columns in the four-component space.
Since $Qn_{\mathrm{sp}}=-\theta/2$, their quotient coordinates are
\[
 Y_2=\begin{pmatrix}
 -2cf_1n_{\mathrm{sp}}&-2cf_2n_{\mathrm{sp}}&\xi_q+2cf_Sn_{\mathrm{sp}}\\
 -R_1&-R_2&R_S+\frac32c\mathbf{m}\cdot q_\theta
 \end{pmatrix}.
\]
Put $B_\theta=(n_{\mathrm{sp}},\xi_q)\in GL(2,\mathbb C)$
and $s_S=R_S+3c\mathbf m\cdot q_\theta/2$. The displayed matrix factors as
\[
 Y_2=\operatorname{diag}(B_\theta,1)
 \begin{pmatrix}
 -2cf_1&-2cf_2&2cf_S\\
 0&0&1\\
 -R_1&-R_2&s_S
 \end{pmatrix}.
\]
Expanding the second determinant in its second row gives
$\det Y_2=-2c\det\Phi\det(n_{\mathrm{sp}},\xi_q)=2c$,
independently of $f_S$ and $R_S$. The resulting frames are independent
of $b$ and $g$.

For each $b$, choose $\chi\in C_c^\infty(I)$ equal to one near $\operatorname{supp}b$, and use $\chi(x_\perp)$ as a cutoff in the transverse variable $x_\perp$. Primes on $\chi$ below denote differentiation with respect to the real transverse variable $x_\perp$. The commutator $[M,\chi]=M\circ\chi-\chi M$ records all source terms created by this cutoff. For the third column, set
\[
 A_{0,\chi}^S=\chi A_0^S,\qquad
 A_{1,\chi}^S=\chi A_{1,W}^S+w_\chi,
 \quad\mathsf T_\theta w_\chi=[M,\chi]A_0^S.
\]
Since $D_{\theta}\chi=0$,
\begin{equation}\label{i2:eq-commutator}
 [M,\chi]A_0^S=
 \begin{pmatrix}\chi''q_\theta-\chi'\mathbf{m}(\nu_\theta\cdot q_\theta)\\
 d\chi'(\nu_\theta\cdot q_\theta)\end{pmatrix}.
\end{equation}

With the fixed local four-component frame from
Lemma~\ref{a31:first-correction}, set
\[
w_\chi=\frac12 C_\theta\mathcal C_\theta
\left(\chi_{\rm leaf}C_\theta^{-1}[M,\chi]A_0^S\right)
\]
on the chosen transverse chart and extend it smoothly by zero.
The source has compact transverse support inside $I$. Since the
primitive preserves $x_\perp$, $w_\chi$ and all its derivatives with respect to $x$
vanish on the open interval where $\chi=1$. The same support argument
as in Lemma~\ref{a31:first-correction} gives
$\mathsf T_\theta w_\chi=[M,\chi]A_0^S$ on all of $U$.
Moreover, $[\mathcal W,\chi]=0$, so the full transport relations
for $(A_{0,\chi}^S,A_{1,\chi}^S)$ hold on all of $U$.

For the potential correction, $\theta\cdot q_\theta=0$ gives
\[
 J_\theta(A_{1,W}^S-A_{1,0}^S)
 =J_\theta(\beta q_\theta,0)^T=0.
\]
The unscaled leading physical column therefore agrees with that for $W=0$. The transition correction $w_\chi$
vanishes with its derivatives on a neighborhood of $\supp b$:
for such transverse parameters its entire Cauchy source is zero,
and the primitive acts separately on each leaf. Thus the leading
physical column there is exactly the previously fixed $Y^S$.
For the first two columns, take $\chi A_{0,j}^F$,
$j=1,2$, and use Lemma~\ref{a31:first-correction} to solve their
first-correction equations on all of $U$. Their scaled leading
columns are $\chi Y_j^F$, equal to the fixed $Y_j^F$ near
$\operatorname{supp}b$.
The amplitudes and sources have compact transverse support and extend
smoothly by zero. The chosen leaf disks contain all leaf sections of
the working domains.
The zero-extended corrections have the fixed $H^2(U)$ bounds
needed in Proposition~\ref{pb27:remainder}.
Lemma~\ref{i2:lem-shell} gives the error estimate throughout the target, including the transition region of $\chi$.
\end{proof}

\begin{lemma}\label{i2:leaf-extraction}
Let $I\subset\mathbb R$ be an interval, $D\subset\mathbb C$ a bounded
disk, and let $\mathsf F\in C(I;L^1(D;\mathbb C^{3\times3}))$.
Suppose that for every fixed $g\in\mathbb C[z]$ and every
$b\in C_c^\infty(I)$,
\begin{equation}\label{i2:leaf-extraction-hyp}
 \int_I b(t)\int_D g(z)\mathsf F(t,z)\,dA(z)\,dt=0.
\end{equation}
Then $\int_Dg(z)\mathsf F(t,z)\,dA(z)=0$ for every $t\in I$ and every
$g\in\mathbb C[z]$.
\end{lemma}
\begin{proof}
Fix $g$ and put $J_g(t)=\int_Dg(z)\mathsf F(t,z)\,dA(z)$.
The boundedness of $g$ on $D$ gives
$\|J_g(t)-J_g(t')\|\le\|g\|_{L^\infty(D)}
\|\mathsf F(t,\cdot)-\mathsf F(t',\cdot)\|_{L^1(D)}$,
so $J_g$ is continuous. Applied entrywise, the fundamental lemma for
distributions and~\eqref{i2:leaf-extraction-hyp} imply that $J_g=0$
almost everywhere in $I$. Continuity makes this identity pointwise.
The argument applies to every polynomial using the same $\mathsf F$.
\end{proof}

\begin{theorem}\label{i2:thm-nine}
Under the notation of Definition~\ref{i2:def-input} and the constructions in Proposition~\ref{a28:moment-inputs}, choose
$\overline\Omega\subset B_{R_0}$, $R_0<R_\perp<R_1$, and
an interval $I$ selected in Lemma~\ref{i2:lem-fixed} such that
$\overline I\subset(-R_\perp,R_\perp)$. Choose $R_1\le L$ with
$D_L$ the working disk of Proposition~\ref{a28:moment-inputs}.
On the corresponding preselected common local parameter chart,
there exist frames
$X_1,Y_2$, selected independently of all subsequent polynomial
and transverse tests, such that
\begin{equation}\label{i2:eq-nine}
 \int_{|z|<\sqrt{R_1^2-x_\perp^2}}
 g(z)Y_2^T(\mathsf{P}_1-\mathsf{P}_2)X_1\,dA=0,
 \quad\text{for every }x_\perp\in I\text{ and }g\in\mathbb C[z].
\end{equation}
This is a full $3\times3$ matrix identity on each leaf, for every polynomial.
\end{theorem}
\begin{proof}

Let $X_1,Y_2$ be the frames on the common parameter chart from
Lemma~\ref{i2:lem-fixed}. Let $Z_1$ denote the four-component
lift constructed in the proof of Lemma~\ref{i2:lem-fixed}.
With $U_\theta=\operatorname{diag}(Q,1)$, it satisfies
\[
\mathsf T_{1,\theta}Z_1=0,
\qquad
U_\theta^TZ_1=X_1.
\]

Now fix $b\in C_c^\infty(I)$ and $g\in\mathbb C[z]$.
Since $b=b(x_\perp)$ and $g=g(z)$ is holomorphic,
\[
D_\theta b=0,\qquad D_\theta g=0.
\]
Consequently, the positive-phase leading amplitude
\[
A=bgZ_1
\]
satisfies
\[
\mathsf T_{1,\theta}A
=2D_\theta(bg)Z_1+bg\mathsf T_{1,\theta}Z_1=0.
\]
Choose the first correction $A_1$ to satisfy
\[
\mathsf T_{1,\theta}A_1=-M_1^{W_1}A.
\]
Applying the selected right inverse at $-\theta$ gives
\[
z_{1,h}=e^{\theta\cdot x/h}
(A+hA_1+r_h),
\]
with the estimates in~\eqref{pb27:two}. Apply the reduced
Betti identity columnwise to the three exact solutions $z_{1,h}$
associated with $Z_1$ and the three physical testing columns. The resulting nine
scalar identities assemble into the full matrix identity
used below.
Write $u_{2,h}^{(1)},u_{2,h}^{(2)},u_{2,h}^{(3)}$ for the three
elastic displacements of Lemma~\ref{i2:lem-fixed}, with the first two
rescaled by $-h$ and the third unscaled, and set
$\mathcal Y_{2,h}=(S_2^Tu_{2,h}^{(1)},
S_2^Tu_{2,h}^{(2)},S_2^Tu_{2,h}^{(3)})$.
The three negative-phase lifts then form
\[
 \mathcal Y_{2,h}=e^{-\theta\cdot x/h}(\mathcal Y_\chi^{(0)}+q_h),
 \quad\|q_h\|_2=O(h),\qquad
 \mathcal Y_\chi^{(0)}=U_\theta Y_2\ \text{near }\operatorname{supp}b.
\]
Here $\mathcal Y_\chi^{(0)}$ denotes the fixed leading matrix amplitude of the negative-phase tests.
Write $D^W=M_1^{W_1}-M_2^{W_2}=\sum_k\Delta A_k\partial_k+\Delta B^W$, where $\Delta A_k=A_{1,k}^{\rm phys}-A_{2,k}^{\rm phys}$, $\Delta A_\theta=\sum_k\theta_k\Delta A_k$, and $\Delta B^W$ is the difference of the zeroth-order matrices in~\eqref{bt:MW}. The notation $D^W$ denotes this difference operator. The factor $h$ in the following expansion isolates the first-order coefficient difference in the limit:
Throughout the target, the exact expansion is
\begin{align*}
 he^{-\theta\cdot x/h}D^Wz_{1,h}
 &=\Delta A_\theta A
 +h(\Delta A_\theta A_1+D^WA)
 +\Delta A_\theta r_h+h^2D^WA_1+hD^Wr_h\\
 &=\Delta A_\theta A+e_h,\qquad\|e_h\|_2=O(h).
\end{align*}
In particular, $\|\Delta A_\theta r_h\|_2=O(h^2)$ and $\|hD^Wr_h\|_2=O(h^2)$.
The potential and potential-gradient differences occur in the zeroth-order term $\Delta B^W$, carrying the factor $h$.
The phases cancel, and Cauchy--Schwarz gives
\begin{align*}
 \|h\mathcal Y_{2,h}^TD^Wz_{1,h}-(\mathcal Y_\chi^{(0)})^T\Delta A_\theta A\|_1
 &\le \|q_h\|_2\|\Delta A_\theta A\|_2
 +(\|\mathcal Y_\chi^{(0)}\|_2+\|q_h\|_2)\|e_h\|_2\\
 &=O(h).
\end{align*}
The first two elastic solutions are multiplied by the constant $-h$; the third is unscaled. Thus every column satisfies the original-domain identity
\eqref{i2:eq-betti}. We take the limit on $\Omega$ using the
preceding $L^1$ estimate. It gives
\[
 \int_\Omega bgY_2^TU_\theta^T\Delta A_\theta Z_1\,dx=0.
\]
The limiting integral extends to the working domain because
$\Delta A_\theta$ depends only on the Lam\'e pair and is supported
in $\overline\Omega$. This support property also holds with the
separate potential extensions of Proposition~\ref{ext:finite-common}.
Using $U_\theta^TZ_1=X_1$ and $A_{j,\theta}^{\rm phys}e_\theta=0$,
$U_\theta^T\Delta A_\theta Z_1=2(\mathsf{P}_1-\mathsf{P}_2)X_1$.
In leaf coordinates $dx=dA\,dx_\perp$. The matrix
$\mathsf F(x_\perp,z)=Y_2^T(\mathsf P_1-\mathsf P_2)X_1$
is continuous in $x_\perp$ with values in $L^1$ on any fixed leaf
disk containing $B_{R_0}$; the coefficient difference is supported in
$B_{R_0}$, and the fixed frames are continuous in the parameters.
Fubini therefore gives precisely hypothesis~\eqref{i2:leaf-extraction-hyp}
of Lemma~\ref{i2:leaf-extraction}, for every fixed $g$ and every
$b\in C_c^\infty(I)$. The lemma yields~\eqref{i2:eq-nine}
for every $x_\perp\in I$ and every polynomial $g$, with one and the
same choice of $X_1,Y_2$.
\end{proof}

\begin{remark}\label{i2:rem-scope}
The frames are fixed on each local parameter chart before the tests $b,g$;
the cutoffs, amplitudes, fields, constants, and threshold $h_0$ may depend
on each fixed test. The same remainder satisfies $L^2=O(h^2)$ and
$H^1=O(h)$ throughout the target, including the transverse cutoff transitions.
\end{remark}

\section{Relative transport and rigidity of the Lam\'e parameters}\label{if:section}

This section turns the complete matrix moment identities of
Theorem~\ref{i2:thm-nine} into a global relative transport on the
characteristic quotient and then uses its holomorphic structure to prove
rigidity of the Lam\'e parameters. The construction has two stages.
We first identify the primal and dual quotient connections and their
transformation laws between normalized leaf coordinates and the two
affine spinor charts of $\mathbb{CP}^1$. The moment identities then yield
a compactly supported primitive through the Cauchy transform. In dual
quotient coordinates this primitive defines a relative matrix $H$;
the chart transition laws show that its local representatives agree
on overlaps. The relative transport equation and the common exterior
normalization give the determinant identity
\[
 \det H=\sqrt{\frac{\mu_2}{\mu_1}}>0.
\]
Thus the resulting global endomorphism is an automorphism of the
characteristic quotient bundle. Its holomorphic dependence on the
projective null direction is established by differentiating the
transport equation while keeping the physical point $x$ fixed.

The second stage derives the rigidity of the material coefficients
from this global holomorphic structure. The holomorphic block
classification, together with the relative transport equations,
identifies the common parameter $\kappa_1=\kappa_2$ and reduces the
remaining Lam\'e parameter difference to a trace-free Hessian equation
for $\ell_{12}=\log(\mu_1/\mu_2)$. Introducing the four-component state
\[
 U=(\nabla\ell_{12},\Delta\ell_{12})^T
\]
converts this overdetermined second-order equation into a homogeneous
first-order system with a $4\times4$ coefficient matrix. The common
exterior normalization gives $U=0$ and $\ell_{12}=0$ on a nonempty
open set. Uniqueness along paths then yields $U=0$ throughout the
connected domain. Thus $\nabla\ell_{12}=0$, and the exterior value
fixes the constant to be zero. Consequently,
$\mu_1=\mu_2$, and the relation between $\lambda$, $\mu$, and
$\kappa$ gives $\lambda_1=\lambda_2$.

\subsection{The induced dual transport connection}
We use $\ell,\mathbf m,c,d,F_\mu$ with the meanings fixed in
\eqref{bt:notation}--\eqref{bt:MW}. In particular, $F_\mu$ is real and
symmetric, and $A_\theta^{\rm phys}$ is the contraction introduced in
Section~\ref{tb27:section}. A nonzero spinor $(s_0,s_1)$ determines

\[
 \theta=(s_0^2-s_1^2,i(s_0^2+s_1^2),2s_0s_1)^T,\quad
 Q=\frac12\begin{pmatrix}-s_1&-s_0\\is_1&-is_0\\s_0&-s_1\end{pmatrix},
 \quad n_{\mathrm{sp}}=(-s_1,s_0)^T,\quad \mathsf U_\theta=\operatorname{diag}(Q,1).
\]
Direct multiplication gives $\theta^TQ=0$, $Qn_{\mathrm{sp}}=-\theta/2$ and
$\operatorname{rank}Q=2$. Thus $\mathsf U_\theta$ identifies $\mathbb C^3$ with the dual
of the characteristic quotient,
\[
 E_\theta=\mathbb C^4/\mathbb C(\theta,0)^T,\qquad
 E_\theta^*=\{(p,\eta):\theta^Tp=0\}.
\]
The duality is the bilinear pairing fixed in Section~\ref{sec:notation}. In particular,
$E_\theta^*$ is the annihilator of the characteristic line.

\begin{proposition}\label{if:connection}
Put $r_{\mathbf{m}}=-2\mathbf{m}^TQ$. The dual connection induced by the physical reduction is
\begin{equation}\label{if:Gamma}
 \Gamma=\begin{pmatrix}n_{\mathrm{sp}}r_{\mathbf{m}}&2dn_{\mathrm{sp}}\\-\theta^TF_\mu Q&0\end{pmatrix},
 \qquad -(A_\theta^{\rm phys})^T\mathsf U_\theta=\mathsf U_\theta\Gamma.
\end{equation}
For every field $Y$,
\begin{equation}\label{if:dualtransport}
 (2D_\theta-(A_\theta^{\rm phys})^T)(\mathsf U_\theta Y)=\mathsf U_\theta(2D_\theta+\Gamma)Y.
\end{equation}
If $V$ is any right inverse of $\mathsf U_\theta^T$ that is independent of $x$, the quotient
connection matrix $\mathcal A=\mathsf U_\theta^TA_\theta^{\rm phys} V$ is independent of $V$ and
satisfies $\Gamma=-\mathcal A^T$.
\end{proposition}
\begin{proof}
The identity $A_\theta^{\rm phys}(\theta,0)^T=0$ shows that $A_\theta^{\rm phys}$ induces a
well-defined operator on the quotient. Using $Qn_{\mathrm{sp}}=-\theta/2$ and $F_\mu^T=F_\mu$,
we compute
\[
 \mathsf U_\theta\Gamma=
 \begin{pmatrix}\theta \mathbf{m}^TQ&-d\theta\\-\theta^TF_\mu Q&0\end{pmatrix}
 =-(A_\theta^{\rm phys})^T\mathsf U_\theta.
\]
Since $\mathsf U_\theta$ is independent of $x$, this also proves
\eqref{if:dualtransport}. If $\widetilde V$ is another right inverse,
then $\operatorname{ran}(V-\widetilde V)$ lies in
$\ker \mathsf U_\theta^T=\mathbb C(\theta,0)^T$, which $A_\theta^{\rm phys}$ annihilates. Hence
$\mathcal A$ is independent of the right inverse. Multiplying
\eqref{if:Gamma} on the left by $V^T$ gives
$\Gamma=-V^T(A_\theta^{\rm phys})^T\mathsf U_\theta=-\mathcal A^T$.
\end{proof}

The matrix $\mathsf{P}=\mathcal A/2$ in~\eqref{i2:eq-frames} is the primal quotient coefficient, whereas $\Gamma$ acts on dual quotient coordinates. Proposition~\ref{if:connection} therefore identifies them by $\Gamma=-2\mathsf{P}^T$. This identity converts the matrix moments into the relative dual transport equation used below.
The primal and dual frame equations are $(2D_\theta+\mathcal A)X=0$
and $(2D_\theta+\Gamma)Y=0$, respectively. In the leaf coordinate
$z$, the identity $D_\theta=2\partial_{\bar z}$ gives
\begin{equation}\label{if:leaf-frame-equations}
 2\partial_{\bar z}X+\mathsf P X=0,\qquad
 2\partial_{\bar z}Y-\mathsf P^TY=0.
\end{equation}

\subsection{Normalization and transition between the two affine charts}
The moment construction uses normalized real leaf coordinates, whereas holomorphic dependence on direction is expressed in spinor coordinates. We record their change of frame before constructing the relative field. For an invertible relative matrix $H$ in dual quotient coordinates, write $G=H^{-T}$ for the corresponding primal matrix; Proposition~\ref{if:actualH} will construct such a matrix $H$ and prove its invertibility.
Let $\theta_{\mathrm{hol}},Q_{\mathrm{hol}},n_{\mathrm{sp,hol}},\mathsf U_{\mathrm{hol}}$ denote the quantities associated with a holomorphic spinor representative. The subscripts $\mathrm{hol}$ and $\mathrm{nor}$ distinguish this representative from its real-normalized version. Direct calculation gives
\[
 \theta_{\mathrm{hol}}\cdot\overline{\theta_{\mathrm{hol}}}=2(|s_0|^2+|s_1|^2)^2.
\]
Set
\[
 a_{\rm nor}=(|s_0|^2+|s_1|^2)^{-1}>0,\quad
 \varsigma=\sqrt{a_{\rm nor}},\quad
 S=\operatorname{diag}(\varsigma I_2,1),\quad
 \theta_{\mathrm{nor}}=a_{\rm nor}\theta_{\mathrm{hol}}.
\]
The normalized spinor is $\varsigma(s_0,s_1)$ and $\mathsf U_{\mathrm{nor}}=\mathsf U_{\mathrm{hol}}S$. Coordinates of
the same ambient dual vector satisfy $v_{\mathrm{hol}}=Sv_{\mathrm{nor}}$, while primal quotient
coordinates satisfy $u_{\mathrm{nor}}=S^Tu_{\mathrm{hol}}$. Accordingly,
\begin{align}
 Y_{\mathrm{hol}}&=SY_{\mathrm{nor}},&X_{\mathrm{hol}}&=S^{-T}X_{\mathrm{nor}},\notag\\
 \Gamma_{\mathrm{nor}}&=a_{\rm nor}S^{-1}\Gamma_{\mathrm{hol}}S,&
 \mathcal A_{\mathrm{nor}}&=a_{\rm nor}S^T\mathcal A_{\mathrm{hol}}S^{-T},\notag\\
 H_{\mathrm{hol}}&=SH_{\mathrm{nor}}S^{-1},&G_{\mathrm{hol}}&=S^{-T}G_{\mathrm{nor}}S^T.
 \label{if:normalization}
\end{align}
These change-of-frame matrices depend only on the direction and are
therefore independent of $x$. Real normalization involves
$\bar\zeta$; holomorphic dependence will be established after passing
back to the holomorphic spinor representatives.

The two affine charts of $\mathbb{CP}^1$ use the spinor representatives $(1,\zeta)$ with $\zeta=s_1/s_0$ and $(\upsilon,1)$ with $\upsilon=s_0/s_1$. On their overlap $\upsilon=\zeta^{-1}$. Put
$S_{\zeta}=\operatorname{diag}(\zeta I_2,1)$. Homogeneity gives
\[
 \theta_0=\zeta^2\theta_\infty,\quad \mathsf U_0=\mathsf U_\infty S_{\zeta},\quad
 \Gamma_0=\zeta^2S_{\zeta}^{-1}\Gamma_\infty S_{\zeta},
 \quad\mathcal A_0=\zeta^2S_{\zeta}^T\mathcal A_\infty S_{\zeta}^{-T}.
\]
The relative matrices therefore obey
\begin{equation}\label{if:chart}
 H_0=S_{\zeta}^{-1}H_\infty S_{\zeta},\qquad G_0=S_{\zeta}^TG_\infty S_{\zeta}^{-T}.
\end{equation}
These are the coordinate expressions of the corresponding quotient-bundle
automorphisms.

\begin{lemma}
\label{if:leaf-uniqueness}
Let $D\subset\mathbb C$ be an arbitrary connected open set, and let $\mathsf P_1,\mathsf P_2$ be continuous
$3\times3$ matrix fields on $D$.
Suppose that invertible $C^1$ frames $X_1,Y_2$ on $D$ satisfy
\[
 2\partial_{\bar z}X_1=-\mathsf P_1X_1,
 \qquad 2\partial_{\bar z}Y_2=\mathsf P_2^TY_2.
\]
If $H$ and $\widetilde H$ are $C^1$ solutions of
\begin{equation}\label{if:leaf-relative}
 2\partial_{\bar z}H-\mathsf P_1^TH+H\mathsf P_2^T=0
\end{equation}
on $D$ and agree on a nonempty open subset, then
$H=\widetilde H$ throughout $D$. The conclusion is unaffected by a
change of quotient or dual-quotient frame that is independent of $x$,
provided both connections and both relative fields
are transformed accordingly.
\end{lemma}
\begin{proof}
Set $V=H-\widetilde H$. By the three displayed transport equations,
\[
 2\partial_{\bar z}(X_1^TVY_2)
 =-X_1^T\mathsf P_1^TVY_2
  +X_1^T(\mathsf P_1^TV-V\mathsf P_2^T)Y_2
  +X_1^TV\mathsf P_2^TY_2=0.
\]
Thus each entry of $X_1^TVY_2$ is holomorphic. It vanishes on the
prescribed open set and hence on all of connected $D$ by the identity
theorem. Since both frames are invertible, $V=0$ on $D$.
An $x$-independent invertible frame change conjugates the transformed
relative equation to~\eqref{if:leaf-relative}; it also takes equality
on the open set to equality there. This proves the last assertion.
\end{proof}

\subsection{Construction from the matrix moment identities}
\begin{lemma}\label{if:compact-primitive}
Let $0<R_0<L$ and let $\mathsf F\in C_c^2(D_L;\mathbb C^{3\times3})$
be supported in $\overline D_{R_0}$, extended by zero to $\mathbb C$.
Define its moments and Cauchy transform by
\begin{equation}\label{if:primitive-data}
 \mathsf M_q=\int_{\mathbb C}w^q\mathsf F(w)\,dA(w),\qquad
 (\mathscr C\mathsf F)(z)=\frac1\pi\int_{\mathbb C}
                          \frac{\mathsf F(w)}{z-w}\,dA(w),\quad q\ge0.
\end{equation}
The following conditions are equivalent: all $\mathsf M_q$ vanish;
$\mathscr C\mathsf F=0$ for $|z|>R_0$; and there exists a $C^2$
matrix $\Psi$ supported in $\overline D_{R_0}$ with
$\partial_{\bar z}\Psi=\mathsf F$. Such a compactly supported
primitive is unique and equals $\mathscr C\mathsf F$.
For any fixed matrix norm and its associated $L^1$ norm,
\begin{equation}\label{if:cauchy-tail}
 \left\|\mathscr C\mathsf F(z)-\frac1\pi
       \sum_{q=0}^N\frac{\mathsf M_q}{z^{q+1}}\right\|
 \le\frac{\|\mathsf F\|_{L^1}}{\pi(|z|-R_0)}
       \left(\frac{R_0}{|z|}\right)^{N+1},
 \qquad |z|>R_0,\quad N\ge0.
\end{equation}
Under the equivalent conditions, the tests extend from polynomials to
all scalar holomorphic functions on $D_L$:
\begin{equation}\label{if:holomorphic-tests}
 \int_{D_L}\varphi(z)\mathsf F(z)\,dA(z)=0,
 \qquad \varphi\in\mathcal O(D_L).
\end{equation}
Here $\mathcal O(D_L)$ denotes the holomorphic functions on the disk.
\end{lemma}
\begin{proof}
The distributional identity~\eqref{f27:cauchy-identities} gives
$\partial_{\bar z}\mathscr C\mathsf F=\mathsf F$.
For $C^2$ compactly supported sources the real derivatives through
order two transfer to the source. Local integrability of
$|z-w|^{-1}$ then gives $C^2$ regularity of the transform. The finite geometric identity
\[
 \frac1{z-w}=\sum_{q=0}^N\frac{w^q}{z^{q+1}}
       +\frac{w^{N+1}}{z^{N+1}(z-w)}
\]
gives~\eqref{if:cauchy-tail} by termwise integration of the
finite sum and estimation of its remainder. If all moments vanish,
letting
$N\to\infty$ at each $|z|>R_0$ proves exterior vanishing and hence
the asserted compact primitive.

Conversely, suppose $\partial_{\bar z}\Psi=\mathsf F$ with the
specified support. Given $\varphi\in\mathcal O(D_L)$, choose
$\chi\in C_c^\infty(D_L)$ equal to one near $\overline D_{R_0}$.
Integration by parts gives
\[
 \int_{D_L}\varphi\mathsf F\,dA
 =-\int_{D_L}\Psi\,\partial_{\bar z}(\chi\varphi)\,dA=0,
\]
since $\partial_{\bar z}\varphi=0$ and $\nabla\chi$ vanishes on
$\supp\Psi$. This proves~\eqref{if:holomorphic-tests} and, in
particular, all polynomial moments. The difference of two $C^2$
compact primitives is an entire holomorphic matrix and vanishes
outside a compact set; the identity theorem makes it zero. Applying
this to $\Psi$ and the already compact $\mathscr C\mathsf F$
proves uniqueness and the formula.
\end{proof}

\begin{proposition}\label{if:intrinsic-relative}
Let $\mathsf P_1,\mathsf P_2$ be $C^2$ matrix fields on $D_L$, with
$\supp(\mathsf P_1-\mathsf P_2)\subset\overline D_{R_0}$, where
$0<R_0<L$. Let $X_1,Y_2$ be $C^2$ invertible frames satisfying
\[
 2\partial_{\bar z}X_1=-\mathsf P_1X_1,\qquad
 2\partial_{\bar z}Y_2=\mathsf P_2^TY_2.
\]
Put $\mathsf F_{\rm mom}=Y_2^T(\mathsf P_1-\mathsf P_2)X_1$ and
$\mathsf W=Y_2^TX_1$. Vanishing of every polynomial moment of
$\mathsf F_{\rm mom}$ is equivalent to the existence of a $C^2$
solution $H$ of~\eqref{if:leaf-relative} with $H=I_3$ on
$D_L\setminus\overline D_{R_0}$. This solution is unique and is
explicitly
\begin{equation}\label{if:intrinsic-formula}
 H=I_3+\tfrac12X_1^{-T}(\mathscr C\mathsf F_{\rm mom})^TY_2^{-1}.
\end{equation}
Invertibility is established after the global construction by
\eqref{if:detH}. The field is independent of the choices of
the two solution frames, while
holding the quotient coordinates fixed.
\end{proposition}
\begin{proof}
The leaf frame equations~\eqref{if:leaf-frame-equations} give
$2\partial_{\bar z}\mathsf W=-\mathsf F_{\rm mom}$.
If all moments vanish, Lemma~\ref{if:compact-primitive} gives the
compact primitive $\mathscr C\mathsf F_{\rm mom}$, so
\[
 \mathsf R_{\rm hol}=\mathsf W+\tfrac12\mathscr C\mathsf F_{\rm mom}
\]
is holomorphic. Then
$H=X_1^{-T}\mathsf R_{\rm hol}^TY_2^{-1}$ equals
\eqref{if:intrinsic-formula}. Differentiating this expression gives
\eqref{if:leaf-relative}; it equals $I_3$ in the outer annulus.
Uniqueness follows from Lemma~\ref{if:leaf-uniqueness}.
Conversely, for any such $H$ the frame equations imply
$\partial_{\bar z}(Y_2^TH^TX_1)=0$. Hence
\[
 \Psi=2Y_2^T(H^T-I_3)X_1
\]
is of class $C^2$, supported in $\overline D_{R_0}$, and satisfies
$\partial_{\bar z}\Psi=\mathsf F_{\rm mom}$ after zero extension.
Lemma~\ref{if:compact-primitive} proves all the moments.

To establish independence directly, let $\widetilde X_1,\widetilde Y_2$
be any other invertible frames for these same equations. Define
\begin{equation}\label{if:solution-frame-change}
 \mathsf G_X=X_1^{-1}\widetilde X_1,\qquad
 \mathsf G_Y=Y_2^{-1}\widetilde Y_2,\qquad
 \partial_{\bar z}\mathsf G_X=\partial_{\bar z}\mathsf G_Y=0.
\end{equation}
The last two equalities follow by differentiating the first two
expressions and canceling the identical connection terms.
For the new source,
$\widetilde{\mathsf F}_{\rm mom}
=\mathsf G_Y^T\mathsf F_{\rm mom}\mathsf G_X$.
The $C^2$ compactly supported matrix
$\mathsf G_Y^T(\mathscr C\mathsf F_{\rm mom})\mathsf G_X$,
extended by zero outside $D_L$, has
$\partial_{\bar z}$ derivative $\widetilde{\mathsf F}_{\rm mom}$.
The lemma therefore gives all its moments and the exact identity
\begin{equation}\label{if:cauchy-frame-covariance}
 \mathscr C\widetilde{\mathsf F}_{\rm mom}
 =\mathsf G_Y^T(\mathscr C\mathsf F_{\rm mom})\mathsf G_X
 \quad\hbox{on }D_L.
\end{equation}
Inserting this into~\eqref{if:intrinsic-formula}, with
$\widetilde X_1=X_1\mathsf G_X$ and
$\widetilde Y_2=Y_2\mathsf G_Y$, cancels both frame factors
in their stated order and gives exactly the same $H$.
The inverse factors are also holomorphic, so the covariance identity
holds for arbitrary holomorphic frame changes arising from the same
transport equations.
\end{proof}

The solution-frame changes in~\eqref{if:solution-frame-change} leave
$H$ unchanged with the quotient coordinates fixed. The
compact-primitive covariance identity~\eqref{if:cauchy-frame-covariance}
transfers the moment identities to any other invertible frames
solving the same transport equations. Spinor coordinate changes
act by the conjugation in~\eqref{if:chart}.

\begin{proposition}\label{if:actualH}
Suppose the two extended Lam\'e pairs are of class $C^4$, agree on
$B_R\setminus\overline\Omega$, and have differences supported in
$B_{R_0}$, with $R_0<R$. For each normalized direction and each centre
of the transverse parameter interval, suppose the moment construction
provides invertible, jointly $C^2$ frames $X_1,Y_2$ on a
fixed leaf disk $D_L$, $L>R$, over a local transverse and direction
parameter chart, such that
\[
 2\partial_{\bar z}X_1=-\mathsf{P}_1X_1,\qquad
 2\partial_{\bar z}Y_2=\mathsf{P}_2^TY_2,
 \qquad \int_{D_L}z^qY_2^T(\mathsf{P}_1-\mathsf{P}_2)X_1\,dA=0\quad(q\ge0).
\]
The frames are fixed before the tests are chosen, and the identities hold
for every leaf parameter in the local chart. Then there is a unique
jointly $C^2$ relative bundle automorphism $H$ over
$B_R\times\mathbb{CP}^1$ with the properties in
Proposition~\ref{if:binding}. Its determinant is
\begin{equation}\label{if:detH}
 \det H=\sqrt{\mu_2/\mu_1}>0.
\end{equation}
\end{proposition}
\begin{proof}
On each fixed leaf of a local parameter cylinder, apply
Proposition~\ref{if:intrinsic-relative} with
$\mathsf F_{\rm mom}=Y_2^T(\mathsf P_1-\mathsf P_2)X_1$.
Its hypotheses hold because the coefficient difference is supported
inside $B_{R_0}$, and hence inside $\overline D_{R_0}$ on every leaf.
Formula~\eqref{if:intrinsic-formula} gives $H$ satisfying
\eqref{if:leaf-relative} and equal to $I_3$ for $|z|>R_0$.
The Cauchy transform $\mathscr C$ in~\eqref{if:primitive-data}
has twice the normalization of $\mathcal C_\theta$ defined in
Definition~\ref{i2:def-input}. Proposition~\ref{if:connection}
converts the leaf equation to
$2D_\theta H+\Gamma_1H-H\Gamma_2=0$.

The determinant identity now establishes invertibility. Formula~\eqref{if:Gamma} gives $\operatorname{tr}\Gamma_j=\mathbf{m}_j\cdot\theta$. Let $\operatorname{adj}(H)$ denote the classical adjugate, the transpose of the cofactor matrix. It is defined for every square matrix and satisfies $H\operatorname{adj}(H)=\det(H)I$. The columnwise determinant identity reads
\[
 D_\theta\det H=\operatorname{tr}(\operatorname{adj}(H)D_\theta H).
\]
Together with
$H\operatorname{adj}(H)=\operatorname{adj}(H)H=\det(H)I$, the relative
equation yields
\[
 2D_\theta\det H+(D_\theta\log\mu_1-D_\theta\log\mu_2)\det H=0.
\]
Thus $\sqrt{\mu_1/\mu_2}\det H$ is holomorphic on each leaf and equals one
on an exterior open subset. The identity theorem on the connected leaf
disk proves \eqref{if:detH}. In particular, $H$ is invertible.

The source is supported in a fixed compact subset of the leaf disk
on a smaller parameter chart. In the fixed leaf coordinates,
integration by parts transfers leaf derivatives to the source.
Parameter derivatives act on the source and the pullback on the
same disk. Local integrability of $|z-w|^{-1}$ bounds all derivatives
through order two, and the chain rule gives joint $C^2$ regularity
in physical coordinates.

For a fixed transverse coordinate $t$, the exterior region
\[
 E_t^{\mathrm{ext}}=\{z\in D_L:|z|^2+t^2>R_0^2\}
\]
is connected: it is an annulus with inner radius
$\sqrt{R_0^2-t^2}$ if $|t|<R_0$, a punctured disk if $|t|=R_0$,
and the full disk if $|t|>R_0$.
Since $L>R>R_0$, every $E_t^{\mathrm{ext}}$ meets $\{|z|>R_0\}$.
The material connections agree on $E_t^{\mathrm{ext}}$, so $I$ solves
the relative transport equation there. The frames on the full disk
restrict to invertible frames on $E_t^{\mathrm{ext}}$.
Lemma~\ref{if:leaf-uniqueness} therefore gives $H=I$ throughout that region.
For $|t|\ge R_0$, the moment source is supported in the empty
intersection of the leaf with the \emph{open} ball $B_{R_0}$, so
it vanishes on the entire leaf and the construction gives $H=I$
directly, including at equality.

On an overlap of two local direction and transverse parameter charts,
express both constructions in the same dual quotient coordinates. The changes from
real-normalized coordinates and between the two spinor charts are given
by~\eqref{if:normalization} and~\eqref{if:chart}; their factors are
independent of $x$. In particular, on a chart
overlap, if $S_\zeta=\operatorname{diag}(\zeta I_2,1)$ and
$\widetilde H_0=S_\zeta^{-1}H_\infty S_\zeta$, then
\[
 2D_{\theta_0}\widetilde H_0+
   \Gamma_{1,0}\widetilde H_0-
   \widetilde H_0\Gamma_{2,0}
 =\zeta^2S_\zeta^{-1}
  (2D_{\theta_\infty}H_\infty+
    \Gamma_{1,\infty}H_\infty-
    H_\infty\Gamma_{2,\infty})S_\zeta=0.
\]
Both local constructions therefore solve the same equation in these
coordinates and equal $I$ on the same exterior shell. Each connected
section $B_R\cap\{\nu_\theta\cdot x=t\}$ is a disk and meets
the common exterior normalization shell
$E_{\mathrm{ext}}=\{R_0<|x|<R\}$ in a nonempty relatively open set. Applying
Lemma~\ref{if:leaf-uniqueness} on that section shows that the two
constructions agree wherever both are defined. Independence of the
solution frames on a common disk is also the explicit covariance
\eqref{if:cauchy-frame-covariance}; this is distinct from the spinor
coordinate conjugation used above.

Choose $R_0<R_\perp<R$, use the local construction for
$|\nu_\theta\cdot x|<R_\perp$, and set $H=I$ where
$|\nu_\theta\cdot x|>R_0$. These sets cover $B_R$, and on their
overlap $H=I$ by the exterior argument. The local constructions are jointly $C^2$ in $x$ and the real
direction parameters. Agreement on overlaps patches them to a
jointly $C^2$ bundle field over
$B_R\times\mathbb{CP}^1$. Formula~\eqref{if:detH} is invariant
under these coordinate conjugations and makes the patched field an
automorphism. Proposition~\ref{if:binding} gives holomorphic dependence on the
direction at each fixed physical point and the transition
law~\eqref{if:chart}.
\end{proof}

For the application of this proposition, Theorem~\ref{i2:thm-nine}
supplies the frames and moment identities for smooth inputs.
Subsection~\ref{fr:interior-assembly} establishes these inputs
at the regularity of Theorem~\ref{a28:main}. The resulting field
is constructed from the identities on the original physical domain.

\subsection{Uniqueness and holomorphic dependence on direction}
\begin{proposition}\label{if:binding}
Suppose $H$ is a jointly $C^2$ invertible relative field on the fixed
physical ball $B_R$, with
\[
 2D_\theta H+\Gamma_1H-H\Gamma_2=0,
 \qquad H=I\quad\text{on }E_{\mathrm{ext}}=\{R_0<|x|<R\}
\]
for every direction. If $\Gamma_j$ are given by \eqref{if:Gamma}, then $H$
is holomorphic in the two affine spinor charts for each fixed physical
point $x$ and satisfies \eqref{if:chart}.
\end{proposition}
\begin{proof}
After the changes of frame in \eqref{if:normalization} and
\eqref{if:chart}, all candidates satisfy the same relative equation in $x$ and agree with $I$ on the same shell. Their difference satisfies a homogeneous relative $\bar\partial$
system on each connected leaf section. The frames required in
Lemma~\ref{if:leaf-uniqueness} exist by the planar frame construction.
Every such section meets $E_{\mathrm{ext}}$ in a nonempty relatively open set, so
Lemma~\ref{if:leaf-uniqueness} proves uniqueness and the stated
transition law.

In a holomorphic spinor chart, the three nonzero blocks of
\eqref{if:Gamma} have respective direction degrees $2,1,3$, and the
material coefficients depend only on $x$. Hence
$\partial_{\bar\zeta}\Gamma_j=0$ and
$\partial_{\bar\zeta}\theta=0$ at fixed $x$. Use
\eqref{if:normalization} to express the transport equation in the
holomorphic spinor representative, where the normalization factor
is absent. Joint $C^2$ regularity permits differentiation at fixed $x$
and commutation of the derivatives in $x$ and $\bar\zeta$. Thus
$V=\partial_{\bar\zeta}H|_x$ satisfies
\[
 2D_\theta V+\Gamma_1V-V\Gamma_2=0,\qquad V|_{E_{\mathrm{ext}}}=0.
\]
Let $X_1$ and $Y_2$ be invertible primal and dual frames on a leaf, expressed in the current holomorphic spinor representative. Their equations and the equation for $V$ give $D_\theta(X_1^TVY_2)=0$. This matrix is therefore holomorphic in the leaf variable and vanishes on the nonempty relatively open intersection with $E_{\mathrm{ext}}$. The identity theorem gives $V=0$ on the connected leaf section, in each full affine chart.
\end{proof}

\subsection{Finite-type rigidity for the Lam\'e parameters}
\label{i4:section}

Propositions~\ref{if:actualH} and~\ref{if:binding} give a relative transport matrix with a common exterior normalization and holomorphic dependence on direction at fixed physical points. We now prove that these properties determine the Lam\'e parameters. Lemma~\ref{i4:classification} first restricts the matrix blocks by their chart transition laws. The transport equation then determines all blocks in terms of $\ell_{12}=\log(\mu_1/\mu_2)$ and yields the trace-free Hessian identity~\eqref{i4:tf-Hessian}. Differentiating this identity once gives a closed homogeneous first-order system for the three components of $\nabla\ell_{12}$ and the scalar $\Delta\ell_{12}$. Zero data on the normalization set propagate through this system.

Let $\mathcal{O}\subset\mathbb R^3$ be connected and open and let $E_{\mathrm{eq}}\subset \mathcal{O}$ be
nonempty and open. For $j=1,2$, let $\mu_j,\kappa_j\in C^3(\mathcal O;\mathbb R)$,
with $\mu_j>0$ and $1/4<\kappa_j<1$, and assume
$(\mu_1,\kappa_1)=(\mu_2,\kappa_2)$ on $E_{\mathrm{eq}}$.
Thus $E_{\mathrm{eq}}$ is the open equality set used in the rigidity argument.
Use $\ell_j,\mathbf m_j,c_j,F_{\mu_j}$ as in
\eqref{bt:notation} and~\eqref{bt:MW}, and put
$d_j=\kappa_j\sqrt{\mu_j}$. We retain the tensor and gradient
conventions used in the reduction and the bilinear pairing convention
of Section~\ref{sec:notation}.
For every nonzero spinor $(s_0,s_1)\in\mathbb C^2$, define
\begin{equation}\label{i4:spinor}
 \begin{aligned}
 \theta&=(s_0^2-s_1^2,i(s_0^2+s_1^2),2s_0s_1)^T,\\
 Q&=\frac12\begin{pmatrix}-s_1&-s_0\\is_1&-is_0\\s_0&-s_1\end{pmatrix},\\
 n_{\mathrm{sp}}&=(-s_1,s_0)^T,\qquad R_a=-2a^TQ.
 \end{aligned}
\end{equation}
Direct multiplication gives $\theta^T\theta=0$, $\theta^TQ=0$,
$Qn_{\mathrm{sp}}=-\theta/2$ and $R_an_{\mathrm{sp}}=a\cdot\theta$.
The three minors of $Q$ are $is_0s_1/2,(s_0^2+s_1^2)/4,i(s_0^2-s_1^2)/4$;
they do not vanish simultaneously. Thus $Q$ has rank two, including at
$s_0=0$ and at $s_1=0$.

\begin{theorem}\label{i4:main}
Suppose there is an invertible dual quotient automorphism $H$ on
$\mathcal{O}\times\mathbb{CP}^1$, jointly $C^2$ in $x$ and the real direction
coordinates, holomorphic in projective direction for every fixed physical $x$,
and equal to $I_3$ on $E_{\mathrm{eq}}$ for every direction. Suppose its matrices in the two
full affine spinor charts $(1,\zeta)$ and $(\upsilon,1)$ obey
\begin{equation}\label{i4:input}
\begin{aligned}
 H_0(\zeta)=S_\zeta^{-1}H_\infty(\zeta^{-1})S_\zeta,
 \qquad S_\zeta=\operatorname{diag}(\zeta I_2,1),\\
 2D_\theta H+\Gamma_1H-H\Gamma_2=0,\qquad
 \Gamma_j=\begin{pmatrix}n_{\mathrm{sp}}R_{\mathbf{m}_j}&2d_jn_{\mathrm{sp}}\\-\theta^TF_{\mu_j}Q&0\end{pmatrix},
 \quad D_\theta=\theta\cdot\nabla_x.
\end{aligned}
\end{equation}
Then $\mu_1=\mu_2$, $\kappa_1=\kappa_2$ and $H=I_3$ on $\mathcal{O}$.
Consequently, if
$\lambda_j=\mu_j(2\kappa_j-1)/(1-\kappa_j)$, then
$\lambda_1=\lambda_2$.
\end{theorem}

\begin{lemma}\label{i4:classification}
Under the direction and chart hypotheses of Theorem~\ref{i4:main},
\begin{equation}\label{i4:blocks}
 H=\begin{pmatrix}A&0\\ C&\delta_{H}\end{pmatrix},\qquad
 C=s_0C_{s_0}+s_1C_{s_1},
\end{equation}
Here $A$ is a $2\times2$ invertible matrix, $\delta_H$ is the nonzero
lower-right scalar block of $H$, and $C_{s_0},C_{s_1}$ are $1\times2$ row vectors. All are of class $C^2$ in $x$; $A$ and $\delta_H$ are independent of direction, while the row $C$ is linear in the homogeneous spinor coordinates.
\end{lemma}
\begin{proof}
Write the chart blocks as $(A,B;C,\delta_{H})$. Their transition laws are
$A_0=A_\infty$, $\delta_{H,0}=\delta_{H,\infty}$,
$B_0=\zeta^{-1}B_\infty$, $C_0=\zeta C_\infty$, with infinity entries
evaluated at $\zeta^{-1}$. Holomorphy on each full affine chart implies, for large
$|\zeta|$, bounds $A_0,\delta_{H,0}=O(1)$, $B_0=O(|\zeta|^{-1})$ and
$C_0=O(1+|\zeta|)$. For an entire scalar entry $f$ satisfying
$|f(\zeta)|\le C(1+|\zeta|)^m$ with $m=0$ or $1$, the Cauchy formula
on the circle of radius $R$ centered at any fixed $\zeta_*$ gives
\[
 |f^{(m+1)}(\zeta_*)|
 \le (m+1)!\,C(1+|\zeta_*|+R)^mR^{-m-1}\longrightarrow0.
\]
Thus $A_0,\delta_{H,0}$ are constant and $C_0$ is affine. The entire
entries of $B_0$ are bounded and tend to zero at infinity, so the same
argument with $m=0$ makes them identically zero. Consequently
$B_0=0$ and $C_0=C_{s_0}+\zeta C_{s_1}$. The infinity expression is
$C_\infty=\upsilon C_{s_0}+C_{s_1}$. The coefficients are of class $C^2$ in $x$, by evaluation at fixed
chart points. The determinant
is $(\det A)\delta_{H}$, so invertibility gives the asserted nonvanishing.

The dual ambient frame is $\mathsf U_\theta=\operatorname{diag}(Q,1)$.
For $b\in\mathbb C\setminus\{0\}$, put $S_b=\operatorname{diag}(bI_2,1)$. Replacing the spinor $(s_0,s_1)$ by $b(s_0,s_1)$ changes the frame to $\mathsf U_\theta S_b$, so the matrix of the same dual quotient map changes to $S_b^{-1}HS_b$. The resulting homogeneities are
\begin{equation}\label{i4:weights}
\begin{aligned}
 &(Q,n_{\mathrm{sp}},C)\mapsto b(Q,n_{\mathrm{sp}},C),\quad\theta\mapsto b^2\theta,
 \quad(A,\delta_{H})\mapsto(A,\delta_{H}),\\
 &(\Gamma_j^{11},\Gamma_j^{12},\Gamma_j^{21})
       \mapsto(b^2\Gamma_j^{11},b\Gamma_j^{12},b^3\Gamma_j^{21}).
\end{aligned}
\end{equation}
where $\Gamma_j^{11}=n_{\mathrm{sp}}R_{\mathbf{m}_j}$, $\Gamma_j^{12}=2d_jn_{\mathrm{sp}}$ and $\Gamma_j^{21}=-\theta^TF_{\mu_j}Q$.
Thus the block equations have the induced spinor homogeneities.
\end{proof}

\begin{lemma}\label{i4:determined}
Set $\ell_{12}=\log(\mu_1/\mu_2)$ and $r=\nabla\ell_{12}$. Then
\begin{equation}\label{i4:determined-equations}
 A=I_2,\qquad \delta_{H}=e^{-\ell_{12}/2},\qquad
 \kappa_1=\kappa_2=:\kappa,\qquad
 C=-\frac{R_r}{2d_1}=\mathbf{g}^TQ,\qquad
 \mathbf{g}=\frac{c_1}{\kappa}r.
\end{equation}
\end{lemma}
\begin{proof}
Multiplication of the matrices in their given material order yields
\begin{align}\label{i4:four-equations}
 2D_\theta A+\Gamma_1^{11}A+\Gamma_1^{12}C-A\Gamma_2^{11}&=0,\nonumber\\
 \Gamma_1^{12}\delta_{H}-A\Gamma_2^{12}&=0,\nonumber\\
 2D_\theta C+\Gamma_1^{21}A-C\Gamma_2^{11}-\delta_{H} \Gamma_2^{21}&=0,\\
 2D_\theta\delta_{H}-C\Gamma_2^{12}&=0.\nonumber
\end{align}
Their homogeneous spinor degrees are respectively $2,1,3,2$, so the equations
in both affine charts hold for every homogeneous nonzero spinor. These
spinors cover every complex null ray: when $\theta_1-i\theta_2\ne0$ the
chart parameter is $\theta_3/(\theta_1-i\theta_2)$, and the exceptional ray
is represented by $(0,1)$.

The second equation gives $d_1\delta_{H} n_{\mathrm{sp}}=d_2An_{\mathrm{sp}}$. Since $n_{\mathrm{sp}}(1,0)=(0,1)^T$
and $n_{\mathrm{sp}}(0,1)=(-1,0)^T$, it follows that $A=\alpha I_2$ with
$\alpha=d_1\delta_{H}/d_2$. The first equation becomes
\[
 2(D_\theta\alpha)I_2+
 n_{\mathrm{sp}}\{\alpha(R_{\mathbf{m}_1}-R_{\mathbf{m}_2})+2d_1C\}=0.
\]
The second term has rank at most one, whereas a nonzero multiple of $I_2$
has rank two. Hence $D_\theta\alpha=0$ for every null direction. Since
$2e_i=(e_i+ie_j)+(e_i-ie_j)$ for $j\ne i$, $\nabla\alpha=0$.
Connectedness and the data on $E_{\mathrm{eq}}$ give $\alpha=1$. The nonzero column $n_{\mathrm{sp}}$
then gives $C=-R_r/(2d_1)$ and
$\delta_{H}=d_2/d_1>0$.

The fourth equation, using $R_rn_{\mathrm{sp}}=D_\theta\ell_{12}$, now reads
$2D_\theta\delta_{H}+\delta_{H} D_\theta\ell_{12}=0$.
Therefore $\delta_{H} e^{\ell_{12}/2}$ is constant on $\mathcal{O}$, equal to one on $E_{\mathrm{eq}}$.
Thus $\delta_{H}=e^{-\ell_{12}/2}$. Since
$d_2/d_1=(\kappa_2/\kappa_1)e^{-\ell_{12}/2}$, cancellation gives the common
$\kappa$. Finally $-R_r/(2d_1)=(c_1/\kappa)r^TQ$.
\end{proof}

\begin{lemma}\label{i4:null-lines}
The third equation of \eqref{i4:four-equations} is equivalent to
$Q^TK_{\mathrm{rig}}\theta=0$, where
\begin{equation}\label{i4:K}
 K_{\mathrm{rig}}=2\nabla\mathbf{g}-F_{\mu_1}+e^{-\ell_{12}/2}F_{\mu_2}-\mathbf{m}_2\mathbf{g}^T.
\end{equation}
Here $K_{\mathrm{rig}}$ is the $3\times3$ matrix used in the rigidity
argument, distinct from the four-component differential operator
$K_{\mathrm{lift}}$. For every $x$, $K_{\mathrm{rig}}(x)$ is a scalar multiple of $I_3$.
\end{lemma}
\begin{proof}
With $C=\mathbf{g}^TQ$, the third equation becomes
\[
 \theta^T\{2(\nabla\mathbf{g})^T-F_{\mu_1}+e^{-\ell_{12}/2}F_{\mu_2}
                         -\mathbf{g} \mathbf{m}_2^T\}Q=0.
\]
Indeed, $D_\theta C=\theta^T(\nabla\mathbf{g})^TQ$ and
$-Cn_{\mathrm{sp}}R_{\mathbf{m}_2}=-(\mathbf{g}\cdot\theta)\mathbf{m}_2^TQ$. Transposition gives
\eqref{i4:K}.
Since $\operatorname{ran}Q=\theta^\perp$ and the ambient complex bilinear
form is nondegenerate, $\ker Q^T=\mathbb C\theta$, even for a null
$\theta$. Thus $K_{\mathrm{rig}}\theta\in\mathbb C\theta$ for all null $\theta$.
At each fixed $x$, $K_{\mathrm{rig}}$ is independent of the direction variable. For distinct $i,j,k$,
the $k$ components of $K_{\mathrm{rig}}(e_i\pm ie_j)$ give
$(K_{\mathrm{rig}})_{ki}\pm i(K_{\mathrm{rig}})_{kj}=0$, hence both entries vanish. Choosing each remaining
index removes all six off-diagonal entries. Preservation of
$\mathbb C(e_i+ie_j)$ then makes the diagonal entries equal. Therefore
$K_{\mathrm{rig}}=qI_3$, with $q=\operatorname{tr}K_{\mathrm{rig}}/3$ real-valued and of class $C^1$.
\end{proof}

\begin{proposition}\label{i4:beta}
Let $\mathbf{k}_{\kappa}=\nabla\log\kappa$,
$S=\mathbf{m}_1+\mathbf{m}_2$ and
$\Xi=\Delta\ell_{12}-S\cdot r$. Then the full matrix identity is
\begin{equation}\label{i4:full-K}
 \frac{K_{\mathrm{rig}}}{c_1}=2(\kappa^{-1}-1)\nabla^2\ell_{12}
 +(1-(2\kappa)^{-1})(S\otimes r+r\otimes S)
 -2\kappa^{-1}r\otimes \mathbf{k}_{\kappa}+2\Xi I_3.
\end{equation}
Consequently, with $\operatorname{tf}B=B-(\operatorname{tr}B)I_3/3$,
\begin{equation}\label{i4:tf-Hessian}
 \operatorname{tf}\nabla^2\ell_{12}
 =\operatorname{tf}(\beta\otimes r+r\otimes\beta),\qquad
 \beta=\frac{2\mathbf{k}_{\kappa}-(2\kappa-1)(\mathbf{m}_1+\mathbf{m}_2)}{4(1-\kappa)}.
\end{equation}
\end{proposition}
\begin{proof}
Since $\nabla\log(c_1/\kappa)=-\mathbf{m}_1/2-\mathbf{k}_{\kappa}$ and
$e^{-\ell_{12}/2}c_2=c_1$, direct expansion gives
\begin{align*}
 2\nabla\mathbf{g}&=\frac{2c_1}{\kappa}\nabla^2\ell_{12}
                -\frac{c_1}{\kappa}r\otimes \mathbf{m}_1
                -\frac{2c_1}{\kappa}r\otimes \mathbf{k}_{\kappa},\\
 -F_{\mu_1}+e^{-\ell_{12}/2}F_{\mu_2}
   &=2c_1\{-\nabla^2\ell_{12}+\mathbf{m}_1\otimes \mathbf{m}_1-\mathbf{m}_2\otimes \mathbf{m}_2+\Xi I_3\}.
\end{align*}
Because $r=\mathbf{m}_1-\mathbf{m}_2$, both
$\mathbf{m}_1\otimes \mathbf{m}_1-\mathbf{m}_2\otimes \mathbf{m}_2$ and
$r\otimes \mathbf{m}_1+\mathbf{m}_2\otimes r$ equal
$(S\otimes r+r\otimes S)/2$. Substitution in \eqref{i4:K} proves
\eqref{i4:full-K}, including its scalar and skew terms. Its skew part is
$-(c_1/\kappa)(r\otimes \mathbf{k}_{\kappa}-\mathbf{k}_{\kappa}\otimes r)$ and is zero because $K_{\mathrm{rig}}=qI_3$.
The symmetric trace-free part is
\[
 0=2(\kappa^{-1}-1)\operatorname{tf}\nabla^2\ell_{12}
 +(1-(2\kappa)^{-1})\operatorname{tf}(S\otimes r+r\otimes S)
 -\kappa^{-1}\operatorname{tf}(r\otimes \mathbf{k}_{\kappa}+\mathbf{k}_{\kappa}\otimes r).
\]
Division by $2(1-\kappa)/\kappa>0$ gives \eqref{i4:tf-Hessian}.
\end{proof}

\begin{lemma}\label{i4:closure}
Suppose $\ell_{12}\in C^3(\mathcal{O};\mathbb R)$, $\beta\in C^1(\mathcal{O};\mathbb R^3)$
and \eqref{i4:tf-Hessian} holds. Put $r=\nabla\ell_{12}$ and $w=\Delta\ell_{12}$. The column $U=(r_1,r_2,r_3,w)^T$ collects the four scalar unknowns. In this lemma $\partial_j=\partial/\partial x_j$, $\delta_{ij}$ is the Kronecker delta, and a repeated index is summed from $1$ to $3$. Define
\[
 T^k_{ij}=\beta_i\delta_{jk}+\beta_j\delta_{ik}
                    -\frac23\delta_{ij}\beta_k.
\]
The trace-free Hessian identity determines the derivatives of $r$ up to the scalar $w$. Its differentiated equation determines the derivatives of $w$, giving $U\in C^1$ and $\partial_jU=\mathcal C_jU$ with continuous $4\times4$ coefficient matrices
\begin{equation}\label{i4:closed-system}
 \mathcal C_j=\begin{pmatrix}(T^k_{ij})_{i,k=1}^3&(\delta_{ij}/3)_{i=1}^3\\
                         (B_j^k)_{k=1}^3&5\beta_j/3\end{pmatrix},
 \qquad
 B_j^k=\frac32\left[(\operatorname{div}\beta)\delta_{jk}
 +\partial_k\beta_j-\frac23\partial_j\beta_k
 +\frac13|\beta|^2\delta_{jk}+\frac19\beta_j\beta_k\right].
\end{equation}
\end{lemma}
\begin{proof}
Restoring the trace in \eqref{i4:tf-Hessian} gives
$\partial_jr_i=\delta_{ij}w/3+T^k_{ij}r_k$. Differentiate in $x_i$
and sum. Commutation of third derivatives gives
\[
 \frac23\partial_jw=(\partial_iT^k_{ij})r_k+
 T^k_{ij}\left(\frac13\delta_{ki}w+T^\ell_{ki}r_\ell\right).
\]
The contractions are
\begin{align*}
 \sum_iT^i_{ij}&=\frac{10}3\beta_j,&
 \sum_i\partial_iT^k_{ij}
   &=(\operatorname{div}\beta)\delta_{jk}
                      +\partial_k\beta_j-\frac23\partial_j\beta_k,\\
 \sum_{i,\ell}T^\ell_{ij}T^k_{\ell i}
   &=\frac13\sum_i\beta_iT^k_{ji}
     =\frac13|\beta|^2\delta_{jk}+\frac19\beta_j\beta_k.
\end{align*}
For the last equality, insert the three terms in $T^\ell_{ij}$;
the middle term vanishes by $\sum_iT^k_{ii}=0$, and the first and last
combine with coefficients $1-2/3=1/3$. Multiplication by $3/2$ yields
\eqref{i4:closed-system}. The coefficient matrices involve only
$\beta$ and its first derivatives, and the system is homogeneous.
\end{proof}

\begin{lemma}\label{i4:propagation}
Let $\mathcal O\subset\mathbb R^3$ be connected and open. Suppose
$U\in C^1(\mathcal O;\mathbb R^4)$ satisfies
\[
 \partial_jU=\mathcal C_jU,\qquad j=1,2,3,
\]
where $\mathcal C_j\in C^0(\mathcal O;\mathbb R^{4\times4})$.
If $U(x_0)=0$ at some $x_0\in\mathcal O$, then $U=0$ on
$\mathcal O$. For the state
$U=(r_1,r_2,r_3,\Delta\ell_{12})^T$ in Lemma~\ref{i4:closure},
the initial condition requires both $r(x_0)=0$ and
$\Delta\ell_{12}(x_0)=0$.
\end{lemma}
\begin{proof}
A connected open subset of Euclidean space is polygonally connected.
Fix $x\in\mathcal O$ and a polygonal path
$\mathfrak p:[0,1]\to\mathcal O$ from $x_0$ to $x$.
On each segment, $V=U\circ\mathfrak p$ satisfies
\[
 V'(t)=\left(\sum_{j=1}^3\dot{\mathfrak p}_j(t)
             \mathcal C_j(\mathfrak p(t))\right)V(t).
\]
The coefficient matrix is bounded on each compact segment.
Gronwall's inequality, starting from $V(0)=0$ and then proceeding
through the segment endpoints, gives $V=0$ along the path.
Thus $U(x)=0$, and $x$ was arbitrary.
\end{proof}

\begin{proof}[Proof of Theorem~\ref{i4:main}]
Lemma~\ref{i4:determined} gives the common $\kappa$ and expresses
all remaining blocks in terms of $\ell_{12}$.
Lemma~\ref{i4:null-lines} and Proposition~\ref{i4:beta} yield
\eqref{i4:tf-Hessian}. By Lemma~\ref{i4:closure}, the full state
\[
 U=(r_1,r_2,r_3,\Delta\ell_{12})^T,\qquad r=\nabla\ell_{12},
\]
belongs to $C^1$ and satisfies $\partial_jU=\mathcal C_jU$ with
the matrices in~\eqref{i4:closed-system}. For this fixed pair of
materials, $\beta$ and $\mathcal C_j$ are fixed coefficient fields
that depend continuously on $x$.

Choose $x_0\in E_{\mathrm{eq}}$. Since $\ell_{12}=0$ on the
open set $E_{\mathrm{eq}}$, its first and second derivatives
vanish there. In particular,
\[
 r(x_0)=0,\qquad \Delta\ell_{12}(x_0)=0,
 \qquad \ell_{12}(x_0)=0.
\]
Lemma~\ref{i4:propagation} gives $U=0$ on the connected set
$\mathcal O$. Thus $\nabla\ell_{12}=0$, and connectedness and
$\ell_{12}(x_0)=0$ give $\ell_{12}=0$ throughout $\mathcal O$.
Lemma~\ref{i4:determined} now yields $\mu_1=\mu_2$, the common
$\kappa$, and $H=I_3$; the inverse formula for the Lam\'e parameters
gives $\lambda_1=\lambda_2$.
\end{proof}

Propositions~\ref{if:actualH} and~\ref{if:binding} provide the relative
matrix on the common ball with the hypotheses of Theorem~\ref{i4:main}.
In primal quotient coordinates the corresponding matrix $G$ satisfies
$2D_\theta G+\mathcal A_1G-G\mathcal A_2=0$. Since
$\Gamma_j=-\mathcal A_j^T$, differentiation gives
$2D_\theta(G^{-T})=-\Gamma_1G^{-T}+G^{-T}\Gamma_2$, so $H=G^{-T}$
has the material order used throughout this section.

\section{Density recovery and completion of the proof}\label{sec:potential}

This section completes the uniqueness argument after the Lam\'e
parameters have been determined in Section~\ref{i4:section}.
With the common Lam\'e operator fixed, the physical Green--Betti
identity~\eqref{pt:betti} isolates the difference of the two
zeroth-order coefficients $W_j=\omega^2\rho_j$ in a bilinear
integral of products of elastic displacements. The remaining task
is to construct exact physical solutions whose products have a
controlled high-frequency limit and detect this difference.

For each fixed nonzero Fourier vector $\xi$, choose null phases
whose sum is $i\xi$. The local frames, first corrections and exact
solutions provided at finite regularity in
Subsection~\ref{fr:interior-assembly} give an $L^1$ limit for the
products of the physical displacements. This yields a Fourier identity
with a localized transverse weight. A finite family of local tests
whose squares sum to one on the projection of $\overline\Omega$
removes the weight and determines the Fourier transform of
$(W_1-W_2)/\mu$. Fourier uniqueness gives $W_1=W_2$, and the fixed
known frequency $\omega>0$ then gives $\rho_1=\rho_2$.

\subsection{The physical integral identity and the potential reduction}
Since the Lam\'e coefficients have been recovered, let $L=L_1=L_2$ and
$t=t_1=t_2$ be the common elastic operator and traction from~\eqref{bt:stress}.
Let $\nu$ be the outward unit normal. In the bilinear setting of
Section~\ref{sec:notation}, the trace $\gamma_0$ is understood as the bounded map
$H^1(\Omega;\mathbb C^3)\to H^{1/2}(\partial\Omega;\mathbb C^3)$.
In this section $(L+W_j)u_j=0$
abbreviates $(L+W_jI_3)u_j=0$. For $f_j=\gamma_0u_j$, Green's formula
and symmetry of the stress give
\begin{equation}\label{pt:betti}
 \int_\Omega(W_1-W_2)u_1\cdot u_2\,dx
 =\langle(\Lambda_2^\Omega-\Lambda_1^\Omega)f_1,f_2\rangle=0.
\end{equation}
We use the bilinear pairing convention fixed in Section~\ref{sec:notation}. Here $\Lambda_j^\Omega=\Lambda_{\lambda_j,\mu_j,\rho_j}^\Omega(\omega)$ denotes the physical Dirichlet-to-Neumann map on $\Omega$. To obtain the displayed sign, apply Green's formula first to two solutions of the same
system, which gives the bilinear symmetry of $\Lambda_j^\Omega$, and then to
$u_1,u_2$.

Use the common reduction coefficients and maps $a,c,\mathbf m,\kappa,T,S$
from~\eqref{bt:notation}--\eqref{bt:TS}. For each $j=1,2$, write
$M^{W_j}=M_j^{W_j}$ for the reduced operator with these common Lam\'e
coefficients and potential $W_j$. Lemma~\ref{bt:scalarlemma} gives
$(L+W_jI_3)T=SM^{W_j}$. The potential matrix $\mathcal W$ in
\eqref{bt:MW} acts by multiplication, including the $\nabla W_j$
term in its upper-right block. It contributes at order zero
in $M^{W_j}$, so the leading transport operator is common to
the two reduced systems.

\subsection{Null phases and the leading displacement amplitudes}
Fix $\xi\ne0$, put $k_\xi=|\xi|$ and $\widehat\xi=\xi/k_\xi$, and choose a positively
oriented real orthonormal basis $(e_1,e_2,\widehat\xi)$. Let $\theta=e_1+ie_2$ and
set
\[
 d_h=\frac{h^2k_\xi^2}{16},\quad \alpha=1+d_h,\quad \epsilon=\frac h\alpha,
 \qquad
 \zeta_{1,2}=\pm\frac\theta h+\frac{i\xi}{2}
                  \pm\frac{hk_\xi^2}{16}\bar\theta,
 \qquad \vartheta_j=-\epsilon\zeta_j.
\]
Here $\widehat\xi$ is the unit Fourier direction; $k_\xi$ is distinct
from the boundary-layer oscillation parameter $k$, and $d_h$ is a scalar
phase correction distinct from the reduced coefficient $d$.
The vectors $\zeta_j$ are the exponential phases, with $\zeta_j^2=\zeta_j\cdot\zeta_j=0$ and $\zeta_1+\zeta_2=i\xi$. For fixed $\xi$, the auxiliary parameter $\epsilon=h/\alpha$ satisfies $\epsilon\sim h$ as $h\downarrow0$. The vectors $\vartheta_j=-\epsilon\zeta_j$ are normalized transport directions, so $e^{\zeta_j\cdot x}=e^{-\vartheta_j\cdot x/\epsilon}$ has the form used in Proposition~\ref{pb27:remainder}. Define the leading polarization vectors
\[
 \mathbf q_1=\widehat\xi-\frac{ihk_\xi}{4}\bar\theta,\qquad
 \mathbf q_2=\widehat\xi+\frac{ihk_\xi}{4}\bar\theta.
\]
Their bilinear products are
\begin{equation}\label{pt:phase}
 \vartheta_j\cdot \mathbf q_j=0,\quad \mathbf q_1\cdot\mathbf q_2=1,\quad
 \mathbf q_1\cdot\vartheta_2=\vartheta_1\cdot\mathbf q_2=-i\epsilon k_\xi,
 \quad \vartheta_1\cdot\vartheta_2=-\epsilon^2k_\xi^2/2.
\end{equation}
The real and imaginary parts of $\vartheta_j$ are orthogonal unit vectors.
The corresponding leaf normals are
\[
 \nu_{1,2}=\frac{(1-d_h)\widehat\xi\mp(hk_\xi/2)e_2}{1+d_h}.
\]
Thus $\nu_j\to \widehat\xi$ and $|\nu_j-\widehat\xi|=O_\xi(\epsilon)$. For fixed $\xi$,
each direction path lies in a fixed compact chart for sufficiently small
$h$.

The limiting directions are $\vartheta_1\to-\theta$ and
$\vartheta_2\to\theta$. At each prescribed transverse centre, choose
local frame families for these two directions separately and intersect
their transverse parameter intervals. Shrinking this common interval
$I$ and the two direction neighborhoods leaves a fixed compact
parameter neighborhood for each path when $h$ is sufficiently small.
Choose a fixed real $b\in C_c^\infty(I)$. Since $\vartheta_j\cdot \mathbf q_j=0$ and $\vartheta_j\cdot\nu_j=0$, the leading amplitude $a_{0,j}$ below solves $\mathsf T_{\vartheta_j}a_{0,j}=0$, where $\mathsf T_\vartheta=2D_\vartheta+A_\vartheta^{\rm phys}$ uses the common Lam\'e coefficients.
For the fixed materials, $\xi$, and $b$, the leading amplitudes
have transverse support in the fixed compact set
$\operatorname{supp}b\Subset I$. Choose the auxiliary $\chi_\perp$
of Lemma~\ref{a31:first-correction} to equal one near this set.
Then $\widetilde a_{0,j}=a_{0,j}$. With the finite-order
frames of Subsection~\ref{fr:interior-assembly}, the formula
\eqref{a31:first-correction-formula} gives first corrections that
extend by zero in the transverse variable, with the transport
equations holding on $U$. The positive support margins and the
finite-order frame bounds give uniform $H^2(U)$ bounds along the
two compact parameter paths. Applying
\eqref{pb27:two-sign-correction} with the finite-order right inverse
from the same subsection gives exact solutions on
$\Omega_{\rm CGO}$:
\[
 Z_j=e^{\zeta_j\cdot x}(a_{0,j}+\epsilon a_{1,j}+r_j),\qquad
 a_{0,j}=(b(\nu_j\cdot x)\mathbf q_j,0)^T,
 \quad M^{W_j}Z_j=0,
\]
with $a_{1,j}$ uniformly bounded in $H^2$ and
\begin{equation}\label{pt:remainder}
 \|r_j\|_2=O(\epsilon^2),\qquad \|r_j\|_{H^1}=O(\epsilon).
\end{equation}
\begin{lemma}\label{pt:product-lemma}
Fix $\xi\in\mathbb R^3\setminus\{0\}$ and a real test
$b\in C_c^\infty(I)$ on an interval common to the two limiting
direction charts. Let $Z_j$ be the exact reduced solutions constructed
above, with phases $\zeta_j$, leading vectors $\mathbf q_j$, uniformly bounded
first corrections $a_{1,j}$, and remainders satisfying
\eqref{pt:remainder}. The physical elastic displacements
$u_j=TZ_j$ satisfy
\begin{equation}\label{pt:limit}
 e^{-i\xi\cdot x}u_1\cdot u_2\longrightarrow
 \mu^{-1}b(\widehat\xi\cdot x)^2\quad\hbox{in }L^1(\Omega).
\end{equation}
For each fixed $\xi\ne0$ and $b$, the convergence holds along
the two compact direction paths selected above.
\end{lemma}
\begin{proof}
Write $a_{1,j}=(v_{1,j},f_{1,j})$ and $r_j=(r_{v,j},r_{f,j})$,
separating their vector and scalar components. Set $u_j=TZ_j$. The factorization $(L+W_j)T=SM^{W_j}$ proved above gives $(L+W_j)u_j=0$, so these displacements may be used in~\eqref{pt:betti}. Expanding the lift gives
\[
 e^{-\zeta_j\cdot x}u_j=U_j+E_{j,\mathrm{err}},\qquad
 U_j=c b(\nu_j\cdot x)\mathbf q_j-a\vartheta_jf_{1,j}.
\]
The longitudinal term in $U_j$ is of order one. The remaining terms
satisfy
\[
 \|E_{j,\mathrm{err}}\|_2\le C\left[\epsilon\|a_{1,j}\|_{H^1}
              +\|r_j\|_{H^1}+\epsilon^{-1}\|r_j\|_2\right]
 =O(\epsilon).
\]
Using \eqref{pt:phase}, we obtain
\begin{align*}
 U_1\cdot U_2={}&\mu^{-1}b_1b_2
 +i\epsilon k_\xi c a(b_1f_{1,2}+b_2f_{1,1})
 -\tfrac12\epsilon^2k_\xi^2a^2f_{1,1}f_{1,2},
 \qquad b_j=b(\nu_j\cdot x).
\end{align*}
Hence $\|U_1\cdot U_2-\mu^{-1}b_1b_2\|_{L^1}=O(\epsilon)$.
The $H^2$ bounds for $a_{1,j}$ give $\|U_j\|_2=O(1)$. Therefore Cauchy--Schwarz and $\|E_{j,\mathrm{err}}\|_2=O(\epsilon)$ bound $E_{1,\mathrm{err}}\cdot U_2$, $U_1\cdot E_{2,\mathrm{err}}$, and $E_{1,\mathrm{err}}\cdot E_{2,\mathrm{err}}$ in $L^1$ by $O(\epsilon)$. Since $\nu_j\to \widehat\xi$, smoothness of the fixed compactly supported $b$ gives $b_1b_2\to b(\widehat\xi\cdot x)^2$ uniformly on the bounded working region. Finally $\zeta_1+\zeta_2=i\xi$, so the product of the exponential
phases is exactly $e^{i\xi\cdot x}$. These facts prove
\eqref{pt:limit}. Each potential has its own first correction; the
estimates apply to this pair of exact elastic displacement fields.
\end{proof}

\subsection{Removal of the local weights and Fourier uniqueness}
\begin{proposition}\label{pt:recovery}
Suppose the Lam\'e parameters are common, continuous and strongly
convex, and $W_j\in C(\overline\Omega;\mathbb R)$ for $j=1,2$.
Assume that both Dirichlet problems on $\Omega$ are nonresonant
and that the physical DN maps agree. Suppose that for each fixed
$\xi\ne0$ and each transverse centre there is an open interval $I$
containing that centre such that for every fixed real-valued
$b\in C_c^\infty(I)$, exact solutions of the two physical systems
satisfy the product limit~\eqref{pt:limit}. Then $W_1=W_2$ on
$\overline\Omega$.
\end{proposition}
\begin{proof}
Since $W_1-W_2$ is bounded, the $L^1(\Omega)$ convergence
in~\eqref{pt:limit} permits passage to the limit
in~\eqref{pt:betti}. This gives
\[
 \int_\Omega e^{i\xi\cdot x}\frac{W_1-W_2}{\mu}
                      b(\widehat\xi\cdot x)^2\,dx=0.
\]
The compact set $J=\widehat\xi\cdot\overline\Omega$ has a finite cover by local
intervals $I_1,\ldots,I_N$. Choose
$\chi_\ell\in C_c^\infty(I_\ell;\mathbb R)$ so that
$S=\sum\chi_\ell^2>0$ near $J$, and take
$\eta\in C_c^\infty(\{S>0\})$ equal to one near $J$. Put
$b_\ell=\eta\chi_\ell/\sqrt S$ on $\{S>0\}$ and extend by zero.
The denominator is bounded below on $\operatorname{supp}\eta$, so these
are smooth local tests and $\sum b_\ell^2=1$ on $J$. Applying the local
identity to each $b_\ell$ and adding the resulting integrals gives
\[
 \int_\Omega e^{i\xi\cdot x}(W_1-W_2)/\mu\,dx=0
 \quad(\xi\ne0).
\]
Let $\mathbf1_\Omega$ denote the indicator of $\Omega$ and set $F=\mathbf1_\Omega(W_1-W_2)/\mu\in L_c^1(\mathbb R^3)\cap L^2(\mathbb R^3)$. In the bilinear convention of Section~\ref{sec:notation}, use the Fourier transform
\[
 \widehat F(\xi)=\int_{\mathbb R^3}e^{-i\xi\cdot x}F(x)\,dx.
\]
Replacing $\xi$ by $-\xi$ in the preceding identity gives $\widehat F(\xi)=0$ for $\xi\ne0$. Compact support and integrability make this transform continuous, so it vanishes at $\xi=0$ as well. For
$G_t(x)=(4\pi t)^{-3/2}e^{-|x|^2/(4t)}$, absolute integrability and
Fubini's theorem give
\[
 G_t*F(x)=(2\pi)^{-3}\int e^{i\xi\cdot x-t|\xi|^2}
                                  \widehat F(\xi)\,d\xi=0.
\]
Since the Gaussian kernels form an approximate identity,
$G_t*F\to F$ in $L^1$, and hence $F=0$ almost everywhere. Continuity of
the physical coefficients gives $W_1=W_2$ on $\overline\Omega$.
\end{proof}

Lemma~\ref{pt:product-lemma}, using the fixed-frame construction,
first corrections, and the same two-norm remainder estimate, supplies
the local product limits required in Proposition~\ref{pt:recovery}. Since the
common measurement frequency is known and nonzero,
$\rho_1-\rho_2=(W_1-W_2)/\omega^2=0$.

\subsection{Interior argument at finite regularity}
\label{fr:interior-assembly}

The interior step takes the boundary matching conditions
\eqref{bt:matching} and the existence of exterior coefficient extensions
as hypotheses. In the proof of Theorem~\ref{a28:main},
Proposition~\ref{bl:finite-jets} recovers the boundary traces from the
data, and Proposition~\ref{ext:finite-common} constructs the extensions.

\begin{proposition}
\label{fr:interior}
Let $\Omega\subset\mathbb R^3$ be a bounded $C^2$ domain. For $j=1,2$,
let $(\lambda_j,\mu_j,W_j)$ be real coefficients on $\overline\Omega$
that admit extensions
\[
 (\widetilde\lambda_j,\widetilde\mu_j)
     \in C^{16}(\mathbb R^3;\mathbb R)^2,
 \qquad \widetilde W_j\in C_c^{15}(\mathbb R^3;\mathbb R).
\]
Assume that both extended Lam\'e pairs are strongly convex on
$\mathbb R^3$, agree on $\mathbb R^3\setminus\overline\Omega$, and equal
the same strongly convex constant pair $(\lambda_*,\mu_*)$ outside
a compact set. Suppose that the physical Dirichlet realizations are
invertible on $\Omega$, their Dirichlet-to-Neumann maps agree, and the
boundary traces satisfy~\eqref{bt:matching}. Then
\[
 \lambda_1=\lambda_2,\qquad \mu_1=\mu_2,\qquad W_1=W_2
 \quad\text{on }\overline\Omega .
\]
\end{proposition}

\begin{proof}
\emph{Coefficient bounds and right inverse.}
We use the notation $(\lambda_j,\mu_j,W_j)$ also for the extensions
when evaluating coefficients outside $\overline\Omega$; the boundary
maps always refer to $\Omega$. In \eqref{pb27:MW} the first-order
matrices use at most two derivatives of $\mu$ and depend only on the
Lam\'e pair. The zeroth-order matrix uses at most two derivatives of
$\mu$ and one of $W$. Thus $A_k^{\rm phys},V^W\in C^{14}$, with
$(A_k^{\rm phys},V^W)=(A_{k,*},0)$ outside a fixed ball. The product
rule in Lemma~\ref{pb27:background} gives
\begin{equation}\label{fr:interior-coefficients}
 \widetilde A_k,\widetilde V^W\in C_c^{14}(\mathbb R^3),
 \qquad k=1,2,3 .
\end{equation}
Here and below material indices are suppressed in formulas for a single
medium.

Take $\overline\Omega\subset B_{R_0}$, $R_0<R_\perp<R<L$ as in the
construction of Proposition~\ref{if:actualH}. Choose the balls and cutoff
$\chi_0$ in \eqref{r27:geometry}, with a compact set $K_{\rm coef}$
containing the supports of $\widetilde A_k$ and $\widetilde V^W$
for both media and with $K_{\rm coef}\Subset\operatorname{int}K_{\rm src}$.
Take the auxiliary balls large enough that an open ball $U$ centered
at the origin satisfies
\begin{equation}\label{fr:interior-geometry}
 \operatorname{supp}\chi_0\cup\overline\Omega_{\rm CGO}
 \cup \overline B_{\sqrt{L^2+R_\perp^2}}
 \Subset U\Subset\Omega_1 .
\end{equation}
Indeed, for every normalized null direction $\theta$ the real leaf
coordinates obey $|x|^2=|z|^2+x_\perp^2$. Thus \eqref{fr:interior-geometry}
includes the physical cylinder $|z|\le L$, $|x_\perp|\le R_\perp$ for
every such $\theta$. Fix radii $L_{\rm src}<L_{\rm fr}$ with
$\pi_\theta(\overline U)\Subset D_{L_{\rm src}}\Subset D_{L_{\rm fr}}$
for every $\theta$, as in Lemma~\ref{a31:first-correction}.

In normalized coordinates the coefficients in
\eqref{i2:eq-fr-system}, \eqref{i2:eq-frames}, and
\eqref{pb27:two-sign-transport} are jointly $C^{14}$: the
two-component system uses $D_\theta^2a$, $a=\mu^{-1}$, and the quotient
systems use $A_\theta^{\rm phys}$ and the spinor matrices. Since
$D_\theta=2\partial_{\bar z}$, a four-component frame $C_\theta$
satisfying $\mathsf T_\theta C_\theta=0$ solves
\begin{equation}\label{fr:interior-frame}
 \partial_{\bar z}C_\theta=-\tfrac14 A_\theta^{\rm phys}C_\theta .
\end{equation}
The pulled-back coefficients of \eqref{fr:interior-frame} and the
two smaller systems are bounded in $C^{14}$ on the larger charts and
have compact closure in $C^{13,\alpha}$ on smaller charts
for $0<\alpha<1$. Apply Theorem~3.1 of~\cite{DDStab2026} with $r=12$
to the finite product of coefficient families. The common
refinement supplies invertible frames and inverses jointly $C^{12}$
on the full leaf disks. Choose a finite cover of the compact space
of ordered orthonormal pairs and transverse parameters ranging over
$\pi_\perp(\overline U)$ by relatively compact subcharts. On each
subchart the source interval $I_\perp$ and the intervals carrying
the transverse partition of unity are chosen as in
Lemma~\ref{a31:first-correction}. For the moment tests at each center
$|x_\perp|<R_\perp$, use a local interval
$I\Subset I_\perp\cap(-R_\perp,R_\perp)$ as in Lemma~\ref{i2:lem-fixed}.
The charts, radii, and cutoffs are fixed independently of $h$ and the
tests.

By \eqref{fr:interior-coefficients}, Theorem~3.2 of~\cite{DDStab2026}
at Sobolev order zero and system dimension four applies to
$\widetilde M^W$ in \eqref{pb27:compact} for the normalized null
phases. Write $\widetilde{\mathcal R}^{(0)}_{h,\theta}$ for its right
inverse of $e^{\theta\cdot x/h}\widetilde M^We^{-\theta\cdot x/h}$.
For $\supp F\subset K_{\rm src}$, put
\[
 \mathcal R^{(0)}_{h,\theta}F
 =\mathcal B_{\rm bg}\widetilde{\mathcal R}^{(0)}_{h,\theta}
       (\mathcal B_{\rm bg}^{-1}F).
\]
The matrices $\mathcal B_{\rm bg}^{\pm1}$ preserve source support,
act boundedly on $L^2$ and $H^1$ on the fixed balls, and commute with
the scalar exponential. Therefore
\begin{equation}\label{fr:interior-inverse}
\begin{aligned}
 e^{\theta\cdot x/h}M^We^{-\theta\cdot x/h}
       \mathcal R^{(0)}_{h,\theta}F&=F\quad\text{in }\Omega_1,\\
 h^{-1}\|\mathcal R^{(0)}_{h,\theta}F\|_{L^2(\Omega_0)}
       +\|\mathcal R^{(0)}_{h,\theta}F\|_{H^1(\Omega_0)}
 &\le C\|F\|_{L^2(\Omega_0)} .
\end{aligned}
\end{equation}
The inverse is fixed for all sources at the given direction, and
the estimates are uniform on the compact subcharts, including those
for $-\theta$.

\emph{Amplitudes and exact solutions.}
Write $\partial^\alpha$ for joint derivatives in the real leaf
coordinates $(\Re z,\Im z)$, the transverse coordinate $x_\perp$,
and the local direction coordinates. For a source $F$ compactly
supported in the leaf variable, the kernel of $\mathcal C_\theta$
in the pulled-back coordinates is independent of the parameters.
Integration by parts in the leaf variable and differentiation in the
parameters give
\begin{equation}\label{fr:interior-cauchy}
 \partial^\alpha\mathcal C_\theta F
 =\mathcal C_\theta(\partial^\alpha F),\qquad
 \|\partial^\alpha\mathcal C_\theta F\|_{L^\infty}
 \le C\|\partial^\alpha F\|_{L^\infty}.
\end{equation}
The suprema are taken on the fixed target and source coordinate sets,
respectively. The bound follows from local integrability of $|z-w|^{-1}$
and the bounded source disk.

In Lemma~\ref{i2:lem-fixed}, the term
$\mathsf T_{1,\theta}(V_\theta X_1)$ uses one derivative of the
$C^{12}$ frame $X_1$. In Lemma~\ref{i2:lem-f}, each $H_j$ uses one
derivative of the $C^{12}$ two-component frame $\Phi$. The
amplitude formulas in these lemmas and \eqref{fr:interior-cauchy} give
\begin{equation}\label{fr:interior-amplitudes}
 Z_1,\ A_{0,1}^F,\ A_{0,2}^F\in C^{11}.
\end{equation}
Here $j$ indexes the columns of $\Phi$ in the formula for $H_j$.
The remaining leading amplitude $A_0^S=(q_\theta,0)^T$ is independent
of $x$ and smooth in the spinor chart. For fixed tests, multiplication
by the cutoffs and polynomials preserves the orders in
\eqref{fr:interior-amplitudes}. The sources $M^Wa_0$ therefore have
at least nine joint derivatives. Using \eqref{fr:interior-cauchy}
in \eqref{a31:first-correction-formula}, we obtain
\begin{equation}\label{fr:interior-correction}
 \|a_1\|_{H^2(U)}\le C,\qquad
 \supp(\chi_0M^Wa_1)\subset K_{\rm src},\qquad
 \|\chi_0M^Wa_1\|_{L^2(\Omega_0)}\le C .
\end{equation}
The same bounds hold for the third-column correction of
Proposition~\ref{i2:prop-s}, including the correction generated by the
commutator source in \eqref{i2:eq-commutator}.
The transverse cutoffs are supported strictly inside their
intervals, and the partition functions are constant on each leaf.
Transverse zero extension and summation thus preserve the transport
equations on $U$.

Use $Z_{h,\mathfrak s}$ and $r_{h,\mathfrak s}$ from
\eqref{pb27:two-sign-correction}, with $\mathfrak s=\pm1$. The
identities \eqref{pb27:two-sign-conjugation}--\eqref{pb27:two-sign-residual}
and the bounds \eqref{fr:interior-inverse} and
\eqref{fr:interior-correction} give
\begin{equation}\label{fr:interior-two-phase}
\begin{aligned}
 e^{\mathfrak s\theta\cdot x/h}M^W Z_{h,\mathfrak s}
   &=h(1-\chi_0)M^Wa_1=0\quad\text{in }\Omega_{\rm CGO},\\
 h^{-1}\|r_{h,\mathfrak s}\|_{L^2(\Omega_0)}
     +\|r_{h,\mathfrak s}\|_{H^1(\Omega_0)}
   &\le Ch\|\chi_0M^Wa_1\|_{L^2(\Omega_0)}\le Ch .
\end{aligned}
\end{equation}
For each fixed $h$, \eqref{fr:interior-two-phase} gives
\begin{equation}\label{fr:interior-bootstrap}
 \Delta Z_{h,\mathfrak s}
 =-\sum_k A_k^{\rm phys}\partial_k Z_{h,\mathfrak s}
    -V^WZ_{h,\mathfrak s}\quad\text{in }\Omega_{\rm CGO}.
\end{equation}
If $Z_{h,\mathfrak s}\in H^s_{\rm loc}(\Omega_{\rm CGO})$ for
$1\le s\le15$, the right-hand side of \eqref{fr:interior-bootstrap}
lies in $H^{s-1}_{\rm loc}(\Omega_{\rm CGO})$ by the $C^{14}$
coefficient bounds. Interior elliptic regularity therefore gives
$Z_{h,\mathfrak s}\in H^{s+1}_{\rm loc}(\Omega_{\rm CGO})$. Starting
with $s=1$, we obtain $Z_{h,\mathfrak s}\in
H^{16}_{\rm loc}(\Omega_{\rm CGO})$. Since
$\overline\Omega\Subset\Omega_{\rm CGO}$, these fields have the
classical derivatives needed in \eqref{bt:factorization} and
Propositions~\ref{bt:totalprop} and~\ref{bt:weakprop}. The $C^2$
boundary permits the tangential differentiation and Hessian identity
of Lemma~\ref{bt:hessian}. The classical tractions agree with the
weak Dirichlet-to-Neumann pairings, so the matching conditions
\eqref{bt:matching} give \eqref{i2:eq-betti}.

For the third column, Proposition~\ref{i2:prop-s} and
Lemma~\ref{i2:lem-fixed} give
\[
 \det Y_2=2c_2\ne0,\qquad
 A_{1,W}^S-A_{1,0}^S=(\beta q_\theta,0)^T\in\ker J_\theta
 \quad(\theta\cdot q_\theta=0).
\]
The commutator source in \eqref{i2:eq-commutator} vanishes on the
leaves near $\supp b$, so its Cauchy correction $w_\chi$ vanishes
there as well. Throughout $\Omega_{\rm CGO}$,
\eqref{i2:eq-k-exact} reduces the physical lift on the exact solution
to $\widehat K_{\mathrm{lift}}^W$. With $r_h=r_{h,+1}$ for this
column, the remainder terms in Lemma~\ref{i2:lem-shell} satisfy
\[
 \|\widehat K_{\mathrm{lift}}^W r_h\|_{L^2(\Omega_{\rm CGO})}
   +h^{-1}\|N_\theta r_h\|_{L^2(\Omega_{\rm CGO})}
 \le C\bigl(\|r_h\|_{H^1(\Omega_0)}
                 +h^{-1}\|r_h\|_{L^2(\Omega_0)}\bigr)=O(h)
\]
by \eqref{fr:interior-two-phase}.

\emph{Matrix moments.}
Fix $b\in C_c^\infty(I)$ and $g\in\mathbb C[z]$ on a moment chart
chosen above. Use the positive-phase matrix $z_{1,h}$ and the
negative-phase test matrix $\mathcal Y_{2,h}$ of
Theorem~\ref{i2:thm-nine}. Applying \eqref{i2:eq-betti} on $\Omega$
and using the expansions in that theorem, with
$A=bgZ_1$, gives
\begin{equation}\label{fr:interior-moment-estimates}
\begin{aligned}
 0&=\int_\Omega h\mathcal Y_{2,h}^{T}D^Wz_{1,h}\,dx,\\
 h e^{-\theta\cdot x/h}D^Wz_{1,h}
    &=\Delta A_\theta A+O_{L^2(\Omega)}(h),\\
 e^{\theta\cdot x/h}\mathcal Y_{2,h}
    &=\mathcal Y_\chi^{(0)}+O_{L^2(\Omega)}(h),\\
 \mathcal Y_\chi^{(0)}&=U_\theta Y_2\quad\text{near }\supp b.
\end{aligned}
\end{equation}
The product in \eqref{fr:interior-moment-estimates} converges in
$L^1(\Omega)$ by Cauchy--Schwarz. The zeroth-order difference
$\Delta B^W$ may be nonzero outside $\overline\Omega$, so we take
this limit on $\Omega$. The first-order difference
$\Delta A_\theta$ is supported in $\overline\Omega$. Extending only
the resulting integral, using
$U_\theta^T\Delta A_\theta Z_1=2(\mathsf P_1-\mathsf P_2)X_1$
and applying Fubini gives
\begin{equation}\label{fr:interior-tested-moments}
\begin{aligned}
 0&=\int_\Omega bgY_2^TU_\theta^T\Delta A_\theta Z_1\,dx\\
  &=2\int_I b(x_\perp)\int_{D_L}
       g(z)Y_2^T(\mathsf P_1-\mathsf P_2)X_1\,dA\,dx_\perp .
\end{aligned}
\end{equation}
The finite-order frames and compact support give
$Y_2^T(\mathsf P_1-\mathsf P_2)X_1\in C(I;L^1(D_L))$.
Applying Lemma~\ref{i2:leaf-extraction} to
\eqref{fr:interior-tested-moments} gives the moments for every leaf
and every polynomial, with the same frames for all tests.
Equivalently, using the support in $B_{R_0}$, these are the moments
\eqref{i2:eq-nine} with $R_1=R$.

\emph{Relative transport and recovery.}
Put $\mathsf F_{\rm mom}=Y_2^T(\mathsf P_1-\mathsf P_2)X_1$,
as in Proposition~\ref{if:intrinsic-relative}. This source is jointly
$C^2$ and has compact leaf support inside the frame disks.
By \eqref{fr:interior-cauchy}, the same differentiation applies to
the Cauchy transform $\mathscr C$ of \eqref{if:primitive-data}.
The moments and \eqref{if:cauchy-tail} give
\[
 \mathscr C\mathsf F_{\rm mom}\in C^2,\qquad
 \mathscr C\mathsf F_{\rm mom}=0\quad (|z|>R_0).
\]
The compactly supported primitive transforms by
\eqref{if:cauchy-frame-covariance} under the frame changes
\eqref{if:solution-frame-change}. The local formula
\eqref{if:intrinsic-formula} and Proposition~\ref{if:actualH}
therefore give a field jointly $C^2$ in $x$ and the real direction
parameters, with
\[
 H|_{E_{\mathrm{ext}}}=I_3,\qquad
 \det H=\sqrt{\mu_2/\mu_1}>0 .
\]
In the holomorphic representative from \eqref{if:normalization},
differentiation at fixed $x$ gives
\[
 V=\partial_{\bar\zeta}H|_x,\qquad
 2D_\theta V+\Gamma_1V-V\Gamma_2=0,\qquad
 V|_{E_{\mathrm{ext}}}=0.
\]
Proposition~\ref{if:binding} gives $V=0$ and the transition law
\eqref{if:chart}, with the corresponding conclusion for
$\partial_{\bar\upsilon}H|_x$ on the other chart.
The finite-type system \eqref{i4:closed-system} in
Theorem~\ref{i4:main} uses derivatives of $\ell_{12}$ through order
three. The extended Lam\'e pairs agree on $E_{\mathrm{ext}}$, so
that theorem gives $\lambda_1=\lambda_2$ and $\mu_1=\mu_2$ on $B_R$.

For the common Lam\'e pair, use the phases of \eqref{pt:phase}
for each fixed $\xi\ne0$. The construction above supplies the
amplitudes for both phase paths on a common interval $I$ near any
transverse center, as in Lemma~\ref{pt:product-lemma}. Its
product limit \eqref{pt:limit}, together with the weak identity
\eqref{pt:betti} and $W_1-W_2\in L^\infty(\Omega)$, yields
\[
 \int_\Omega e^{i\xi\cdot x}(W_1-W_2)
                 \mu^{-1}b(\widehat\xi\cdot x)^2\,dx=0
 \qquad\text{for every real }b\in C_c^\infty(I).
\]
Proposition~\ref{pt:recovery} removes these transverse weights and
applies Fourier uniqueness. Hence $W_1=W_2$ on $\Omega$ and, by
continuity, on $\overline\Omega$.
\end{proof}

\subsection{Completion of the uniqueness proof of Theorem~\ref{a28:main}}\label{a28:assembly}
\begin{proof}[Proof of Theorem~\ref{a28:main}]
For the $C^9$ boundary of Theorem~\ref{a28:main},
Proposition~\ref{bl:finite-jets} applied to the physical
Dirichlet-to-Neumann maps gives
\[
 \left.\partial_x^\alpha\delta_{12}\lambda\right|_{\partial\Omega}
 =\left.\partial_x^\alpha\delta_{12}\mu\right|_{\partial\Omega}=0
 \quad (|\alpha|\le16),\qquad
 \left.\partial_x^\alpha\delta_{12}W\right|_{\partial\Omega}=0
 \quad (|\alpha|\le14).
\]
These equalities supply the traces in \eqref{bt:matching} and
the boundary derivatives of the Lam\'e parameters used in Proposition~\ref{ext:finite-common}.
That proposition supplies the exterior extensions, so all
the hypotheses of Proposition~\ref{fr:interior} are satisfied.
It follows that $\lambda_1=\lambda_2$, $\mu_1=\mu_2$, and
$W_1=W_2$ on $\overline\Omega$. By \eqref{en:W-definition},
\[
 \rho_1-\rho_2=\omega^{-2}(W_1-W_2)=0\quad\text{on }\overline\Omega .
\]
\end{proof}

\section{Finite-time uniqueness for dynamical elasticity}
\label{sec:dynamical}

We prove Theorem~\ref{td:uniqueness} in two steps. First we extend the
equality of the response maps in \eqref{td:data-equality} to the positive
time axis under the shear-wave filling-time condition of
Belishev and Lasiecka~\cite{BL2002}. We then use the Laplace transform
to identify the two displacement-to-traction maps at a common positive
frequency to which Theorem~\ref{a28:main} applies.

\subsection{Dense terminal states at the shear filling time}

A separated control is a finite sum of smooth time factors times
traces of ambient smooth spatial fields. Such spatial traces are dense
in $H^{1/2}(\partial\Omega;\mathbb C^3)$ by the trace lifting and
$H^1$ extension theorems. Choose an orthonormal basis from the real
traces and use its complexification. Let $\Pi_N$ be the orthogonal
projection onto the first $N$ basis vectors. The projections
preserve real data, are contractions and converge strongly to the
identity.
For $f\in C_c^\infty((0,T);H^{1/2}(\partial\Omega;\mathbb C^3))$,
compactness of the time images of the derivatives gives
\begin{equation}\label{td:control-approximation}
 \sum_{k=0}^3\|\partial_t^k(\Pi_N f-f)\|_{L^\infty(0,T;H^{1/2})}
 \longrightarrow0 .
\end{equation}
The time support of $\Pi_N f$ is contained in that of $f$.
By \eqref{td:response-bound} in Proposition~\ref{an:dynamic-wellposedness},
the response on separated controls determines the response on the full
input class.

We suppress the medium index in \eqref{td:filling-time} and
\eqref{td:dirichlet-generator}.

\begin{lemma}\label{td:dense-terminal}
If $T>\tau_S$, the set of terminal states $u^f(T)$ for real separated
controls compactly supported in time in $(0,T)$ is dense in the real
state space $\mathcal H$.
\end{lemma}

\begin{proof}
Let $y\in\mathcal H$ be real and orthogonal to the terminal states.
For $0\le t\le2T$ use the homogeneous Dirichlet wave
\begin{equation}\label{td:adjoint-wave}
 w(t)=\frac{\sin((t-T)\sqrt{\mathsf A_D})}{\sqrt{\mathsf A_D}}y,\qquad
 w(T)=0,\qquad \partial_tw(T)=y.
\end{equation}
It lies in $C([0,2T];V)\cap C^1([0,2T];\mathcal H)$ and is odd about
$T$. Given a real spatial trace $g$ and $\phi\in C_c^\infty(0,T;\mathbb R)$,
use \eqref{td:lifted-wave} with $f=\phi g$. Duhamel's formula yields
\begin{equation}\label{td:terminal-duality}
 \begin{aligned}
 u^{\phi g}(T)
  &=-\int_0^T\frac{\sin((T-t)\sqrt{\mathsf A_D})}{\sqrt{\mathsf A_D}}
                   \phi''(t)Gg\,dt,\\
 0=(u^{\phi g}(T),y)_{\mathcal H}
  &=\int_0^T\phi''(t)(w(t),Gg)_{\mathcal H}\,dt .
 \end{aligned}
\end{equation}
The weak conormal of $w$ paired with $g$ is
$\partial_t^2(w,Gg)_{\mathcal H}$ in time distributions, since
$\mathfrak a(w(t),Gg)=0$. Thus $(w,Gg)_{\mathcal H}$ is affine on
$(0,T)$ by \eqref{td:terminal-duality}, and on $(T,2T)$ by oddness.
Its $C^1$ regularity at $T$ gives
\begin{equation}\label{td:zero-conormal}
 \gamma_0 w=0,\qquad \sigma_{\lambda,\mu}(w)\nu=0
 \quad\text{on }(0,2T)\times\partial\Omega .
\end{equation}
The second equality extends from the separated tests to general
smooth-in-time boundary tests by density of the spatial traces.

Extend the coefficients to an open neighborhood $\Omega_e\supset
\overline\Omega$ by the Whitney construction in the proof of
Proposition~\ref{ext:finite-common}. Shrink $\Omega_e$ to preserve the
positive bounds. Extend $w$ by zero to $\Omega_e\setminus\overline\Omega$.
The first condition in \eqref{td:zero-conormal} gives an $H^1$ spatial
extension. For a test restricted to $\Omega$, split it into a static
lift of its trace and a zero-trace part. The conormal equality in
\eqref{td:zero-conormal} treats the first part and the Dirichlet
equation treats the second. Consequently
$(\rho\partial_t^2-L_{\lambda,\mu})w=0$ on
$(0,2T)\times\Omega_e$ in distributions.

Fix $T'\in(\tau_S,T)$ and mollify $w$ in time by a smooth even
kernel of radius $\delta<T-T'$. The resulting $u=w_\delta$ solves the
same equation on a neighborhood of $[T-T',T+T']\times\Omega_e$
and vanishes off $\overline\Omega$. It is smooth in time with values
in $H^1_{\rm loc}(\Omega_e)$. On every compact time interval $J$ in
this neighborhood and for $O_1\Subset O_2\Subset\Omega_e$, interior
elliptic regularity gives
\begin{equation}\label{td:mollified-regularity}
 \|u\|_{L^2(J;H^2(O_1))}
 \le C\left(\|\rho\,u_{tt}\|_{L^2(J;L^2(O_2))}
             +\|u\|_{L^2(J;H^1(O_2))}\right).
\end{equation}
Together with $u_t\in L^2(J;H^1(O_2))$, this yields
$U=(u,v,\operatorname{curl}u)\in H^1_{\rm loc}$ in space-time, where
$v=\operatorname{div}u$.

Use the reduction of \cite[Appendix~A, Lemma~A.1]{dHLLO2023}.
With $b=(\lambda+\mu)/\rho$ and $q_k=(\partial_k\mu)/\rho$, we have
\begin{equation}\label{td:scalar-principal-parts}
 \begin{aligned}
 u_{tt}-c_S^2\Delta u-b\nabla v
    &=\rho^{-1}(\nabla\lambda)v
       +\sum_{k=1}^3 q_k(\partial_k u+\nabla u_k),\\
 \Delta u&=\nabla v-\operatorname{curl}\operatorname{curl}u,\qquad
 c_P^2=c_S^2+b=(\lambda+2\mu)/\rho .
 \end{aligned}
\end{equation}
In taking divergence and curl, the second derivatives of $u$ apart
from the principal parts reduce by
\begin{align*}
 \operatorname{div}(c_S^2\Delta u)-c_S^2\Delta v
   &=\nabla(c_S^2)\cdot\Delta u,\\
 \operatorname{curl}(c_S^2\Delta u)-c_S^2\Delta\operatorname{curl}u
   &=\nabla(c_S^2)\times\Delta u,\\
 \operatorname{div}(\partial_k u+\nabla u_k)
   &=\partial_k v+\Delta u_k,\qquad
 \operatorname{curl}(\partial_k u+\nabla u_k)
   =\partial_k\operatorname{curl}u .
\end{align*}
The second line of \eqref{td:scalar-principal-parts} expresses
$\Delta u$ by first derivatives of $U$. Thus, in distributions,
\begin{equation}\label{td:diagonal-wave-system}
 (P+\mathcal L)U=0,\qquad
 P=\diag(P_S I_3,P_P,P_S I_3),\qquad
 P_j=\partial_t^2-c_j^2\Delta\quad(j=S,P),
\end{equation}
where $\mathcal L$ has order at most one and locally bounded,
time-independent coefficients. The speeds $c_S,c_P$ are positive and
time-independent of class $C^1$.

We show the local continuation property for \eqref{td:diagonal-wave-system}.
Write $z=(t,x)$, let $\psi\in C^{2,\alpha}$ satisfy $\psi(z_*)=0$
and $d\psi(z_*)\ne0$, and suppose $U$ vanishes on $\{\psi>0\}$
near $z_*$. Assume the surface $\psi=0$ is noncharacteristic for
both speeds there. The principal symbols are
$p_j(z,\xi)=-\xi_t^2+c_j^2\xi_x\cdot\xi_x$. On the time-conormal set
$\xi_t=0$, the real characteristic set has no nonzero covectors.
For the complex surface condition, with $\varsigma>0$,
\[
 \operatorname{Im}\left[
   \partial_\xi p_j(z_*,\xi+i\varsigma d\psi(z_*))\cdot d\psi(z_*)
 \right]=2\varsigma\,p_j(z_*,d\psi(z_*))\ne0 .
\]
These are the surface conditions of
\cite[Definition~2.3 and Lemma~2.5]{Tataru1999}. Choose a single
$m>0$ sufficiently large that $\varphi=e^{m\psi}$ is strongly
pseudoconvex in the time-conormal set for both symbols. If
$T^2_{z_*}\varphi$ denotes the quadratic Taylor polynomial at $z_*$,
choose $\sigma>0$ small and then shrink the ball about $z_*$. The common
quadratic weight
\begin{equation}\label{td:shared-continuation-weight}
 \phi(z)=T^2_{z_*}\varphi(z)-\varphi(z_*)-\sigma|z-z_*|^2
\end{equation}
is strongly pseudoconvex in that set for both symbols, by stability
under small perturbations of the second-order boundary derivatives as in
\cite[Lemma~2.5]{Tataru1999} and
\cite[Section~3.1, Step~3]{BKL2018}. Taylor's formula gives
$\phi(z)\le\varphi(z)-1-(\sigma/2)|z-z_*|^2$ on a sufficiently small
ball. In particular $\phi(z_*)=0$ and $\phi<0$ on $\{\psi\le0\}$
there away from $z_*$.

Take $\chi\in C_c^\infty$ in a common ball for the scalar estimates,
equal to one near $z_*$. Put $V=\chi U$ and
$F=(P+\mathcal L)V=[P+\mathcal L,\chi]U$. We have $V\in H^1$,
$PV\in L^2$ and $\phi\le-b_*$ on $\supp F$ for some $b_*>0$,
since the derivatives of $\chi$ vanish near $z_*$. For $\tau\ge1$ set
\[
 \|Y\|_{1,\tau}^2=\|\nabla_{t,x}Y\|_{L^2}^2+\tau^2\|Y\|_{L^2}^2,
 \qquad
 \mathcal G_\tau=e^{-\eta D_t^2/(2\tau)},\quad D_t=-i\partial_t,
 \qquad
 V_\tau=e^{\tau\phi}V,\quad Z_\tau=\mathcal G_\tau V_\tau .
\]
The scalar Carleman estimate of \cite[Theorem~2.1 and
Section~3.2]{BKL2018}, applied to the diagonal blocks with the common
weight \eqref{td:shared-continuation-weight}, gives, for one fixed
$\eta>0$ small enough and all sufficiently large $\tau$,
\begin{equation}\label{td:scalar-carleman-estimate}
 \|Z_\tau\|_{1,\tau}
 \le C\tau^{-1/2}\|\mathcal G_\tau e^{\tau\phi}PV\|_{L^2}
        +Ce^{-a\tau}\|V_\tau\|_{1,\tau},
 \qquad a>0.
\end{equation}
All fields inside the Gaussian are extended by zero beyond the
cutoff region. Multiplication by the coefficients of $\mathcal L$
commutes with the time Gaussian. For each space-time derivative
$\partial_\ell$, the affine dependence of $\partial_\ell\phi$ on $t$
gives
\begin{align*}
 [\mathcal G_\tau,\tau\partial_\ell\phi]V_\tau
    &=\eta(\partial_t\partial_\ell\phi)\partial_t Z_\tau,\\
 \mathcal G_\tau e^{\tau\phi}\partial_\ell V
    &=\partial_\ell Z_\tau-\tau(\partial_\ell\phi)Z_\tau
         -\eta(\partial_t\partial_\ell\phi)\partial_t Z_\tau .
\end{align*}
Indeed, for the Gaussian kernel
$g_\tau(r)=(\tau/(2\pi\eta))^{1/2}e^{-\tau r^2/(2\eta)}$,
$-r g_\tau(r)=(\eta/\tau)g_\tau'(r)$.
On a fixed time slab containing the support of $V$ with a gap
$d_*>0$ to its boundary, $\partial_\ell\phi$ is bounded on the
spatial projection of $\supp V$. Outside this slab, the Gaussian
convolution uses only $|r|\ge d_*$ and
$|\partial_\ell\phi(t,x)|\le C(1+|r|)$ when $(t-r,x)\in\supp V$.
Young's inequality and
\[
 \|(1+|r|)\mathbf 1_{\{|r|\ge d_*\}}g_\tau(r)\|_{L^1(dr)}
       \le Ce^{-\tau d_*^2/(4\eta)}
\]
then yield, for some $a'>0$,
\begin{equation}\label{td:gaussian-lower-order}
 \|\mathcal G_\tau e^{\tau\phi}\mathcal LV\|_{L^2}
 \le C\|Z_\tau\|_{1,\tau}+Ce^{-a'\tau}\|V_\tau\|_{1,\tau}.
\end{equation}
Substitute $PV=F-\mathcal LV$ in
\eqref{td:scalar-carleman-estimate} and absorb the first term of
\eqref{td:gaussian-lower-order} after multiplication by $C\tau^{-1/2}$.
Since $\|V_\tau\|_{1,\tau}\le C(1+\tau)$ and
$\|e^{\tau\phi}F\|_{L^2}\le Ce^{-b_*\tau}$, this gives, for some
$\gamma>0$,
\begin{equation}\label{td:local-continuation}
 \|\mathcal G_\tau e^{\tau\phi}V\|_{L^2}\le Ce^{-\gamma\tau},
 \qquad
 \supp V\subset\{\phi\le-\gamma\}.
\end{equation}
The support inclusion follows from the first inequality by
\cite[Theorem~4]{Tataru1999}, component by component. As $\chi=1$
near $z_*$ and $\phi(z_*)=0$, this excludes $z_*$ from $\supp U$.

Fix $x_0\in\Omega$. The strict bound
$d_{g_S}(x_0,\partial\Omega)<T'$ gives a smooth embedded regular arc in
$\Omega_e$ from an exterior point to $x_0$ of shear length $T_0<T'$,
continued slightly past $x_0$. To obtain the arc, extend a
near-minimizing path by a short exterior segment in $\Omega_e$.
Polygonal approximation, small perturbations of vertices in
$\mathbb R^3$ to separate nonadjacent edges, and rounding of corners
preserve the strict length bound and the two endpoints. In a thin
Euclidean normal tube around the arc, let $s(x)$ be the shear
arclength of the axial foot point measured from the exterior
endpoint, and $r(x)$ the Euclidean distance to the axis in the normal
disk. Then $c_S|\nabla s|=1$ on the axis and
$\nabla s\cdot\nabla(r^2)=0$ throughout the tube.

Choose $T_0<T_1<T'$ and $1<\kappa_*<T'/T_1$ such that the distal cap
has $s>\kappa_*T_0$. Narrow the tube so that its initial cap is
exterior and $c_S|\nabla s|/T_0>1/T_1$. Take $r_0>0$ with
$r^2/r_0^2>\kappa_*$ on the outer tube wall, and
$0<\varepsilon<T_1(\kappa_*-1)$. On this tube set
\begin{equation}\label{td:continuation-surface}
 \begin{aligned}
 \Psi(t,x)&=\frac{\sqrt{(t-T)^2+\varepsilon^2}}{T_1}
                +\frac{s(x)}{T_0}+\frac{r(x)^2}{r_0^2},\\
 |\partial_t\Psi|&\le \frac1{T_1}
       <\frac{c_S|\nabla s|}{T_0}
       \le c_S|\nabla_x\Psi|\le c_P|\nabla_x\Psi|,\\
 \Psi(T,x_0)&=1+\varepsilon/T_1<\kappa_* .
 \end{aligned}
\end{equation}
By the chosen bounds, the closed sublevel $\Psi\le\kappa_*$ avoids
the time faces $|t-T|=T'$, the distal cap, and the outer tube wall.
The initial cap lies in the exterior zero region. If the support of
$(w_\delta,\operatorname{div}w_\delta,\operatorname{curl}w_\delta)$
meets this closed sublevel, $\Psi$ attains its minimum $\kappa$ on
the intersection at an interior point. The field vanishes on the
local side $\Psi<\kappa$. Let $z_*$ be the first contact point and set $\psi=\kappa-\Psi$.
By \eqref{td:continuation-surface},
\[
 \min_{j\in\{S,P\}}p_j(z_*,d\psi(z_*))
 =\left.\bigl(c_S^2|\nabla_x\Psi|^2-|\partial_t\Psi|^2\bigr)
   \right|_{z=z_*}>0 .
\]
Thus $\psi=0$ is noncharacteristic for both speeds, and $U=0$ locally
on $\{\psi>0\}$. The local continuation \eqref{td:local-continuation}
excludes $z_*$ from $\supp U$, a contradiction. The third line of
\eqref{td:continuation-surface} now gives vanishing of $w_\delta$ near
$(T,x_0)$. Since $x_0\in\Omega$ was arbitrary,
\[
 \partial_tw_\delta(T)=0\quad\text{in }\mathcal H,\qquad
 \partial_tw_\delta(T)\longrightarrow\partial_tw(T)=y
 \quad(\delta\to0)
\]
by the time continuity in \eqref{td:adjoint-wave}. This proves the density.
\end{proof}

\subsection{Continuation of the full-boundary response}

\begin{proposition}\label{td:response-continuation}
For the domains and media in Theorem~\ref{td:uniqueness}, let
$T>\max\{\tau_{S,1},\tau_{S,2}\}$. If $R_1^{2T}=R_2^{2T}$, then the
causal response maps agree on the positive time axis.
\end{proposition}

\begin{proof}
We first use real controls and real state spaces. Let $\mathcal M_T$
consist of separated controls on $[0,T]$ whose time factors vanish
near $0$ and extend smoothly by odd reflection about $T$. For
$f\in\mathcal M_T$, write $f_-$ for that reflection on $[0,2T]$.
It vanishes near both endpoints. Set
\[
 \mathcal W_jf=u_j^{f_-}(T),\qquad f\in\mathcal M_T .
\]
Each time derivative $\partial_t^k u_j^{f_-}$ solves the causal wave
equation with boundary value $\partial_t^k f_-$ by
\eqref{td:lifted-wave} and \eqref{td:energy-estimate}. In particular
$\partial_t^2u_j^{f_-}(T)\in\mathcal H_j$. Since $f(T)=0$ and
$\partial_t^2f\in\mathcal M_T$, we have
\begin{equation}\label{td:terminal-generator}
 \mathcal W_jf\in D(\mathsf A_{D,j}),\qquad
 \mathsf A_{D,j}\mathcal W_j f=-\mathcal W_j(\partial_t^2 f).
\end{equation}
For these controls, subtracting an $H^2$ spatial lift of the boundary
trace and applying Dirichlet elliptic regularity to
$L_j u_j^{f_-}(t)=\rho_j\partial_t^2u_j^{f_-}(t)$ gives
$u_j^{f_-}(t)\in H^2(\Omega)$. The corresponding tractions agree with
the weak tractions of \eqref{td:weak-response}.

\smallskip\noindent\emph{Density in the graph norm.}
Fix $j$, put $\mathsf A=\mathsf A_{D,j}$, and choose
$0<\delta<T-\max_j\tau_{S,j}$. The set of $\mathcal W_j f$ for
separated $f$ compactly supported in $(\delta,T)$ is dense in
$\mathcal H_j$ by Lemma~\ref{td:dense-terminal} at time $T-\delta$.
Choose a nonnegative even $\eta\in C_c^\infty(-1,1)$ with integral
one. For $0<\varepsilon<\delta$ set
$\eta_\varepsilon(s)=\varepsilon^{-1}\eta(s/\varepsilon)$ and
$f_\varepsilon=(\eta_\varepsilon*f_-)|_{[0,T]}$, using zero extension
of $f_-$ to $\mathbb R$. The convolution preserves odd reflection
and vanishing near $0$, so $f_\varepsilon\in\mathcal M_T$.
With $u=u_j^{f_-}$ extended by zero before $t=0$,
$(\eta_\varepsilon*u)(T)=\mathcal W_jf_\varepsilon$ by causality:
the convolved wave has the prescribed trace and vanishes near $0$.
Also $u(T+s)+u(T-s)$ has zero trace for $0\le s\le\varepsilon$,
and has displacement $2\mathcal W_j f$ and velocity zero at $s=0$.
Hence
\begin{equation}\label{td:terminal-smoothing}
 \begin{aligned}
 u(T+s)+u(T-s)&=2\cos(s\sqrt{\mathsf A})\mathcal W_jf,\\
 \mathcal W_j f_\varepsilon&=Q_\varepsilon(\mathsf A)\mathcal W_j f,
 \qquad
 Q_\varepsilon(\zeta)
 =\int_\mathbb R\eta_\varepsilon(s)\cos(s\sqrt{\zeta})\,ds
 \quad(\zeta\ge0).
 \end{aligned}
\end{equation}
Since $\eta_\varepsilon$ is smooth and compactly supported, its
Fourier transform gives
\[
 \sup_{\zeta\ge0}(1+\zeta)|Q_\varepsilon(\zeta)|<\infty,\qquad
 |Q_\varepsilon(\zeta)|\le1,\qquad
 Q_\varepsilon(\zeta)\longrightarrow1
 \quad(\varepsilon\to0,\ \zeta\ge0).
\]
Thus $Q_\varepsilon(\mathsf A):\mathcal H_j\to D(\mathsf A)$ is bounded
in the graph norm, and $Q_\varepsilon(\mathsf A)y\to y$ in that norm
for $y\in D(\mathsf A)$ by spectral dominated convergence. At each
fixed $\varepsilon$, approximate $y$ in $\mathcal H_j$ by terminal
states supported in time in $(\delta,T)$ and use
\eqref{td:terminal-smoothing}. The graph-norm closure of
$\mathcal W_j\mathcal M_T$ contains $Q_\varepsilon(\mathsf A)y$, hence
also $y$. This proves that $\mathcal W_j\mathcal M_T$ is a core.

\smallskip\noindent\emph{The terminal Gram form.}
For $f,g\in\mathcal M_T$ and $(s,t)\in[0,2T]\times[0,T]$, put
$B_j(s,t)=(u_j^{f_-}(s),u_j^{g_-}(t))_{\mathcal H_j}$. Green's
formula and the causal initial data give
\begin{equation}\label{td:gram-wave}
 \begin{aligned}
 (\partial_t^2-\partial_s^2)B_j(s,t)
   &=\langle(R_j^{2T}g_-)(t),f_-(s)\rangle
        -\langle(R_j^{2T}f_-)(s),g(t)\rangle,\\
 B_j(s,0)&=\partial_tB_j(s,0)=0 .
 \end{aligned}
\end{equation}
The d'Alembert formula on the backward triangle with vertex $(T,T)$
gives
\begin{equation}\label{td:terminal-gram}
 (\mathcal W_j f,\mathcal W_j g)_{\mathcal H_j}
 =-\frac12\int_0^T\int_t^{2T-t}
        \langle (R_j^{2T}f_-)(s),g(t)\rangle\,ds\,dt.
\end{equation}
The first term of \eqref{td:gram-wave} has zero integral in $s$ by
odd symmetry of $f_-$ about $T$. The weak traction bound makes all
pairings integrable. Equality of $R_j^{2T}$ therefore gives equality
of the Gram forms. Since both ranges are dense, the assignment
\[
 \mathcal U(\mathcal W_1f)=\mathcal W_2f\quad(f\in\mathcal M_T)
\]
extends to an orthogonal isomorphism on the real spaces and to a
unitary $\mathcal U:\mathcal H_1\to\mathcal H_2$. On the core of the
first generator, \eqref{td:terminal-generator} gives
$\mathcal U\mathsf A_{D,1}\mathcal W_1f
=\mathsf A_{D,2}\mathcal U\mathcal W_1f$. Taking graph closures
gives
\begin{equation}\label{td:generator-intertwining}
 \mathcal U D(\mathsf A_{D,1})=D(\mathsf A_{D,2}),\qquad
 \mathcal U\mathsf A_{D,1}=\mathsf A_{D,2}\mathcal U
 \quad\text{on }D(\mathsf A_{D,1}).
\end{equation}

\smallskip\noindent\emph{The response after the source has stopped.}
Let $f$ be a real separated control compactly supported in $(0,T)$,
and extend it by zero after $T$. The terminal data for its free
evolution are $(\mathcal W_jf,\mathcal W_j\partial_t f)$.
Under \eqref{td:generator-intertwining}, the spectral formulas give
\[
 \mathcal U u_1^f(t)=u_2^f(t),\qquad
 \mathcal U\partial_tu_1^f(t)=\partial_tu_2^f(t)\qquad(t\ge T).
\]
Fix a free interval of length $T$ and call the corresponding time
translate of the free wave $w_j$. For a real separated probe $h$
compactly supported in $(0,T)$, Green's formula with $u_j^h$ gives
\begin{equation}\label{td:free-traction}
 \begin{aligned}
 \int_0^T\langle\sigma_{\lambda_j,\mu_j}(w_j(t))\nu,h(t)\rangle\,dt
   &=\int_0^T\bigl((\partial_t^2w_j,u_j^h)_{\mathcal H_j}
                 -(w_j,\partial_t^2u_j^h)_{\mathcal H_j}\bigr)\,dt\\
   &=(\partial_t w_j(T),\mathcal W_jh)_{\mathcal H_j}
       -(w_j(T),\mathcal W_j(\partial_t h))_{\mathcal H_j}.
 \end{aligned}
\end{equation}
The first line uses the zero trace of $w_j$ and
$\rho_j\partial_t^2w_j=L_j w_j$; the initial terms in the second line
vanish by causality of $u_j^h$. The endpoint pairings coincide by the
unitary correspondence. Density of the separated probes gives
equality of the tractions on every such interval. The measured
response gives equality for times preceding $T$ and in a neighborhood
of $T$.

A smooth time partition of unity expresses each compactly supported
separated control on $(0,\infty)$ as a finite sum of time translates
of controls compactly supported in $(0,T)$. Time translation, causality
and linearity extend the equality to their sums. Applying
\eqref{td:control-approximation} and \eqref{td:response-bound} on
each finite interval, and then complexifying the controls, gives
\begin{equation}\label{td:full-response}
 R_1^\infty=R_2^\infty
\end{equation}
for the causal response maps on the positive time axis.
\end{proof}

\subsection{The Laplace transform and proof of Theorem~\ref{td:uniqueness}}

\begin{proof}[Proof of Theorem~\ref{td:uniqueness}]
Choose $\tau_{\rm c}$ with
$\max\{\tau_{S,1},\tau_{S,2}\}<\tau_{\rm c}<T_{\rm obs}/2$.
Zero extension of inputs and restriction of outputs in
\eqref{td:data-equality} give $R_1^{2\tau_{\rm c}}=R_2^{2\tau_{\rm c}}$
by causality. Proposition~\ref{td:response-continuation} gives
\eqref{td:full-response}.

\smallskip\noindent\emph{The full-time Laplace transform.}
Fix a trace $g$ of an ambient smooth field and a nonzero
$\phi\in C_c^\infty((0,\infty);\mathbb R)$ with $\phi\ge0$.
Let $u_j=u_j^{\phi g}$. Given $s>0$, set
\[
 \widehat u_j(s)=\int_0^\infty e^{-st}u_j(t)\,dt\in H^1(\Omega;\mathbb C^3),
 \qquad
 \widehat\phi(s)=\int_0^\infty e^{-st}\phi(t)\,dt>0.
\]
By the zero initial data and \eqref{td:time-bound}, two integrations
by parts give
$\int_0^\infty e^{-st}\partial_t^2u_j(t)\,dt=s^2\widehat u_j(s)$ in
$\mathcal H_j$. The weak response formula \eqref{td:weak-response}
and a bounded trace lift give
\[
 \int_0^\infty e^{-st}\|R_j^\infty(\phi g)(t)\|_{H^{-1/2}}\,dt
 \le C_j\int_0^\infty e^{-st}
       \bigl(\|u_j(t)\|_{H^1(\Omega)}
                 +\|\partial_t^2u_j(t)\|_{\mathcal H_j}\bigr)\,dt<\infty.
\]
Taking the transform in \eqref{td:wave} and in the weak traction pairing
yields
\begin{equation}\label{td:laplace-system}
 \begin{aligned}
 (L_j-s^2\rho_jI_3)\widehat u_j(s)&=0,\\
 \gamma_0\widehat u_j(s)&=\widehat\phi(s)g,\\
 \left.\sigma_{\lambda_j,\mu_j}(\widehat u_j(s))\nu\right|_{\partial\Omega}
   &=\int_0^\infty e^{-st}R_j^\infty(\phi g)(t)\,dt .
 \end{aligned}
\end{equation}
The equation and traction in \eqref{td:laplace-system} are the weak
stationary ones at spectral parameter $z=-s^2$.

\smallskip\noindent\emph{The stationary family.}
For $z\in\mathbb C\setminus\operatorname{spec}(\mathsf A_{D,j})$
and $f\in H^{1/2}(\partial\Omega;\mathbb C^3)$, the weak solution of
$(L_j+z\rho_jI_3)u_j(z,f)=0$, $\gamma_0u_j(z,f)=f$, is
\begin{equation}\label{td:poisson-resolvent}
 u_j(z,f)=G_j f+z(\mathsf A_{D,j}-z)^{-1}G_j f.
\end{equation}
Indeed, the second term has zero trace and
\[
 (-\rho_j^{-1}L_j-z)G_j f=-zG_j f,\qquad
 (\mathsf A_{D,j}-z)\bigl[z(\mathsf A_{D,j}-z)^{-1}G_j f\bigr]=zG_j f .
\]
Uniqueness follows by applying the resolvent to the zero-trace
difference of two solutions in the domain \eqref{td:dirichlet-generator}.
Define the stationary weak traction map $\widetilde\Lambda_j^\Omega(z):
H^{1/2}(\partial\Omega;\mathbb C^3)\to H^{-1/2}(\partial\Omega;\mathbb C^3)$
by
\begin{equation}\label{td:weak-dn}
 \begin{aligned}
 \langle\widetilde\Lambda_j^\Omega(z)f,g\rangle
  =\int_\Omega\bigl[
    &\lambda_j\operatorname{div}u_j(z,f)\operatorname{div}v_g
      +2\mu_j\varepsilon(u_j(z,f)):\varepsilon(v_g)\\
    &-z\rho_j u_j(z,f)\cdot v_g\bigr]dx,
 \end{aligned}
\end{equation}
where $v_g\in H^1(\Omega;\mathbb C^3)$ is a bounded lifting of $g$.
The equation makes this bilinear pairing independent of the lifting.
The resolvent is holomorphic in the graph-norm domain, which embeds
in $H^1$ by \eqref{td:elliptic-bound}. Thus
\eqref{td:poisson-resolvent} and \eqref{td:weak-dn} give operator
holomorphy of $\widetilde\Lambda_j^\Omega$. In the normalization of
\eqref{en:dn-weak-form},
\begin{equation}\label{td:positive-frequency-map}
 \widetilde\Lambda_j^\Omega(\omega^2)
 =\Lambda_{\lambda_j,\mu_j,\rho_j}^{\Omega}(\omega)
 \qquad
 \bigl(\omega>0,\ \omega^2\notin\operatorname{spec}(\mathsf A_{D,j})\bigr).
\end{equation}

From \eqref{td:full-response} and \eqref{td:laplace-system},
\[
 \widehat\phi(s)
 \bigl(\widetilde\Lambda_1^\Omega(-s^2)-\widetilde\Lambda_2^\Omega(-s^2)\bigr)
 g=0 .
\]
Division by $\widehat\phi(s)$ and density of the spatial traces give
\begin{equation}\label{td:negative-parameter-equality}
 \widetilde\Lambda_1^\Omega(-s^2)=\widetilde\Lambda_2^\Omega(-s^2)
 \qquad(s>0).
\end{equation}
Let $\alpha_D=\min\{\alpha_{D,1},\alpha_{D,2}\}>0$. For fixed
$f,g$, the scalar pairing of
$\widetilde\Lambda_1^\Omega(z)-\widetilde\Lambda_2^\Omega(z)$ is
holomorphic on the common half-plane $\operatorname{Re}z<\alpha_D$
and vanishes at $z=-s^2$ for all $s>0$. The identity theorem yields
equality of the two maps throughout the half-plane. At the choice
$\omega_*=\sqrt{\alpha_D/2}$, \eqref{td:spectral-gap} gives
\[
 \mathfrak a_j(v,v)-\omega_*^2\|v\|_{\mathcal H_j}^2
 \ge \left(1-\frac{\omega_*^2}{\alpha_{D,j}}\right)\mathfrak a_j(v,v)
 \ge \tfrac12\mathfrak a_j(v,v)\qquad(v\in V).
\]
Therefore the physical Dirichlet operators at $\omega_*$ satisfy
\eqref{en:nonresonance}, and \eqref{td:positive-frequency-map} gives
\[
 \Lambda_{\lambda_1,\mu_1,\rho_1}^\Omega(\omega_*)
 =\Lambda_{\lambda_2,\mu_2,\rho_2}^\Omega(\omega_*).
\]
Apply Theorem~\ref{a28:main} at this frequency to obtain
$\lambda_1=\lambda_2$, $\mu_1=\mu_2$, and $\rho_1=\rho_2$ on
$\overline\Omega$.
\end{proof}

\begin{remark}[The terminal term in the finite-time transform]
\label{td:qualitative}
For a solution with zero initial data, the equation for the finite-time
Laplace transform is
\begin{equation}\label{td:finite-laplace-term}
 (L_j-s^2\rho_j I_3)\int_0^{T_{\rm obs}}e^{-st}u_j(t)\,dt
 =e^{-sT_{\rm obs}}\rho_j
        \bigl(\partial_tu_j(T_{\rm obs})+su_j(T_{\rm obs})\bigr).
\end{equation}
The right-hand side involves the interior terminal data.
Response continuation \eqref{td:full-response} permits the
full-time transform. The bound \eqref{td:time-bound} makes the
exponentially weighted terminal term vanish as the upper endpoint
tends to infinity.
\end{remark}

\begin{remark}[A uniform sufficient observation time]\label{td:uniform-time}
Suppose both media satisfy $\mu_j\ge\mu_->0$ and $\rho_j\le\rho_+$ on
$\overline\Omega$. Set
$r_{\rm in}(\Omega)=\sup_{x\in\Omega}\operatorname{dist}_E(x,\partial\Omega)$.
For $x\in\Omega$ let $y$ be a nearest Euclidean boundary point. The segment
$[x,y)$ lies in $\Omega$ by the definition of distance. Therefore
\[
 d_{g_{S,j}}(x,\partial\Omega)
 \le |x-y|\sqrt{\rho_+/\mu_-},\qquad
 \tau_{S,j}\le r_{\rm in}(\Omega)\sqrt{\rho_+/\mu_-}.
\]
Consequently, the condition
\begin{equation}\label{td:uniform-threshold}
 T_{\rm obs}>2r_{\rm in}(\Omega)\sqrt{\rho_+/\mu_-}
\end{equation}
suffices for Theorem~\ref{td:uniqueness} uniformly over the class of
media with these bounds.
\end{remark}

\appendix

\section{The free inverse and Fourier operator estimates}
\label{fr27:section}

We construct a scalar free inverse, acting componentwise on vectors.
We use the bilinear convention of Section~\ref{sec:notation}, the weighted
Sobolev spaces of Section~\ref{r27:inverse-inputs}, and the Fourier convention
of Section~\ref{r27:kernel-estimates}. In this section the weight exponent associated with $\langle x\rangle$ is
written $s_{\mathrm{wt}}$ to distinguish it from the boundary-layer
phase sign $\mathfrak{s}$ used in Section~\ref{bl:section}.
The derivatives are distributional;
equivalent unweighted Bessel norms may be used in the Fourier-kernel arguments below.

\begin{theorem}\label{fr27:free}
Let $\delta=\tau(\mathbf v_0+i{\mathbf v_1})$, where $\tau\geq1$ and
$\mathbf v_0,{\mathbf v_1}\in\mathbb R^3$ are orthonormal, and put
\[
 p_\delta(\xi)=|\xi|^2+2\delta\cdot\xi.
\]
The ordinary locally integrable function $m_\delta=1/p_\delta$ defines on
Schwartz functions the selected operator
$E_\delta f=\mathcal F^{-1}(m_\delta\widehat f)$, where $\mathcal F$
denotes the Fourier transform and $\mathcal F^{-1}$ its inverse. It has a unique bounded
extension
\[
 E_\delta:L^2_1(\mathbb R^3)\longrightarrow H^1_{-1}(\mathbb R^3)
\]
with
\begin{equation}\label{fr27:free-estimate}
 \tau\|E_\delta f\|_{L^2_{-1}}
       +\|E_\delta f\|_{H^1_{-1}}
 \leq C\|f\|_{L^2_1}.
\end{equation}
The same output satisfies both estimates, and
\begin{equation}\label{fr27:right-identity}
 (D+\delta)^2E_\delta f=f\quad\hbox{in }\mathcal S'(\mathbb R^3).
\end{equation}
The constant is independent of $\tau$ and of the orthonormal pair.
\end{theorem}

\begin{proof}
Rotate to $(\mathbf v_0,{\mathbf v_1})=(e_1,e_2)$; radial weighted derivative norms are
unchanged up to fixed dimension constants under this rotation.  Write
$x'=(x_1,x_3)$, $t=x_2$, $u=e^{-i\tau x_1}v$, and
$F=e^{i\tau x_1}f$. A tilde denotes Fourier transformation in
$x'=(x_1,x_3)$ alone, with dual variable $\eta'=(\eta_1,\eta_3)$;
the normal variable $t=x_2$ is retained. Tangential Fourier transformation
of the conjugated equation gives
\begin{equation}\label{fr27:fiber}
 L_b=-\partial_t^2+2\tau\partial_t+b^2-\tau^2
     =-(\partial_t-a_+)(\partial_t-a_-),\qquad
 b=|\eta'|,\quad a_\pm=\tau\pm b.
\end{equation}
For a nonzero real number $a$ set
\[
 G_a(t)=
 \begin{cases}-\mathbf1_{t<0}e^{at},&a>0,\\
                \mathbf1_{t>0}e^{at},&a<0.
 \end{cases}
\]
Its value has jump one at zero, so
$(\partial_t-a)G_a=\delta_0$ distributionally, where $\delta_0$
is the Dirac mass at $t=0$. For $b\ne\tau$ put
$K_b=-G_{a_+}*G_{a_-}$.  Both factors are integrable for these fibers.
Direct convolution gives
\begin{equation}\label{fr27:explicit-kernels}
 K_b(t)=
 \begin{cases}
 \displaystyle\mathbf1_{t<0}
        \frac{e^{(\tau+b)t}-e^{(\tau-b)t}}{2b},&0<b<\tau,\\[5pt]
 \displaystyle\frac1{2b}
  \begin{cases}e^{(\tau+b)t},&t<0,\\e^{(\tau-b)t},&t\ge0,\end{cases}
       &b>\tau,\\[7pt]
 \mathbf1_{t<0}t e^{\tau t},&b=0.
 \end{cases}
\end{equation}
In each case $K_b$ is continuous at zero and
$K_b'(0+)-K_b'(0-)=-1$.  The regular part solves $L_bK_b=0$ away from zero;
the distributional second derivative therefore gives
$L_bK_b=\delta_0$, including the double-root fiber $b=0$.

The characteristic fiber $b=\tau$ has zero measure in the two
tangential frequency variables. The construction uses only the zeroth-order frequency multiplier,
since radial-frequency derivatives may fail to be uniformly bounded
across the critical circle $b=\tau$. On this fiber one may take either
of the following representatives:
\begin{equation}\label{fr27:critical-kernels}
 K_\tau^-(t)=\frac{\mathbf1_{t<0}(e^{2\tau t}-1)}{2\tau},
 \qquad
 K_\tau^+(t)=\frac{\mathbf1_{t<0}e^{2\tau t}
                         +\mathbf1_{t\ge0}}{2\tau}.
\end{equation}
Both are continuous representatives at zero, both have derivative jump
$-1$, and both solve $L_\tau K=\delta_0$.  Their difference is the constant
$1/(2\tau)$.  They are the respective one-sided limits of
\eqref{fr27:explicit-kernels}.  The direct integral with respect to $d\eta'$ depends only on the
almost-everywhere values of $K_b$.

The factorization also proves bounds uniform through the critical circle:
\begin{equation}\label{fr27:kernel-bounds}
 \|G_{a_+}\|_1=(\tau+b)^{-1},\qquad
 \|G_{a_-}\|_\infty\leq1,\qquad
 K_b'=a_-K_b-G_{a_+},
\end{equation}
and hence
\[
 |K_b(t)|\leq(\tau+b)^{-1},\qquad
 b|K_b(t)|\leq1,\qquad |K_b'(t)|\leq2
\]
for almost every tangential frequency.  The same inequalities hold for the
chosen critical representatives.  In particular the apparent factor
$1/(\tau-b)$ never enters these bounds.

Let $w(t)=\langle t\rangle^{-1}$, and for a fiber source $h$ write
$h=w g$, so $\|g\|_2=\|\langle t\rangle h\|_2$.  The kernel from $g$ to
$w(K_b*h)$ is $w(t)K_b(t-s)w(s)$.  Its Hilbert--Schmidt norm is at most
\[
 \frac{\|w\|_2^2}{\tau+b}.
\]
For the differentiated kernels $bK_b$ and $K_b'$ the corresponding bounds
are $\|w\|_2^2$ and $2\|w\|_2^2$.  Indeed these statements follow by
squaring the pointwise kernel bounds and integrating $w(t)^2w(s)^2$.
Tangential Plancherel, applied to
$\widetilde v(\eta',t)=K_{|\eta'|}*\widetilde F(\eta',\cdot)(t)$, now gives
\begin{align*}
 \tau\|wv\|_2&\leq C\|\langle t\rangle F\|_2,\\
 \|w\partial_tv\|_2+\|w\partial_{x_1}v\|_2
                    +\|w\partial_{x_3}v\|_2
       &\leq C\|\langle t\rangle F\|_2.
\end{align*}
Undoing the modulation replaces $\eta_1$ by $\eta_1-\tau$; its extra
$\tau K_b$ term is bounded by one.  Finally
\[
 \langle x\rangle^{-1}\leq\langle t\rangle^{-1},\qquad
 \langle t\rangle\leq\langle x\rangle
\]
prove \eqref{fr27:free-estimate} for this direct-integral operator.

It remains to identify this operator with exactly the stated multiplier.
Put $s=\eta_2$.  On every noncritical fiber, direct Fourier transformation
of the integrable kernels gives
\begin{equation}\label{fr27:fiber-multiplier}
 \widehat K_b(s)=
 \frac1{s^2+2i\tau s+b^2-\tau^2}.
\end{equation}
The joint characteristic set in $(\eta',s)$ is
$\{b=\tau,s=0\}$.  Near it, the map
\[
 (b,s)\longmapsto
 (u,v)=(b^2+s^2-\tau^2,2\tau s)
\]
has Jacobian $4b\tau\ne0$.  In tubular charts the absolute reciprocal is
therefore bounded by a constant times $(u^2+v^2)^{-1/2}$.  Its normal
integral is bounded by
$C\int_0^\varepsilon \rho^{-1}\rho\,d\rho<\infty$.
The tangential polar-coordinate density is $b$, which is positive
near $b=\tau$. At $b=0$, the denominator is $(s+i\tau)^2$ and
has modulus $s^2+\tau^2$ for real $s$.
Outside a ball depending on the fixed $\tau$, the denominator has size at
least $c|\eta|^2$.  Consequently the reciprocal is a locally integrable
tempered function, and its product with any Schwartz test is absolutely
integrable.

The direct-integral distribution pairs against a
Schwartz function $\Phi(\eta',t)$ by an absolutely convergent integral:
\[
 \int\!\int |K_b(t)\Phi(\eta',t)|\,dt\,d\eta'
 \leq \int\frac{1}{\tau+b}
                \int|\Phi(\eta',t)|\,dt\,d\eta'<\infty.
\]
For almost every $\eta'$, apply \eqref{fr27:fiber-multiplier} to the
one-dimensional Schwartz test.  The integral of the resulting frequency
pairing is absolutely convergent by the preceding local-integrability and
large-frequency bounds.  Fubini therefore identifies the full Fourier
transform of this direct integral with the ordinary function in
\eqref{fr27:fiber-multiplier} on the entire frequency space.
Undoing the modulation shifts
$\eta_1=\xi_1+\tau$, giving precisely $m_\delta=1/p_\delta$.
Since $p_\delta m_\delta=1$ almost everywhere, the same identity holds as
distributions.  This proves \eqref{fr27:right-identity} on Schwartz inputs.

Smooth compactly supported functions are dense in $L^2_1$: first truncate
with a smooth cutoff at radius $R$, then mollify on the resulting compact
set.  The bound extends the same operator to all of $L^2_1$.  The limiting
weighted derivative norms define the claimed $H^1_{-1}$ output, and local
$L^2$ convergence implies distributional convergence; hence the right
identity also extends.  One may still interpret the extension as the
multiplier product: $f\in L^2_1$ implies $\widehat f\in H^1$, so locally
$\widehat f\in L^6$, while $m_\delta\in L^{6/5}_{\rm loc}$ by its two
normal-dimensional singularity.  H\"older's inequality controls the product
and its approximation by Schwartz inputs.  At infinity the quadratic
denominator and Cauchy--Schwarz suffice.  This also identifies the extension with the multiplier product.
\end{proof}

\subsection{Fourier kernel estimates with finite derivatives}
\label{fr27:fourier-section}
For a fixed compact set $K\subset\mathbb R^3$ and an $m\times m$ matrix
symbol supported in $K$ in the $x$ variable, define for nonnegative
integers $N,r$
\[
 s_{N,r}(a)=\max_{|\alpha|\leq N,\,|\beta|\leq r}
        \sup_{x,\xi}\|\partial_x^\alpha\partial_\xi^\beta a(x,\xi)\|.
\]
We derive the exact identities for smooth symbols on Schwartz
inputs. The estimates and extensions use the finite seminorms
specified below. Their constants depend on $K$, the dimension, and
the matrix norm, independently of the input support and any auxiliary
frequency cutoff.
Left quantization is
\[
 \operatorname{Op}(a)u(x)
 =(2\pi)^{-3}\int e^{ix\cdot\xi}a(x,\xi)\widehat u(\xi)\,d\xi.
\]

\begin{lemma}\label{fr27:schur}
For $j=0,1$,
\begin{equation}\label{fr27:ordinary-symbol}
 \|\operatorname{Op}(a)\|_{H^j\to H^j}
 \leq C_K s_{6,0}(a).
\end{equation}
\end{lemma}

\begin{proof}
As in Section~\ref{r27:kernel-estimates}, $\widehat a(v,\xi)$ denotes Fourier transformation in $x$ alone.
Applying $(1-\Delta_x)^3$ before integration gives
\begin{equation}\label{fr27:symbol-transform}
 \|\widehat a(v,\xi)\|
 \leq C_Ks_{6,0}(a)\langle v\rangle^{-6}.
\end{equation}
Here compact support bounds the $L^1_x$ norms of all derivatives in that
sixth-order expression.  The full Fourier kernel is
$(2\pi)^{-3}\widehat a(\eta-\xi,\xi)$.
After conjugation by the Bessel weights of $H^j$, its norm is at most
\[
 C_Ks_{6,0}(a)\langle\eta-\xi\rangle^{j-6},
\]
because $\langle\eta\rangle^j/\langle\xi\rangle^j
\leq C_j\langle\eta-\xi\rangle^j$.
Both Schur integrals are bounded by
$C_Ks_{6,0}(a)\int\langle v\rangle^{j-6}dv<\infty$.  This proves the
claim first for Schwartz inputs and then by density in $H^j$.
\end{proof}

\begin{proposition}
\label{fr27:weight-exchange}
For a symbol $a$ with compact support in $x$ and $s_{6,2}(a)<\infty$, the following holds for $j=0,1$:
\begin{equation}\label{fr27:weighted-local}
 \|\operatorname{Op}(a)u\|_{H^j}
 \leq C_Ks_{6,2}(a)\|u\|_{H^j_{-1}}.
\end{equation}
The extension from Schwartz inputs is unique, has support in $x$ contained in $K$,
and retains all contributions of the input outside $K$.  On Schwartz
inputs the exact identity is
\begin{equation}\label{fr27:exact-weight}
 \operatorname{Op}(a)\langle x\rangle^2
 =\langle x\rangle^2\operatorname{Op}(a)
   -2i\sum_{k=1}^3x_k\operatorname{Op}(\partial_{\xi_k}a)
   -\operatorname{Op}(\Delta_\xi a).
\end{equation}
\end{proposition}

\begin{proof}
Since $\widehat{x_ku}=i\partial_{\xi_k}\widehat u$, integration by parts
gives
$\operatorname{Op}(a)x_k=x_k\operatorname{Op}(a)
-i\operatorname{Op}(\partial_{\xi_k}a)$.
Applying the identity twice and summing over $k$ gives
\eqref{fr27:exact-weight}. Schwartz decay and the bounded symbol
derivatives justify the integrations by parts.

For $u\in H^j_{-1}$ set $v=\langle x\rangle^{-2}u$.  For $j=0$,
$\|v\|_2\leq\|\langle x\rangle^{-1}u\|_2$.  For $j=1$ use
\[
 \partial_kv=\langle x\rangle^{-2}\partial_ku
                   -2x_k\langle x\rangle^{-4}u
\]
to get $\|v\|_{H^j}\leq C\|u\|_{H^j_{-1}}$.
Apply Lemma~\ref{fr27:schur} to $a$, each
$\partial_{\xi_k}a$, and $\Delta_\xi a$ in
\eqref{fr27:exact-weight}.  Their outputs are all supported in $K$.
The factors $x_k$ and $\langle x\rangle^2$ are therefore bounded
multipliers on these $H^j$ outputs, including their first derivatives.
This proves \eqref{fr27:weighted-local} using at most six
$x$-derivatives of the indicated symbols.

To justify the extension, $C_c^\infty$ is dense also in $H^1_{-1}$.
For a smooth radial cutoff $\chi(x/R)$, the derivative error
$R^{-1}(\partial\chi)(x/R)u$ has weighted $L^2$ norm bounded by
$C R^{-1}\|\mathbf1_{R\lesssim|x|\lesssim2R}
\langle x\rangle^{-1}u\|_2\to0$; the other terms follow by tail
convergence.  Mollification after truncation completes the density proof.
The weighted estimate makes the operator limit independent of the
approximating sequence.  Distributional limits of fields supported in the
fixed closed set $K$ remain supported there.  The estimate thus applies to the full weighted input.
\end{proof}

\paragraph{Localized Fourier multipliers.}
For $\psi\in C_c^\infty$ and a Fourier multiplier $a_0(\xi)$,
$\psi\operatorname{Op}(a_0)=\operatorname{Op}(\psi(x)a_0(\xi))$.
Proposition~\ref{fr27:weight-exchange} therefore applies whenever the
stated seminorms of this symbol with compact support in $x$ are finite.

\begin{proposition}
\label{fr27:composition}
Let $a,c$ be matrix symbols whose supports in $x$ lie in one fixed
compact set $K$, and let
\[
 A_1=\max_{k=1,2,3}s_{6,0}(\partial_{\xi_k}a),\qquad C_0=s_{6,0}(c).
\]
For $j=0,1$,
\begin{equation}\label{fr27:composition-bound}
 \|\operatorname{Op}(a)\operatorname{Op}(c)-\operatorname{Op}(ac)
          \|_{H^j\to H^j}
 \leq C_K A_1 C_0.
\end{equation}
In particular a uniform bound $A_1\leq C\varepsilon$ gives a defect
$O(\varepsilon)$ when $C_0$ is uniformly bounded.  The same conclusion
holds if the left symbol is a Fourier multiplier $a_0(\xi)$, with $A_1$
replaced by $\|\nabla_\xi a_0\|_\infty$.  A right Fourier multiplier
has zero composition defect in left quantization.
\end{proposition}

\begin{proof}
On Schwartz inputs the exact full Fourier kernel of the defect, including
the Fourier normalization, is
\begin{equation}\label{fr27:exact-defect-kernel}
 (2\pi)^{-6}\int
 [\widehat a(\eta-\zeta,\zeta)-\widehat a(\eta-\zeta,\xi)]
                  \widehat c(\zeta-\xi,\xi)\,d\zeta.
\end{equation}
The product $ac$ retains the matrix order in
\eqref{fr27:exact-defect-kernel}. The identity is first established
with compact frequency cutoffs, which are removed using the
integrable bounds below.
The mean value theorem at fixed first Fourier variable, followed by
\eqref{fr27:symbol-transform} for $\partial_{\xi_k}a$, bounds
\eqref{fr27:exact-defect-kernel} by
\[
 C_K A_1C_0\int |\zeta-\xi|
        \langle\eta-\zeta\rangle^{-6}
        \langle\zeta-\xi\rangle^{-6}\,d\zeta.
\]
For the $H^j$ conjugation use
\[
 \frac{\langle\eta\rangle^j}{\langle\xi\rangle^j}
 \leq C_j\langle\eta-\zeta\rangle^j
                \langle\zeta-\xi\rangle^j.
\]
The weighted kernel is thus bounded by
$C_KA_1C_0(F_j*G_j)(\eta-\xi)$, where
\[
 F_j(v)=\langle v\rangle^{j-6},\qquad
 G_j(v)=|v|\langle v\rangle^{j-6}.
\]
Both functions are integrable in dimension three for $j=0,1$; for
$j=1$ their large-radius powers are $-5$ and $-4$.  The two Schur
integrals are at most
$C_KA_1C_0\|F_j\|_1\|G_j\|_1$.  This proves
\eqref{fr27:composition-bound}.

If the left symbol is a multiplier, the kernel is instead
\[
 (2\pi)^{-3}[a_0(\eta)-a_0(\xi)]
                                  \widehat c(\eta-\xi,\xi).
\]
Its weighted bound is
$C_K\|\nabla a_0\|_\infty C_0
 |\eta-\xi|\langle\eta-\xi\rangle^{j-6}$,
which has finite Schur integrals by the same powers.  If the right symbol
is a multiplier $c_0(\xi)$, multiplication of $\widehat u$ by $c_0$
before left quantization gives
$\operatorname{Op}(a)\operatorname{Op}(c_0)=\operatorname{Op}(ac_0)$
exactly.  These statements also handle symbols that are finite sums of
symbols with compact support in $x$ and Fourier multipliers, by bilinearity; only
the left-multiplier/right-compact and compact/compact defects remain.
\end{proof}

\begin{remark}\label{fr27:budget}
Six Cartesian derivatives suffice in all ordinary, weighted-to-local, and
composition bounds at $j=0,1$.  Two frequency derivatives suffice for
the weight exchange.  The composition estimate uses six Cartesian
derivatives of one frequency derivative of its left symbol and six
Cartesian derivatives of its right symbol. The singular multiplier
$1/p_\delta$ is used by multiplication.  Applying the six-derivative bound to
a residual containing $-\Delta_xb$ requires eight Cartesian derivatives
of $b$; mixed derivatives such as
$\partial_x^\alpha\partial_\xi^\beta b$, $|\alpha|\leq8$,
$|\beta|\leq2$, are the frame bounds used in the localized construction.
\end{remark}

\section*{Acknowledgments}
The work of H. Diao was supported by the National Natural Science
Foundation of China (Grant No. 12371422) and the Fundamental Research
Funds for the Central Universities, Jilin University.

\end{document}